\def\VersionDateTime{6/April/2017}

\documentclass[12pt, leqno]{amsart}
\usepackage{amscd, amsmath, amssymb, amsfonts}
\usepackage{latexsym}
\usepackage{mathrsfs}
\usepackage{color}

\newtheorem{Theorem}{Theorem}[section]
\newtheorem{Proposition}[Theorem]{Proposition}
\newtheorem{Lemma}[Theorem]{Lemma}
\newtheorem{Corollary}[Theorem]{Corollary}
\newtheorem{Claim}{Claim}[Theorem]
\newtheorem{TheoremNoNum}{Theorem}[Theorem]
\newtheorem{PropositionNoNum}{Proposition}[Theorem]

\theoremstyle{definition}
\newtheorem{Definition}[Theorem]{Definition}
\newtheorem{Remark}[Theorem]{Remark}
\newtheorem{Example}[Theorem]{Example}

\newtheorem{Question}[Theorem]{Question}

\newcommand{\TT}{{\mathbb{T}}}
\newcommand{\GG}{{\mathbb{G}}}

\newcommand{\ZZ}{{\mathbb{Z}}}

\newcommand{\RR}{{\mathbb{R}}}

\newcommand{\PP}{{\mathbb{P}}}
\newcommand{\OO}{{\mathcal{O}}}

\newcommand{\Fscr}{{\mathscr{F}}}

\newcommand{\Lscr}{{\mathscr{L}}}
\newcommand{\Mscr}{{\mathscr{M}}}
\newcommand{\Nscr}{{\mathscr{N}}}
\newcommand{\Sscr}{{\mathscr{S}}}
\newcommand{\Uscr}{{\mathscr{U}}}
\newcommand{\Xscr}{{\mathscr{X}}}

\newcommand{\mfrak}{{\mathfrak{m}}}

\newcommand{\Hom}{\operatorname{\operatorname{Hom}}}

\newcommand{\Irr}{\operatorname{\operatorname{Irr}}}

\newcommand{\Div}{\operatorname{Div}}

\newcommand{\GL}{\operatorname{GL}}

\newcommand{\ord}{\operatorname{ord}}

\newcommand{\rest}[2]{\left.{#1}\right\vert_{{#2}}}  

\newcommand{\Spec}{{\operatorname{Spec}}}

\newcommand{\trop}{\operatorname{trop}}

\newcommand{\an}{{\operatorname{an}}}

\newcommand{\Aff}{\mathbb{A}}

\newcommand{\id}{\operatorname{id}}

\newcommand{\Supp}{\operatorname{Supp}}

\newcommand{\Dscr}{\mathscr{D}}

\newcommand{\st}{\operatorname{st}}
\newcommand{\val}{\operatorname{val}}
\newcommand{\Sing}{\operatorname{Sing}}

\newcommand{\Trop}{\operatorname{Trop}}

\usepackage{mathrsfs}

\newcommand{\Proof}{{\sl Proof.}\quad}
\newcommand{\QED}{{\unskip\nobreak\hfil\penalty50\quad\null\nobreak\hfil
{$\Box$}\parfillskip0pt\finalhyphendemerits0\par\medskip}}

\begin{document}

\title[tropicalization associated to a linear system]{Effective faithful tropicalizations associated to linear systems on curves}
\author{Shu Kawaguchi}
\address{Department of Mathematical Sciences, 
Doshisha University, Kyoto 610-0394, Japan}
\email{kawaguch@mail.doshisha.ac.jp}
\author{Kazuhiko Yamaki}
\address{Institute for Liberal Arts and Sciences, 
Kyoto University, Kyoto, 606-8501, Japan}
\email{yamaki.kazuhiko.6r@kyoto-u.ac.jp}
\date{\VersionDateTime}
\thanks{The first named author partially supported by KAKENHI 15K04817 and 25220701,  and the second named author partially supported by KAKENHI 26800012.}


\newcommand{\Proj}{\operatorname{\operatorname{Proj}}}
\newcommand{\Prin}{\operatorname{Prin}}
\newcommand{\Rat}{\operatorname{Rat}}
\newcommand{\zero}{\operatorname{div}}
\newcommand{\Func}{\operatorname{Func}}
\newcommand{\red}{\operatorname{red}}
\newcommand{\pr}{\operatorname{pr}}
\newcommand{\can}{\operatorname{can}}
\newcommand{\relin}{\operatorname{relin}}

\begin{abstract}
For a connected smooth projective curve $X$ of genus
 $g$, global sections of 
any line bundle $L$ with $\deg(L) \geq 2g+ 1$ give an embedding 
of the curve into projective space.  We consider an analogous 
statement for a Berkovich skeleton in nonarchimedean geometry, 
in which projective space is replaced by tropical projective space, and 
an embedding is replaced by a homeomorphism onto its image preserving integral structures (called a faithful tropicalization).  
Let $K$ be an algebraically closed field which is complete with respect to a non-trivial nonarchimedean value. 
Suppose that $X$ is defined over $K$ and has genus $g \geq 2$ and that $\Gamma$ is a skeleton (that is allowed to have ends) of the analytification $X^{\mathrm{an}}$ of 
 $X$ in the sense of Berkovich. 
We show that if $\deg(L) \geq 3g-1$, then global sections of $L$ give 
a faithful tropicalization of $\Gamma$ into tropical projective space. 
As an application, 
when $Y$ is a suitable affine curve, we describe the analytification $Y^{\mathrm{an}}$ as the limit 
of tropicalizations of an effectively bounded degree. 
\end{abstract}

\maketitle


\setcounter{equation}{0}
\section{Introduction}
\label{sec:intro}

Let $X$ be a connected smooth projective curve 
over a field $K$ of genus $g$, and 
let $L$ be an ample line bundle over $X$. 
If $\deg(L) \geq 2g+1$, then 
there exist nonzero global sections $s_0, \ldots, s_N \in H^0(X, L)$ such that 
the associated morphism 
\[
  X \longrightarrow \PP^N, \quad 
  p \mapsto (s_0(p): \cdots: s_N(p))
\]
is an embedding into projective space. 
In this paper, we consider an analogous statement for a Berkovich skeleton in nonarchimedean geometry, in which projective space is replaced by tropical projective space, and 
an embedding is replaced by a homeomorphism onto its image preserving integral structures,
called a faithful tropicalization. 

In the following, we assume that $K$ is an algebraically closed field 
which is complete with respect to a non-trivial nonarchimedean absolute value. 
Let $X^\an$ denote the analytification of $X$ in the sense of Berkovich, and let $\Gamma$ be a skeleton (which is allowed to have ends) of $X^{\an}$;
$\Gamma$ is a closed subset of $X^\an \setminus X(K)$ and has an  integral structure  (see \S\ref{subsec:Berko:skeleta} for details). 
The integral structures play important roles
in many aspects.
For example, they are used to describe
measures on the Berkovich spaces (cf. \cite{gubler1, gubler4,
yamaki5, yamaki6});
they are also of interest in mirror symmetry through the work of Kontsevich--Soibelman \cite{KS}. 

We set $\TT:= \RR \cup \{+\infty\}$. The $N$-dimensional tropical projective space is defined to be 
$\TT\PP^N = (\TT^{N+1}\setminus\{(+\infty, \ldots, +\infty)\})/\sim$, 
where $(x_0, \ldots, x_N) \sim (y_0, \ldots, y_N)$ if 
there exists $c \in \RR$ such that $y_i = x_i + c$ for any $i$ (see \cite{MZ}). 
We remark that $\TT\PP^{N}$ has a natural integral structure;
see \S\ref{subsec:tropical:geometry}.

A faithful tropicalization is a tropicalization that is faithful
with respect to the integral structures.
Nonzero global sections $s_0, \ldots, s_N \in H^0(X, L)$ define a
map 
\begin{equation}
\label{eqn:trop:map:L}
  \varphi: X^\an \longrightarrow \TT\PP^N, 
  \quad p = (p, |\cdot|) \mapsto \left(-\log|s_0(p)|: \cdots: -\log|s_N(p)|\right), 
\end{equation}
called a tropicalization map associated to $L$.
We say that \emph{$\varphi$ is a faithful tropicalization of $\Gamma$}
if its restriction to $\Gamma$ is a 
homeomorphism onto its image preserving integral structures.
We say that \emph{$L$ admits a faithful tropicalization of $\Gamma$}
or that
$\Gamma$ has a {\em faithful tropicalization associated to the linear system $|L|$} 
if 
for some nonzero $s_0, \ldots, s_N \in H^0(X, L)$,
$\varphi$ is a faithful tropicalization of $\Gamma$.
See \S\ref{subsec:unimodular:faithful} for details.

Then,
in view of the classical result of 
embedding of connected smooth projective curves
by linear systems,
it will be natural to ask the following question.

\begin{Question}
\label{question:first}
Does there exist a constant $d(g)$ depending only on $g \geq 0$ 
with the following property?
For any connected smooth projective curve $X$ over $K$ of genus $g$, for any
skeleton $\Gamma$ of $X^{\an}$,
and for any line bundle $L$ over $X$,
if $\deg (L) \geq d(g)$, then
$L$ admits a faithful tropicalization of $\Gamma$.
Furthermore, if there exists such one,
can we
give a concrete $d(g)$?
\end{Question}

The following is our main result, which answers Question~\ref{question:first}.

\begin{Theorem}
\label{thm:main:unimodular:faithful}
Let $K$ be an algebraically closed field that is complete with respect to a non-trivial nonarchimedean absolute value. 
Set
\[
t (g) :=
\begin{cases}
1 & \text{if $g=0$,} 
\\
3 & \text{if $g=1$,} 
\\
3g-1 & \text{if $g\geq 2$.} 
\end{cases}
\]
Let $X$ be a connected smooth projective curve over $K$
of genus $g$, let 
$\Gamma$ be a skeleton of $X^{\an}$,
and let
$L$ be a line bundle over $X$.
Then
if
$\deg (L) \geq t (g) $,
then 
$L$ admits 
a faithful tropicalization of $\Gamma$.
\end{Theorem}

As a consequence of Theorem~\ref{thm:main:unimodular:faithful},
we obtain
the following result on tropicalization associated to the pluri-canonical systems, because  
$g \geq 2$
implies
$3\deg(\omega_X) 
\geq 
t (g)$.

\begin{Corollary}
\label{cor:faithful:canonical}
Let $K$, $X$, and $\Gamma$ be as in Theorem~\ref{thm:main:unimodular:faithful}. 
Suppose that $g \geq 2$.
Then  $\Gamma$ has a faithful tropicalization associated to the pluri-canonical linear system $|\omega_X^{\otimes m}|$ for any $m \geq 3$. 
\end{Corollary}

Once one obtains an affirmative answer to Question~\ref{question:first},
the next natural question would be the following.

\begin{Question}
\label{Question:optimalbound}
What is the best lower bound $d_{{\rm best}}(g)$, which depends only on $g \geq 0$, with the following property? 
For any connected smooth projective curve $X$ over $K$ of genus $g$, for any
skeleton $\Gamma$ of $X^{\an}$,
and for any line bundle $L$ over $X$,
if $\deg (L) \geq d_{{\rm best}}(g)$, then
$L$ admits a faithful tropicalization of $\Gamma$.\end{Question} 

By Theorem~\ref{thm:main:unimodular:faithful}, 
$d_{{\rm best}}(g) \leq t (g)$.
Although we do not have the complete answer to 
Question~\ref{Question:optimalbound},
we have the answer for small $g$;
for $g = 0,1,2$, the bound in Theorem~\ref{thm:main:unimodular:faithful}
is optimal, i.e., $d_{{\rm best}}(g) = t (g)$.

\begin{Theorem}[cf. Propositions~\ref{prop:complement:0}, 
\ref{prop:complement:1}, \ref{prop:complement:2}]
\label{thm:best:possible}
Let $K$ be as in Theorem~\ref{thm:main:unimodular:faithful}. Then 
there exist a connected smooth projective curve $X$ 
over $K$ of genus $2$ \textup{(}resp. $1$, $0$\textup{)}, 
a skeleton $\Gamma$ of $X^{\an}$, and 
a line bundle $L$ of degree $4$ \textup{(}resp. $2$, $0$\textup{)} such that 
 $\Gamma$ does {\em not} have a faithful tropicalization associated to $|L|$.  
\end{Theorem} 

As one sees from the proofs of Propositions~\ref{prop:complement:0}, 
\ref{prop:complement:1}, and \ref{prop:complement:2},
line bundles that
we consider for the proof of Theorem~\ref{thm:best:possible} are not very ample. 
One might expect that if $L$ is very ample, then it should always
admit a faithful tropicalization of a given skeleton. 
In fact, this expectation is not true, which makes 
Question~\ref{Question:optimalbound} more interesting. 
We show the existence of a connected smooth projective curve 
with a {\em very ample} line bundle 
that does not admit a faithful tropicalization of any skeleton of the curve. 

\begin{Proposition}
\label{prop:countereg:very:ample}
Let $K$ be as in Theorem~\ref{thm:main:unimodular:faithful}. Let $d \geq 4$ be an integer. 
Then there exists a connected smooth projective plane curve $X$ of $\PP^2$ of degree $d$ over $K$ 
such that 
no skeleton of $X^{\an}$
has a faithful tropicalization
associated to $|\OO_X(1)|$, 
where $\OO_X(1)$ is the very ample line bundle
over $X$ given by restriction of $\OO_{\PP^2}(1)$. 
\end{Proposition}

\medskip
Let us explain the relation to previous
works on faithful tropicalization of skeleta.
Katz--Markwig--Markwig \cite{KMM, KMM2}, Baker--Payne--Rabinoff \cite{BPR}, 
Gubler--Rabinoff--Werner \cite{GRW}, Baker--Rabinoff \cite{BR} showed that for a smooth projective variety $X$
and any skeleton $\Gamma$ of $X^{\an}$, 
there exist nonzero rational functions $f_1, \ldots, f_N$ of $X$ such that
\[
  \psi: X^{\an} \dasharrow \RR^N, \quad p = (p, |\cdot|)\mapsto (-\log|f_1(p)|, \ldots, -\log|f_N(p)|)
\] 
gives a faithful tropicalization of $\Gamma$. 
Note that there exists a positive integer $d$ such that
for any line bundle $L$ over $X$
with
$\deg (L) \geq d$,
$1, f_1, \ldots, f_N$ can be regarded as elements of $H^0(X, L)$.
This means that
for any $X$ and $\Gamma$ as above,
there exists a positive integer 
$d (X, \Gamma)$, which may depend on $X$ and $\Gamma$, 
such that any line bundle $L$
with $\deg (L) \geq d ( X, \Gamma)$
admits a faithful tropicalization of $\Gamma$.
Theorem~\ref{thm:main:unimodular:faithful}, on the other hand,
tells us that once $g$ is fixed,
then one has
a uniform and effective lower bound for the degree of a line bundle
such that the line bundle
admits a faithful tropicalization
of any skeleton of $X^{\an}$, as far as 
$X$ is a connected smooth projective curve of genus $g$.

We remark that the tropicalization map \eqref{eqn:trop:map:L} is not injective on the whole $X^{\an}$.  We also remark that, from a different perspective, Hrushovski--Loeser--Poonen \cite{HLP} studies if $X^{\an}$ is embedded in the Euclidean space $\RR^N$ 
by a homeomorphism onto its image (where one does not consider 
integral structures).  
In another direction, Haase--Musiker--Yu \cite{HMY} and Amini \cite{Amini} study embeddings of {\em tropical curves} associated to linear systems. See also Cueto--Markwig \cite{CM} for an algorithmic side of tropicalizations of plane curves using tropical modifications.

\medskip
After we wrote up the first draft of this paper, 
through 
communications with Sam Payne, we 
realize that 
Theorem~\ref{thm:main:unimodular:faithful} has an 
unexpected application to 
limit tropicalizations of affine space curves.
We begin by  recalling Payne's result \cite{payne}
that the Berkovich analytification is the limit of the tropicalizations.
Let $Y$ be an affine variety over $K$.
We call a
closed embedding of $Y$ into affine space
an \emph{affine embedding of $Y$}.
Any affine embedding $\iota: Y \hookrightarrow \Aff^n$
induces
a standard tropicalization map
$Y (K) \to \TT^n$
given by 
$p \mapsto (-\log |x_1 (\iota (p) )|_K, \ldots, -\log |x_n (\iota (p))|_K)$,
where $x_1 , \ldots , x_n$ are the standard coordinates of $\Aff^n$ 
and $|\cdot|_K$ is the absolute value of $K$. 
The \emph{tropicalization
$\Trop(Y, \iota)$ of $Y$ with respect to $\iota$} 
is defined to be 
the closure of the image of the standard tropicalization map.
By \cite[Proposition~2.2]{payne},
this map extends to a unique continuous map
$\pi_\iota: Y^{\an} \to \Trop(Y, \iota)$. 
Let $I$ be the set of all affine  embeddings of $Y$.
One makes $I$ a directed set by declaring $\iota_1 \leq \iota_2$ 
for two elements $\iota_1: Y \hookrightarrow \Aff^{n_1}$ and $\iota_2: Y \hookrightarrow \Aff^{n_2}$ of $I$
if there exists a torus-equivariant morphism $\varphi: \Aff^{n_2} \to \Aff^{n_1}$ with $\varphi\circ \iota_2 = \iota_1$ (cf. \S\ref{sec:limit:trop}). 
If $\iota_1 \leq \iota_2$, then one has a natural
continuous map 
$\Trop(Y, \iota_2) \to \Trop(Y, \iota_1)$,
and thus $\left(\Trop(Y, \iota)\right)_{\iota \in I}$ constitutes an inverse system in the category of topological spaces. 
Then Payne's theorem
\cite[Theorem~1.1]{payne} asserts that the map 
\[
  {\rm LTrop}: Y^{\an} \to \varprojlim_{\iota \in I} \Trop(Y, \iota) 
\]
induced from $( \pi_{\iota} )_{\iota \in I}$ is a homeomorphism. This homeomorphism says that various tropicalizations of $Y$ give rise to 
the Berkovich space $Y^{\an}$ associated to $Y$.

Now, fix 
a distinguished closed embedding
$Y \subset \Aff^N = \Spec(K[z_1, \ldots, z_N])$.
Then for any positive integer $D$,
we have the notion of affine embeddings of degree at most $D$ with
respect to $z_1, \ldots, z_N$,
so that
we consider a bounded degree version of 
Payne's construction. 
To be precise,
we say that an affine embedding $\iota : Y \hookrightarrow \mathbb{A}^n$ has 
\emph{degree at most $D$} if there exist $f_1 , \ldots , f_n \in K [z_1 , \ldots , z_N]$ such that $\deg (f_j) \leq D$ for $j = 1, \ldots, n$ and $\iota = (f_1, \ldots, f_n)$. 
Let $I_{\leq D}$ denote the set of all affine embeddings of $Y$ whose degree is at most $D$. 
One sees that $I_{\leq D}$ is a directed set,
and 
thus $\left( \Trop(Y, \iota ) \right)_{\iota \in I_{\leq D}}$ 
again constitutes an inverse system. 
We denote by
\begin{equation}
\label{eqn:limit:tropical:map:d:intro}
\mathrm{LTrop}_{\leq D} : Y^{\an} \to \varprojlim_{\iota \in I_{\leq D}} \Trop(Y, \iota)
\end{equation}
the map induced from $( \pi_{\iota} )_{\iota \in I_{\leq D}}$. 

As an application of Theorem~\ref{thm:main:unimodular:faithful}, we show that 
when $Y$ is a suitable affine curve,
$Y^{\an}$ is homeomorphic to the inverse limit 
of tropicalizations of an {\em effectively bounded degree},
i.e., $\mathrm{LTrop}_{\leq D}$ is a homeomorphism for some
effectively computed $D$.
Let $Y \subset \Aff^N = \Spec(K[z_1, \ldots, z_N])$ be a connected smooth affine curve. We say that 
\emph{$Y$ has smooth compactification in $\PP^N$} if the 
closure of $Y$ in $\PP^N$ is smooth,
where we regard $\Aff^N$ as an open subset of $\PP^N$. 
The \emph{degree $d$ of $Y$} is defined as the degree of $X$. 
For a real number $x$, $\lceil x \rceil$ 
denotes the smallest integer with $\lceil x \rceil \geq x$. 
We set 
\begin{equation}
\label{eqn:D:def}
D := 
\begin{cases}
\max
\left\{
\left\lceil
\frac{3d^2 - 9 d + 4}{2d}
\right\rceil
,
1
\right\}
& 
\parbox{8cm}{%
\flushleft %
if 
$N \leq 2$, or
if
$N \geq 3$ and
$Y$ is contained in
some affine plane in $\Aff^N$,}
\\
\max \{ d-2 , 1 \}
& \text{otherwise.}
\end{cases}
\end{equation}
Then we have the following;
see also Remark~\ref{remark:morethanD}. 

\begin{Theorem}
\label{thm:limit:trop:d:intro}
Let $Y \subset\Aff^N$ be a connected smooth affine curve of degree $d \geq 1$.
Assume that $Y$ has smooth compactification in 
$\PP^N$. Let $D$ be as in \eqref{eqn:D:def}. 
Then the map $\mathrm{LTrop}_{\leq D}$ in
\eqref{eqn:limit:tropical:map:d:intro} is a homeomorphism. 
\end{Theorem}

As sample examples, 
let $Y_3$ be a connected affine plane curve of degree $3$ that has smooth compactification in 
$\PP^2$ (an affine part of an elliptic curve). 
Then $D = 1$, 
and Theorem~\ref{thm:limit:trop:d:intro} says that  
$Y_3^{\an}$ is the limit of 
the tropicalizations of 
embeddings of degree at most $1$,
namely, the limit of 
{\em linear} tropicalizations. 
Next, 
let $Y_4$ be a connected affine plane curve of degree $4$ that has smooth compactification in 
$\PP^2$
(an affine part of a connected smooth non-hyperelliptic curve of genus $3$). 
Then $D = 2$, and Theorem~\ref{thm:limit:trop:d:intro} says that  
$Y_4^{\an}$ is the limit of tropicalization of 
embeddings of degree $2$,
namely, the limit of
{\em quadratic} tropicalizations.
We note that in general, 
$Y_4^{\an}$ is not the limit of linear tropicalizations.
For details, see Examples~\ref{eg:limit:trop:linear}, \ref{eg:limit:trop:linear:impossible}, and \ref{eg:limit:trop:quad}. 

\medskip
Let us explain the ideas of our proof of 
Theorem~\ref{thm:main:unimodular:faithful}. 
For simplicity, assume here that $X$ has genus $g \geq 2$ 
and the given skeleton is the minimal
skeleton $\Gamma_{\min}$,
which is
associated to a Deligne--Mumford strictly semistable
curve with generic fiber $X$.
It
is well-defined for $X$ and does not depend
on the choice of such a Deligne--Mumford semistable model.
It is known that $\Gamma_{\min}$ has a canonical structure of finite graph
such that the set $V$ of vertices corresponds to the set of the
irreducible components
of the special fiber
of the stable model of $X$.
With this finite graph structure, the set 
$E$
of edges
corresponds to the set of nodes of the special fiber of the stable model.

To realize a faithful tropicalization of $\Gamma_{\min}$ by $L$,
we construct global sections of $L$ according to the following
aims:
for a given $e \in E$, to 
give a unimodular tropicalization of $e$, where
a map is said to be unimodular if it is a piecewise 
linear map preserving the integral structures,
and to 
separate
points in the relative interior $\mathrm{relin} (e)$ of $e$;
for given distinct $e , f \in E$,
to 
separate two points $x$ and $y$ with
$x \in \mathrm{relin} (e)$ and $y \in f$;
for given two distinct points of $V$,
to 
separate the two points.
Then the collection of those global sections will give 
a faithful tropicalization of $\Gamma_{\min}$.

Global sections as above
are constructed as global sections of models.
Here, a model $(\mathscr{X} , \mathscr{L})$ means a pair
consisting of a semistable curve over the ring of integers $R$ of $K$ equipped with an isomorphism $\Xscr\times_{\Spec(R)}\Spec(K) \cong X$
and a line bundle $\mathscr{L}$ over $\mathscr{X}$ with 
$\Lscr\otimes_R K \cong L$.
For each aim as above,
we 
will 
define ``suitable'' global sections of a model. 

Now the issue 
is divided
into two parts:
the first part is to construct a ``good'' model;
the second part is to construct 
``suitable'' global sections of the model.

In the first part, we use the theory of $\Lambda$-divisors on $\Gamma_{\min}$, 
which has been recently developed (see, for example, \cite{AC, BN, GK, MZ}),
together with Raynaud's type theorem
on the surjectivity of the principal divisors under the specialization
map
(see, for example, \cite{AC,  B, BR, Ra}). 
We take a divisor $\tilde{D}$ on $X$ such that $L \otimes \omega_X^{\otimes -1} \cong \OO_X(\tilde{D})$. 
The retraction map $X^{\an} \to \Gamma_{\min}$ induces
the specialization map 
$
\tau_*: \Div(X) \to \Div_\Lambda(\Gamma_{\min}) 
$
from the divisor group on $X$ to the $\Lambda$-divisor group on $\Gamma_{\min}$. 
Using Riemann's inequality, the notion of reduced divisors on metric graphs, and so forth, we construct,
according to each aim, a $\Lambda$-divisor 
$D$ on $\Gamma_{\min}$ that are linearly equivalent to $\tau_*(\tilde{D})$. 
Then using a Raynaud's type theorem,
we construct
a ``good'' Deligne--Mumford model $(\Xscr, \Lscr)$ of $(X, L)$ such that 
$\sum_{C} \deg
\left(\rest{\Lscr \otimes \omega_{\mathscr{X}/R}^{\otimes -1}}{C}
\right) [C] = D$, where $C$ runs through all the irreducible component of the special fiber 
$\Xscr_s$ and $[C] \in \Gamma_{\min}$ is the point in the Berkovich space 
$X^{\an}$ associated to $C$,
called the Shilov point.  

In the second part,
we do algebraic geometry on semistable curves.
We construct (1) global sections of the restriction
of (a modification of) $\Lscr$
to the special fiber $\mathscr{X}_s$, 
and then (2) we lift those sections to global sections
of $\Lscr$.
In (1),
we use some basepoint freeness over nodal curves,
and in (2),
we use base-change theorem.
In those arguments, the key
is vanishing of cohomologies of
line bundles over the special fiber.
By the choice of the above $D \in \Div_{\Lambda} ( \Gamma_{\min} )$,
we can show desired vanishing of cohomologies.
Thus we find desired global sections and
obtain a faithful tropicalization of the minimal skeleton.

Let us explain here why $3g-1$ appears when $g \geq 2$,
even in a very rough form. 
In the second part above,
to show vanishing of cohomologies,
we need
$\Lscr \otimes \omega_{\Xscr/R}^{\otimes -1}$
to be nef over the special fiber, 
which requires 
$\deg (L) \geq 2g-2$, at least.
Further, depending on our each aim, 
we often want $\Lscr \otimes \omega_{\Xscr/R}^{\otimes -1}$ to have positive degree
over a specified  irreducible component $C$.  To obtain such positivity,
we want a $\Lambda$-divisor $D$ in the first part above to be effective and positive at $[C] \in \Gamma_{\min}$.
To find such a $D$,
we use Riemann's inequality 
($|E| \neq \emptyset$ if $\deg(E) \geq g$) 
on $\Gamma_{\min} \subset X^{\an}$,
and this requires
$\deg \left( \tau_{\ast} ( \tilde{D}) - [v] \right) \geq g$
and hence $\deg \left( \tilde{D} \right) \geq g+1$.
Since
$\deg
\left( \tau_{\ast} ( \tilde{D}) 
\right) = \deg \left( L \otimes \omega_{X}^{\otimes -1} \right)$,
this requires
that
$\deg (L)$ should be at least
$g + 1+ (2g-2) = 3g -1$. 

The above description how $3g-1$ appears
indicates that
by our 
approach, it is difficult to reduce $3g-1$ further.
Indeed, when $g=2$, $t(2)$ is optimal by Theorem~\ref{thm:best:possible}. 
However,
a different approach may give a better lower bound for a faithful tropicalization
and may contribute to Question~\ref{Question:optimalbound}.

\begin{Remark}
As explained as above, for the proof of Theorem~\ref{thm:main:unimodular:faithful}, we show the existence of ``good'' divisors on a metric graph, and then we lift them to divisors on a curve under the specialization map. Albeit from a different perspective, this type of lifting argument is sometimes quite successful in tropical geometry. For example, Cools, Draisma, Payne and  Robeva \cite {CDPR} gave a tropical proof of the Brill-Noether theorem. See also \cite{CLM, CJP, KY0, KY1}, for example.   
\end{Remark}

\begin{Remark}
In \cite{KY}, we consider
faithful tropicalizations associated to line bundles
for projective varieties of any dimension.
We recall what the main result of \cite{KY} asserts for curves:
Let $X$ be a connected smooth projective curve over $K$ of genus $g$, 
and let $L$ be an ample line bundle over $X$;
suppose that there exist a discrete valuation ring $R_0$ dominated by $R$,
a strictly semistable regular scheme $\mathscr{X} \to \Spec (R_0)$ 
over $R_0$ with
$\mathscr{X} \otimes_{R_0} K = X$, 
a \emph{relatively ample} line bundle $\mathscr{N}$
over $\mathscr{X}$,
and an integer $m \geq 2$ 
such that $\rest{\mathscr{N}^{\otimes m} \otimes \omega_{\mathscr{X}/R_0}}{X} = L$; 
then $L$ admits a faithful tropicalization of the skeleton associated to $\Xscr$.

For curves, Theorem~\ref{thm:main:unimodular:faithful} in this paper has an advantage over \cite{KY}. 
For example, 
the result in \cite{KY} cannot contribute to Question~\ref{question:first},
even if we work under the condition that $X$ can be defined over a discrete valued field and we consider only 
a faithful tropicalization of the minimal skeleton $\Gamma_{\min}$ of $X^{\an}$; 
see \S\ref{subsection:comparison} for details.

In \cite{KY}, we work with a given model $(\mathscr{X} , \mathscr{N})$.
On the other hand,
to establish Theorem~\ref{thm:main:unimodular:faithful}
in this paper,
we need much finer analysis for constructing  
various carefully chosen models $(\Xscr, \Lscr)$ of $(X, L)$ as we do in this paper. That is one of the main reasons why this paper has much more pages
than \cite{KY}.
\end{Remark}

\begin{Remark}
Theorem~\ref{thm:limit:trop:d:intro} does not follow automatically from 
Theorem~\ref{thm:main:unimodular:faithful}, because we need to separates type IV points of $Y^{\an}$ (cf. Remark~\ref{rmk:four:types})
and they
are outside skeleta of $Y^{\an}$. To separate a type IV point
from another point, we will use 
global sections constructed in \S\ref{sec:ft:general:model} 
originally for the proof of Theorem~\ref{thm:main:unimodular:faithful}. 
Other ingredients of the proof of Theorem~\ref{thm:limit:trop:d:intro} 
are Castelnuovo's bound for a curve in projective space and 
a criterion for a complete linear system
to be traced out by hypersurfaces on the projective space. See \S\ref{sec:limit:trop} for details. 
\end{Remark}

\medskip
The organization of this paper is as follows. 
Section~\ref{sec:preliminaries} is preliminary. 
We recall semistable models and semistable pairs, some known facts on Berkovich spaces and tropical geometry. In particular, we recall a minimal skeleton $\Gamma_{\min}$ for a connected smooth projective curve. 
We also give the definition of a faithful tropicalization. 
In \S\S\ref{sec:models}-\ref{section:FTcan}, we generally assume that $g \geq 2$. 
In \S\ref{sec:models}, we introduce the notion of a good model. 
Further, for a given effective divisor on the minimal skeleton
of $X^{\an}$, we construct a good model with a line bundle 
which is associated to a divisor that  is a ``lifting'' of the divisor on the skeleton.
In \S\ref{sec:usuful:lemmas}, 
we show several lemmas that are useful for constructing 
required sections. Then we show that the minimal skeleton 
$\Gamma_{\min}$ can be embedded into $\TT\PP^N$ 
piecewisely isometrically (viz. a unimodular tropicalization). 
In \S\ref{section:FTcan},
we construct sections that separate points in $\Gamma_{\min}$, and establish a 
faithful tropicalization for $\Gamma_{\min}$. 
In \S\ref{section:FT:lowgenus}, we discuss the tropicalization
when $g = 0$ or $g=1$. 
In \S\ref{sec:ft:general:model}, 
we establish a faithful tropicalization for an arbitrary skeleton 
possibly 
with ends, 
thus completing the proof of Theorem~\ref{thm:main:unimodular:faithful}. 
In \S\ref{sec:complementary}, we give proofs of complementary results such as 
Theorem~\ref{thm:best:possible} and Proposition~\ref{prop:countereg:very:ample}. 
We also compare Theorem~\ref{thm:main:unimodular:faithful} with \cite[Theorem~1.1]{KY}. 
In \S\ref{sec:limit:trop}, we give a proof of Theorem~\ref{thm:limit:trop:d:intro}. 

\medskip
{\sl Acknowledgments.}\quad
 Originally, we essentially worked only 
over the completion of the algebraic closure of a discretely valued field, 
and the bound of the degree was $5g-4$. 
While we were making gradual improvement of 
our results, we had opportunities to give talks in several conferences and have fruitful discussions.
We especially thank Professor Antoine Ducros for helpful discussions, which lead us to work over any algebraically closed field that is complete with respect to a non-trivial nonarchimedean absolute value. 
Our deep thanks go to Professor Sam Payne for helpful communications,
which has lead us to Proposition~\ref{prop:countereg:very:ample} and Theorem~\ref{thm:limit:trop:d:intro}.
Our deep thanks go to Professor Matthew Baker
 for helpful comments, which improve the exposition. 

\subsection*{Notation and Conventions}
By the convention of this paper,
an inclusion $A \subset B$ of sets allows the case $A =B$. 

A {\em variety} over a field $k$ means a reduced separated scheme 
of finite type over $k$. (We allow a variety to be reducible
unless otherwise specified.) A \emph{curve} is a variety of pure dimension $1$.
For a curve $Y$ over $k$, let $\Irr (Y)$ denote the set of 
irreducible components of $Y$
and let $\Sing (Y)$ denote the set of non-smooth points of 
$Y \to \Spec (k)$.

Let $Y$ be a curve over $k$.
A {\em curve in $Y$}  is a closed subscheme of $Y$ that
is a curve over $k$ in the above sense.
If there is no danger of confusion,
we identify a curve in $Y$ with a 
union of irreducible components of $Y$. 
A \emph{connected curve in $Y$}
is a curve in $Y$ that is connected.
For a curve $D$ in $Y$, 
we denote by $Y - D$  the union of irreducible components of $Y$ 
other than those in $D$, i.e.,
$Y - D := \bigcup_{C \in \Irr (Y) \setminus \Irr (D)} C$. 
By convention, if $D = Y$, then $Y-D = \emptyset$. 
Notice that $Y - D$ is not equal to the (set-theoretic) complement $Y \setminus D$ 
of $D$ unless $D$ is a union of connected components of $Y$.
Let $D_1 , \ldots , D_r$ be curves in $Y$.
Assuming that 
$\Irr (D_i) \cap \Irr (D_j) = \emptyset$ for $i \neq j$, 
we denote by $D_1 + \cdots + D_r = \sum_{i=1}^r D_i$
the reduced curve $\bigcup_{i = 1}^r D_i$
in $Y$. 
The notation $\sum_{i=1}^r D_i$ is used only under this assumption.

Throughout this paper, 
$K$ denotes an algebraically closed field 
that is complete with respect to a non-trivial nonarchimedean absolute value $|\cdot|_K$. 
We denote by $R$ the ring of integers of $K$, by $\mfrak$ the maximal ideal of $R$, and by $k$ the residue field $R/\mfrak$. 
Let $v_K : K^{\times} \to \RR$
be the valuation map given by $v_K := - \log | \cdot |_K$, 
and let $\Lambda := \{v_K(x) \mid x \in K^\times\} \subset \RR$ be the value group of $K$. 

We use notions on graphs.
A {\em finite graph} $G$  means a finite connected graph (that is allowed to have loops and multiple edges), and we denote by  $V(G)$ and $E(G)$  the set of vertices and the set of edges, respectively. 
We assign each edge $e \in E(G)$ length in $\Lambda$. 
Then the graph $G$ assigned with edge lengths gives a $\Lambda$-metric graph $\Gamma$ 
together with the set of vertices $V(G)$ and the set of edges $E(G)$. 
A {\em $\Lambda$-metric graph} $\Gamma$ is a connected metric graph such that there exists 
a finite graph $G$ assigned with edge lengths in $\Lambda$ that induces $\Gamma$ 
as above. A point $x$ in a $\Lambda$-metric graph is a {\em $\Lambda$-valued point} if the distance from 
$x$ to one (hence any) vertex $v$ of $G$ lies in $\Lambda$. 
We denote by $\Gamma_{\Lambda}$ the set of $\Lambda$-valued points of $\Gamma$. 
Notice that the notion of $\Lambda$-valued points is independent of the choice of 
a finite graph $G$ as above.

Let $\Gamma$ be a $\Lambda$-metric graph with a set $V$ of vertices and a set $E$ of edges.
We naturally regard $V$ as a finite subset of $\Gamma$ and each $e \in E$ as a closed
subset of $\Gamma$.
We take any $e \in E$
and set $\partial e := e \cap V$.
We call a point in $\partial e$ an \emph{end vertex} of $e$.
If $e$ is isometric to a closed interval, then $\# \partial e = 2$;
if $e$ is a loop, then $\# \partial e = 1$.
We define the \emph{relative interior of $e$} to be
$
\mathrm{relin}(e) := e \setminus \partial e
$.

\setcounter{equation}{0}
\section{Preliminaries}
\label{sec:preliminaries}
In this section, after recalling semistable models and semistable pairs, we briefly review some facts on Berkovich spaces and tropical geometry, 
which we use later. Then we define a unimodular and a faithful tropicalization associated to a linear system. Our basic references for Berkovich spaces and skeleta are
\cite{BPR2}, \cite[\S5]{gubler4}, and \cite[\S3]{GRW}. 
We use the language of schemes (rather than formal schemes), which might
be more familiar to readers in algebraic geometry, even though actual proofs use formal geometry.

\subsection{Semistable models and semistable pairs}
\label{subsec:semistable:models:pairs}
This subsection is mainly to fix the notation. We recall semistable curves, 
semistable models, and semistable pairs. 

\subsubsection*{Semistable curves}
We first recall the notion of semistable curves over schemes. 
Let $S$ be a scheme. A {\em semistable curve} over $S$ of genus $g \geq 0$ 
is a proper flat morphism $\pi: \Xscr \to S$ whose geometric fibers 
$\Xscr_{\bar{s}}$ are reduced, separated, connected, $1$-dimensional schemes
that has only ordinary double points (called nodes) as singularities
and satisfies $h^1(\OO_{\Xscr_{\bar{s}}}) = g$. 
A semistable curve $\Xscr$ over a scheme $S$ 
is said to be {\em strictly semistable} if for any scheme point $s$ of $S$, any 
singular point of $\Xscr_s$ is defined over the residue field $\kappa(s)$ and any irreducible component of $\Xscr_s$ is geometrically irreducible and smooth. 

Let $\Xscr$ be a strictly semistable curve over a scheme $S$
of genus $g$. 
For $p \in \Sing (\mathscr{X}_s)$,
we say that $p$ is of \emph{connected type}
if the partial normalization
$\widetilde{\mathscr{X}_s}$ of $\mathscr{X}_s$ at $p$ is connected;
we say that $p$ is of \emph{disconnected type} 
otherwise.

A line bundle $\Lscr$ over $\Xscr$ is {\em vertically nef} if $\rest{\Lscr}{\Xscr_{\bar{s}}}$ is nef for 
all geometric fibers $\Xscr_{\bar{s}}$. 
A semistable curve is called a {\em Deligne--Mumford semistable curve} if 
the dualizing sheaf $\omega_{\Xscr/S}$ is vertically nef. A Deligne--Mumford semistable curve is called a Deligne--Mumford {\em stable} curve (or simply a {\em stable} curve) if  $\omega_{\Xscr/S}$ is ample. 

Let $\Xscr$ be a semistable curve over $S$. 
An irreducible component $C$ of $\Xscr_{\bar{s}}$ is called a {\em $(-1)$-curve}  
if $C \cong \PP^1$ and $\#(C \cap (\Xscr_{\bar{s}}-C)) = 1$, and is  
called a {\em $(-2)$-curve}  
if $C \cong \PP^1$ and $\#(C \cap (\Xscr_{\bar{s}}-C)) = 2$. 
By the adjunction formula, an irreducible component 
$C$  of $\Xscr_{\bar{s}}$ is a $(-1)$-curve if and only if $\deg(\rest{\omega_{\Xscr/S}}{C}) = -1$, 
and $C$ is a $(-2)$-curve if and only if 
$\deg(\rest{\omega_{\Xscr/S}}{C}) = 0$.
Suppose that $\Xscr$ is a Deligne--Mumford semistable curve
of genus $g \geq 1$.
Let $E$ be a connected curve in $\Xscr_{\bar{s}}$.
Then $E$ consists of $(-2)$-curves if and only
if $\deg(\rest{\omega_{\Xscr/S}}{E}) = 0$.
If this is the case, then $E$ is a cycle of $(-2)$-curves or a chain of 
$(-2)$-curves. Furthermore, if $g \geq 2$, then $\mathscr{X}_{\bar{s}}$ 
does not contain a cycle of $(-2)$-curves, and hence $E$ is a chain.

A chain consisting of $(-2)$-curves are called a \emph{$(-2)$-chain}.
Let $E$ be a $(-2)$-chain in $\Xscr_{\bar{s}}$.
We say that 
$E$ is {\em of connected type} if $\Xscr_{\bar{s}} - E$ is connected,
and 
$E$ is {\em of disconnected type} if $\Xscr_{\bar{s}} - E$ is 
disconnected. 
Note that if  
$\mathscr{X}$ is a Deligne--Mumford semistable curve and 
$E$ is of disconnected type,
then $\Xscr_{\bar{s}} - E$ consists of two connected components 
each of which has positive arithmetic genus. 

A  $(-2)$-chain $E$ is called a {\em maximal} $(-2)$-chain if $F$ is a $(-2)$-chain with 
$E \subset F$, then $E = F$. 
When $g \geq 2$, any connected curve in $\mathscr{X}_{\bar{s}}$
consisting of $(-2)$-curves is a chain,
and hence any two maximal $(-2)$-chains are disjoint.

\subsubsection*{Semistable models}
\label{subsec:semistabe:models}
Let $X$ be a variety over $K$.
A \emph{model of $X$ over $R$} is a flat morphism $
\pi : \mathscr{X} \to \Spec (R)$ 
of finite type equipped with an identification of the generic fiber of $\pi$
with $X$. When $X$ is proper over $K$, we require that 
a model be also proper over $R$. The special fiber of $\Xscr$ is denoted by 
$\Xscr_s$. 

For two models $\mathscr{X}'$ and $\mathscr{X}$ of $X$ over $R$, 
we say that \emph{$\mathscr{X}'$ dominates $\mathscr{X}$} 
if the identity morphism on $X$ extends to a morphism 
$\mathscr{X}' \to \mathscr{X}$. 

Now, let $X$ be a connected smooth projective curve over $K$
of genus $g$. 
A {\em semistable model} of $X$ over $R$ is a model 
$\pi : \mathscr{X} \to \Spec (R)$ that 
is a semistable curve. 
If a semistable model $\Xscr$ is strictly semistable, 
then $\Xscr$ is called a {\em strictly} semistable model. 
If $\Xscr$ is Deligne--Mumford semistable (resp. stable), 
then $\Xscr$ is called a {\em Deligne--Mumford semistable} (resp. {\em stable}) model of $X$. 
By \cite[Theorem~7.1]{BL1} or  \cite[Corollary~4.23 and Remark~4.24]{BPR2},
there exists a Deligne--Mumford semistable model of $X$ if $g \geq 1$,
and there exists
a unique stable model of $X$
if $g \geq 2$.
When $g \geq 2$ and we denote a model of $X$ by
$\mathscr{X}^{\st} \to \Spec (R)$, this stands for the stable model.
It is known that the stable model is ``minimal'' in the sense that 
if $\mathscr{X}$ is any semistable model of $X$, 
then $\mathscr{X}$ dominates $\mathscr{X}^{\st}$.

Let $L$ be a line bundle over $X$. Let $\mathscr{X}$ be a model of $X$ over $R$.  
We call a line bundle $\mathscr{L}$ over $\mathscr{X}$ 
a {\em model} of $L$ if $\rest{\mathscr{L}}{X} = L$. In this case, we 
call 
the pair $(\mathscr{X} , \mathscr{L})$ a \emph{model} of $(X, L)$. 
If $\Xscr$ is semistable (resp. strictly semistable, Deligne--Mumford semistable, 
stable), we say that $(\mathscr{X} , \mathscr{L})$ is a 
{\em semistable}  (resp. {\em strictly semistable}, {\em Deligne--Mumford semistable}, 
{\em stable}) model of $(X, L)$. 
  
\subsubsection*{Strictly semistable pairs}
\label{subsec:semistabe:pairs}
Let $X$ be a connected smooth projective curve over $K$. 
Let $\pi: \mathscr{X} \to \Spec (R)$ be a strictly semistable model of $X$.  
For a section $\sigma: \Spec(R) \to \mathscr{X}$ of $\pi$, 
we denote the image of $\sigma$ by $\sigma (R)$,
the image of the special (resp. generic) point of $\Spec (R)$ by $\sigma (k)$
(resp. $\sigma(K)$).

Let $\sigma_1 , \ldots , \sigma_r$ be sections of $\pi$.
The pair $( \mathscr{X} ; \sigma_1 , \ldots , \sigma_r)$
is called a \emph{strictly semistable pair}
if $\sigma_i (R)$ is contained in the smooth locus of $\pi$
and $\sigma_i (k) \neq \sigma_j (k)$ for $i \neq j$.
We note that the definition of a strictly semistable pair is the same as in \cite{GRW}.
Indeed, since $\sigma_i(R)$ is contained in the smooth locus of $\pi$, $\sigma_1 (R) , \ldots , \sigma_r (R)$ are prime Cartier divisors on $\mathscr{X}$, and by the condition $\sigma_i (k) \neq \sigma_j (k)$ for $i \neq j$, these Cartier divisors have disjoint supports.

\subsection{Berkovich spaces}
\label{subsec:Berko:skeleta}
Let $X$ be a variety over $K$. For each $x \in X$, let $\kappa(x)$ denote the residue field. 
Let $X^{\an}$ be the {\em Berkovich (analytic) space} associated to $X$. (For details, 
see fundamental papers  by Berkovich \cite{Be1, Be2, Be3, Be4}.)  Here we recall 
$X^{\an}$ as a set and the underlying topology of $X^{\an}$. 
As a set, $X^{\an}$ is defined by 
\[
  X^{\an} := 
  \left\{
  (x, |\cdot|) \mid \text{$x \in X$ and 
$|\cdot|$ is an absolute value of $\kappa(x)$ extending $|\cdot|_{K}$}
  \right\}. 
\]
The  \emph{Berkovich topology} on $X^{\an}$ is 
the weakest topology 
such that for any Zariski open subset $U$ of $X$ 
and for any regular function $g \in \OO_X(U)$, 
the map $\iota^{-1}(U) \to \RR, (x, |\cdot|) \mapsto |g(x)|$ is continuous, 
where $\iota: X^{\an} \to X$ is the natural map given by 
$(x, |\cdot|) \mapsto x$. 

The set $X(K)$ of $K$-valued points are naturally regarded as points of $X^{\an}$ via 
$x \mapsto (x, |\cdot|_K)$ for $x \in X(K)$,
where $| \cdot |_{K}$ is the original absolute value on $K = \kappa (x)$. 
Those points are called \emph{classical points}. 
It is known that $X(K)$ is dense in $X^\an$.

There is another class of points of $X^{\an}$, called the \emph{Shilov points}.
Let $\mathscr{X} \to \Spec (R)$ be a 
model
of $X$.
Assume that any irreducible component 
of the special fiber $\Xscr_s$
is reduced.
Let $C$ be an irreducible component of $\mathscr{X}_s$
and let $\xi$ be the generic point of $C$.
Then one associates a unique absolute value $|\cdot|_{C}$
on the function field $K ( X )$ of $X$, which is characterized
by the property that for any $f \in K ( X )$,
\[
| f |_{C}
=
\inf
\{
|a| \in \RR
\mid
a \in K^{\times} , 
a^{-1} f \in \OO_{\mathscr{X},\xi}
\}
.
\]
Thus
$C$ gives a point $[C]$ of $X^{\an}$,
which we call the \emph{Shilov point} associated to
$(\mathscr{X} , C)$
or simply to $C$.
A Shilov point associated to $(\mathscr{X} , C)$
for some irreducible component $C$ of $\mathscr{X}_s$
is called a \emph{Shilov point with respect to $\mathscr{X}$}.
In this paper,
let $V( \mathscr{X} )$ denote the set of Shilov
points with respect to $\mathscr{X}$.
If 
$\mathscr{X}^\prime \to \Spec (R)$ is another model of $X$ with reduced
special fiber 
and with a surjective morphism $\mathscr{X}^\prime \to \mathscr{X}$,
then we have 
$V( \mathscr{X} ) \subset V( \mathscr{X}^\prime )$. 
A point which is a Shilov point with respect to some model
is called a \emph{Shilov point}.

\begin{Remark} \label{rem:type2points}
Let $x \in X^{\an}$ be a Shilov point.
Then for any $f \in K(X)^{\times}$,
we have $- \log |f ( x )| \in \Lambda$.
One sees this fact from 
the maximum modulus principle 
\cite[Proposition 6.2.1/4]{BGR}
with the description
\cite[\S2.3]{BPR2}
or from 
the description of $\mathrm{Val}$ in \cite[\S5.3]{gubler4}.
\end{Remark}

\subsubsection*{Reduction map}
We recall the reduction map associated to a model $\Xscr$ of $X$.
Here, we assume that $X$ is proper over $K$.
Let $(x, |\cdot|)$ be a point of $X^\an$. 
Then we have a natural morphism
$\Spec ( \kappa (x) ) \to \mathscr{X}$.
Let $V$ be the ring of integers of $\kappa (x)$ with respect to $| \cdot |$.
Since $V$ is a valuation ring and $\mathscr{X} \to \Spec (R)$ is proper,
the morphism $\Spec ( \kappa (x) ) \to \mathscr{X}$
extends to a unique morphism
$\Spec (V) \to \mathscr{X}$.
Let $\red_{\mathscr{X}} ( x, |\cdot| ) \in \mathscr{X}_s$ 
be the image of the special point of $\Spec (V)$ 
by this morphism.
Then the assignment $( x, |\cdot| ) \to \red_{\mathscr{X}} ( x, |\cdot| )$ defines
a map 
$\red_{\mathscr{X}} : X^{\an} \to \mathscr{X}_s$.
We call this map
the \emph{reduction map associated to $\mathscr{X}$}.

If $x = (x, |\cdot|_K) \in X^{\an}$ is a classical point
and if $\sigma$ is the corresponding section of the
model $\mathscr{X} \to \Spec (R)$ guaranteed by the valuative criterion, 
then $\red_{\mathscr{X}} (x) = \sigma (k)$.

\subsection{Skeleta associated to strictly semistable models}
\label{subsection:skeleta}
We recall the notion of skeleta.  
For simplicity, in the sequel, we  assume 
that $X$ is a connected smooth projective curve,
although skeleta can be defined in a more general setting.

\subsubsection*{Standard model and its skeleton}
Let $\varpi \in R$ be an element with $0 < |\varpi| < 1$. 
We set $
\mathscr{S} :=
\Spec ( R 
[x,y]/(xy - \varpi))
$, and call $\Sscr$ a  \emph{standard model}.

Let $\mathscr{S}_K = \Spec(K[x, y]/(xy-\varpi))$ be the generic fiber of $\mathscr{S} \to \Spec (R)$.
We first define the skeleton
$S ( \mathscr{S}) \subset \mathscr{S}_K^{\an}$.
Since
$K [x,y]/(x y - \varpi)
\cong K [y^{\pm}]$,
we regard an element of $K [x,y]/(x y - \varpi)$ as a Laurent 
polynomial on $y$.
For any 
$v \in [ 0 , - \log |\varpi| ]$,
one constructs a unique absolute value $|\cdot|_{v}$ on 
the function field of $K [x,y]/(x y - \varpi)$
such that
for any Laurent polynomial 
$f= \sum_{m} a_{m} y^m 
$, we have 
$
|f|_{v} =
\max_{m}
\left\{ 
|a_{m}|
\exp(- vm
)
\right\}
.
$

Let $\eta$ be the generic point of $\mathscr{S}$. 
The assignment $v \to (\eta, |\cdot|_{v})$ gives a map
$[ 0 , - \log |\varpi| ] \to \mathscr{S}_K^{\an}$. 
It is known that this map is continuous and injective. 
The {\em skeleton} $S(\Sscr)$ of $\Sscr$ is defined to be 
the image of this map. Then the map 
$[ 0 , - \log |\varpi| ] \to  S(\Sscr)$ is a homeomorphism 
with inverse 
\addtocounter{Claim}{1}
\begin{align} \label{eq:inverseisometry}
(\eta, | \cdot |) \mapsto - \log |y|
.
\end{align}
Thus $S ( \mathscr{S} )$
has a structure of $1$-simplex via the homeomorphism
$[ 0 , - \log |\varpi| ]  \cong S ( \mathscr{S} )$.
One sees that
the end points of $S ( \mathscr{S} )$
are the Shilov points corresponding to the irreducible components
of $\mathscr{S}_s$.

\subsubsection*{Local \'etale atlas that distinguishes a node}
Let $\mathscr{X} \to \Spec (R)$ be a
strictly semistable model of $X$.
We take any $p \in \Sing ( \mathscr{X}_s )$. 
Then there exist $C, C' \in
\Irr ( \mathscr{X}_s)$ such that $C \neq C'$ 
and $p \in C \cap C'$.
By the same argument as in the proof of 
\cite[Proposition~4.3]{BPR2},
one can take an open neighborhood of $p$ in $\Xscr$
and an \'etale morphism from this open subset to a standard model.
Furthermore,
replacing the open neighborhood by a smaller one,
we can take an affine open neighborhood $\mathscr{U}_p
\subset \mathscr{X}$ of $p$,
an element $\varpi_p \in R$ with $0 < |\varpi_p| < 1$,
and an \'etale morphism
$\psi_p : \mathscr{U}_p \to \mathscr{S}:=\Spec ( R [x,y]/(x y - \varpi_p))$
such that
the special fiber
$\mathscr{U}_p \cap \Xscr_s$ consists of exactly two irreducible components
and $p$ is a unique node of $\mathscr{U}_p \cap \Xscr_s$
(cf. \cite[Proposition~5.2]{gubler4}).
We call $\psi_p$ a \emph{local \'etale atlas that distinguishes $p$}.

Note that, with the above notation, 
$v_K(\varpi_p)$ depends only on $p \in \mathscr{X}$
and is
independent of the choice of $\mathscr{U}_p$. 
Thus the following definition makes sense.

\begin{Definition}[multiplicity] \label{def:multiplicityatanodes}
We call $v_K (\varpi_p)$ the \emph{multiplicity }
(of $\mathscr{X}$) at $p \in \Sing ( \mathscr{X}_s)$. 
\end{Definition}

\subsubsection*{Skeleton $S ( \mathscr{X} )$ associated to $\Xscr$} 
Let $\mathscr{X} \to \Spec (R)$ be a
strictly semistable model of $X$, 
and let $p \in \Sing(\Xscr_s)$.

\begin{Definition}[canonical $1$-simplex $\Delta_p$ corresponding to a node $p$]
\label{def:canonical:1:simplex:node}
The \emph{canonical $1$-simplex corresponding to $p$}
is a closed subset $\Delta_p \subset X^{\an}$ 
characterized by the
following conditions.
\begin{enumerate}
\renewcommand{\labelenumi}{(\alph{enumi})}
\item
For any local \'etale atlas
$\psi_p : \mathscr{U}_p \to \mathscr{S}:=\Spec ( R [x,y]/(x y - \varpi_p))$
that distinguishes $p$,
we have $\Delta_p \subset U_p^{\an}$,
where $U_p$ is the generic fiber of $\mathscr{U}_p \to \Spec (R)$.
\item
The restriction to $\Delta_p$ of the induced map $\psi_p^{\an} : U^{\an} \to \mathscr{S}_K^{\an}$ gives a homeomorphism $\Delta_p \to S ( \mathscr{S} )$.
\end{enumerate}
\end{Definition}

Note that
the end points $\Delta_p$ are the Shilov points
which correspond to the irreducible components of 
the special fiber of $\mathscr{U}$.

\begin{Definition}[skeleton associated to a model]
\label{def:skeleton:model}
Let $\Xscr$ be a strictly semistable model of $X$. 
The {\em skeleton} $S(\Xscr)$ associated to $\Xscr$ is defined by 
\[
  S(\Xscr) = \bigcup_{p \in 
\Sing ( \mathscr{X}_s )} \Delta_{p}, 
\]
which is a compact subset of $X^{\an}$. 
\end{Definition}

One sees that for any $p_1 , p_2 \in \Sing ( \mathscr{X}_s )$,
$\Delta_{p_1} = \Delta_{p_2}$ if and only if $p_1 = p_2$;
if
$p_1 \neq p_2$,
then
$\mathrm{relin} ( \Delta_{p_1}) \cap \mathrm{relin} ( \Delta_{p_2}) = \emptyset$, 
where $\mathrm{relin}$ means the relative interior, 
and
the
points in $\Delta_{p_1} \cap \Delta_{p_2}$ 
are the Shilov points
associated to
irreducible components of $\mathscr{X}_s$
that contain
$p_1$ and $p_2$.
Thus
$S ( \mathscr{X} )$ has a canonical structure of simplicial set
and is regarded as the dual graph of $\mathscr{X}_s$. 
Note that the set $V (\mathscr{X}) =
\{ [C] \mid C \in \Irr (\mathscr{X}_a )\  \}$ 
of Shilov points with respect to 
$\mathscr{X}$ is actually a set of vertices of this simplicial set.
Each canonical $1$-simplex is called an edge,
and we denote by $E ( \mathscr{X} )$ the set of edges.

\medskip
It is known that $X^{\an}\setminus X(K)$ has a canonical metric structure (cf. \cite[Corollary~5.7]{BPR2}). We describe this canonical metric structure on 
$S(\Xscr)$. We take any $p \in \Sing ( \mathscr{X}_s )$ and denote by 
$\lambda_p$ the multiplicity at $p$
(cf. Definition~\ref{def:multiplicityatanodes}). 
The composition of the above homeomorphisms 
$\Delta_p \to S(\Sscr)$ and $S(\Sscr) \to [0, \lambda_p]$ (see \eqref{eq:inverseisometry}) gives a 
homeomorphism $\Delta_p \to [0, \lambda_p]$. 
Then, with respect to the canonical metric structure on $\Delta_p$ 
induced from $X^{\an}\setminus X(K)$, 
this homeomorphism is an isometry (\cite[Theorem~5.15]{BPR2} and \cite[\S4]{GRW}). 
In particular, 
if $C$ and $C'$ are the distinct irreducible components with $p
\in C \cap C'$
and if $v$ and $v'$ are the Shilov points of $C$ and $C'$, 
then the distance between $v$ and $v'$ 
in $\Delta_p$
equals the multiplicity at $p$.

Thus the skeleton $S ( \mathscr{X} )$ has a structure of a $\Lambda$-metric
graph.
For a node $p \in \mathscr{X}_s$,
a point of $\Delta_p$ that lies in 
$[0 , \lambda_p] \cap \Lambda$ 
via the isometry
$\Delta_p \to [0 , \lambda_p]$
is a $\Lambda$-rational point of $S ( \mathscr{X} )$.
We denote the set of $\Lambda$-rational points by
$S_{\Lambda} (\mathscr{X} )$.
When we write the skeleton as $\Gamma$, we 
write the subset of those points as $\Gamma_\Lambda$,
which is compatible with ``Notation and convention'' in the introduction. 

We end this subsection by defining a canonical $1$-simplex 
corresponding to a connected curve.

\begin{Definition}[canonical $1$-simplex $\Delta_E$ for a connected curve]
\label{def:canonical:1:simplex:curve}
Let $\Xscr$ be a strictly semistable model of $X$. 
Let $E$ be a connected curve in $\mathscr{X}_s$. 
We set 
\[
\Delta_E := \bigcup_{p \in E \cap \Sing ( \mathscr{X}_s)} \Delta_p,
\]
which we call the 
\emph{canonical $1$-simplex corresponding to $E$}.
\end{Definition}

Note that if $E$ is a chain in $\mathscr{X}_s$,
then $\Delta_E$ is isometric to a closed interval
of length $\sum_{p \in E \cap \Sing ( \mathscr{X}_s)} \lambda_p$,
where $\lambda_p$ denotes the multiplicity at $p$.

\subsection{Skeleta associated to strictly semistable pairs}
\label{subsec:skeleta:ass:ssp}
We recall
the notion of skeleta associated to strictly semistable pairs, which generalize 
skeleta associated to strictly semistable models in \S\ref{subsection:skeleta}.  
See also \cite{Tyo}. 

\subsubsection*{Standard pair for ends and its skeleton}
We set $\mathscr{T} := \Spec ( R [y] )$
and denote by $\mathscr{T}_K = \Spec (K [y] )$ the generic fiber 
of $\mathscr{T} \to \Spec(R)$. 
Let $\tau$ be a section of $\mathscr{T} \to \Spec (R)$
given by $R[y] \to R [y]/(y) = R$.
We call
$(\mathscr{T} ; \tau )$ the \emph{standard pair for an end}.

One defines the skeleton $S(\mathscr{T} ; \tau ) \subset \mathscr{T}_K^{\an}$ as follows. 
For any $v \in \RR_{\geq 0}$, one constructs a unique absolute value $|\cdot|_v$ on 
$K (y)$ such that for any $f = \sum_{m} a_m y^m \in K [y]$, 
we have
$
|f|_{v} := \max_{m} 
\{ |a_m| 
\exp (-vm)
\}
.
$

Let $\eta$ be the generic point of $\mathscr{T}$. 
The assignment $v \to (\eta, | \cdot |_v)$ gives a
map $\RR_{\geq 0} \to \mathscr{T}_K^{\an}$. It is known that 
this map is continuous and injective. The {\em skeleton} $S(\mathscr{T} ; \tau )$ of 
$(\mathscr{T} ; \tau )$ is defined to be the image of this map.
The map $\RR_{\geq 0} \to S ( \mathscr{T} ; \tau )$ is a homeomorphism
 with inverse given by $(\eta, | \cdot |) \to - \log |y|$. 

The point $\xi$ in $S ( \mathscr{T} ; \tau )$ 
corresponding to $0 \in \RR_{\geq 0}$
is the Shilov point associated to
the special fiber of $\mathscr{T}$.
We write $\mathrm{relin} (S ( \mathscr{T} ; \tau ))$ for the relative interior
of $S ( \mathscr{T} ; \tau )$.
We remark that $\mathrm{relin} ( S ( \mathscr{T} ; \tau )) =
S ( \mathscr{T} ; \tau )\setminus \{ \xi \}$.

\subsubsection*{Local \'etale atlas that distinguishes a section}
Let $\mathscr{X} \to \Spec (R)$ be a strictly semistable
curve
and let $\sigma$ be a section. 
Assume that
$(\mathscr{X} ; \sigma)$ is a strictly semistable pair.
Then there exist
an affine open neighborhood $\mathscr{U}_{\sigma}$ of 
$\sigma$
and an \'etale morphism
\[
\psi_{\sigma} : \mathscr{U}_{\sigma} \to 
\Spec ( R [ y ] ) = \mathscr{T}
\]
such that
the special fiber of
$\mathscr{U}$ is irreducible and
$\zero ( \psi^{\ast} (y)) = \sigma (R)$.
We call such an \'etale morphism a
\emph{local \'etale atlas that distinguishes $\sigma$}.

\subsubsection*{Skeleton with ends}

Let $\mathscr{X}  \to \Spec(R)$ be a strictly semistable
model of $X$, and 
let $(\mathscr{X} ; \sigma)$ be a strictly semistable pair. 

\begin{Definition}[canonical end $\Delta(\sigma)$ corresponding to a section $\sigma$]
The {\em canonical end} 
$\Delta ( \sigma)$
associated
to $ \sigma $ is a closed subspace of $X^{\an} \setminus X (K)$ 
which is characterized by the following conditions.
\begin{enumerate}
\renewcommand{\labelenumi}{(\alph{enumi})}
\item
For any local \'etale atlas 
$\psi_{\sigma} : \mathscr{U}_{\sigma} \to 
\Spec ( R [ y ] ) = \mathscr{T}$
that distinguishes $\sigma$,
we have $\Delta (  \sigma ) \subset U_{\sigma}^{\an}$,
where $U_{\sigma}$ is the generic fiber of $\mathscr{U}_{\sigma}$.
\item
The restriction to $\Delta ( \sigma )$ of the 
induced map $\psi_{\sigma}^{\an} : U_{\sigma}^{\an} \to \mathscr{T}_K^{an}$
gives a homeomorphism $\Delta ( \sigma ) \to S ( \mathscr{T} ; \tau )$.
\end{enumerate}
\end{Definition}

The canonical metric structure on $X^{\an}\setminus X(K)$ gives a metric
on
$\Delta ( \sigma )$, and the map
$\Delta ( \sigma ) \to \RR_{\geq 0}$
given by $(\eta' , | \cdot |) \mapsto - \log |\psi_p^{\ast} (y)|$
is
 an isometry,
where $\eta'$ is the generic point of $X$ and
$\psi_p^{\ast} (y)$ is regarded as a rational function on $X$.
The endpoint of  $\Delta ( \sigma)$,
i.e, the point corresponding to $0 \in \RR_{\geq 0}$ via the isometry,
is the Shilov point associated to the irreducible component
of $\mathscr{X}_s$ that contains the point $\sigma (k)$.

\begin{Definition}[skeleton associated to a pair]
\label{def:definitionofskeleton}
Let $\mathscr{X}$ be a strictly semistable model of $X$, and 
let $( \mathscr{X} ; \sigma_1 , \ldots , \sigma_r)$ be
a strictly semistable pair.
The \emph{skeleton}  $S (\mathscr{X} ; \sigma_1 , \ldots , \sigma_r)$ associated to
$( \mathscr{X} ; \sigma_1 , \ldots , \sigma_r)$
is defined by
\begin{align}
\label{align:definitionofskeleton}
S (\mathscr{X} ; \sigma_1 , \ldots , \sigma_r)
:= S ( \mathscr{X} ) \cup \bigcup_{i = 1}^{r} \Delta ( \sigma_i).
\end{align}
\end{Definition}

One shows that $\mathrm{relin} ( \Delta ( \sigma_i) )
\cap \mathrm{relin} ( \Delta ( \sigma_j)) = \emptyset$
for $i \neq j$
and that
$\mathrm{relin} ( \Delta ( \sigma_i) )
\cap S ( \mathscr{X} ) = \emptyset$
for any $i$. 
On the other hand,
for each $i$,
since
the boundary point $\xi_i$ 
of
$\Delta ( \sigma_i)$
is the Shilov point associated to the
irreducible component of $\mathscr{X}_s$ that meets $\sigma_i$,
$\xi_i$ belongs to $S ( \mathscr{X} )$ as well.
Thus 
$S (\mathscr{X} ; \sigma_1 , \ldots , \sigma_r)$ is 
the space obtained by successive one-point sum
of $S (\mathscr{X} )$
with $\Delta (\sigma_1) , \ldots , \Delta (\sigma_r)$
at Shilov points associated to the
irreducible components that contain $\sigma_1 (k) , 
\ldots , \sigma_r (k)$, respectively.

Since $X^{\an} \setminus X(K)$ has a canonical metric,
the skeleton $S (\mathscr{X} ; \sigma_1 , \ldots , \sigma_r)$
is equipped with a metric.
For each $\sigma_i$,
we have a unique isometry $\Delta ( \sigma_i) \cong \RR_{\geq 0}$.
We also consider $\Lambda$-rational points of
$\Delta ( \sigma_i)$ via the isometry.
Thus
one has the notion of
$\Lambda$-rational points of 
$S (\mathscr{X} ; \sigma_1 , \ldots , \sigma_r)$,
and we denote the union of
the set of those points
and
$S_{\Lambda} ( \mathscr{X})$
by
$S_{\Lambda} (\mathscr{X} ; \sigma_1 , \ldots , \sigma_r)$.
When we write $\Gamma$ for the skeleton, we 
write $\Gamma_\Lambda$  for
the subset of $\Lambda$-rational points.

Note that Definition~\ref{def:definitionofskeleton} 
makes sense even if $r = 0$, 
which is the skeleton $S(\Xscr)$ associated to the model $\Xscr$ 
in the sense of Definition~\ref{def:skeleton:model}.
Then we make the following definition.

\begin{Definition}[skeleton of $X^\an$]
\label{def:skeletongeneral}
A  subset $\Gamma$ of $X^{\an}$ is called a \emph{skeleton of $X^\an$}
if there exists a strictly semistable pair 
$(\mathscr{X} ; \sigma_1 , \ldots , \sigma_r)$ 
(possibly $r=0$) such that 
$\Gamma = S(\mathscr{X} ; \sigma_1 , \ldots , \sigma_r)$. 
A skeleton $\Gamma$ of $X^\an$ is called a {\em compact skeleton} if 
$\Gamma$ is compact. 
\end{Definition}

Remark that a skeleton $\Gamma$ is a compact skeleton if and only 
if $\Gamma = S(\Xscr)$ for some strictly semistable model of $X$.
We remark also that the definition of skeleta in this paper
and that in \cite{BPR2} are the same by
\cite[Theorem~4.11]{BPR2}.

\subsection{Some properties of skeleta}
\label{subsec:properties:skeleta}
In this subsection, we recall some properties of skeleta.
Basic references are \cite{BPR2}.

\subsubsection*{Subdivision of a skeleton}
As we have seen, a skeleton associated to a strictly semistable model 
has a structure of a simplicial set such that  the set of irreducible components of the special fiber of the model coincides with the set of vertices.
Here, we give remarks on the relationship between the subdivision of this simplicial structure
and strictly semistable models dominating the given model.

The following propositions
are essentially parts of \cite[Theorem~4.11]{BPR2}.
Indeed, they follow from \cite[Theorem~4.11]{BPR2}
and the fact that an admissible formal model of a smooth projective
curve is algebraizable (cf. \cite[Proposition 10.3.2]{FK}).

\begin{Proposition} \label{prop:subdivision1}
Let $X$ be a connected smooth projective curve over $K$, and let
$\pi : \mathscr{X} \to \Spec (R)$ be a strictly semistable model of $X$.
Let $\pi' : \mathscr{X}' \to \Spec (R)$
be a strictly semistable model of $X$ that dominates $\mathscr{X}$, 
and let $\mu : \mathscr{X}' \to \mathscr{X}$ be the morphism
extending the identity morphism on $X$. 
Let $\{E_1 , \ldots , E_m\}$ be the set of irreducible components of 
$\mathscr{X}'_s$ such that each $\mu (E_i)$ is a singleton of $\mathscr{X}_s$, 
and we assume that $\{E_1 , \ldots , E_m\}$ consists of $(-2)$-curves. 
\textup{(}Namely, we assume that $\mu$ is given by contracting only $(-2)$-curves.\textup{)}
For $i = 1 , \ldots , m$, let $[E_i]$ be the Shilov point associated to $E_i$.
\begin{enumerate}
\item
Then $S(\Xscr') = S(\Xscr)$ as subsets of $X^{\an}$, and 
$V ( \mathscr{X}) = V ( \mathscr{X}') \setminus 
\{ [E_1] , \ldots ,  [E_m] \}$.
\item
Let $p \in \Sing ( \mathscr{X}_s )$, and we set 
$\{q_0 , \ldots , q_r\} := \{q  \in \Sing (\mathscr{X}'_s) \mid \mu(q)= p\}$. 
Then $\Delta_{p} = \bigcup_{i=1}^{r} \Delta_{q_i}$.
\end{enumerate}
\end{Proposition}

The next proposition is a converse of Proposition~\ref{prop:subdivision1}. 

\begin{Proposition} \label{prop:subdivision1:b}
Let $X$ be a connected smooth projective curve over $K$ and let
$\pi : \mathscr{X} \to \Spec (R)$ be a strictly semistable model of $X$. 
Let $V'$ be a finite subset of $S_{\Lambda} (\mathscr{X})$
with $V' \supset V (\mathscr{X})$.
Then there exists a unique
\textup{(}up to a canonical isomorphism\textup{)}
strictly semistable model 
$\pi' : \mathscr{X}' \to \Spec (R)$
of $X$
such that $V ( \mathscr{X}') = V'$.
Further,
$\mathscr{X}'$
has the following properties.
\begin{enumerate}
\item[(i)]
There exists
a unique morphism $\mu : \mathscr{X}' \to \mathscr{X}$ 
that extends the identity morphism on $X$.
\item[(ii)]
Any $E \in \Irr  ( \mathscr{X}'_s)$ such that $\mu (E)$ is a singleton
is a $(-2)$-curve.
\item[(iii)]
We have
$S ( \mathscr{X}') = S (\mathscr{X})$
as subsets of $X^\an$.
\end{enumerate}
\end{Proposition}

\subsubsection*{Retraction map} \label{subsect:retraction}
Let $X$ be a connected smooth projective curve over $K$.
Let $\mathscr{X} \to \Spec (R)$ be a strictly semistable model of $X$,
and set $\Gamma : = S ( \mathscr{X} )$.
Then we have a unique strongly deformation retraction
$\tau : X^{\an} \to \Gamma$ that is the identity on $\Gamma$.
Indeed,
one shows that for each connected component $B$ of 
$X^{\an} \setminus \Gamma$, 
the boundary $\partial B$ consists of a unique point $v$ in $\Gamma_{\Lambda}$ 
and 
the closure $B \cup \{ v \}$ of $B$ is retracted to $v$.
We call $\tau$ the \emph{retraction map} with respect to $\Gamma$.
The retraction map is also described in terms of valuation map;
see \cite[\S5.3]{gubler4}.

\begin{Lemma}
\label{lemma:retraction:rational}
Let $\Gamma$ be a compact skeleton of $X^{\an}$,
and let $\tau_{\Gamma} : X^{\an} \to \Gamma$ be the retraction map.
Let $A$ be a connected component of $X^{\an} \setminus \Gamma$.
Then the following hold.
\begin{enumerate}
\item
There exists $v \in \Gamma$ such that 
$\tau_{\Gamma} (A) = \{ v \}$. Further, $v \in \Gamma_\Lambda$. 
\item
Let $v$ be as in \textup{(1)}.
Let $\mathscr{X}$ be a strictly semistable model of $X$
such that $S( \mathscr{X} ) = \Gamma$ and $v \in V( \mathscr{X})$.
Let $C \in \Irr ( \mathscr{X}_s )$ with $[C] = v$.
Then there exists a unique $q \in C(k) \setminus \Sing ( \mathscr{X}_s )$
such that $\red_{\mathscr{X}} (A) = \{ q \}$.
\end{enumerate}
\end{Lemma}

\Proof
The first assertion 
in (1)
follows from the definition of the retraction map. 
We show the second assertion, i.e., that $v$ is a $\Lambda$-rational point. 
Since $\Gamma$ is a compact skeleton,
there exists a strictly semistable model $\mathscr{X}$
of $X$
such that $S(\mathscr{X}) = \Gamma$.
We take $p \in \Sing ( \mathscr{X}_s )$ such that
$v \in \Delta_p$.
Let
$\psi_p : \mathscr{U}_p \to \mathscr{S}:=\Spec ( R [x,y]/(x y - \varpi_p))$
be a
local \'etale atlas that distinguishes $p$.
Then $\varphi := - \log \left| \rest{\psi_p^{\ast} (y)}{X} \right| $
gives an isometry $\Delta_p \to [0 , - \log | \varpi_p |_K]$.
Therefore,
(1) is reduced to showing that $\varphi (v) \in \Lambda$.
We see from
\cite[\S5.3]{gubler4}
or \cite[Lemma~3.8]{BPR2}
that 
\begin{align}
\label{align:retraction:val}
\varphi \circ \rest{\tau_{\Gamma}}{A} = \rest{\varphi}{A}.
\end{align}
Since $X(K)$ is dense in $X^{\an}$,
there exists $P \in A \cap X(K)$.
Since $\tau_{\Gamma} (A) = \{ v \}$, we have $v = \tau_{\Gamma} (P)$.
By (\ref{align:retraction:val}), 
$\varphi (v) = \varphi (\tau_{\Gamma} (P) ) = \varphi (P)$,
and this belongs to $\Lambda$ by the definition of $\varphi$.
Thus 
we obtain (1).

Let us prove (2).
Set $C^{\circ} := C \setminus \Sing ( \mathscr{X}_s )$.
By \cite[Proposition~5.7]{gubler4},
$\red_{\mathscr{X}} ( \tau_{\Gamma}^{-1} (v) ) = C^{\circ}$.
Recall that $v$ is the unique point in $X^{\an}$ such
that $\red_{\mathscr{X}} (v)$ is the generic point of $C$.
Then we have 
$\tau_{\Gamma}^{-1} (v) \setminus \{ v \} \subset 
\red_{\mathscr{X}}^{-1} (C^{\circ} (k) )$,
so that
$A \subset \red_{\mathscr{X}}^{-1} (C^{\circ} (k) )$.
Note that
$
\red_{\mathscr{X}}^{-1} (C^{\circ} (k) ) = \bigcup_{q \in C^{\circ} (k)}
\red_{\mathscr{X}}^{-1} ( q )$,
Since $\red_{\mathscr{X}}$ is anti-continuous,
each $\red_{\mathscr{X}}^{-1} ( q )$ is an open subset of $X^{\an}$,
and thus 
$\bigcup_{q \in C^{\circ} (k)}
\red_{\mathscr{X}}^{-1} ( q )$ a disjoint union of open subsets
of $\red_{\mathscr{X}}^{-1} (C^{\circ} (k) )$.
Since $A$ is a connected open subset of $\red_{\mathscr{X}}^{-1} (C^{\circ} (k) )$,
it follows that
there exists a unique $q \in C^{\circ} (k)$
such that $A \subset \red_{\mathscr{X}}^{-1} (q)$.
This proves (2).
\QED

\subsubsection*{Minimal skeleton and Deligne--Mumford stable model}

Let $X$ be a connected smooth projective curve of genus $g$.
Let $\mathscr{X}$ be a strictly semistable model of $X$.
We say that $\mathscr{X}$ is 
\emph{minimal} if $\mathscr{X}_s$ does not have a $(-1)$-curve.
Remark that in the case where $g \geq 1$,
$\mathscr{X}$ is minimal if and only if $\mathscr{X}$ is 
Deligne--Mumford semistable.

\begin{Definition}[minimal skeleton]
\label{def:minimal:skeleton}
We say that a skeleton $\Gamma$ of $X^{\an}$ is a \emph{minimal skeleton} 
if there exists a minimal model $\mathscr{X}$ of $X$ with $\Gamma = S ( \mathscr{X} )$.
\end{Definition}

We note that the notion of a minimal skeleton is the same as 
that in \cite{BPR2}. 
Remark that if $\Gamma$ is minimal,
then any model $\mathscr{X}'$ with $\Gamma = S ( \mathscr{X}')$ is minimal.

Assume that $g \geq 1$. 
Then the minimal skeleton is unique. 
Indeed,
let $\mathscr{X} \to \Spec (R)$ be any  
Deligne--Mumford strictly semistable model of $X$;
then the minimal skeleton equals
$S ( \mathscr{X} )$ as a subspace of $X^{\an}$,
and it
does not depend on the choice of $\Xscr$. 
We denote by $\Gamma_{\min}$ the minimal skeleton.

Furthermore, assume that $g \geq 2$.
Then the minimal skeleton 
$\Gamma_{\min}$ of $X^{\an}$
has a {\em canonical} finite graph
 structure determined by the 
{\em stable} model. 
Let $\mathscr{X}^{\st}$ be the stable model of $X$.
As is noted in \S\ref{subsec:semistabe:models}, there exists a 
Deligne--Mumford strictly semistable model $\mathscr{X}$
that dominates $\mathscr{X}^{\st}$.
It follows that the set of Shilov points $V (\mathscr{X}^{\st})$ with
respect to 
$\mathscr{X}^{\st}$ is a subset of $V ( \mathscr{X})$
(in fact, $V (\mathscr{X}^{\st})$ equals 
to the subset of $V ( \mathscr{X}) $ consisting of the Shilov points 
associated to non-$(-2)$-curves),
and hence 
$V (\mathscr{X}^{\st})$ 
is also a subset of $\Gamma_{\min}$. 
We give 
$\Gamma_{\min}$ a finite graph structure so that 
$V (\mathscr{X}^{\st})$ is the set of vertices. 
Let $E ( \mathscr{X}^{\st})$ denote the set of edges with respect to this
graph structure.
Note that this finite graph may have self-loops in general.

We remark that for any $p \in \Sing (\mathscr{X}^{\st}_s)$,
$\Delta_p$ makes sense.
Indeed, let $\mathscr{X}$ be a Deligne--Mumford strictly semistable model 
of $X$ and let
$\mu : \mathscr{X} \to \mathscr{X}^{\st}_s$
be the morphism extending the identity morphism on $X$. 
Let $\{q_1 , \ldots , q_r\}$ be the set of nodes of $\mathscr{X}_s$ 
such that $\mu (q_1) = \cdots = \mu (q_r) = p$. 
Then $\bigcup_{i = 1}^{r} \Delta_{q_i}$ is a subset of $\Gamma_{\min}$
and does not depend on the choice of $\mathscr{X}$. 
We set $\Delta_p := \bigcup_{i = 1}^{r} \Delta_{q_i}$. 
Note that $\Delta_p \in E(\mathscr{X}^{\st})$. 

When $g \geq 2$,
we use the following convention  unless otherwise specified:
When we say that $v$ is a vertex of $\Gamma_{\min}$,
this means that $v \in V ( \mathscr{X}^{\st})$;
when we say that $e$ is an edge,
this means that $e \in E (\mathscr{X}^{\st})$.

\subsection{Tropical geometry}
\label{subsec:tropical:geometry}
Let $\GG_{m}^N$ be the algebraic torus of dimension $N$
over $K$
with coordinates $z_1 , \ldots , z_N$
and let $\GG_m^{N , \an}$ be the associated analytic space.
The map
\begin{align}
\label{eqn:trop:map:torus}
\trop_{\GG_m^N} : \GG_{m}^{N ,\an} \to \RR^N
,
\quad
p = (p , | \cdot |)
\mapsto
( - \log |z_1 (p) | , \ldots , - \log |z_N (p)|)
\end{align}
is
called the \emph{tropicalization map of algebraic torus}.
We extend this tropicalization map from
the analytic space
associated to projective space to tropical projective space.
As in the introduction, 
we set  $\TT := \RR\cup\{+ \infty\}$. The 
\emph{$N$-dimensional tropical projective space} is defined to be 
\[
\TT\PP^N := 
(\TT^{N+1}\setminus\{(+ \infty, \ldots, +\infty)\})/\sim, 
\]
where $x := (x_0, \ldots, x_N), y := (y_0, \ldots, y_N) \in \TT^{N+1}\setminus\{(+\infty, \ldots, +\infty)\}$ satisfy $x \sim y$ if there exists 
$c \in \RR$ such that $y_i = x_i + c$ for all $i = 0, \ldots, N$ (see \cite{MZ}). The equivalence class of $x$ in $\TT\PP^N$ is written as 
$(x_0: \cdots: x_N)$. 
Now, let $\PP^N$ be the $N$-dimensional projective space over $K$ with
homogeneous coordinates
$X_0, \ldots, X_N$
and let $\PP^{N, \an}$ denote the Berkovich analytification of $\PP^N$. 
Then we define the tropicalization map to be
\begin{equation}
\label{eqn:trop:map}
\trop: \PP^{N, \an} \to \TT\PP^N, \qquad p = (p, |\cdot|)  \mapsto (-\log|X_0(p)|
: \cdots: -\log|X_N(p)|). 
\end{equation}

The tropicalization map on (\ref{eqn:trop:map})
extends that on (\ref{eqn:trop:map:torus}), which we make precise now.
The tropical projective space is equipped with $(N+1)$ charts 
$U_i := \{x = (x_0: \cdots: x_N) \in \TT\PP^N\mid x_i \neq + \infty\}$. 
We set $E:=\bigcap_{i = 0}^N U_i$.
For each $i= 0 , \ldots , N$,
we have a natural homeomorphism 
\[
\phi_i : \RR^{N} \to E, 
\quad 
(u_1 , \ldots , u_N) \mapsto ( u_1 : \cdots : u_{i-1} : 0 : u_{i}
: \cdots : u_N ).
\]
Thus $E$ is an $N$-dimensional Euclidean space embedded in $\TT\PP^N$.
Further,
for each $i=0,\ldots ,N$, let 
\[
\psi_i : \GG_m^N \hookrightarrow
\PP^N, \quad 
(z_1 , \ldots , z_N)
\mapsto (z_1 : \cdots : z_{i-1} : 1 : z_{i} \cdots : z_N )
\]
denote the open immersion. Then we have
$\trop \circ \psi_i = \phi_i \circ \trop_{\GG_m^N}$.
Remark that $\trop^{-1} ( E ) = \psi_i ( \GG_m^N)$
for any $i=0 , \ldots , N$.

\begin{Remark}
\label{remark:integralstructure}
Let $i , j \in \{ 0 , 1 , \ldots , N \}$.
Then we have an homeomorphism
$\phi_j^{-1} \circ \phi_i : \RR^N \to \RR^N$.
It is clear from the definition of $\phi_i$
and $\phi_j$ that
this restricts an isomorphism from $\ZZ^N$ to $\ZZ^N$
as $\ZZ$-modules,
i.e., $\phi_j^{-1} \circ \phi_i$ is given by
the multiplication of an element of $\GL_N (\ZZ)$.
\end{Remark}

Let $Y$ be an irreducible closed subvariety of $\PP^N$.
Fix an $i = 0 , \ldots , N$.
Then $Y_i^{\circ} := \psi_i^{-1} (Y)$
is an irreducible closed subvariety of $\GG_m^N$.
The results of Bieri--Groves and Speyer--Sturmfels (cf.~\cite{MS} and \cite[Theorem~3.3]{GuTool})
describe the polytopal structure of 
$\trop_{\GG_m^N}\left(Y_i^{\circ, \an}\right)$.
Although we do not recall their results in full generality,
we describe the polytopal structure
when $\dim(Y_i^\circ) = 1$, because that is the case we are concerned with.
In this case,
there exists a finite subset $V \subset 
\trop_{\GG_m^N}\left(Y_i^{\circ, \an}\right) \cap \Lambda^N$
such that if $\Sigma$ is the set of closures of the connected components
of
of $\trop_{\GG_m^N}\left(Y_i^{\circ, \an}\right) \setminus V$,
then each $\Delta \in \Sigma$
is of form
\begin{enumerate}
\item
$\Delta = \{y + t z \mid t \in [0, \ell]\}$
or
\item
$\Delta = \{y + t z \mid t \in [0, +\infty)\}$, 
\end{enumerate}
where $y \in \Lambda^N$, 
$z = (z_1, \ldots, z_N) \in \ZZ^N$ with ${\rm GCD}(z_1, \ldots, z_N) = 1$,
and $\ell \in \Lambda \cap \RR_{> 0}$.

Using the above illustration,
we endow $\trop_{\GG_m^N}\left(Y_i^{\circ, \an}\right)$ with a metric structure.
Indeed, we put a metric
by lattice length:
If $\Delta \in \Sigma$ is as in (1)
(resp. (2)),
then  we identify $\Delta$ with $[0, \ell]$ 
(resp. $[0, + \infty)$)
via 
$y + t z \mapsto t$, and with this identification,
$\Delta$ is a metric space;
since $\trop_{\GG_m^N}\left(Y_i^{\circ, \an}\right) = \bigcup_{\Delta} \Delta$, 
this metric structure does not depend on the choice of $\Sigma$
and is well-defined for $\trop_{\GG_m^N}\left(Y_i^{\circ, \an}\right)$.
Thus $\trop_{\GG_m^N}\left(Y_i^{\circ, \an}\right)$ is
a metric space.

Furthermore,
the above metric structure on each $\trop_{\GG_m^N}\left(Y_i^{\circ, \an}\right)$
gives a well-defined metric on
$\trop (Y^{\an}) \cap E$.
Indeed,
for each $i = 0 , \ldots , N$, 
since
$\trop^{-1} ( E ) = \psi_i ( \GG_m^N)$,
we have a homeomorphism
$\rest{\phi_i}{\trop_{\GG_m^N}\left(Y_i^{\circ, \an}\right)} : 
\trop_{\GG_m^N}\left(Y_i^{\circ, \an}\right) \to \trop (Y^{\an}) \cap E$,
and since $\trop_{\GG_m^N}\left(Y_i^{\circ, \an}\right)$
has a metric structure, this homeomorphism induces a metric structure
on $\trop (Y^{\an}) \cap E$. 
Since the metric on 
$\trop_{\GG_m^N}\left(Y_i^{\circ, \an}\right)$ is 
given by the lattice length,
one sees from Remark~\ref{remark:integralstructure}
that this metric structure on $\trop (Y^{\an}) \cap E$
does not depend on $i$.
Thus $\trop (Y^{\an}) \cap E$ is a metric space.

Remark that
tropical geometry near the boundary $\TT\PP^N \setminus 
E$ is rather subtle (see \cite{MZ}).
However, for our purposes, namely, for faithful tropicalizations, we will not need detailed analysis on the boundary 
because the tropicalization of a skeleton is contained in $E$,
as we will see in the next subsection.

\subsection{Faithful tropicalization}
\label{subsec:unimodular:faithful}
Let $X$ be a connected smooth projective curve over $K$.
Let $L$ be a line bundle over $X$.
Suppose that we are given nonzero global sections 
$s_0, s_1, \ldots, s_N \in H^0(X, L)$.  
Associated to those sections, we define a map
$\varphi : X^\an \to \TT\PP^N$ as follows:
Let $\varphi^\prime: X^\an \longrightarrow \PP^{N, \an}$ be the morphism 
induced by $p  \mapsto \left(s_0(p): \cdots: s_N(p)\right)$;
we define $\varphi : X^\an \to \TT\PP^N$
to be $\trop \circ \varphi^\prime$,
where $\trop$ is the map in (\ref{eqn:trop:map}).
We write
\begin{equation}
\label{eqn:map:tropical}
  \varphi: X^\an \longrightarrow \TT\PP^N, \quad p = (p, |\cdot|) \mapsto \left(-\log|s_0(p)|: \cdots: -\log|s_N(p)|\right).
\end{equation}

Let $E$ be the $N$-dimensional Euclidean space embedded in $\TT\RR^N$
as in the previous subsection.
Note that $\varphi ( X^{\an} \setminus X (K) )
\subset \varphi (X^{\an} ) \cap E$.
Indeed,
since each $s_i$ has zero only at points in $X (K)$,
$- \log |s_i (p)| \neq  + \infty$ for any $p \in X^\an \setminus X (K)$,
and hence
$\varphi ( X^\an \setminus X (K) ) \subset E$.

We describe $\rest{\varphi}{X^\an \setminus X (K) }$ 
in terms of the coordinates given by 
the homeomorphism $\phi_i : \RR^N \to E$ for $i = 0, \ldots, N$. 
We explain only the case $i = 0$, since the other cases are the same. 
We identify $E = \RR^N$ via $\phi_0$.
Then
\begin{equation}
\label{eq:trop:affine}
\rest{\varphi}{X^\an \setminus X (K) } = 
\left( 
- 
\log \left| \frac{s_1}{s_0}
\right|
,
\ldots 
,
- \log \left| \frac{s_N}{s_0}
\right|
\right)
.
\end{equation}

Let $\Gamma = S(\Xscr; \sigma_1, \ldots, \sigma_r)$ 
be a skeleton (possibly with ends)
and write $\Gamma = \bigcup_{q \in \Sing ( \mathscr{X}_s)} \Delta_q \cup  \bigcup_{\sigma_i} \Delta(\sigma_i)$ be the decomposition 
of $\Gamma$ as in \eqref{align:definitionofskeleton}. 
Recall that  $\Gamma$ is a metric space;
$\Delta_q$ is canonically identified with 
$[0 , \lambda_q]$, where $\lambda_q$ is the multiplicity at $q$,
and
$\Delta (\sigma_i)$ is canonically identified with $\RR_{\geq 0}$. 
By the definition of $\Gamma$,
we have
$\Gamma \subset X^\an \setminus X (K)$.
Since $\varphi ( X^{\an} \setminus X (K) )
\subset \varphi (X^{\an} ) \cap E$,
we have a map
$\Gamma \to \varphi (X^{\an} ) \cap E$
between metric spaces by restriction. 
Thus the following definition makes sense.

\begin{Definition}[unimodular tropicalization]
Let $X$, $L$, and $\Gamma$ be as above.
Let $s_0, s_1, \ldots, s_N$ be nonzero 
global sections of $L$, and let $\varphi: X^\an \longrightarrow \TT\PP^N$ denote the  associated morphism \eqref{eqn:map:tropical}. We call $\varphi$  a {\em unimodular tropicalization} of $\Gamma$ if there exists a finite subset $V \subset \Gamma_\Lambda$ 
such that
for the closure 
$\Delta$
of any connected component of $\Gamma \setminus V$,
the restriction map $\rest{\varphi}{\Delta}: \Delta\to\varphi(\Delta)$ is an isometry. 
\end{Definition}

Namely, a unimodular tropicalization is a piecewise isometry.

\begin{Remark}
\label{remark:unimodular}
Let $X$, $\Gamma$, and $L$ be as above.
Let $s_0 , \ldots , s_M , s_{M+1} , \ldots , s_N$ be global sections
of $L$
($N \geq M+1$).
Let $\varphi_1 : X^{\an} \to \TT\PP^M$ be the  associated morphism \eqref{eqn:map:tropical}
to $s_0 , \ldots , s_M$,
and let
$\varphi_2 : X^{\an} \to \TT\PP^N$ be the  associated morphism \eqref{eqn:map:tropical}
to $s_0 , \ldots , s_N$.
Let $e$ be a connected subspace of $\Gamma$.
Suppose that $\rest{\varphi_1}{e}$ is an isometry.
Then $\rest{\varphi_2}{e}$ is also an isometry.
Indeed, 
with the description on (\ref{eq:trop:affine}),
we have
$\rest{\varphi_1}{e} = 
\left( 
- 
\log \left| \frac{s_1}{s_0}
\right|
,
\ldots 
,
- \log \left| \frac{s_M}{s_0}
\right|
\right)$
and
$\rest{\varphi_2}{e } = 
\left( 
\rest{\varphi_1}{e}
,
- 
\log \left| \frac{s_{M+1}}{s_0}
\right|
,
\ldots 
,
- \log \left| \frac{s_N}{s_0}
\right|
\right)$.
Since 
 the metric on the tropicalization is given by the lattice length,
one then sees that
if $\rest{\varphi_1}{e}$ is an isometry, then so is
$\rest{\varphi_2}{e}$.
\end{Remark}

\begin{Definition}[faithful tropicalization]
Let $X$, $L$, and $\Gamma$ be as above.
Let $s_0, s_1, \ldots, s_N$ be nonzero 
global sections of $L$ and let $\varphi: X^\an \longrightarrow \TT\PP^N$ denote the associated morphism \eqref{eqn:map:tropical}. We 
call $\varphi$ a {\em faithful tropicalization} of $\Gamma$ if 
it is an injective unimodular tropicalization;
in other words,
the restriction of $\varphi$ to $\Gamma$ is 
a homeomorphism of $\Gamma$ onto its image 
preserving the metric. 
\end{Definition}

\begin{Remark}
In relation to faithful tropicalizations, we recall the notion of $\Lambda$-rational polyhedral complexes and their integral structures, and remark that  faithful tropicalizations are exactly homeomorphisms preserving 
the integral structures. 

Here we describe $\Lambda$-rational polyhedral complexes of dimension at most one. 
A $\Lambda$-rational  polyhedron $\sigma$ in $\RR$  of dimension at most one is a subset defined by 
$\sigma = \{x \in \RR \mid A x \leq b\}$ for some $A \in M_{r, 1}(\ZZ)$ and 
$b \in \Lambda^r$. A $\Lambda$-rational polyhedral complex $\Sigma$ of dimension at most one on a Hausdorff topological space $Y$ is a finite collection of closed topological spaces $\Delta \subset Y$ such that $\Delta$ is homeomorphic to a $\Lambda$-rational  polyhedron $\sigma$ in $\RR$  of dimension at most one with the following properties: Any face of $\Delta \in \Sigma$ belongs to $\Sigma$;  if $\Delta$ and $\Delta'$ are in $\Sigma$, then $\Delta \cap \Delta'$ (if non-empty) are in $\Sigma$.  
Such a  polyhedral complex $\Sigma$ is said to be of pure dimension one if $|\Sigma| := \{y \in Y \mid \text{$y \in \Delta$ for some $\Delta \in \Sigma$} \}$ is connected and not a singleton. 
A {\em refinement} of $\Sigma$ is a $\Lambda$-rational polyhedral complex $\widetilde{\Sigma}$ with 
$|\widetilde{\Sigma}| = |\Sigma|$ such that for any $\widetilde{\Delta} \in \widetilde{\Sigma}$, there exists $\Delta \in \Sigma$ with $\widetilde{\Delta} \subset \Delta$. 

Let $\Sigma$ be a $\Lambda$-rational polyhedral complex of pure dimension one on $Y$. 
A member $\Delta \in \Sigma$ is said to be zero-dimensional (resp. one-dimensional) if the 
corresponding  $\Lambda$-rational  polyhedron $\sigma$ in $\RR$ is zero-dimensional (resp. one-dimensional). Let $|\Sigma|_0$ be the union of all $\Delta \in \Sigma$ such that $\Delta$ is zero-dimensional, and let $|\Sigma|_1$ be the union of all $\Delta \in \Sigma$ such that $\Delta$ is one-dimensional.  Then $|\Sigma|_0$ (resp. $|\Sigma|_1$) is naturally a manifold of dimension $0$ (resp. $1$) with an integral structure. Namely, for $n = 0$, this means nothing but 
that each member of $|\Sigma|_0$ is identified with a point in $\Lambda \subset \RR$. For $n=1$, this means that there exist an open covering $\{U_i\}$ of $|\Sigma|_1$ and homeomorphisms $\phi_i: U_i \to V_i \subset \RR$ such that 
for each $i$ and $j$, the transition map $\rest{\phi_i \circ \phi_j^{-1}}{\phi_j(U_i \cap U_j)}: \phi_j(U_i \cap U_j) \to 
\phi_i(U_i \cap U_j)$ is the restriction of the map $x \mapsto c x + d$, where $c = \pm1 $ and 
$d \in \Lambda$.  

Let $\Sigma$ (resp. $\Sigma'$) be a $\Lambda$-rational polyhedral complex of pure dimension one on 
$Y$ (resp. $Y'$). Let $f:Y \to Y'$ be a continuous map. 
We say that the restriction of $f$ to $|\Sigma|$ is a homeomorphism onto $|\Sigma'|$ preserving 
the integral structures if there exist refinements $\widetilde{\Sigma}$ of $\Sigma$ 
and 
$\widetilde{\Sigma'}$ of $\Sigma'$ with the following properties: 
$|\widetilde{\Sigma}|_0$ and $|\widetilde{\Sigma'}|_0$ are bijective under $f$;  
there exist open coverings $\{U_i\}$ of $|\widetilde{\Sigma}|_1$, $\{U'_i\}$ of 
$|\widetilde{\Sigma'}|_1$, and homeomorphisms $\phi_i: U_i \to V_i \subset \RR$
and 
$\phi'_i: U'_i \to V'_i \subset \RR$ that give the integral structures of $|\widetilde{\Sigma}|_1$ and 
$|\widetilde{\Sigma'}|_1$ such that 
$\phi'_i \circ f \circ\phi_i^{-1}:  V_i \to V'_i$ is a homeomorphism of the form $x \mapsto c x + d$, 
where $c = \pm 1$ and $d \in \Lambda$. 

Let $X$, $L$, and $\Gamma$ be as above.
Let $s_0, s_1, \ldots, s_N$ be nonzero 
global sections of $L$ and let $\varphi: X^\an \longrightarrow \TT\PP^N$ denote the associated morphism \eqref{eqn:map:tropical}. Then each of
$\Gamma$ and $\varphi(\Gamma)$ is equipped with 
a natural structure of a $\Lambda$-rational polyhedral complex of dimension at most one, and 
$\varphi$ is a faithful tropicalization of $\Gamma$ in the above definition if and only if 
$\rest{\varphi}{\Gamma}: \Gamma \to \varphi(\Gamma)$ is 
a homeomorphism of $\Gamma$ onto its image preserving 
the integral structures. 
See \cite[\S3]{MZ},  \cite{BPR2}, and \cite[\S\S2.2--2.3]{GRW} for more details. 
\end{Remark}

\setcounter{equation}{0}
\section{Good models}
\label{sec:models}
Let $X$ be a connected smooth projective curve over $K$ of genus $g$,
and let $L$ be a line bundle over $X$.
In the proof of the main Theorem~\ref{thm:main:unimodular:faithful},
we need to construct many suitable global sections of
a line bundle $L$ over $X$ 
which will give a faithful tropicalization for a skeleton.
In our strategy,
we will construct those global sections as global sections
of $\mathscr{L}$
for some models $( \mathscr{X} , \mathscr{L} )$ 
of $(X,L)$.
Thus we will need to construct various models of $(X,L)$,
and to do that,
we use the theory of divisors on graphs. 

In this section, we firstly define the notion of good models $\mathscr{X}$ of 
$X$. In subsequent sections, over a good model $\Xscr$ of $X$, 
we will construct a line bundle $\Lscr$ with 
$\Lscr \otimes_R K \cong L$ and 
suitable global sections of $\Lscr$. 
We secondly recall some recent development on the divisor theory on 
$\Lambda$-metric graphs.
We thirdly introduce the notion of {\em islands} of weighted $\Lambda$-metric
graphs
and prove a positivity result on divisors on 
$\Lambda$-metric graphs.
In the last subsection of this section,
we
introduce the machinery that produces
various models $( \mathscr{X} , \mathscr{L} )$
which will be crucially used in the subsequent arguments.
Further,
we use it to construct a model 
of $(X,L)$ which will be frequently used.

\subsection{Good models of $X$}
\label{subsec:def:good:model}
In this subsection, we define the notion that a model $\Xscr$ of $X$ is 
good. 
Let 
$\pi : \mathscr{X} \to \Spec (R)$ be a strictly semistable model of $X$. 
Recall from \S\ref{subsec:semistable:models:pairs} that 
a $(-2)$-chain is a connected curve in $\mathscr{X}_s$ consisting of 
$(-2)$-curves, and it is a {\em maximal} $(-2)$-chain if it is maximal with respect to 
the inclusion among $(-2)$-chains. 
Recall also from \S\ref{subsec:semistable:models:pairs} that 
a $(-2)$-chain $E$ is of connected type if $\mathscr{X}_s - E$ is connected, and 
is of disconnected type otherwise.

\begin{Definition}[symmetric multiplicities]
\label{def:symmetricmultiplicity}
Let $E$ be a maximal $(-2)$-chain of connected type
in $\mathscr{X}_s$.
We say that $E$ \emph{has symmetric multiplicities} if it satisfies
the following conditions.
\begin{enumerate}
\item[(i)]
The number $\# \Irr(E)$ is odd and at least $3$. We write 
$\# \Irr(E) = 2 \ell - 1$ for $\ell \geq 2$.
\item[(ii)]
Put $D := \mathscr{X}_s - E$.
Let $E_1 , \ldots , E_{2 \ell - 1}$ be the irreducible components of $E$
such that
$E_{i} \cap E_{i+1} \neq \emptyset$ for $i = 1 , \ldots , 2 \ell -2$.
Let $p_0 , \ldots , p_{2 \ell -1}$ be the nodes of $\mathscr{X}_s$
such that $\{ p_i \} = E_i \cap E_{i+1}$ for $i = 1 , \ldots , 2 \ell -2$,
$\{ p_0 \} = D \cap E_1$, and $\{ p_{2 \ell-1} \} = D \cap E_{2 \ell -1}$.
For $i = 0 , \ldots , 2 \ell - 1$,
let $\lambda_i$ be the multiplicity at $p_i$ (cf. Definition~\ref{def:multiplicityatanodes}).
Then
$\lambda_j = \lambda_{2\ell - 1 -j}$ for $j = 0 , \ldots , \ell -1$.
\item[(iii)]
We have $\lambda_0 = \lambda_{\ell-1}$. (Thus $\lambda_0 = \lambda_{\ell-1} 
= \lambda_{\ell} = \lambda_{2\ell-1}$.) 
\end{enumerate}
\end{Definition}

From here on to the end of this subsection,
we assume that $g \geq 2$.
Let $\mathscr{X}$ be a strictly semistable
model of $X$.
Then two distinct maximal $(-2)$-chains 
in $\mathscr{X}_s$
are disjoint.

\begin{Definition}[good model] \label{def:goodmodel:new}
Let $X$ be a connected smooth projective curve of genus $g \geq 2$, and 
let $\mathscr{X}$ be a 
model of $X$.
We say that $\mathscr{X}$ is a \emph{good model} if 
it satisfies the following conditions.
\begin{enumerate}
\item[(i)]
The model $\mathscr{X}$ is Deligne--Mumford strictly semistable.
\item[(ii)]
For any $p \in \Sing ( \mathscr{X}_s)$, there exists a 
$(-2)$-chain $E$ such that $p \in E$.
\item[(iii)]
Any maximal $(-2)$-chain $E$ in $\mathscr{X}_s$ 
of connected type
has symmetric multiplicities.
\end{enumerate}
\end{Definition}

We make clear the relationship between a good model $\mathscr{X}$ and the canonical finite graph structure of the minimal skeleton $\Gamma_{\min}$ with the set $V(\Xscr^{\st})$ of vertices and the set $E(\Xscr^{\st})$ of edges. 
Recall that $V(\Xscr^{\st}) =
\{ [C] \mid C \in \Irr (\mathscr{X}^{\st}_s) \}$ by definition,
where $[C]$ is the Shilov point associated to $C$,
and an element of $E(\Xscr^{\st})$ is characterized by the closure of a connected component 
of $\Gamma_{\min} \setminus V(\Xscr^{\st})$. We say that $e \in E(\Xscr^{\st})$ is {\em of connected type} if $\Gamma_{\min} \setminus 
\mathrm{relin}(e)$ is connected, and $e \in E(\Xscr^{\st})$ is {\em of disconnected type} otherwise. 

\begin{Lemma}
\label{lemma:make:clear}
Let $\mathscr{X} \to \Spec (R)$ be a good model of $X$, 
and let $\mu: \mathscr{X} \to \Xscr^{\st}$ be the morphism extending the identity morphism on $X$. 
\begin{enumerate}
\item
There exist natural one-to-one correspondences between the following three sets: 
\[
   \Sing(\Xscr_s^{\st}) \longleftrightarrow E(\Xscr^{\st}) 
   \longleftrightarrow \{\text{maximal $(-2)$-chain in $\Xscr_s$}\}. 
\]
\item
Let $e \in E(\Xscr^{\st})$ correspond to a maximal $(-2)$-chain $E$ in $\Xscr_s$ under the natural correspondence in \textup{(1)}. Then $e$ is  of connected type \textup{(}resp. of disconnected type\textup{)} if and only if $E$  is of connected type \textup{(}resp. of disconnected type\textup{)}.  
\end{enumerate}
\end{Lemma}

\Proof
(1) Since $\Xscr$ is Deligne--Mumford semistable, 
the irreducible components of $\Xscr_s$ contracted by $\mu$ 
are
the $(-2)$-curves. 
For $p \in \Sing(\Xscr_s^{\st})$, we set $E := \mu^{-1}(p)$. Since $\Xscr$ satisfies  condition (ii) of Definition~\ref{def:goodmodel:new}, $E$ is not a singleton, and thus 
is a maximal $(-2)$-chain in $\Xscr_s$. This gives a natural one-to-one correspondence between the first and the third sets. For a maximal $(-2)$-chain $E$ in $\Xscr_s$, we set $\Delta_E:= \bigcup_{q \in E \cap \Sing ( \mathscr{X}_s)} \Delta_q$ as in Definition~\ref{def:canonical:1:simplex:curve}. 
Then $\Delta_E \in E(\Xscr^{\st})$, and this gives a natural one-to-one correspondence between the third and the second sets. Finally, we remark that $\Delta_p$ at the end of \S\ref{subsec:properties:skeleta} is 
exactly $\Delta_E$. We have thus natural correspondences: 
\begin{equation}
\label{eqn:correspondences}
   p = \mu(E) \quad\longleftrightarrow\quad \Delta_p = \Delta_E 
   \quad\longleftrightarrow\quad E:= \mu^{-1}(p). 
\end{equation}

(2) The assertion follows from the definitions of being of connected type and of disconnected type. 
\QED

Let $E$ be a maximal $(-2)$-chain in $\mathscr{X}_s$.
Then $\Delta_{E}$
is a circle or a closed line segment in $S ( \mathscr{X}_s)$.
The irreducible components $C \in \Irr(\Xscr_s)$  
having the properties 
that $C \notin \Irr (E)$ and 
$C \cap E \neq \emptyset$ give points in $\Delta_E$,
and
those points are in $V ( \mathscr{X}^{\st})$. 

We finish this subsection by introducing the notation,
 which will  be frequently used.
Let $\mathscr{X}$ be a good
model of $X$ and let 
$\mu : \mathscr{X} \to \mathscr{X}^{\st}$ be the morphism
extending the identity morphism on $X$. 
For each $q \in \Sing ( \mathscr{X}_s )$,
let $\Delta_q$ denote the corresponding $1$-simplex in the minimal skeleton
$\Gamma_{\min}$.
Remark that if $E$ is the maximal $(-2)$-chain with $q \in E$,
then $\Delta_q \subset \Delta_E = \Delta_{\mu (q)}$
with the notation in (\ref{eqn:correspondences}).
For a given compact subset $\Gamma_1 \subset \Gamma_{\min}$,
we set
\begin{align}
\label{align:nodes:subset}
\Sing (\mathscr{X}_s )_{\subset \Gamma_1} :=
\left\{
q \in \Sing (\mathscr{X}_s )
\mid
\Delta_q \subset \Gamma_1
\right\}
.
\end{align}
Remark that
for each $e \in E ( \mathscr{X}^{\st})$,
\begin{align*}
\Sing (\mathscr{X}_s )_{\subset e} =
\left\{
q \in \Sing (\mathscr{X}_s )
\mid
\Delta_{\mu (q)} = e
\right\}
=
\Sing (\mathscr{X}_s ) \cap E
,
\end{align*}
where $E$ is the maximal $(-2)$-chain with $\Delta_E = e$.

\subsection{Theory of divisors on $\Lambda$-metric graphs}
\label{subsec:theory:divisors}
As in ``Notation and conventions,'' 
let $\Lambda \subset \RR$ denote the value group of $K$.  
We review the theory of divisors on $\Lambda$-metric graphs, which we use later. 
Our basic references are Amini, Baker, Brugall\'e, and Rabinoff~\cite{ABBR}, 
Baker and Norine~\cite{BN}, Amini and Caporaso~\cite{AC}, and 
Baker and Rabinoff~\cite{BR}.  See also \cite{KY0, KY1}.

Let $\Gamma$ be a $\Lambda$-metric graph (see ``Notation and conventions''). 
Let $\Div (\Gamma)$ 
be the free abelian group generated 
by the points in $\Gamma$. 
An element of $\Div (\Gamma)$ 
is called a {\em divisor} on $\Gamma$. 
Any divisor $D \in \Div (\Gamma)$ is uniquely written as a finite sum 
$D = \sum_{x \in \Gamma} n_x [x]$ for $n_x \in \ZZ$. We set $D(x) := n_x$. 
Remark that $D = \sum_{x \in \Gamma} D(x) [x]$. 
The {\em support} $\Supp(D)$ of $D$ is the set of points with 
$D(x) \neq 0$. The {\em degree} of $D$ is defined to be  $\deg(D) = 
\sum_{x \in \Gamma} D(x)$. 
We say that a divisor $D$ is {\em effective} and write $D \geq 0$
if $D(x) \geq 0$ for any $x \in \Gamma$. 
For a subset $S$ of $\Gamma$, we set
$\rest{D}{S} := \sum_{x \in S} D(x)[x]$
and
\begin{equation}
\label{eqn:deg:extended}
\deg\left(\rest{D}{S}\right)
=  \sum_{x \in S} D(x). 
\end{equation}

We recall the notion of linear equivalence of divisors on a metric graph.
For a point $x \in \Gamma$,  $\val(x)$ denotes the valency at $x$,
i.e.,
the number of branches emanating from $x$. 
A {\em rational function} on $\Gamma$ is 
a continuous function $f: \Gamma \to \RR$ such that there exists 
a finite subset $\{x_1, \ldots, x_n\}$ of  $\Gamma$ 
containing $\{x \in \Gamma \mid \val(x) \neq 2\}$ such that
$f$ is an affine map with integer slopes on $\Gamma \setminus \{x_1, \ldots, x_n\}$. 
We denote by $\Rat(\Gamma)$ the set of 
rational functions on $\Gamma$. 
For 
$f \in \Rat(\Gamma)$ and $x \in \Gamma$, we define 
$\ord_x(f)$ to be the sum of outgoing slopes of $f$ at $x$. 
Then $\zero(f) := \sum_{x \in \Gamma} \ord_x(f) [x]$ is a 
divisor on $\Gamma$, which is 
called the \emph{principal divisor} of $f$.
We set $ \Prin (\Gamma) := \{\zero(f) \mid f \in \Rat (\Gamma)\}$,
the set of {\em principal divisors} on $\Gamma$.
For $D_1 , D_2 \in \Div  (\Gamma)$, we say that $D_1$ is \emph{linearly
equivalent to $D_2$}, denoted by $D_1 \sim D_2$, 
if $D_1 - D_2 \in \Prin (\Gamma)$.
Linear equivalence is an equivalence relation.

In studying divisors on graphs, it is often convenient
to consider reduced divisors, which we are recalling now. (The notion of reduced divisors 
was considered in \cite{BN} to prove the Riemann--Roch formula for a finite graph.)  
For any closed subset $A$ of $\Gamma$ and $v \in \Gamma$, the {\em out-degree} of $v$ from $A$, denoted by $\mathrm{outdeg}^\Gamma_A(v)$, is defined to be the maximum number of internally disjoint segments of $\Gamma\setminus A$ with an open end $v$. For $D \in \Div(\Gamma)$, 
a point $v \in \partial A$ is  {\em saturated} for $D$ with respect to $A$ 
if $D(v) \geq \mathrm{outdeg}^\Gamma_A(v)$, 
and {\em non-saturated} otherwise. 
We fix a point $v_0\in\Gamma$.  
A divisor $E \in \Div (\Gamma)$ is said to be 
{\em $v_0$-reduced} if 
$E(x) \geq 0$ for any $x \in \Gamma\setminus\{v_0\}$ and every connected compact  
subset $A$ of $\Gamma\setminus\{v_0\}$ contains a non-saturated point 
$v \in \partial A$ for $E$ with respect to $A$. 
It is known that
for any $D \in \Div ( \Gamma )$,
there exists a unique $v_0$-reduced 
$D_{v_0} \in \Div ( \Gamma )$ that is linearly equivalent to $D$ (cf. \cite{MZ}, 
\cite[Theorem~2.3]{Luo}). 

In this paper, we mainly 
consider divisors supported in $\Gamma_{\Lambda}$.
We set $\Div_{\Lambda}(\Gamma) :=
\{ D \in \Div ( \Gamma) \mid \Supp (D) \subset \Gamma_{\Lambda} \}$.
The elements of $\Div_{\Lambda}(\Gamma)$ are called 
{\em $\Lambda$-divisors} on $\Gamma$. 
For a $\Lambda$-divisor $D$, the complete linear system 
$|D|$ is defined by 
\[
  |D| := \{D^\prime \in \Div_{\Lambda}(\Gamma) \mid D^\prime \geq 0, D^\prime - D \in \Prin (\Gamma)\}. 
\]

\begin{Proposition}
\label{prop:reduced:Lambda-div}
For any $D \in \Div_{\Lambda} ( \Gamma )$
and for any $v_0 \in \Gamma_\Lambda$,
the $v_0$-reduced divisor 
$D_{v_0}$ that is linearly equivalent to $D$ is a $\Lambda$-divisor.
Further,  $D_{v_0} \geq 0$ if and only if $|D| \neq \emptyset$.
\end{Proposition}

\Proof
We take any $D \in \Div_\Lambda(\Gamma)$ and $v_0 \in \Gamma_\Lambda$.
Let $D_{v_0} \in \Div ( \Gamma )$ be the $v_0$-reduced divisor
that is linearly equivalent to $D$.
The first assertion $D_{v_0} \in \Div_{\Lambda} ( \Gamma )$ follows from, for example, 
Luo's algorithm  to construct $v_0$-reduced divisors in \cite{Luo}.

For the second assertion,
it is obvious that if $D_{v_0} \geq 0$, then $|D| \neq \emptyset$.
The other implication follows from
\cite[Corollary~2.18]{Luo}.
\QED

Let $g(\Gamma)$ denote the first Betti number of $\Gamma$. 

\begin{Proposition} [Riemann's inequality on a metric graph] \label{prop:RIforMG}
If $\deg(D) \geq g(\Gamma)$, then $|D| \neq \emptyset$. 
\end{Proposition}

\Proof
Suppose that $\deg(D) \geq g(\Gamma)$.
Take any $v_0 \in \Gamma_{\Lambda}$,
and let $D_{v_0}$ be the $v_0$-reduced divisor with $D_{v_0} \sim D$.
Since $\deg(D) \geq g(\Gamma)$,
the Riemann--Roch formula on metric graphs 
(see \cite[Proposition 3.1]{GK}, \cite[Theorem~7.4]{MZ}, 
or \cite[Theorem~1.2]{HKN})
shows that
there exists an effective divisor
that is linearly equivalent to $D$. 
By \cite[Corollary~2.18]{Luo}, it follows that $D_{v_0} \geq 0$.
Further,
by Proposition~\ref{prop:reduced:Lambda-div}, we have
$D_{v_0} \in \Div_{\Lambda} ( \Gamma )$.
Thus $D_{v_0} \in |D|$, which shows the proposition.
\QED

\begin{Remark}
\label{rmk:on:reduced:divisors}
Let $D \in \Div_\Lambda(\Gamma)$ and $v_0 \in \Gamma_\Lambda$. 
Let $D_{v_0}$ be the $v_0$-reduced divisor that is linearly equivalent to $D$. 
Then $D_{v_0}(v_0) \geq \deg(D) - g(\Gamma)$. Indeed, 
we set $E := D - (\deg(D) - g(\Gamma)[v_0]$. 
Since $D$ is a $v_0$-reduced divisor,
it follows from the definition of $v_0$-reduced divisor
that $E$ is also a $v_0$-reduced divisor.
Since $\deg(E) \geq g(\Gamma)$, 
Propositions~\ref{prop:reduced:Lambda-div}
and \ref{prop:RIforMG}
show that
$E \geq 0$.
Thus $D_{v_0}(v_0) \geq \deg(D) - g(\Gamma)$.
\end{Remark}

\subsection{Weighted $\Lambda$-metric graphs}
\label{subsec:weighted:Lambda:metric:graph}
In this subsection, we show in Proposition~\ref{prop:good:eff:div} the existence of a certain effective divisor on a 
$\Lambda$-metric graph, which will be used to 
construct models $(\Xscr, \Lscr)$ of $(X, L)$. Here we use the theory of {\em weighted} $\Lambda$-metric graphs. 
For a general account of the theory of weighted metric graphs, we refer to 
Amini and Caporaso~\cite{AC}. 

A {\em weighted} $\Lambda$-metric graph $\bar{\Gamma} = (\Gamma, \omega)$ is 
a pair of a $\Lambda$-metric  graph $\Gamma$ and a function $\omega: \Gamma \to \ZZ_{\geq 0}$ such that 
the set $\{x \in \Gamma \mid \omega(x) \neq 0\}$ is finite and contained in $\Gamma_\Lambda$. 
The function $\omega$ is called a weight function. 
The {\em genus} of $\bar{\Gamma}$ is defined to be $g(\bar{\Gamma}) := g(\Gamma) + \sum_{x \in \Gamma_\Lambda} \omega(x)$. 

Let $\bar{\Gamma} = (\Gamma, \omega)$ be a weighted $\Lambda$-metric graph.
Since $g(\bar{\Gamma}) \geq g(\Gamma)$,
the following follows from
Proposition~\ref{prop:RIforMG}.
 
\begin{Proposition}[Riemann's inequality on a weighted $\Lambda$-metric graph]
\label{prop:riemann:inequality}
For $D \in \Div_\Lambda(\Gamma )$,  
if $\deg(D) \geq g(\bar{\Gamma})$, then $|D| \neq \emptyset$. 
\end{Proposition}

Let $\bar{\Gamma} = (\Gamma, \omega)$ be a weighted $\Lambda$-metric graph. 
Assume now that $g(\bar{\Gamma}) \geq 2$ and that 
there does not exist $x \in \Gamma_\Lambda$ with 
$\val(x) = 1$ and $\omega(x) = 0$. 
Then $\Gamma$ has a finite graph structure with the set of vertices
$V(\bar{\Gamma}):= \{x \in \Gamma_\Lambda \mid \val(x) \neq 2\} 
\cup \{x \in \Gamma_\Lambda \mid \omega(x) \neq 0\}$. 
Let $E(\bar{\Gamma})$ denote the set of edges. 

An edge $e \in E(\bar{\Gamma})$ is called an edge of connected type if 
$\Gamma \setminus {\mathrm{relin}(e)}$ is disconnected;
$e$ is called an edge of disconnected type otherwise. 
(An edge of disconnected type is also called a bridge, 
but in this paper, we do not use the terminology 
``bridge.'') 
Let $\{e_1, \ldots, e_r\}$ be the set of 
the edges of disconnected type, and we decompose $\Gamma \setminus \left({\mathrm{relin}(e_1)}\cup\cdots\cup {\mathrm{relin}(e_r)}\right)$ as 
\[
  \Gamma \setminus \left({\mathrm{relin}(e_1)}\cup\cdots\cup {\mathrm{relin}(e_r)}\right)
  = \Gamma_1 \amalg \cdots \amalg \Gamma_{r+1}, 
\] 
where each $\Gamma_i$ is a connected component and thus is a
$\Lambda$-metric graph. 
We set $\omega_i := \rest{\omega}{\Gamma_i}: \Gamma_i \to \ZZ_{\geq 0}$. 
Then $\bar{\Gamma}_i := (\Gamma_i, \omega_i)$ is a weighted $\Lambda$-metric graph, 
and $g(\bar{\Gamma}) = \sum_{i=1}^{r+1} g(\bar{\Gamma}_i)$. 
Since we assume that there does not exist $x \in \Gamma_\Lambda$ such that 
$\val(x) = 1$ and $\omega(x) = 0$, we have $g(\bar{\Gamma}_i) \geq 1$ for each $i$. 

\begin{Definition}[island]
\label{def:island}
We call each $\Gamma_i$ an {\em island} of a weighted $\Lambda$-metric graph $\bar{\Gamma} = (\Gamma, \omega)$. By slight abuse of 
terminology, we also call $\Gamma_i$ an island of $\Gamma$ if there is no confusion 
for the choice of a weight function $\omega$. 
\end{Definition}

\smallskip
The following proposition shows the existence of an effective divisor on $\Gamma$
whose restriction to each island is not trivial. 

\begin{Proposition}
\label{prop:good:eff:div}
Let $\bar{\Gamma} = (\Gamma, \omega)$ be a weighted $\Lambda$-metric graph. 
Assume that $g(\bar{\Gamma}) \geq 2$ and that 
there does not exist $x \in \Gamma_\Lambda$ such that 
$\val(x) = 1$ and $\omega(x) = 0$. 
Let $D \in \Div_{\Lambda}(\Gamma)$. 
If $\deg(D) \geq g(\bar{\Gamma})$, then there exists  $E \in  \Div_{\Lambda}(\Gamma)$ with the following properties\textup{:} 
\begin{enumerate}
\item[(i)]
$E \in |D|$, i.e., $E$ is effective and linearly equivalent to $D$.   
\item[(ii)]
For any island $\Gamma_i$ of $\Gamma$, we have 
$\deg\left(\rest{E}{\Gamma_i}\right) \geq 1$.  
\item[(iii)]
For any edge $e \in E(\bar{\Gamma})$ of disconnected type, 
we have  $\deg\left(\rest{E}{{\mathrm{relin}(e)}}\right) =0$. 
Further, for any edge $e \in E(\bar{\Gamma})$ \textup{(}that is of connected type\textup{)},  
we have  $\deg\left(\rest{E}{{\mathrm{relin}(e)}}\right) \leq 1$.  
\end{enumerate}
\end{Proposition}

\Proof
Suppose that $\deg(D) \geq g(\bar{\Gamma})$.
We construct such 
$E \in \Div_\Lambda(\Gamma)$ step by step. 
Let $\{e_1, \ldots, e_r\} \subset E(\bar{\Gamma})$ be the set of edges of disconnected type. 

\smallskip
{\bf Step 1.}\quad 
 By Riemann's inequality (Proposition~\ref{prop:riemann:inequality}), 
 one has $|D| \neq \emptyset$, and thus there exists 
 an effective divisor $E_1$ with $E_1 \sim D$. 
The divisor $E_1$ satisfies condition (i).

\smallskip
{\bf Step 2.}\quad 
In this step, we replace $E_1$ with another $E_2 \in |D|$ so that 
$E_2$ satisfies conditions (i) and (ii)
and the first condition in (iii). 
We are going to show a stronger statement that 
for any effective divisor $E_1 \in \Div_\Lambda(\Gamma)$ with $\deg(E_1) \geq g(\bar{\Gamma})$, 
there exists an 
effective divisor $E_2 \in \Div_\Lambda(\Gamma)$ such that 
$E_1 \sim E_2$ and $\deg\left(\rest{E_2}{\Gamma_i}\right) \geq g(\bar{\Gamma}_i)$ for any 
island $\Gamma_i$ of $\Gamma$.  
Since any two points on an edge of disconnected type are linearly equivalent to each other, 
we may assume that $\deg\left(\rest{E_1}{{\mathrm{relin}(e_i)}}\right) = 0$ 
for any $1 \leq i \leq r$.  

Suppose that $r = 0$, i.e., $\Gamma$ has no edge of disconnected type. In this case, 
$\Gamma$ is the only island,
and we have $\deg\left(\rest{E_1}{\Gamma}\right) 
= \deg(E_1) \geq g(\bar{\Gamma})$.  Thus it suffices to take $E_2 := E_1$.  

Suppose that $r \geq 1$. Let $v^\prime$ and $ v^{\prime\prime}$ be the end vertices of $e_1$. 
Let $\Gamma^\prime$ and $\Gamma^{\prime\prime}$ be the connected components of 
$\Gamma\setminus{\mathrm{relin}(e_1)}$ with $v^\prime \in \Gamma^\prime$ and 
$v^{\prime\prime} \in \Gamma^{\prime\prime}$. We set $\bar{\Gamma}^\prime = (\Gamma^\prime, \omega^\prime)$
and $\bar{\Gamma}^{\prime\prime} = (\Gamma^{\prime\prime}, \omega^{\prime\prime})$, where 
$\omega^\prime := \rest{\omega}{\Gamma^\prime}$ and 
 $\omega^{\prime\prime} := \rest{\omega}{\Gamma^{\prime\prime}}$. 
Since 
\[
\deg(\rest{E_1}{\Gamma^\prime}) + \deg(\rest{E_1}{\Gamma^{\prime\prime}}) = \deg(E_1) 
\geq g(\bar{\Gamma}) = g(\bar{\Gamma}^\prime) + g(\bar{\Gamma}^{\prime\prime}), 
\]
we may assume without loss of generality that 
$\deg(\rest{E_1}{\Gamma^{\prime\prime}})  \geq  g(\bar{\Gamma}^{\prime\prime})$. 
On $\Gamma'$,
set
\[
 E^{\prime}_1 := \rest{E_1}{\Gamma^{\prime}} 
 + \left(\deg(\rest{E_1}{\Gamma^{\prime\prime}}) -  g(\bar{\Gamma}^{\prime\prime})\right) [v^{\prime}] 
 \in \Div_\Lambda(\Gamma^{\prime})
.
\]
Then $E^{\prime}_1 \geq 0$.
Further, 
$
 \deg(E^{\prime}_1) = \deg(E_1) - g(\bar{\Gamma}^{\prime\prime}) 
 \geq g(\bar{\Gamma}) -  g(\bar{\Gamma}^{\prime\prime})  = g(\bar{\Gamma}^{\prime})
$.
On $\Gamma''$,
since 
\[
\deg
\left( \rest{E_1}{\Gamma^{\prime\prime}} - 
(\deg(\rest{E_1}{\Gamma^{\prime\prime}})  -  g(\bar{\Gamma}^{\prime\prime}))[v''] \right) =
g ( \bar{\Gamma}^{\prime\prime}),
\]
Riemann's inequality 
(cf. Proposition~\ref{prop:riemann:inequality})
gives an effective divisor $E_1^{\prime\prime} \in \Div_{\Lambda} ( \Gamma^{\prime\prime})$
such that 
\[
E_1^{\prime\prime} 
\sim \rest{E_1}{\Gamma^{\prime\prime}}
- (\deg(\rest{E_1}{\Gamma^{\prime\prime}})  -  g(\bar{\Gamma}^{\prime\prime}))[v''].
\]
Since $\Gamma^\prime$ and $\Gamma^{\prime\prime}$ are subgraphs of $\Gamma$, we naturally regard
$E^\prime_1$ and $E^{\prime\prime}_1$ as $\Lambda$-divisors on $\Gamma$. 
Since 
$v'$ and $v''$ are the end vertices of an edge of disconnected type,
we have $[v^{\prime}] \sim [v^{\prime\prime}] $, and hence
$E^{\prime}_1 + E^{\prime\prime}_1 \in |E_1| = |D|$. 
Further,
$
\deg \left( \rest{(E^{\prime}_1 + E^{\prime\prime}_1)}{\Gamma'}\right)
=
 \deg(E^{\prime}_1) \geq g(\bar{\Gamma}^{\prime})
$
and
$
\deg \left( \rest{(E^{\prime}_1 + E^{\prime\prime}_1)}{\Gamma''}\right)
=
 \deg(E^{\prime\prime}_1) = g(\bar{\Gamma}^{\prime\prime})
$.

If $g(\bar{\Gamma}^\prime) =1$ (resp. $g(\bar{\Gamma}^{\prime\prime}) =1$), then 
$\Gamma^{\prime}$ (resp. $\Gamma^{\prime\prime}$) is an island of $\Gamma$, 
and we have $\deg\left(\rest{
E^{\prime}_1 + E^{\prime\prime}_1
}{\Gamma^\prime}\right) \geq g(\Gamma^\prime)$ 
(resp. $\deg\left(\rest{
E^{\prime}_1 + E^{\prime\prime}_1
}{\Gamma^{\prime\prime}}\right) \geq g(\Gamma^{\prime\prime})$).  If $g(\bar{\Gamma}^\prime) \geq 2$ (resp. $g(\bar{\Gamma}^\prime) \geq 2$), 
we repeat this argument to $\bar{\Gamma}^\prime$ and $E^\prime_1$ (resp. 
$\bar{\Gamma}^{\prime\prime}$ and $E^{\prime\prime}_1$). 
Then we obtain $E_2 \in \Div_\Lambda(\Gamma)$ that 
satisfies conditions (i) and (ii).
Further, we may assume that 
$\deg\left(\rest{E_2}{{\mathrm{relin}(e_i)}}\right) = 0$ for any $1 \leq i \leq r$.

\smallskip
{\bf Step 3.}\quad 
In the final step, 
we replace $E_2$ with $E_3$ so that 
$E_3$ satisfies conditions (i)--(iii). 
Suppose that $e$ is an edge that is of connected type and 
that $\deg\left(\rest{E_2}{{\mathrm{relin}(e)}}\right) \geq 2$. 
Let $\gamma : [0,\ell] \to e$ be a parameterization
that restricts to an isometry
$\rest{\gamma}{(0,\ell)} : (0,\ell) \to \mathrm{relin} (e)$.
Then there exist 
$t_1$ and $t_2$ with
$0 < t_1 \leq t_2 < \ell$ such that $E_2 - [\gamma(t_1)] - 
[\gamma(t_2)] \geq 0$. 
We set $m = \min\{t_1, \ell- t_2\}$ and $E_2^\prime:= E_2 - [t_1] - [t_2] + 
[\gamma(t_1 -m)] + [\gamma(t_2 + m)] \geq 0$. 
Then we have 
$\deg\left(\rest{E_2^\prime}{{\mathrm{relin}(e)}}\right) < \deg\left(\rest{E_2}{{\mathrm{relin}(e)}}\right)$ 
and $\deg\left(\rest{E_2}{\Gamma_i}\right) = \deg\left(\rest{E_2^\prime}{\Gamma_i}\right)$ 
for any island $\Gamma_i$, and  $E_2^\prime$ still satisfies conditions (i)(ii) and the first condition 
of (iii).  
Applying this argument for all such $e$ repeatedly, 
we obtain $E_3 \in \Div_\Lambda(\Gamma)$ that 
satisfies conditions (i)--(iii). 
\QED

\subsection{Skeleton as a weighted $\Lambda$-metric graph (with a finite graph structure)}
\label{subsec:skeleton:weighted:metric:graph}

Let $X$ be a connected smooth projective curve over $K$ of genus $g \geq 0$. 
Let $\mathscr{X}$ be a strictly semistable model of $X$.
Then the skeleton $\Gamma := S(\mathscr{X})$ is equipped with 
a natural weight function $\omega$, and thus $(\Gamma, \omega)$ becomes 
a weighted $\Lambda$-metric graph. 
Indeed, we recall that we put 
$V(\Xscr) := \{[C] \mid C \in \Irr(\Xscr_s)\}$,
which is the set of Shilov points with respect to $\mathscr{X}$.
We define a weight function $\omega: \Gamma \to \ZZ_{\geq 0}$ by letting 
$\omega([C])$ be the geometric genus of $C \in \Irr(\Xscr_s)$ and 
letting $\omega(x) = 0$ for any $x \not\in  V(\Xscr)$.  
Then $\bar{\Gamma} = (\Gamma, \omega)$ is a weighted 
$\Lambda$-metric graph. 
One sees that $g (\bar{\Gamma}) = g$.
Further, $g ( \Gamma ) = g$ if and only if $\omega$ is trivial.

\subsubsection*{Minimal weighted skeleton}
Assume that $g \geq 1$. Then 
we have a unique minimal skeleton $\Gamma_{\min}$ of $X^{\an}$ (cf. the paragraph after Definition~\ref{def:minimal:skeleton}).
We take any Deligne--Mumford semistable model 
$\Xscr$ of $X$.
Remark that $\Gamma_{\min} = S ( \mathscr{X} )$,
and
we define a canonical weight function $\omega$ from $\mathscr{X}$ as above.
Note that $\omega$ does not 
depend on the choice of a Deligne--Mumford semistable model $\Xscr$ of $X$. 
Thus we obtain a weighted $\Lambda$-metric graph 
\[
\bar{\Gamma}_{\min} := (\Gamma_{\min} , \omega)
\] 
for $X$.
We call $\bar{\Gamma}_{\min}$ the 
\emph{minimal weighted skeleton of $X^{\an}$}.
Remark that since $\mathscr{X}$ is Deligne--Mumford semistable,
one sees that
there does not exist $x \in \Gamma_{\min , \Lambda}$ such that 
$\val(x) = 1$ and $\omega(x) = 0$.

\subsubsection*{Minimal weighted skeleton with the canonical finite graph structure}
Assume that $g \geq 2$. Recall that we have given the minimal skeleton $\Gamma_{\min}$ with 
the canonical finite graph structure (cf. the second paragraph after Definition~\ref{def:minimal:skeleton})
arising from the stable model:  
the set $V(\Xscr^{\st})$ of vertices is given by the set of Shilov points 
with respect to
the stable model $\Xscr^{\st}$ of $X$, 
and the set $E(\Xscr^{\st})$  of edges is determined by $V(\Xscr^{\st})$. 
On the other hand,
since $\Gamma_{\min}$ has a natural weight function $\omega$,
this
gives a finite graph structure to $\Gamma_{\min}$ as in \S\ref{subsec:weighted:Lambda:metric:graph}:  
the set $V(\bar{\Gamma}_{\min})$ of vertices is 
given by 
\[
  V(\bar{\Gamma}_{\min}):=\{x \in \Gamma_{\min , \Lambda} \mid \val(x) \neq 2\} 
\cup \{x \in \Gamma_{\min , \Lambda} \mid \omega(x) \neq 0\}, 
\]
and the set $E(\bar{\Gamma}_{\min})$  of edges is determined by $V(\bar{\Gamma}_{\min})$. 
Then we have 
\[
  V(\Xscr^{\st}) = V(\bar{\Gamma}_{\min})
  \quad\text{and}\quad
  E(\Xscr^{\st}) = E(\bar{\Gamma}_{\min}). 
\]
This says that 
the finite graph structures on $\Gamma_{\min }$
that are determined by two ways
are in fact the same,
and we have a canonical finite graph structure on $\Gamma_{\min }$.
The minimal weighted skeleta $\bar{\Gamma}_{\min}$ for genus $g \geq 2$ 
with this canonical finite graph structures are important examples on which 
one can apply Proposition~\ref{prop:good:eff:div}. 

\begin{Remark}
\label{remark:valence:stable}
Assume that $g \geq 2$. 
With the notation above, suppose that $g ( \Gamma_{\min} ) = g$. 
Then $\omega$ is trivial. 
Since there does not exist $x \in \Gamma_{\min, \Lambda}$ 
such that $\val(x) = 1$ and $\omega(x) = 0$, it follows that 
$V ( \mathscr{X}^{\st} ) = 
\{ x \in \Gamma_{\min , \Lambda} \mid \val (x) \geq 3  \}$.
\end{Remark}

\subsection{Construction of a model of $(X, L)$}
\label{subsec:starting:point}
Let $X$ be a connected smooth projective curve over $K$.
Let $\mathscr{X}$ be a strictly semistable model of $X$, and let
$\Gamma = S ( \mathscr{X} )$ be the associated skeleton.
Let $\tau : X^{\an} \to \Gamma$ be the retraction map.
Since $X(K) \subset X^{\an}$, 
the restriction of $\tau$ to $X(K)$ gives the map $\tau: X(K) \to \Gamma$ (still denoted by $\tau$). 
Since $K$ is algebraically closed, we have $\Div(X) = \Div(X(K))$. 
By linearity and Lemma~\ref{lemma:retraction:rational}(1), we extend $\tau$ to 
\begin{equation}
\label{eqn:specialization:map}
\tau_*: \Div(X) \to \Div_\Lambda(\Gamma). 
\end{equation}
This map $\tau_*$ is called the {\em specialization map}. Note that 
$\tau_*$ preserves degrees and linearly equivalent classes (cf. \cite[Lemma~2.1]{B})

Let $N$ be a line bundle over $X$. We take a divisor $\tilde{D} \in \Div(X)$ with $N \cong \OO_X(\tilde{D})$, and 
we consider $\tau_*(\tilde{D}) \in \Div_\Lambda(\Gamma)$. Since $\tau_*$ preserves linear equivalent classes, 
the linear system $|\tau_*(\tilde{D})|$ does not depend on the choice of $\tilde{D}$. We denote $|\tau_*(\tilde{D})|$ by $|\tau_*(N)|$. 

For a line bundle $\Nscr$ over $\Xscr$, we set 
\begin{equation}
\label{eqn:def:DM}
  D_{\Nscr} := \sum_{C \in \Irr ( \mathscr{X}_s )}
\deg
(\rest{\mathscr{N}}{C}) [C]
\in \Div_{\Lambda} (\Gamma)
.
\end{equation}
We show the following proposition, 
which will be the starting point to construct 
various good strictly semistable models $(\Xscr, \Lscr)$. 

\begin{Proposition}
\label{prop:pregoodmodel:pre}
Let $X$ be a connected smooth projective curve over $K$
of genus $g \geq 2$. 
Let $N$ be a line bundle over $X$ and let $D \in |\tau_*(N)|$. 
Then there exists a 
 model $(\Xscr, \Nscr)$ of $(X, N)$ with 
the following properties\textup{:}
\begin{enumerate}
\item[(i)]
$\mathscr{X}$ is a 
good model \textup{(}cf. Definition~\textup{\ref{def:goodmodel:new}}\textup{);} 
\item[(ii)]
$D_{\Nscr}  = D$. 
\end{enumerate}
Further, if $\Xscr^0$ is a Deligne--Mumford strictly semistable model of $X$, 
then we may take $\Xscr$ so that $\Xscr$ dominates $\Xscr^0$, i.e., 
there exists a morphism $\mu: \Xscr \to \Xscr^0$ extending the identity morphism on $X$. 
\end{Proposition}

\Proof
{\bf Step 1.}
We take $\tilde{D}^\prime \in \Div(X)$ such that $N \cong \OO_X(\tilde{D}^\prime)$. 
Let $\Gamma_{\min}$ be the minimal skeleton of $X^{\an}$, and
set 
\[
 D^\prime := \tau_{\ast} (\tilde{D}^\prime) \in \Div_\Lambda(\Gamma_{\min}) 
.
\]
Then $D \in |\tau_*(N)| = | \tau_{\ast} (\tilde{D}^\prime)|$. 
We take a
$P \in \Prin_\Lambda(\Gamma_{\min})$ such that
$D  =D^\prime + P$. 
Remark that our $\Prin_\Lambda(\Gamma_{\min})$ coincides with
the group of $\Lambda$-rational principal divisors in \cite{BR}.
Then by \cite[Theorem~1.1]{BR},
which proves Raynaud's type theorem on the surjectivity between 
the groups of
principal divisors under the specialization map, there exists a principal divisor $\tilde{P}$ on $X$ 
such that $\tau_{\ast} (\tilde{P}) = P$. 
We have $D = \tau_*(\tilde{D}^\prime + \tilde{P})$. 
Since $K$ is algebraically closed,
we write $\tilde{D}^\prime + \tilde{P} = \sum_{i= 1}^{m} P_i$ with $P_i \in X(K)$. 

Let $\mathscr{X}^0$ be
a Deligne--Mumford strictly semistable model of $X$. 
Then Proposition~\ref{prop:subdivision1:b} 
gives a strictly semistable model $
\pi^{1} : \mathscr{X}^{1} \to
\Spec (R)$ such that
$V ( \mathscr{X}^{1} ) = V ( \mathscr{X}^0) \cup \{ \tau (P_1) ,
\ldots , \tau (P_m )\}$. 
Since $S(\mathscr{X}^{1}) = S(\mathscr{X}^{0}) = \Gamma_{\min}$, 
$\mathscr{X}^{1}$ is Deligne--Mumford semistable. 

\medskip
{\bf Step 2.}
We modify $\Xscr^1$ to obtain a good model.
Let $p_1 , \ldots , p_\alpha \in \Sing ( \mathscr{X}^1_s)$ be the nodes
such that there does not exist $(-2)$-chain $E$ with $p_i \in E$.
Let $\Delta_{p_i}$ be the canonical $1$-simplex corresponding to $p_i$. 
Take a $\Lambda$-valued points $v_i$ in $\mathrm{relin}( \Delta_{p_i})$
and put $M_1 := \{ v_1 ,  , \ldots , v_\alpha  \}$.
Then by Proposition~\ref{prop:subdivision1:b}, 
there exists a strictly semistable model $\mathscr{X}^2 \to \Spec (R)$ of $X$
such that $V \left( \mathscr{X}^2 \right) =
V \left( \mathscr{X}^1 \right) \cup M_1$.
By the construction, $\mathscr{X}^2$ satisfies condition (ii)
in Definition~\ref{def:goodmodel:new}.
Since $M_1 \subset S(\Xscr^1) = \Gamma_{\min}$, 
$\mathscr{X}^2$ is Deligne--Mumford semistable.
Further,  $\mathscr{X}^2 $ dominates $\mathscr{X}^1$. 

Next, we modify $\mathscr{X}^2$
to make a model 
that satisfies condition (iii) in Definition~\ref{def:goodmodel:new}.
Let $E$ be a maximal $(-2)$-chain of connected type. 
Put $D := \Xscr^2 - E$. 
Let $E_1, \ldots, E_{r-1}$ be the irreducible components of $E$
with a numbering such that $\# (E_i \cap E_{i+1}) = 1$ for $i = 1 , \ldots r-2$.
Let
$p_0 , \ldots , p_{r-1}$ be the distinct nodes such that
$\{ p_0 , \ldots , p_{r-1} \} = \Sing ( \Xscr^2_s ) \cap E
$, where the numbering is given in such a way
that
$\{p_i\} = E_i \cap E_{i+1}$ for $i = 1, \ldots, r-2$, $p_0 \in D \cap E_1$, 
and $p_{r-1} \in D \cap E_{r-1}$. Take $E_0, E_r \in \Irr(D)$ such that 
$p_0 \in E_0$ and $p_{r-1} \in E_r$ (we may have $E_0 = E_r$). 
We set $\Delta_{E} := \bigcup_{i=0}^{r-1} \Delta_{p_i}$
in the skeleton $\Gamma_{\min} = S \left( \mathscr{X}^2 \right)$.
There exists a unique nontrivial isometry $\iota : \Delta_E \to \Delta_E$
such that $\iota ( \partial \Delta_E) = \partial \Delta_E$
and $\iota^2 = \id$,
where $\partial \Delta_E$ is the set of end vertices of $\Delta_E$,
i.e., $\partial \Delta_E = \Delta_E \cap V(\mathscr{X}^{\st})$.
There exists a unique point
$w \in \Delta_E \setminus \partial \Delta_E$ with $\iota (w) = w$.
We can take a finite set 
$M_{E}$
of $\Lambda$-valued points of 
$\Delta_{E} \subset S(\Xscr^2)$ 
such that 
$[E_0], [E_1], \ldots, [E_{r-1}], [E_{r}], w $ are all in  $M_E$
and such that
$\iota (M_E ) = M_{E}$.
Note that $\# \left( M_{E} \setminus \{ [E_0] , [E_r]\} \right)$ is odd,
and we write this number by $2 \ell - 1$.
Furthermore,
adding $\Lambda$-valued points to $M_E$ if necessary,
we may and do assume that $\ell \geq 2$ and 
$M_E$ satisfies the following condition:
Let $w_1 (\neq [E_0]) \in M_E$ be a point that is the nearest to $[E_0]$
and let $w_{\ell -1} (\neq w) \in M_E$ be a point that is the nearest to $w$; 
then the distance between $[E_0]$ and $w_1$ 
equals the distance between $w$ and $w_{\ell - 1}$.

For all maximal $(-2)$-chains $E$ of connected type, 
we take $M_{E}$ as above, and we set $M' := \bigcup_{E} M_E$.
By Proposition~\ref{prop:subdivision1:b},
there exists a strictly semistable model $\mathscr{X} \to \Spec (R)$
such that $V \left( \mathscr{X} \right) = 
V \left( \mathscr{X}^2 \right) \cup M'$ and 
$\mathscr{X}$ dominates $\Xscr^2$. 
Since $M' \subset S(\Xscr^2) = \Gamma_{\min}$, 
$\mathscr{X}$ is Deligne--Mumford semistable. 
By the construction, $\mathscr{X}$ satisfies condition (iii)
in Definition~\ref{def:goodmodel:new}. 
Thus $\Xscr$ is a good model in Definition~\ref{def:goodmodel:new}. 
Further,
$\mathscr{X}$ dominates $\mathscr{X}^2$ and hence dominates $\mathscr{X}^0$.

\medskip

{\bf Step 3.}
Finally, we construct a line bundle $\mathscr{N}$ over $\mathscr{X}$
that satisfies condition (ii) of the proposition.
Let  $P_1, \ldots , P_m \in X(K)$ be as in Step~1. 
We note that $\tau_{\ast} ( P_i ) \in V ( \mathscr{X}^1)
\subset V ( \mathscr{X})$.
Let $\sigma_i: \Spec(R) \to \Xscr$ be the corresponding section to $P_i$
guaranteed by the valuative criterion
of properness.
By Lemma~\ref{lemma:retraction:rational}(2),
there exists,
for any $i = 1 , \ldots , m$,  
a unique $F_i \in \Irr \left( \mathscr{X}_s \right)$ 
such that $\sigma_i (k) \in F_i \setminus \Sing 
\left( \mathscr{X}_s \right)$,
where $\sigma_i (k) = \red_{\mathscr{X}} (P_i)$ by 
the definition of the reduction map $\red_{\mathscr{X}}$.
Since $\sigma_i (k) \notin \Sing \left( \mathscr{X}_s \right)$,
$\sigma_i (R)$ is a Cartier divisor on $\mathscr{X}$.
We define a line bundle $\Nscr$ 
over $\Xscr$
to be
$
\OO_{\Xscr}
\left( \sum_{i=1}^{m} \sigma_i (R)
\right)$.
We have
$\rest{\Nscr}{X} = \OO_X \left(\sum_{i=1}^{m}P_i\right) = \OO_{X} ( \tilde{D}) \cong N$. 
For any $C \in \Irr (\mathscr{X}_s)$
and $i= 1 ,\ldots , m$,
we have
\[
\deg ( \rest{\OO_{\mathscr{X}} (\sigma_i (R) )}{C}) =
\begin{cases}
1
& \text{if $C = F_i$,}
\\
0
&
\text{otherwise.}
\end{cases}
\]
By Lemma~\ref{lemma:retraction:rational},
$\tau_{\ast} (P_i) = [F_i]$.
It follows that for any  $C \in \Irr (\mathscr{X}_s)$,
\[
\deg \left(\rest{\mathscr{N}}{C}\right) = 
\tau_* \left(\tilde{D}^\prime + \tilde{P}\right)([C]) 
= D([C]). 
\]
This shows that $D_{\Nscr} = D$. 
Thus the model $(\Xscr, \Nscr)$ satisfies 
condition (ii).
This 
completes the proof of the proposition.
\QED

Now, via Proposition~\ref{prop:pregoodmodel:pre},
we construct a model $(\Xscr, \Lscr)$ of $(X, L)$, 
which will be frequently used to construct various
global sections of $L$. Recall from the previous subsection that 
we have the minimal weighted skeleton $\bar{\Gamma}_{\min} = (\Gamma_{\min}, \omega)$ of 
$X^{\an}$ for a connected smooth projective curve $X$ of genus $g \geq 2$. 

\begin{Proposition}
\label{prop:pregoodmodel}
Let $X$ be a connected smooth projective curve over $K$ of genus $g \geq 2$,
and let $L$ be a line bundle over $X$. 
Suppose that $\deg (L) \geq 3g-1$. 
Let $x \in \Gamma_{\min}$. Then 
there exist 
a model $(\mathscr{X},\mathscr{L})$
such that
$\mathscr{X}$ is a good model 
\textup{(}cf. Definition~\textup{\ref{def:goodmodel:new}}\textup{)} 
and 
such that
as a divisor on the minimal skeleton $\Gamma_{\min}$,  
\[
D_{\Mscr} - [x] \geq 0,
\]
where 
$\mathscr{M} := \mathscr{L} \otimes \omega_{\mathscr{X}/R}^{\otimes -1}$ 
and $D_{\Mscr}$ is the divisor on $\Gamma_{\min}$ defined in \eqref{eqn:def:DM}. 
Further, we may take $\Lscr$ so that the following properties are also satisfied\textup{:} 
\begin{enumerate}
\item[(i)]
For any island $\Gamma_i$ of $\Gamma_{\min}$, we have 
$\deg\left(\rest{(D_{\Mscr} - [x])}{\Gamma_i}\right) \geq 1$.  
\item[(ii)]
For any edge $e \in E(\Xscr^{\st})$ of $\Gamma_{\min}$, 
we have  $\deg\left(\rest{(D_{\Mscr} - [x])}{{\mathrm{relin}(e)}}\right) \leq 1$.  
\end{enumerate}
Furthermore, 
if $\Xscr^0$ is a Deligne--Mumford strictly semistable model of $X$, 
then we may take $\Xscr$ so that $\Xscr$ dominates $\Xscr^0$.
\end{Proposition}

\Proof
Set $M := L \otimes \omega_X^{\otimes -1}$.
We take $\tilde{D} \in \Div(X)$ with $M \cong \OO_X(\tilde{D})$. 
Since $\deg(M) \geq 3g-1 - (2g-2) = g+1$, 
we have 
\[
\deg \left(\tau_*(\tilde{D}) - 
[x]
\right) \geq 
 g = g(\bar{\Gamma}_{\min}).
\] 
Then 
by Proposition~\ref{prop:good:eff:div}, 
we take $E \in  \left\vert \tau_{\ast} (\tilde{D})  - [x]\right\vert$ 
satisfying conditions (i)--(iii) of 
that proposition.
Now, apply Proposition~\ref{prop:pregoodmodel:pre} to 
$
  E  + [x] \in \left|\tau_{\ast} (\tilde{D})\right| = \left|\tau_{\ast} (M)\right| $,
and we obtain a good model $\Xscr$ of $X$ and a line bundle 
$\Mscr$ with $\rest{\Mscr}{X} \cong M$ such that 
$
  D_{\Mscr} = E  + [x]
$. 
Further,
 if $\Xscr^0$ is a Deligne--Mumford strictly semistable model of $X$, 
then we may take $\Xscr$ so that $\Xscr$ dominates $\Xscr^0$.
Set $\mathscr{L} := \mathscr{M} \otimes \omega_{\mathscr{X}/R}$.
Then, since $E$ satisfies conditions
(ii) and (iii) in Proposition~\ref{prop:good:eff:div},
it is straightforward to see that
$\mathscr{L}$ has the required properties.
\QED

\setcounter{equation}{0}
\section{Unimodular tropicalization of minimal skeleta for $g \geq 2$}
\label{sec:usuful:lemmas}

Let $X$ be a connected smooth projective curve
over $K$ of genus $g$, and let $L$ be a line bundle over $X$.
In this section, 
we prove that  
if $g \geq 2$ and $\deg (L) \geq 3 g - 1$, then 
the minimal skeleton $\Gamma_{\min}$
of $X^{\an}$ has
a unimodular tropicalization
associated to $|L|$. 

\subsection{Useful lemmas}
We begin by showing some lemmas which will be frequently used
not only in this section but also in the following sections. 

We say that a line bundle $L$ is {\em free} at a point $p$ if $p$ is not a base point of $L$. 
We would often like to show that a line bundle is free at a point. 
The following simple lemma serves for this purpose. 
Recall that $k$ denotes the residue field of the ring of integers $R$ of $K$. 

\begin{Lemma} 
\label{lemma:nonbasepoint}
Let $D$ be a 
semistable
curve over $k$, let 
$L$ be a line bundle over $D$, and  let 
$p \in D (k)$. 
We set 
$M := L \otimes \omega_D^{\otimes -1}$. 
\begin{enumerate}
\item
Assume that $p$ is a regular point of $D$. 
If $h^0\left((M(-p))^{\otimes -1}\right) = 0$, then 
$L$ is free at $p$.
\item
Assume that 
$p$ is a node of $D$. 
Let $\nu : \widetilde{D} \to D$ be the partial normalization of $D$ at~$p$. 
If $h^0\left(\nu^*\left(
M^{\otimes -1}\right)\right) = 0$, then 
$L$ is free at $p$. 
\end{enumerate}
\end{Lemma}

\Proof
(1) By the Serre duality, $h^0\left((M(-p))^{\otimes -1}\right) = 0$ is equivalent to $h^{1} ( L ( -p ) ) = 0$. From the exact sequence $0 \to L(-p) \to L \to k_p \to 0$, we then have $h^{0} ( L ( -p ) ) = h^{0} ( L  )  - 1$. Thus $L$ is free at $p$.

(2) 
Let $\frak{m}_p$ be the maximal ideal of $\OO_D$ at $p$. 
We claim that if $\Hom\left(\frak{m}_p, 
 \left(L \otimes \omega_D^{\otimes -1} \right)^{\otimes -1}\right) = 0$, 
then $L$ is free at $p$. Indeed, by the Serre duality, $\Hom\left(\frak{m}_p, 
 \left(L \otimes \omega_D^{\otimes -1} \right)^{\otimes -1}\right) = 0$ is equivalent to 
$h^1(\mfrak_p L) = 0$ (cf. \cite[Theorem~III.7.6]{Ha}). By the same argument as (1) using 
the exact sequence $0 \to \mfrak_p L \to L \to k_p \to 0$,  we obtain the claim.

As is noted in the proof of \cite[Theorem~1.2]{DM},
we have
\begin{equation}
\label{eqn:basepointfreeforpositive-ft:2}
\Hom \left( \mathfrak{m}_p , \left(L \otimes \omega_D^{\otimes -1} \right)^{\otimes -1}\right)
\cong
H^{0}
\left(
\widetilde{D} ,
\nu^*\left((L \otimes \omega_D^{\otimes -1})^{\otimes -1}\right)
\right).
\end{equation}
Thus we get the assertion. 
\QED

In view of Lemma~\ref{lemma:nonbasepoint}, we would 
often like to show 
the vanishing of global sections of 
a certain line bundle. The following lemmas serve for this purpose. 

\begin{Lemma}
\label{lemma:vanishing0}
Let $D$ be a connected proper reduced curve over $k$
and let $M$ be a line bundle over $D$. If $M$ is nef and 
$\deg(M) > 0$, then $h^0(M^{\otimes -1}) = 0$. 
\end{Lemma}

\Proof
Let $C \in \Irr(D)$. Since $\deg(\rest{M}{C}) \geq 0$, we see that 
$H^0(\rest{M^{\otimes -1}}{C}) \neq 0$ if and only if $\rest{M}{C} \cong \OO_C$. 
We take any $\sigma \in H^0(M^{\otimes -1})$. 
Then $\rest{\sigma}{C} \in H^0(\rest{M^{\otimes -1}}{C})$ is either zero or nowhere vanishing. 
On the other hand, there exists $C_1 \in  \Irr(D)$ with $\deg(\rest{M}{C_1}) > 0$, and 
we have $\rest{\sigma}{C_1}  = 0$. Since $D$ is connected, we have $\sigma = 0$. 
\QED

\begin{Lemma}
\label{lemma:vanishing:a}
Let $D$ be  a connected proper reduced curve over $k$ 
and let $M$ be a line bundle over $D$. 
Let $E$ be a curve in $D$
and set $F := D - E$.
Assume that $F \neq \emptyset$, 
$\rest{M}{E}$ is nef and $h^0(\rest{M^{\otimes -1}}{F}) = 0$. 
Then $h^0(M^{\otimes -1}) = 0$. 
\end{Lemma}

\Proof
Our argument is similar to the proof of Lemma~\ref{lemma:vanishing0}. 
We take any $\sigma \in H^0(M^{\otimes -1})$. 
For any $C \in \Irr(E)$, since $\deg(\rest{M}{C}) \geq 0$, 
$\rest{\sigma}{C} \in H^0(\rest{M^{\otimes -1}}{C})$ is either zero or nowhere vanishing. 
On the other hand, $\rest{\sigma}{F} = 0$. Since $D$ is connected, we have $\sigma = 0$. 
\QED

\begin{Lemma}
\label{lemma:vanishing7}
Let $D$ be a proper reduced curve over $k$
with at most nodes as singularities,
and let $M$ be a line bundle over $D$. 
Let $E$ be a curve in $D$
and set $F := D - E$.
Regard $\Sigma := F \cap E$ as a 
Cartier divisor on $F$.
Suppose that $h^0 \left( \rest{M^{\otimes -1}}{E} \right) = 0$
and 
$h^0 \left(
\left( \rest{M}{F} ( \Sigma )
\right)^{\otimes -1}
\right) = 0$.
Then $h^0(M^{\otimes -1}) = 0$. 
\end{Lemma}

\Proof
Take any
$\eta \in H^0 \left( M^{\otimes -1}
\right)$.
Since $h^0 \left( \rest{M^{\otimes -1}}{E} \right) = 0$, we have 
$\rest{\eta}{E} = 0$,
and thus $\rest{\eta}{E \cap F} = 0$.
This means that via the natural inclusion $\rest{M^{\otimes -1}}{F} 
(- \Sigma) \hookrightarrow
\rest{M^{\otimes -1}}{F}$, $\rest{\eta}{F} \in H^0
\left(
M^{\otimes -1} ( - \Sigma )
\right)$.
Since $h^0 \left(
\rest{M^{\otimes -1}}{F} ( - \Sigma)
\right) = 0$, $\rest{\eta}{F} = 0$. 
It follows that $\eta = 0$.
This proves $h^0(M^{\otimes -1}) = 0$.
\QED

Now, let $X$ be a connected smooth projective curve over $K$, 
and let $\Xscr$ be a strictly semistable model of $X$. 
Using the above lemmas,
we give a sufficient condition
in terms of $\Lambda$-metric graphs
for a line bundle over $\Xscr$ 
to be free at a node $p \in \Xscr_s$.  

\begin{Lemma} \label{lemma:freeatnode:graphversion:emptyset}
Let $X$ be a connected smooth projective curve over $K$, 
and let $\Xscr$ be a strictly semistable model of $X$.  
Let $\Gamma := S(\Xscr)$ be the associated skeleton. 
Let $\mathscr{L}$ be a line bundle over $\Xscr$, and we set 
$\Mscr := \Lscr \otimes \omega_{\Xscr/R}^{\otimes -1}$. 
Let $p \in \Sing ( \mathscr{X}_s)$, and denote by $\Delta_p$ the 
corresponding $1$-simplex in $\Gamma$ \textup{(}cf. Definition~\textup{\ref{def:canonical:1:simplex:node}}\textup{)}. 
Assume that that $D_{\mathscr{M}} \geq 0$ \textup{(}cf. \eqref{eqn:def:DM}\textup{)}. 
Suppose that for any connected component $\Gamma^\prime$ of 
$\Gamma \setminus \mathrm{relin}(\Delta_p)$, we have 
$\deg
\left(
\rest{D_{\mathscr{M}}}{\Gamma^\prime} \right) \geq 1$. 
Then $\mathscr{L}$ is free at $p$.
\end{Lemma}

\Proof
To ease notation, we set $M := \rest{\Mscr}{\Xscr_s}$. 
Note that $D_{\mathscr{M}} \geq 0$ means that  $M$ is nef. 
Let $\nu : \widetilde{\mathscr{X}_s} \to \mathscr{X}_s$
be the partial normalization at $p$. 
Let $F$ be any connected component of $\widetilde{\mathscr{X}_s}$. 
Then by the assumptions, $\rest{\nu^*\left(M\right)}{F}$ is nef and 
$\deg\left(\rest{\nu^*\left(M\right)}{F}\right) >0$. 
By Lemma~\ref{lemma:vanishing0}, 
we have
$h^{0} 
\left( 
\nu^*\left(M^{\otimes -1}\right)
\right) = 0$.
By Lemma~\ref{lemma:nonbasepoint}(2),
$\rest{\mathscr{L}}{\mathscr{X}_s}$
is free at $p$.

Since $H^0 \left(
M^{\otimes -1} 
\right) \subset H^0
\left( 
\nu^\ast \left(
M^{\otimes -1}
\right)
\right)$, we have 
$h^{0} 
\left(
M^{\otimes -1}
\right) = 0$. 
By the Serre duality,
$h^1
\left(
\rest{\mathscr{L}}{\mathscr{X}_s}
\right) = 0$.
By the base-change theorem,
it follows that the natural homomorphism
$H^0 ( \mathscr{L}) \to H^{0}
\left(
\rest{\mathscr{L}}{\mathscr{X}_s}
\right)$ is surjective.
Since $\rest{\mathscr{L}}{\mathscr{X}_s}$
is free at $p$,
this concludes that
$\mathscr{L}$ is free at $p$.
\QED

In the next lemma, we consider a case where $D_{\Mscr}$ is not effective. 

\begin{Lemma} \label{lemma:freeatnode:graphversion:nonemptyset}
Let $X$, $\mathscr{X}$, $\Gamma$, $\mathscr{L}$, $\mathscr{M}$, $p$,
and $\Delta_p$ be as in Lemma~\ref{lemma:freeatnode:graphversion:emptyset}.
Let $\Gamma_1$ be a connected component of $\Gamma \setminus \mathrm{relin}(\Delta_p)$.
Set 
$
V := \{
v \in V(\Xscr) \mid 
D_{\mathscr{M}}(v) <0\}
$.
Assume that $V \neq \emptyset$ 
and 
$V \subset \Gamma_1$. 
Suppose the following. 
\begin{enumerate}
\item[(i)]
For any connected component $\Gamma^\prime$ of 
$\Gamma \setminus \mathrm{relin}(\Delta_p)$ with $\Gamma^\prime \neq \Gamma_1$, 
$\deg
\left(\rest{D_{\mathscr{M}}}{\Gamma^\prime} \right) \geq 1$. 
\item[(ii)]
There exist $s \geq 1$ and distinct connected components $\Gamma_{11}^\circ, \ldots, \Gamma_{1s}^\circ$ 
of $\Gamma_1 \setminus V$ such that 
\begin{enumerate}
\item[(a)]
$\deg\left(\rest{D_{\mathscr{M}}}{\Gamma_{1j}^\circ} \right) \geq 1$ for any $1 \leq j \leq s$, and 
\item[(b)]
for any $v \in V$, 
the valence of $\Gamma_{11} \cup \cdots \cup \Gamma_{1s}$ at $v$ is at least $-D_{\Mscr}(v) + 1$, 
where $\Gamma_{1j}$ is the closure of $\Gamma_{1j}^\circ$ in $\Gamma_1$. 
\end{enumerate}
\end{enumerate}
Then $\mathscr{L}$ is free at $p$.
\end{Lemma}

\Proof
Let $\nu : \widetilde{\mathscr{X}_s} \to \mathscr{X}_s$
be the partial normalization at $p$. 
Via $\nu$ we identify $\Irr \left( \widetilde{\mathscr{X}_s} \right)$
with $\Irr \left( \mathscr{X}_s \right)$. 
For $C \in \Irr \left( \widetilde{\mathscr{X}_s} \right)$, 
let $[C]$ denote the corresponding point in $\Gamma \setminus {\rm relin}(\Delta_p)$. 
To ease notation, we set $M := \rest{\Mscr}{\Xscr_s}$. 

\begin{Claim}
\label{claim:in:lemma:freeatnode:graphversion:nonemptyset}
We have $h^{0}\left(\nu^*(M)^{\otimes -1}\right) = 0$.
\end{Claim}

Indeed, we set $F_1 := 
\sum_{C \in \Irr \left(\widetilde{\mathscr{X}_s} \right) , [C] \in \Gamma_1} C$ 
and set
$F_{1j} := \sum_{C \in \Irr \left(\widetilde{\mathscr{X}_s} \right) , [C] \in \Gamma^{\circ}_{1j}} 
C$ for $1 \leq j \leq s$.  
We take $C_1 , \ldots , C_{r} \in \widetilde{\mathscr{X}_s}$ such that 
$V' = \{[C_1], \ldots, [C_r]\}$. By the assumption,
$F_1$ is the connected component of $\Xscr_s \setminus \{p\}$ such that 
$\{C_1 , \ldots , C_{r}\} \subset F_1$. Further, by 
(ii)(a)
in the lemma, 
each $F_{1j}$ is a connected component 
of $F_1 - (C_1 + \cdots + C_{r})$, 
$\rest{\nu^*(M)}{F_{1j}}$ is nef, and 
$\deg\left(\rest{\nu^*(M)}{F_{1j}}\right) > 0$.  
Thus Lemma~\ref{lemma:vanishing0}  gives 
\begin{equation}
\label{eqn:freeatnode:graphversion:nonemptyset}
h^0\left( 
\rest{\nu^*(M)^{\otimes -1}}{\sum_{j=1}^s F_{1j}}
\right) = 0.
\end{equation}

Take any $i = 1 , \ldots , r$. 
Set $\Sigma_i := C_i \cap (\sum_{j=1}^s F_{1j})$,
which we regard as a Cartier divisor on $C_i$.
Then, by condition (ii)(b) in the lemma,
$\deg \left(
\rest{\nu^*(M)}{C_i} (\Sigma_i) 
\right)
\geq 1$. 
Since $C_i$ is irreducible, 
Lemma~\ref{lemma:vanishing0} tells us that 
\begin{equation}
\label{eqn:freeatnode:graphversion:nonemptyset:2}
h^0
\left( 
\left(
\rest{\nu^*(M)}{C_i} (\Sigma_i) 
\right)^{\otimes -1}
\right) = 0.
\end{equation}

Note
that 
$h^{0}
\left(
\rest{\nu^*(M)^{\otimes -1}}{F_1}
\right) = 0$.
Indeed,
we
take any 
$\sigma \in H^{0}
\left(\rest{\nu^*(M)^{\otimes -1}}{\sum_{i=1}^r C_i + \sum_{j=1}^s F_{1j}}
\right)$.
By \eqref{eqn:freeatnode:graphversion:nonemptyset} 
and \eqref{eqn:freeatnode:graphversion:nonemptyset:2}, 
Lemma~\ref{lemma:vanishing7} gives $h^{0}
\left(
\rest{\nu^*(M)^{\otimes -1}}{C_i + \sum_{j=1}^s F_{1j}}
\right) = 0$.
Thus 
$\rest{\sigma}{C_i + \sum_{j=1}^s F_{1j}} = 0$
for each $i=1, \ldots , r$.
Since the union of $\left\{ C_i + \sum_{j=1}^s F_{1j} \right\}_{i=1}^r$
equals $\sum_{i=1}^r C_i + \sum_{j=1}^s F_{1j}$,
this proves $\sigma = 0$.
Thus $h^0
\left( 
\rest{\nu^*(M)^{\otimes -1}}{\sum_{i=1}^r C_i + \sum_{j=1}^s F_{1j}}
\right) = 0$.
Since $\rest{\nu^*(M)}{F_1 - (\sum_{i=1}^r C_i + \sum_{j=1}^s F_{1j})}$
is nef,
Lemma~\ref{lemma:vanishing:a} concludes that
$h^{0}
\left(
\rest{\nu^*(M)^{\otimes -1}}{F_1}
\right) = 0$.

Let $\Gamma^\prime$ be any connected component  of 
$\Gamma \setminus \mathrm{relin}(\Delta_p)$ with $\Gamma^\prime \neq \Gamma_1$, 
and we set $
F^\prime := 
\sum_{C \in \Irr \left(\widetilde{\mathscr{X}_s} \right) , [C] \in \Gamma^\prime} C$. 
It follows from $V \subset \Gamma_1$ that $\rest{\nu^*(M)}{F^\prime}$ is nef. 
It also follows from condition (i) 
in the lemma
that  
$\deg(\rest{\nu^*(M)}{F^\prime}) \geq 1$. 
By Lemma~\ref{lemma:vanishing0}, 
we have $h^{0}
\left(
\rest{\nu^*(M)}{F^\prime}^{\otimes -1}
\right) = 0$. 
Thus
$h^{0}\left(
\nu^*(M)^{\otimes -1}\right) = 0$.
This completes the proof of
Claim~\ref{claim:in:lemma:freeatnode:graphversion:nonemptyset}.

By Claim~\ref{claim:in:lemma:freeatnode:graphversion:nonemptyset} and Lemma~\ref{lemma:nonbasepoint},
it follows that
$\rest{\Lscr}{\Xscr_s}$ is free at $p$. 
Then the same argument as in the proof of Lemma~\ref{lemma:freeatnode:graphversion:emptyset}
concludes that $\Lscr$ is free at $p$. 
\QED

While
Lemmas~\ref{lemma:freeatnode:graphversion:emptyset}
and \ref{lemma:freeatnode:graphversion:nonemptyset} consider
freeness of $\Lscr$ at a node,
the following two lemmas
consider freeness of $\Lscr$ 
at a regular point. 
The reader may skip these two lemmas 
and return when reading \S\ref{section:FTcan},
because
we will 
use them not in this section
but
in \S\ref{section:FTcan}.

\begin{Lemma}
\label{lemma:freeatregular:graphversion:effective}
Let $X$ be a connected smooth projective curve over $K$, 
and let $\Xscr$ be a strictly semistable model of $X$.  
Let $\Gamma := S(\Xscr)$ be the associated skeleton. 
Let $\mathscr{L}$ be a line bundle over $\Xscr$, and we set 
$\Mscr := \Lscr \otimes \omega_{\Xscr/R}^{\otimes -1}$. 
Let $C$ be an irreducible component of $\mathscr{X}_s$, 
and let $v := [C] \in V ( \mathscr{X} )$ be the Shilov point 
in $\Gamma$ corresponding to $C$. 
Suppose that 
$D_{\Mscr} - [v]$ is effective on $\Gamma$ and $\deg(D_{\Mscr} - [v]) \geq 1$ .
Then $\mathscr{L}$ is free at any $p \in C(k)  
\setminus \Sing ( \mathscr{X}_s)$.
\end{Lemma}

\Proof
We take any $p \in C(k)  
\setminus \Sing ( \mathscr{X}_s)$. 
To ease notation, we set $M := \rest{\mathscr{M}}{\mathscr{X}_s}$.
Since $D_{\Mscr} - [v]$ is effective  on $\Gamma$ and has positive degree,
$M (-p)$
is nef and has positive degree. 
By Lemma~\ref{lemma:vanishing0} we have 
$h^0
\left( 
M (-p)^{\otimes -1}
\right) = 0$, and then Lemma~\ref{lemma:nonbasepoint}(1) implies that 
$\rest{\mathscr{L}}{\mathscr{X}_s}$ is free at $p$.
Since $h^0
\left( 
M (-p)^{\otimes -1}
\right) = 0$, we have 
$h^0
\left( 
M^{\otimes -1}
\right) = 0$. By the Serre duality, we have 
$
h^1
\left( 
\rest{\mathscr{L}}{\mathscr{X}_s}
\right)
=
h^0
\left( 
M^{\otimes -1}
\right) = 0$. 
By the base-change theorem, 
$H^0 ( \mathscr{L} ) \to H^0 \left( 
\rest{\mathscr{L}}{\mathscr{X}_s} \right)$
is surjective. 
Thus $\mathscr{L}$ is free at $p$.
\QED

In the next lemma, we consider a case where $D_\Mscr - [v]$ is not effective. 

\begin{Lemma} \label{lemma:freeatregular:graphversion}
Let $X$, $\Xscr$, $\Gamma$, $\mathscr{L}$, $\Mscr$, $C$, and $v$ 
be as in Lemma~\textup{\ref{lemma:freeatregular:graphversion:effective}}.
Suppose that 
$D_{\mathscr{M}}$ is effective
and that
$D_{\Mscr} (v) = 0$.
Further, suppose that
there exist connected components $\Gamma_1^\circ$ 
and $\Gamma_2^\circ$ 
of $\Gamma \setminus \{ v \}$
\textup{(}possibly $\Gamma_1^\circ = \Gamma_2^\circ$\textup{)}
such that
\begin{enumerate}
\item[(i)]
$\deg (\rest{D_{\mathscr{M}}}{\Gamma_j^\circ}) \geq 1$
for any $j = 1 , 2$\textup{;}
\item[(ii)]
the valence of 
$\Gamma_1 \cup \Gamma_2$ at $v$
is at least $2$, where $\Gamma_j:= \Gamma_j^\circ \cup \{v\}$ is the closure 
of $ \Gamma_j^\circ$ in~$\Gamma$. 
\end{enumerate}
Then $\mathscr{L}$ is free at any $p \in C(k)  \setminus \Sing ( \mathscr{X}_s)$.
\end{Lemma}

\Proof
To ease notation, we set $M := \rest{\mathscr{M}}{\mathscr{X}_s}$.
We take any 
$p \in C(k)
\setminus \Sing ( \mathscr{X}_s)$.
We prove that
$h^{0}
\left(
M(-p)^{\otimes -1}
\right) = 0$.
Put 
$F_j := \sum_{C^\prime \in \Irr \left( 
\mathscr{X}_s \right) , [C^\prime] \in \Gamma_j^\circ} C^\prime$
for $j=1,2$.
Then 
each $F_j$ ($j=1,2$)
is a connected component of 
$\mathscr{X}_s -  C_1
$.
By condition (i), 
$\rest{M (-p)}{F_j}$ is nef and has positive degree for $j = 
1 ,2$.
By Lemma~\ref{lemma:vanishing0},
we have 
$h^0
\left( 
\rest{M(-p)}{F_j}^{\otimes -1}
\right) = 0$. 
We set 
\[
F_{3} := 
\begin{cases}
F_1
& \text{if $F_1 = F_2$} 
\\
F_1 + F_2
&
\text{if $F_1 \neq F_2$.} 
\end{cases}
\]
Then we have $h^0(\rest{M(-p)}{F_3}^{\otimes -1}) = 0$. 
We set
$\Sigma :=
C \cap F_3$,
which we regard as a Cartier divisor on $C$.
Then by condition (ii)
and $D_{\mathscr{M}} (v) = 0$, we have 
$\deg \left(
\rest{M(-p)}{C} ( \Sigma ) 
\right)
\geq 1$, 
and then Lemma~\ref{lemma:vanishing0} gives 
$h^0
\left( 
\left(
\rest{M(-p)}{C}  ( \Sigma ) 
\right)^{\otimes -1}
\right) = 0$.
By Lemma~\ref{lemma:vanishing7},
we have $h^{0}
\left(
\rest{M(-p)}{C + F_3}^{\otimes -1}
\right) = 0$.
Since $D_{\mathscr{M}}$ is effective,
$\rest{M(-p)}{\mathscr{X}_s - (C + F_3)}$
is nef. Thus 
Lemma~\ref{lemma:vanishing:a} gives 
$h^{0}
\left(
M(-p)^{\otimes -1}
\right) = 0$.

We have shown that $h^{0}
\left(M(-p)^{\otimes -1}\right) = 0$. Then 
Lemma~\ref{lemma:nonbasepoint}(1) implies that 
$\rest{\mathscr{L}}{\mathscr{X}_s}$
is free at $p$. By the same argument as in the proof 
of Lemma~\ref{lemma:freeatregular:graphversion:effective}, 
we obtain that $\mathscr{L}$ is free at $p$.
\QED

The following remark  
will not be used in this section or in \S\ref{section:FTcan}
but
in \S\ref{subsec:good:model:classical:point}.

\begin{Remark}
\label{remark:proofs:lemma:freeatnode:graphversion:emptyset}
In the proofs of Lemma~\ref{lemma:freeatnode:graphversion:emptyset} and Lemma~\ref{lemma:freeatregular:graphversion}, 
we have actually shown the following statements. 
\begin{enumerate}
\item
Let $X, \mathscr{X}, \Gamma$, $p$ and $\Delta_p$ 
be as in Lemma~\ref{lemma:freeatnode:graphversion:emptyset}.
Let $M$ be a line bundle over $\mathscr{X}_s$, and we set 
$D_{M} := \sum_{C \in \Irr ( \mathscr{X}_s )}
\deg (\rest{M}{C}) [C]
\in \Div_{\Lambda} (\Gamma)$.
Assume that $D_M \geq 0$ 
and that for any connected component $\Gamma'$ of $\Gamma \setminus 
\mathrm{relin} ( \Delta_p )$,  we have $\deg ( \rest{D_M}{\Gamma'}) \geq 1$.
Then $h^{0}\left(
\nu^*(M)^{\otimes -1}\right) = 0$, where $\nu$ is the partial normalization of $\Xscr_s$ at $p$. 
\item
Let $X$, $\Xscr$, $\Gamma$, $C$, and $v$ 
be as in Lemma~\ref{lemma:freeatregular:graphversion:effective}.
Let $M$ be a line bundle over $\mathscr{X}_s$.
Set $D_{M} 
:= \sum_{C \in \Irr ( \mathscr{X}_s )}
\deg (\rest{M}{C}) [C]
\in \Div_{\Lambda} (\Gamma)$.
Suppose that 
$D_{M}$ is effective
and
$D_{M} (v) = 0$.
Further, suppose that there exist connected components $\Gamma_1^\circ$ 
and $\Gamma_2^\circ$ 
of $\Gamma \setminus \{ v \}$
\textup{(}possibly $\Gamma_1^\circ = \Gamma_2^\circ$\textup{)}
such that conditions (i) and (ii) of Lemma~\textup{\ref{lemma:freeatregular:graphversion}} 
are satisfied with $D_{\Mscr}$ replaced by $D_M$. Then 
$h^{0}
\left(
M(-p)^{\otimes -1}
\right) = 0$
for any $p \in C(k) \setminus \Sing ( \mathscr{X}_s )$.
\end{enumerate}
\end{Remark}

\subsection{Fundamental vertical divisors}
\label{subsec:fvd}

Let  $X$ be a connected smooth projective curve over $K$,
and let $\Xscr \to \Spec (R)$ be a strictly semistable model of $X$. 
Since $\mathscr{X}$ is a local-ringed space,
we have the notion of Cartier divisors on $\mathscr{X}$.
In this paper, divisors on $\mathscr{X}$ always mean
Cartier divisors. Let $\Div ( \mathscr{X})$ denote the group of divisors
on $\mathscr{X}$.

Let $C \in \Irr (\mathscr{X}_s)$,
and
let $|\cdot|_C$ denote the absolute value on 
the function field
$K (X)$
that gives the Shilov point associated to $C$.
For any $f \in K (X)^{\times}$, $\ord_C (f) := - \log |f|_C$
is called the \emph{order of $f$ along $C$.}
Let $\Dscr$ be a divisor of $\Xscr$. 
Let $\xi \in \mathscr{X}$ be the generic point of $C$, and
let $g$ be a local equation of $\Dscr$ at $\xi$,
i.e., 
a rational function on $\Xscr$ 
such that $\zero (g) = \Dscr$
on some open neighborhood of $\xi$ in $\mathscr{X}$.
We define 
\emph{order of $\mathscr{D}$ along $C$} to be $\ord_C (\Dscr) := - \log |g|_C$. 
Remark that $\ord_C (\Dscr)$ does not depend on the choice of $g$ and that 
$\ord_C (\Dscr) \in \Lambda$ (cf. Remark~\ref{rem:type2points}).

For a divisor $\Dscr$, its support $\Supp (\Dscr)$
is the locus on $\mathscr{X}$ at which the local equation
of $\Dscr$ is not a locally invertible regular function.
We remark that $\Supp (\Dscr)$ is a closed subset of $\mathscr{X}$.
A divisor $\Dscr$ on $\mathscr{X}$ is said to be \emph{vertical} if 
the restriction of $\Dscr$ to $X$ is trivial, 
which is equivalent to $\Supp (\Dscr) \subset \mathscr{X}_s$.

We recall from Definition~\ref{def:multiplicityatanodes}
the multiplicity at a node.
Let $F$ be a connected curve in $\mathscr{X}_s$. 
For each $p \in F \cap ( \mathscr{X}_s - F)$,
take
a local \'etale atlas
\begin{equation}
\label{eqn:local:atlas:psi:p}
\psi_p :
\mathscr{U}_p \to
\Spec ( R[x, y]/(xy - \varpi_p))
\end{equation}
that
distinguishes $p$,
where $\varpi_p \in R$. 
Then
$v_K ( \varpi_p)$ is called the multiplicity 
of $\Xscr$ at $p$.

Assume that 
$F$ satisfies the following conditions (F):
\begin{itemize}
\item[(F)]
There exists a $\lambda \in \Lambda$
such that for any $p \in F \cap ( \mathscr{X}_s - F)$,
the multiplicity of $\mathscr{X}$ at $p$ equals $\lambda$.
\end{itemize}
Then
we define a vertical divisor $\Fscr$ as follows. 
For each $p \in F \cap ( \mathscr{X}_s - F)$, we take 
a local \'etale atlas \eqref{eqn:local:atlas:psi:p} such that   
$\psi_{p}^{\ast} (x)$ defines
the branch $C$ at $p$ with $C \subset F$.
Then there exists a unique divisor $\Fscr$
characterized by the properties that 
\begin{enumerate}
\item[(i)]
for any $p \in F \cap ( \mathscr{X}_s - F)$,
$\Fscr$ has a local equation $\psi_p^{\ast} (x)$ on $\mathscr{U}_p$;
\item[(ii)]
$\Fscr = \zero ( \varpi_p )$
on some open neighborhood of
$F \setminus ( \mathscr{X}_s - F)$;
\item[(iii)]
$\Fscr$ is trivial on some open neighborhood of $( \mathscr{X}_s - F) \setminus F$.
\end{enumerate}
One sees that 
$\Fscr$ does not dependent on the choice of a local atlas 
$\psi_p$
at
$p$ nor of $\varpi_p$.
Further,  
$\Fscr$ has support $F$.
We call $\Fscr$ the 
\emph{fundamental vertical divisor with support $F$}.
The number $\lambda$ is called the \emph{multiplicity of 
$\Fscr$}.
The fundamental vertical divisor with support 
$F$ for 
some connected curve $F$
is simply called a \emph{fundamental vertical divisor}.

\begin{Remark} \label{rem:restrictionofcartierdivisor}
Let $\Fscr$ be a fundamental vertical divisor on $\mathscr{X}$.
Let $D$ be a connected curve with $\Irr (D) \cap \Irr ( \Supp (\Fscr)) = \emptyset$.
Then $\Fscr|_{D}$ is a reduced Cartier divisor on $D$ with
support $\Supp ( \Fscr ) \cap D$.
This can be checked by using a local equation of $\Fscr$
at each  point of $\Supp ( \Fscr ) \cap D$.
\end{Remark}

Let $X$ be a connected smooth projective curve over $K$, and let
$\pi : \mathscr{X} \to \Spec (R)$ be a strictly semistable
model of $X$.

\begin{Lemma} \label{lem:forinjectivity}
Let $X$ and $\pi : \mathscr{X} \to \Spec (R)$ be as above.
Let $g$ be a non-zero rational function on $\mathscr{X}$, 
and let $h:= \rest{g}{X}$ denote the restriction of $g$ to the generic fiber $X$. 
Let $\mathscr{W}
\in \Div (\mathscr{X})
$ be
a $\ZZ$-linear combination of fundamental vertical divisors.
Then the following hold,
where
for any $C \in \Irr (\mathscr{X}_s)$,
$[C]$ denotes the Shilov point associated to $C$.
\begin{enumerate}
\item
Let $C \in \Irr ( \mathscr{X}_s)$.
Suppose that there exists an open neighborhood $\mathscr{U} \subset 
\mathscr{X}$
of the generic point $\xi_C$ such that
$
\rest{
\left(
\zero (g) - \mathscr{W}
\right)}{\mathscr{U}}
$
is effective.
Then
$- \log |h ([C])| \geq \ord_C (\mathscr{W})$. Furthermore, 
if $\rest{
\left(
\zero (g) - \mathscr{W}
\right)}{\mathscr{U}}$ is trivial, 
then $- \log |h ( [C] )| = \ord_C (\mathscr{W})$.
\item
Let $p \in \Sing ( \mathscr{X}_s )$,
and let $C_1, C_2 \in \Irr (\mathscr{X}_s)$ be 
distinct irreducible components with $p \in C_{1} \cap C_{2}$.
Let $\varphi$ be the affine function on $\Delta_p$
with $\varphi ([C_i]) = \ord_{C_i} (\mathscr{W})$ for $i=1,2$.
Suppose that there exists 
an open neighborhood $\mathscr{U} \subset 
\mathscr{X}$
of $p$ such that
$
\rest{
\left(
\zero (g) - \mathscr{W}
\right)}{\mathscr{U}}
$
is effective.
Then
$- \log |h (u) | \geq \varphi (u)$
for any $u \in \Delta_p$.
Furthermore, if $\rest{
\left(
\zero (g) - \mathscr{W}
\right)}{\mathscr{U}}$ is trivial,
then $- \log |h (u) | = \varphi (u)$ for any $u \in \Delta_p$. 
\end{enumerate}
\end{Lemma}

\Proof
Assertion (1) follows from the definition of $\ord$.
Let us prove (2).
We take a distinguished local \'etale morphism
$\psi_p : \mathscr{U}_p \to \Spec (R [x,y]/(xy - \varpi_p))$
such that $\mathscr{U}_p \subset \mathscr{U}$.
We may assume that, in the special fiber,
$C_1$ is defined by $\psi_p^{\ast} (x)$
and $C_2$  by $\psi_p^{\ast} (y)$.
Since $\mathscr{W}$ is a $\ZZ$-linear combination of fundamental vertical
divisors,
we have a local
equation of $\mathscr{W}$ of the form
$\varpi \psi_p^{\ast} (x)^{m_1} \psi_p^{\ast} (y)^{m_2}$
for some $m_1 , m_2 \in \ZZ$
and $\varpi \in K^{\times}$,
and then 
$\varphi = - \log \left|\rest{\varpi \psi_p^{\ast} (x)^{m_1} \psi_p^{\ast} (y)^{m_2}}{X}\right|$.
Put $g_1:= g \cdot 
\left(\varpi \psi_p^{\ast} (x)^{m_1} \psi_p^{\ast} (y)^{m_2} \right)^{-1}$.
Since
$
\rest{
\left(
\zero (g) - \mathscr{W}
\right)}{\mathscr{U}}
$
is effective,
$g_1$
is regular on $\mathscr{U}_p$
and thus $- \log \left|\rest{g_1}{X}\right| \geq 0$ on $\Delta_p$.
Thus we obtain $- \log  \left|\rest{g}{X}\right| \geq \varphi$.
If $
\rest{
\left(
\zero (g) - \mathscr{W}
\right)}{\mathscr{U}}
$
is trivial on some open neighborhood of $p$, then 
$g_1$ is an invertible regular function around $p$,
and hence $- \log \left|\rest{g_1}{X}\right| = 0$. 
In this case, we have $- \log \left|\rest{g}{X}\right| = \varphi$.
\QED

\subsection{Stepwise vertical divisors}
\label{subsection:stepwise}

Let $X$ be a connected smooth projective curve over $K$ of genus $g$. 
In this subsection, we assume that $g \geq 2$. 
Let $\mathscr{X}$ be a strictly semistable model of $X$.

In this subsection, we define the stepwise vertical divisors on $\Xscr$ 
associated to an oriented edge of disconnected type
and those associated to an edge of connected type.
These divisors is defined as the sum of 
suitable fundamental vertical divisors. 
They will be mainly used to construct a tropicalization
map that is unimodular over the edge.

\medskip
Let $e$ be an edge of $\Gamma_{\min}$, i.e., 
$e \in E ( \mathscr{X}^{\st} )$, where $\mathscr{X}^{\st}$ is the stable model of $X$; 
see the last paragraph in \S\ref{subsec:properties:skeleta}.
Let $E$ be the corresponding maximal $(-2)$-chain in $\mathscr{X}_s$
(cf. Lemma~\ref{lemma:make:clear}).

\subsubsection*{Stepwise vertical divisor for an edge of disconnected type}
We consider the case where $e \in E ( \mathscr{X}^{\st} )$ is of disconnected type.

Since we use orientations on $e$, 
we first fix the notation.
Let $\lambda_e$ be the length of $e$.
Then
there are exactly two isometries from
$[0 , \lambda_e]$ to  $e$,
and an \emph{orientation} means one 
of the choices 
of an isometry. 
Once an orientation on $e$ is fixed,
we usually denote by $\vec{e}$ the edge 
$e$
with the orientation,
and we call $\vec{e}$ 
an \emph{oriented edge}.
For an oriented edge $\vec{e}$ with orientation 
$\gamma : [0 , \lambda_e]
\to e$, we call $\gamma (0)$ the \emph{head of $\vec{e}$}
and $\gamma (\lambda_e )$ the \emph{tail of $\vec{e}$}.
We remark that giving an orientation on $e$ is equivalent to
specifying an end vertex of $e$ to be the head.

\begin{figure}[!h]
\[
\setlength\unitlength{0.08truecm}
\begin{picture}(100, 65)(0,0)
  \qbezier(0, 5)(12.5, 15)(25, 5)
  \put(-20,7){\line(1,0){25}}
  \put(10,7){\line(3,4){17}}
  \put(15,20){\line(1,0){25}}
  \put(30,17){\line(3,4){14}}
  \put(55,28){\line(1,0){25}}
  \put(65,25){\line(3,4){17}}
  \qbezier(70, 35)(82.5, 45)(95, 35)
   \put(85,38){\line(1,0){25}}
  \put(45, 27){$\cdots$}
  \put(0, 10){$E_0$}
  \put(40, 18){$E_2$}
  \put(18, 27){$E_1$}
  \put(70, 47){$E_{r-1}$}
  \put(86, 42){$E_r$}
  \put(81, 26){$E_{r-2}$}
  \put(12.2,10){\circle*{1.5}}
  \put(20,20){\circle*{1.5}}
  \put(67.2,28){\circle*{1.5}}
  \put(75, 38){\circle*{1.5}}
  \put(12, 5){$p_0$}
  \put(20, 15){$p_1$}
   \put(67, 23){$p_{r-2}$}
   \put(75, 33){$p_{r-1}$}
    \put(-22, 12){$\overbrace{\phantom{AAAAAAAAA}}$} 
    \put(-10, 20){$D$}
    \put(85,45){$\overbrace{\phantom{AAAAAAAA}}$} 
    \put(97,53){$T$}
    \put(15, 50){$\overbrace{\phantom{AAAAAAAAAAAAAAAA}}$} 
    \put(42,60){$E$}
\end{picture}
\]
\caption{Maximal $(-2)$-chains of disconnected type}
\label{figure:disconn:type}
\end{figure}
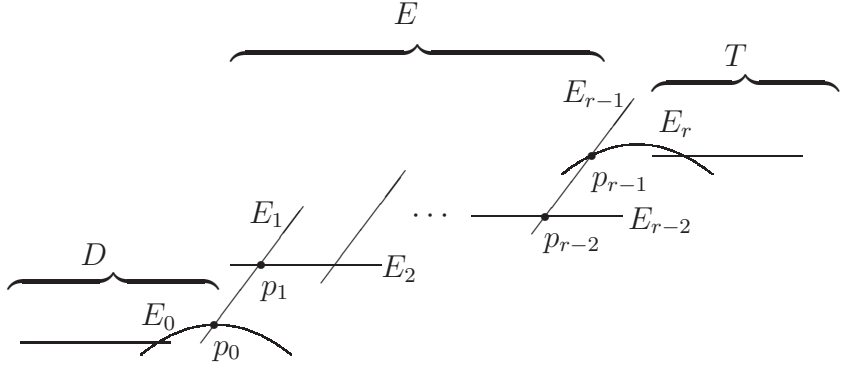

Fix an orientation on $e$, and let $\vec{e}$ denote the 
oriented edge. 
Let $v_0$ be 
the head of $\vec{e}$.
Let $E_0$ be the irreducible component of $\mathscr{X}_s$ with $[E_0] = v_0$.
Let $E_{1} , \ldots , E_{r-1}$ be the irreducible components of $E$, 
where we put the numbering in such a way that
$\# (E_i \cap E_{i+1} ) = 1$ for $i=0 , \ldots , r-2$.
We remark that this numbering is unique for $\vec{e}$.
Since $e$ is of disconnected type, $\mathscr{X}_s - E$ consists of two
connected components;
let $D$ be the connected component with $E_{0} \subset D$
and let $T$ be the other connected component.
Let $p_i$ denote the node with $\{ p_i \} = E_i \cap E_{i+1}$ for $i = 0
, \ldots , r-2$
and let $p_{r-1}$ be the node with $\{ p_{r-1 }\} = E_{r-1} \cap T$.
Let $E_r \in \Irr (T)$ with $p_{r-1} \in E_r$.
See Figure~\ref{figure:disconn:type}. 

For each $j =0, \ldots , r-1$,
the connected curve
$\sum_{i = j+1}^{r-1} E_i + T$
satisfies condition (F) in \S\ref{subsec:fvd},
and 
let $\mathscr{F}_j$ be the fundamental vertical divisor
associated to $\sum_{i = j+1}^{r-1} E_i + T$. 
(Note that $\mathscr{F}_{r-1}$ is the fundamental vertical divisor
associated to $T$.) 
Further, we set
\[
\mathscr{V}_{\vec{e}} :=
\sum_{i=0}^{r-1} \mathscr{F}_i,
\]
which we call the
\emph{stepwise vertical divisor
associated to $\vec{e}$.}
For the opposite orientation $\vec{e}^{\, opp}$, we also have the stepwise
vertical divisor $\mathscr{V}_{\vec{e}^{\, opp}}$.
When we say a \emph{stepwise vertical divisor for $e$},
we mean $\mathscr{V}_{\vec{e}}$
or $\mathscr{V}_{\vec{e}^{\, opp}}$.

We compute the degree of the stepwise vertical divisor
over each $C \in \Irr ( \mathscr{X}_s)$.
Noting Remark~\ref{rem:restrictionofcartierdivisor}, 
we obtain 
\begin{equation}
\label{eqn:degreepositivetype}
\deg
\left(
\rest{\OO_{\mathscr{X}} ( \mathscr{V}_{\vec{e}})}{C}
\right)
=
\begin{cases}
1
& \text{if $C = E_0$}
\\
-1
& \text{if $C = E_{r}$}
\\
0
&
\text{otherwise.}
\end{cases}
\end{equation}

For the order at
any $C \in \Irr ( \mathscr{X}_s )$,
we have
\begin{equation}
\label{eqn:ordVdisconnected}
\ord_{C} \left( \mathscr{V}_{\vec{e}} \right)
=
\begin{cases}
0
&
\text{if $C \in \Irr (D)$,}
\\
\sum_{j=0}^{i-1} \lambda_{j} & 
\text{if $C = E_{i}$ ($i = 1 , \ldots , r-1$ ),}
\\
\sum_{j=0}^{r-1} \lambda_{j}
&
\text{if $C \in \Irr (T)$,}
\end{cases}
\end{equation}
where $\lambda_j$ is the multiplicity at $p_j$, which also
equals the length of the canonical $1$-simplex $\Delta_{p_j}$ corresponding to $p_j$ (cf. \S\ref{subsection:skeleta}).

\begin{Remark}
We explain how a stepwise vertical divisor looks
by pretending that $R$ is a discrete valuation ring
and $\mathscr{X}$ is a regular scheme.
(This is not rigorous and only intuitive
because we are working over a non-Noetherian ring $R$.)
In this situation,
the stepwise vertical divisor $\mathscr{V}_{\vec{e}}$ for an oriented edge of disconnected type is 
\[
\text{
``$E_1 + 2 E_2 + \cdots + (r-1) E_{r-1} + r T$.''}
\]  
Similarly, for an edge $e$ of connected type, the stepwise vertical divisor $\mathscr{V}_e$ defined below 
is 
\[
\text{
``$E_1 + 2 E_2 + \cdots + (\ell -1) E_{\ell -1} + \ell  E_{\ell} 
+ (\ell -1) E_{\ell +1} + \cdots + 2 E_{2\ell-2} + E_{2\ell-1}$.''}
\]
\end{Remark}

\subsubsection*{Stepwise vertical divisor for an edge of connected type}
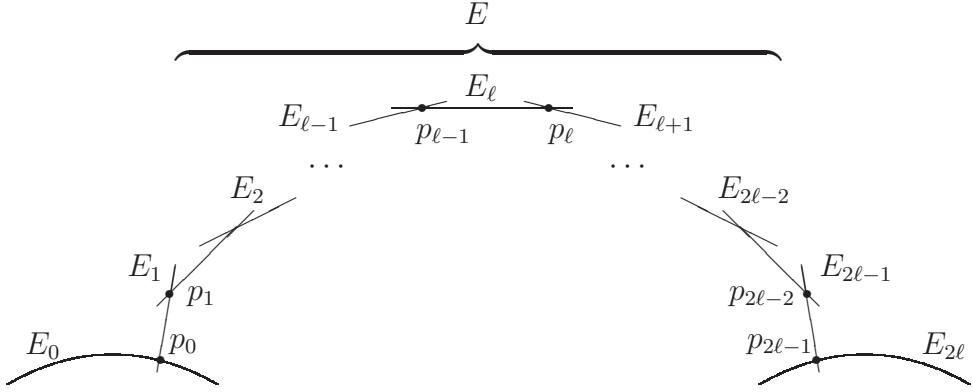
\begin{figure}[!h]
\[
\setlength\unitlength{0.08truecm}
\begin{picture}(180, 70)(0,0)
  \put(57,48){\line(4,1){16}}
  \put(102,48){\line(-4,1){16}}
  \put(50, 40){$\cdots$}
  \put(100, 40){$\cdots$}
  \put(64,51){\line(1,0){30}}
  \put(76, 53){$E_\ell$}
  \put(45, 48){$E_{\ell-1}$}
  \put(104, 48){$E_{\ell+1}$}
  \put(69,51){\circle*{1.5}}
  \put(90,51){\circle*{1.5}}
  \put(68, 46){$p_{\ell-1}$}
  \put(90, 46){$p_{\ell}$}
  \qbezier(0, 5)(17.5, 15)(35, 5)
  \put(25,7){\line(1,6){3}}
  \put(25,18){\line(1,1){16}}
  \put(32,28){\line(2,1){16}}
  \put(27,20){\circle*{1.5}}
  \put(133,20){\circle*{1.5}}
  \put(30, 19){$p_1$}
  \put(120, 19){$p_{2\ell-2}$}
  \put(25.5,9){\circle*{1.5}}
  \put(27, 11){$p_0$}
  \put(134.5,9){\circle*{1.5}}
  \put(123, 11){$p_{2 \ell-1}$}
  \qbezier(125, 5)(142.5, 15)(160, 5)
  \put(135,7){\line(-1,6){3}}
  \put(135,18){\line(-1,1){16}}
  \put(128,28){\line(-2,1){16}}  
  \put(3, 10){$E_0$}
  \put(20,23){$E_1$}
  \put(37,36){$E_2$}
  \put(118,36){$E_{2\ell-2}$}
  \put(135,23){$E_{2\ell-1}$}
  \put(152, 10){$E_{2 \ell}$}        
  \put(28, 55){$\overbrace{\phantom{AAAAAAAAAAAAAAAAAAAAAAAAAA}}$} 
  \put(76, 65){$E$}      
\end{picture}
\]
\caption{Maximal $(-2)$-chain of connected type}
\label{figure:conn:type}
\end{figure}

Next, we consider the case where $e$ is of connected type.
Here we assume that $E$ has symmetric multiplicities
(cf. Definition~\ref{def:symmetricmultiplicity}),
and let $E_{1} , \ldots , E_{2 \ell - 1}$ be the irreducible components of $E$
with the numbering such that
$\# (E_{i} \cap E_{i+1}) = 1$ for $i=1 , \ldots , 2 \ell -2$.
Let $E_0$ be the irreducible component of $\mathscr{X}_s - E$
with $E_{0} \cap E_1 \neq \emptyset$
and let 
$E_{2 \ell}$ be the irreducible component of $\mathscr{X}_s - E$
with $E_{2 \ell} \cap E_{2 \ell - 1} \neq \emptyset$.
We remark that $E_0 = E_{2 \ell}$ if and only if $e$ is a loop.
Let $p_i$ denote the node with $p_i \in E_{i} \cap E_{i+1}$ for $i=0 ,
\ldots , 2 \ell -1$.
See Figure~\ref{figure:conn:type}.

Since $E$ has symmetric multiplicities,
the connected curve
$\sum_{i = j+1}^{2\ell - 1 - j} E_i$ satisfies condition (F) in \S\ref{subsec:fvd} 
for each $j = 0 , \ldots \ell-1$,
and hence we have
the fundamental vertical divisor $\mathscr{F}_j$ with support
$\sum_{i = j+1}^{2\ell - 1 - j} E_i$.
We remark that
the divisor $\mathscr{F}_i$ does not depend on such numbering of the
irreducible components of $E$,
and hence well-defined for $e$.
We set 
\[
\mathscr{V}_{e} := 
\sum_{i = 0}^{\ell - 1} 
\mathscr{F}_i,
\]
which we
call the
\emph{stepwise vertical divisor associated to $e$.} 

We compute the degree of the stepwise vertical divisor
over each $C \in \Irr ( \mathscr{X}_s)$.
If 
$e$ is a loop,
then
$E_0 = E_{2\ell}$,
and 
noting Remark~\ref{rem:restrictionofcartierdivisor}, 
we obtain
\begin{equation}
\label{eqn:degreetype0:pre}
\deg
\left(
\rest{\OO_{\mathscr{X}} ( \mathscr{V}_{e})}{C}
\right)
=
\begin{cases}
2
& \text{if  $C = E_0 \;(=E_{2 \ell})$,}
\\
-2
& \text{if $C = E_{\ell}$,}
\\
0
&
\text{otherwise.}
\end{cases}
\end{equation}
If $e$ is not a loop,
then
$E_0 \neq E_{2\ell}$, and we have 
\begin{equation}
\label{eqn:degreetype0}
\deg
\left(
\rest{\OO_{\mathscr{X}} ( \mathscr{V}_{e})}{C}
\right)
=
\begin{cases}
1
& \text{if  $C = E_0$ or $C = E_{2 \ell}$,}
\\
-2
& \text{if $C = E_{\ell}$,}
\\
0
&
\text{otherwise.}
\end{cases}
\end{equation}

For the orders, regardless $e$ being a loop or not, 
we have 
\[
\ord_{E_i} \left( \mathscr{V}_{e} \right)
=
\ord_{E_{2 \ell - i}} \left( \mathscr{V}_{e} \right)
=
\sum_{j=0}^{i-1} \lambda_{j}
\] 
for $i = 1, \ldots, \ell$, where $\lambda_j$ is the multiplicity at the node $p_j$,
which equals the length of the canonical $1$-simplex $\Delta_{p_j}$ corresponding to $p_j$ (cf. \S\ref{subsection:skeleta}). 
Further, $\ord_{C} \left( \mathscr{V}_{e} \right) = 0$
for any $C \in \Irr (\mathscr{X}_s) \setminus \Irr( E)$.

\subsection{Edge-base sections and edge-unimodularity sections}
Let $X$ be a connected smooth projective curve over $K$ of genus $g$.
From here on to the end of this section,
assume that
$g \geq 2$.
Let $L$ be a line bundle over $X$.
In this subsection, we fix an edge $e \in E(\mathscr{X}^{\st})$. 
Let $(\Xscr, \Lscr)$ be a model of $(X, L)$ such that 
$\Xscr$ is a good model in Definition~\ref{def:goodmodel:new}. 
Recall that $\Sing ( \mathscr{X}_s)_{\subset e} := 
\{ p \in \Sing ( \mathscr{X}_s ) \mid
\Delta_p \subset  e \}$
(cf. (\ref{align:nodes:subset})), where $\Delta_p$ denotes the canonical $1$-simplex corresponding to $p$ 
in $S(\Xscr) = \Gamma_{\min}$. 

\begin{Definition}[edge-base section, $e$-base section]
\label{def:edge-base:section}
A nonzero global section $\widetilde{s} \in H^0(\mathscr{L})$ is 
called an \emph{$e$-base section} if 
$\widetilde{s} (p) \neq 0$
for any
$p \in \Sing ( \mathscr{X}_s)_{\subset e}$. 
\end{Definition}

In the following, 
we define the notion of edge-unimodularity sections
for an edge according to the type of the edge.

\subsubsection*{Edge of disconnected type}
Assume that $e \in E(\Xscr^{\st})$ is of disconnected type. 

\begin{Definition}[edge-unimodularity section, $e$-unimodularity section]
A nonzero global section $\widetilde{s} \in H^0(\mathscr{L})$ is 
called an \emph{$e$-unimodularity section} if 
there exists an orientation on $e$ such that
$\zero(\widetilde{s}) - \mathscr{V}_{\vec{e}}$ is effective on $\Xscr$ 
and is trivial 
on some open neighborhood of $\Sing ( \mathscr{X}_s)_{\subset e}$,
where $\vec{e}$ denotes the oriented edge. 
\end{Definition}

The next lemma asserts that
an $e$-base section and an $e$-unimodularity section
give a tropicalization map which is unimodular over the edge $e$.

\begin{Lemma}
\label{lemma:actuallyunimodular:discon}
Let $e \in E(\Xscr^{\st})$ be an edge of disconnected
type.
Suppose that there exists a model $(\Xscr, \Lscr)$ of $(X, L)$ such that 
$\Xscr$ is a good model 
\textup{(}cf. Definition~\textup{\ref{def:goodmodel:new}}\textup{)} and 
such that 
there exist 
an $e$-base section 
$\widetilde{s}^{(e)}_{0} \in H^0(\Lscr)$
and 
an $e$-unimodularity section
$\widetilde{s}^{(e)}_{1} \in H^0(\Lscr)$.  
Set $s^{(e)}_0 := \rest{\widetilde{s}^{(e)}_0}{X}$ 
and $s^{(e)}_{1} := \rest{\widetilde{s}^{(e)}_{1}}{X}$,
which are non-zero global sections of $L$.
We define a function $\varphi^{(e)} :
X^\an \setminus X(K) \to \RR$ by
$\varphi^{(e)} := - \log \left| s^{(e)}_{1}
/ s^{(e)}_{0} \right|$.
Let $\lambda_e$ denote the length of $e$.
Then 
$\varphi^{(e)}$
gives an isometry $e \to [0 , \lambda_e ]$.
\end{Lemma}

\Proof
To ease notation, we omit the superscript $(e)$ in the proof. 
We define the rational function $g$ on $\mathscr{X}$ by
$
g:=\widetilde{s}_{1} / \widetilde{s}_0
$.
Since $\widetilde{s}_{1}$
is an $e$-unimodularity section,
there exists an oriented edge $\vec{e}$ such that
$\zero \left(
\widetilde{s}_{1}
\right) - \mathscr{V}_{\vec{e}}$
is effective on $\Xscr$ and is trivial on some open neighborhood of 
$\Sing ( \mathscr{X}_s )_{\subset e}$.
Since
$\widetilde{s}_0$ does not have zero on some open neighborhood
of $\Sing ( \mathscr{X}_s )_{\subset e}$,
there exists an open neighborhood
$\mathscr{U}$ of $\Sing ( \mathscr{X}_s )_{\subset e}$
such that
$\rest{\zero \left( g \right)
- \mathscr{V}_{\vec{e}}}{\mathscr{U}}$
is trivial.

We use the notation 
used in \S\ref{subsection:stepwise}; in particular, we have 
$E = \sum_{i=1}^{r-1} E_i$, 
$[E_0]$ is the head of $\vec{e}$, 
$E_0 \subset D$, $E_r \subset T$,
and $\{ p_i \} = E_i \cap E_{i+1}$ for $i = 0, \ldots r-1$. 
For each $i = 0 , \ldots , r-1$,
we take a local \'etale atlas 
$\psi_{i} : \mathscr{U}_{i} \to \Spec ( R [x,y]/ (xy - \varpi_{i}))$
that distinguishes $p_i$ such that $\mathscr{U}_i
\subset \mathscr{U}$.
Exchanging $x$ and $y$ if necessary,
we may assume that in the special fiber $\left(
\mathscr{U}_{i}
\right)_s$ of $\mathscr{U}_i$,
the branch $E_{i+1}$ at $p_i$ is defined by $\psi_i^\ast (y)$.
Since $\rest{\zero \left( g \right)
- \mathscr{V}_{\vec{e}}}{\mathscr{U}}$
is trivial,
$g$ is a local equation of $\mathscr{V}_{\vec{e}}$ at each $p_i$.
By (\ref{eqn:ordVdisconnected}),
it follows that
for any $i = 0 , \ldots , r-1$,
there exists a unit $u_i \in \OO_{\mathscr{X}} ( \mathscr{U}_i)^{\times}$
such that
$g = u_i \prod_{j=0}^{i-1} \varpi_{j} 
\psi_i^\ast (y)$ over $\mathscr{U}_i$;
indeed, see the properties (i)--(iii) in the definition of 
fundamental vertical divisor and the definition of $\mathscr{V}_{\vec{e}}$
in \S\ref{subsec:fvd}. 
Note that $v_K 
\left(
\prod_{j=0}^{i-1} \varpi_{j}
\right) = \sum_{j = 0}^{i-1} \lambda_j$,
where $v_K = - \log |\cdot|_{K}$.
Since $\varphi = - \log \left| \rest{g}{X^{\an}} \right|$,
$\varphi =
\sum_{j = 0}^{i-1} \lambda_j 
- \log \left| \rest{\psi_1^{\ast} (y)}{X^{\an}} \right|$.
It follows from Lemma~\ref{lem:forinjectivity}(1)
that $\varphi$
maps $\Delta_{p_i}$ 
to $\left[ \sum_{j = 0}^{i-1} \lambda_j , 
\sum_{j = 0}^{i} \lambda_j
\right]$
isometrically.
Since $\lambda_e = \sum_{j = 0}^{r-1} \lambda_j$,
this shows that
$\varphi$
maps $e = \Delta_{E} = \bigcup_{i=0}^{r-1} \Delta_{p_i}$ to $\left[ 0,
\lambda_e  \right]$
isometrically.
\QED

\subsubsection*{Edge of connected type}
Assume that $e \in E(\Xscr^{\st})$ is of connected type.
In this case, we divide $e$ at the middle point of $e$ into two $1$-simplices,
and let $e^\prime$ denote one of them.

\begin{Definition}[$e^\prime$-unimodularity section]
A nonzero global section $\widetilde{s} \in H^0(\mathscr{L})$ is 
called an \emph{$e^\prime$-unimodularity section} if 
$\zero(\widetilde{s}) - \mathscr{V}_{e}$ is effective on $\Xscr$ and 
is trivial on some open neighborhood of $\Sing ( \mathscr{X}_s)_{\subset e^\prime}$.
\end{Definition}

\begin{Lemma}
\label{lemma:actuallyunimodular:conn}
Let $e \in E(\Xscr^{\st})$ be of connected type
and let $e^\prime$ be as above.
Suppose that there exists a model $(\Xscr, \Lscr)$ of $(X, L)$ such that 
$\Xscr$ is a good model 
and that there exist 
an $e$-base section 
$\widetilde{s}^{(e^\prime)}_{0} \in H^0(\Lscr)$
and 
an $e^\prime$-unimodularity section
$\widetilde{s}^{(e^\prime)}_1 \in H^0(\Lscr)$.  
Set $s^{(e^\prime)}_0 := \rest{\widetilde{s}^{(e^\prime)}_0}{X}$ 
and $s^{(e^\prime)}_{1} := \rest{\widetilde{s}^{(e^\prime)}_{1}}{X}$,
which are non-zero global sections of $L$.
We define a function $\varphi^{(e^\prime)} :
X^\an \setminus X(K) \to \RR$ by
$\varphi^{(e^\prime)} := - \log \left| s^{(e^\prime)}_{1}
/ s^{(e^\prime)}_{0} \right|$.
Let $\lambda_{e^\prime}$ denote the length of $e^\prime$.
Then 
$\varphi^{(e^\prime)}$
gives an isometry $e^\prime \to [0 , \lambda_{e^\prime} ]$. 
\end{Lemma}

Here we place a superscript $(e^\prime)$ to an edge-base section $\widetilde{s}^{(e^\prime)}_{0}$ 
in (1), because a model $(\Xscr, \Lscr)$ depends on $e^\prime$. 

\medskip
\Proof
We can prove the lemma 
by almost the same idea as
in Lemma~\ref{lemma:actuallyunimodular:discon}, so 
that
we omit the details. 
\QED

\subsection{Unimodular tropicalization}
\label{subsec:edge:good:model:disconnected}
In view of Lemmas~\ref{lemma:actuallyunimodular:discon}
and \ref{lemma:actuallyunimodular:conn},
we would like to prove 
the following Proposition~\ref{prop:existence:sections:stronger}
to complete the proof of a unimodular tropicalization of the minimal skeleton for $g \geq 2$. 
Recall that a good model is defined in Definition~\ref{def:goodmodel:new} 
and an edge-base section is defined in Definition~\ref{def:edge-base:section}. 
Given a model $(\Xscr, \Lscr)$ of $(X, L)$ such that $\Xscr$ is a good model, 
we put $\Mscr := \Lscr \otimes \omega_{\Xscr/R}^{\otimes -1}$ as before, and 
for a line bundle $\Nscr$ over $\Xscr$, 
let $D_{\Nscr}$ be the divisor on $\Gamma_{\min}$ defined in \eqref{eqn:def:DM}.

\begin{Proposition}
\label{prop:existence:sections:stronger}
Let $X$ be a connected smooth projective curve over $K$ of genus $g \geq 2$,
and let $L$ be a line bundle over $X$. Suppose that $\deg(L) \geq 3g-1$. 
Let $e \in E(\Xscr^{\st})$. 
\begin{enumerate}
\item
Suppose that $e$ is of disconnected type.
Fix an orientation on $e$ and let $\vec{e}$ denote the oriented
edge.
Then there exists a model $(\Xscr, \Lscr)$ of $(X, L)$ 
such that $\mathscr{X}$ is a good model 
with the following
properties\textup{:}  
there exist sections 
$\widetilde{s}^{(e)}_{0} , \widetilde{s}^{(e)}_{1} \in H^0(\Lscr)$ 
such that $\widetilde{s}^{(e)}_{0}(p) \neq 0$ for any 
$p \in \Sing(\Xscr_s)$ 
and such that
$\zero \left( \widetilde{s}^{(e)}_{1} \right)
- \mathscr{V}_{\vec{e}}$ is effective on $\Xscr$ and is trivial on
some neighborhood of $\Sing ( \mathscr{X}_s )$.
\textup{(}In particular, 
$\widetilde{s}^{(e)}_{0}$ is an $e$-base section,
and
$\widetilde{s}^{(e)}_{1}$ is an 
$e$-unimodularity section.\textup{)}
\item
Suppose that $e$ is of connected type.
Let $e^\prime$ be a $1$-simplex that appears by subdividing $e$ at the middle point, 
i.e., let $e^\prime$ be a half of $e$. 
Then  there exist a model $(\Xscr, \Lscr)$ of $(X, L)$ such that 
$\Xscr$ is a good model 
and such that there exist an $e$-base section 
$\widetilde{s}^{(e^\prime)}_{0} \in H^0(\Lscr)$ and 
an $e^\prime$-unimodularity section $\widetilde{s}^{(e^\prime)}  \in H^0(\Lscr)$. 
\end{enumerate}
\end{Proposition}

To prove Proposition~\ref{prop:existence:sections:stronger}, 
we show the following lemma based on Proposition~\ref{prop:pregoodmodel}.

\begin{Lemma}
\label{lemma:forbasesectionUT:dc}
Let $X$ be a connected smooth projective curve over $K$ of genus $g \geq 2$,
and let $L$ be a line bundle over $X$.
Assume that $\deg (L) \geq 3g-1$.
Take an edge $e \in E(\Xscr^{\st})$. 
Let $v$ be an end vertex of $e$.
Let $(\Xscr, \Lscr)$ be a model of $(X,L)$.
Assume that
$\mathscr{X}$ is a good model, 
$D_{\Mscr} - [v] \geq 0$ on the minimal skeleton $\Gamma_{\min}$,
and $(\Xscr, \Lscr)$ satisfies 
conditions \textup{(i)} and \textup{(ii)} in Proposition~\textup{\ref{prop:pregoodmodel}}.  
Then the following hold.
\begin{enumerate}
\item
The line bundle
$\mathscr{L}$ is free at
any $q \in \Sing ( \mathscr{X}_s )$.
\item
Suppose that $e$ is of disconnected type.
Put the orientation on $e$ with head $v$,
and let $\vec{e}$ denote the oriented edge.
Then
$\mathscr{L} ( - \mathscr{V}_{\vec{e}})$
is free at any $q \in \Sing ( \mathscr{X}_s )$.
\item
Suppose that $e$ is of connected type.
Divide $e$ into two $1$-simplices at the middle point of $e$,
and let $e^\prime$ be one of them with $v \in e^\prime$.
Then
$\mathscr{L} ( - \mathscr{V}_{e})$
is free at any $q \in \Sing ( \mathscr{X}_s )_{\subset e^\prime}$.
\end{enumerate}
\end{Lemma}

\Proof
To show (1),
we use Lemma~\ref{lemma:freeatnode:graphversion:emptyset}.
Note that $\deg ( D_{\Mscr} ) \geq g+1$
and
$D_{\Mscr} \geq 0$;
indeed, the inequality on the degree follows from $\deg (L) \geq 3g-1$,
and the effectivity follows from $D_{\Mscr} - [v] \geq 0$.
We take any $q \in \Sing ( \mathscr{X}_s )$.
If $q$ is of connected type,
then $\Gamma_{\min}\setminus \mathrm{relin}(\Delta_q)$
is connected
and $\deg(\rest{D_{\Mscr}}{\Gamma_{\min}\setminus \mathrm{relin}(\Delta_q)}) 
\geq (g+1) - \deg(\rest{D_{\Mscr}}{\mathrm{relin}(\Delta_q)})
\geq g  > 1$,
and thus Lemma~\ref{lemma:freeatnode:graphversion:emptyset} 
concludes that 
$\Lscr$ is free at $q$.
Suppose that $q$ is of disconnected type.
Then $\Delta_q$ is a part of an edge of disconnected type. 
Then 
any connected component $\Gamma^\prime$ 
of $\Gamma_{\min}\setminus \mathrm{relin}(\Delta_q)$ contains an island 
$\Gamma_i$ of $\Gamma_{\min}$, and Proposition~\ref{prop:pregoodmodel}(i) 
shows that $\deg(\rest{D_{\Mscr}}{\Gamma^\prime}) \geq 1$. 
It follows from Lemma~\ref{lemma:freeatnode:graphversion:emptyset}  that 
$\Lscr$ is free at $q$. 
Thus we have (1).

Assertion (2) can be shown similarly.
Let $v^\prime$ the end vertex of $e$ with $v \neq v^\prime$.
By \eqref{eqn:degreepositivetype}, we have 
$
D_{\mathscr{M} ( - \mathscr{V}_{\vec{e}})} 
= D_{\mathscr{M}} - [v] + [v^\prime]
$. 
Further,
since $D_{\Mscr} - [v] \geq 0$, we have 
$D_{\mathscr{M} ( - \mathscr{V}_{\vec{e}})} \geq 0$. 
Then the same argument as in (1) proves that
$\Lscr( - \mathscr{V}_{\vec{e}})$ is free at any
$q \in \Sing ( \mathscr{X}_s )$. 
Thus we have~(2).

Let us prove (3).
If $e$ is not a loop, then 
let
$v^\prime$ be the end vertex of $e$ other than $v$;
if $e$ is a loop, 
then set $v = v^\prime$. 
Let $w$ denote the middle point of $e$.
By \eqref{eqn:degreetype0:pre} and \eqref{eqn:degreetype0}, we have 
$
D_{\mathscr{M} ( - \mathscr{V}_{e})} 
= D_{\mathscr{M}} - [v] - [v^\prime] + 2[w]
$. 
We take any $q \in \Sing ( \mathscr{X}_s )_{\subset e^\prime}$.

Suppose that $D_{\mathscr{M} ( - \mathscr{V}_{e})}(v^\prime) \geq 0$. 
Then, since $D_{\Mscr} - [v] \geq 0$, we have $D_{\mathscr{M} ( - \mathscr{V}_{e})}
\geq 0$. 
Since $\Delta_q$ is a part of the edge $e$ of connected type, 
$\Gamma_{\min}\setminus \mathrm{relin}(\Delta_q)$ is connected, and 
Proposition~\ref{prop:pregoodmodel}(i) 
shows that 
\[
\deg(\rest{D_{\Mscr( - \mathscr{V}_{e})}}{\Gamma_{\min}\setminus \mathrm{relin}(\Delta_q)}) 
= \deg(\rest{D_{\Mscr}}{\Gamma_{\min}\setminus \mathrm{relin}(\Delta_q)}) 
\geq (g+1) - \deg\left(\rest{D_{\Mscr}}{ \mathrm{relin}(\Delta_q)}\right)
\geq g \geq 1. 
\] 
It follows from Lemma~\ref{lemma:freeatnode:graphversion:emptyset}  that 
$\Lscr( - \mathscr{V}_{e})$ is free at $q$.

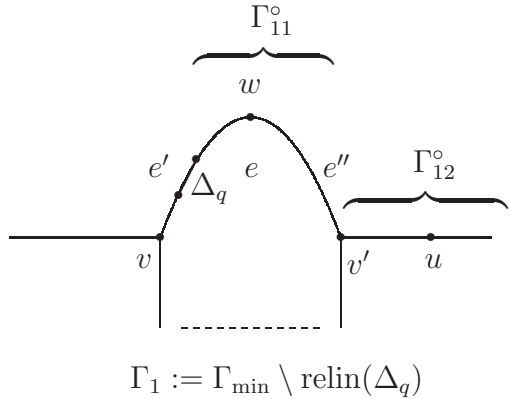
\begin{figure}[!h]
\[
\setlength\unitlength{0.08truecm}
\begin{picture}(100, 65)(0,0)
 \put(5, 25){\line(1,0){25}}
 \qbezier(30, 25)(45, 65)(60, 25)
 \put(60, 25){\line(1,0){25}}
 \put(60, 25){\line(0,-1){15}}
 \put(30, 25){\line(0,-1){15}}
  \put(30, 25){\circle*{1.5}}
  \put(60, 25){\circle*{1.5}}
  \put(45, 45){\circle*{1.5}}
  \put(75, 25){\circle*{1.5}}
  \put(33, 32){\circle*{1.5}}
  \put(36, 38){\circle*{1.5}}
  \multiput(33.5, 10)(2,0){12}{\line(1,0){1}}
  \put(25, 0){$\Gamma_1 := \Gamma_{\min}\setminus\relin(\Delta_q)$}
  \put(26, 20){$v$}
  \put(43, 49){$w$}
  \put(61, 19){$v^\prime$}
  \put(74, 20){$u$}
  \put(35, 32){$\Delta_q$}
  \put(28, 35){$e^\prime$}
  \put(57, 35){$e^{\prime\prime}$}
  \put(44, 35){$e$}
  \put(45, 60){$\Gamma_{11}^\circ$}
  \put(72, 36){$\Gamma_{12}^\circ$}
  \put(35.5, 50){$\overbrace{\phantom{AAAAAA}}$} 
  \put(60.5, 26){$\overbrace{\phantom{AAAAAAA}}$} 
  \thicklines
  \put(33, 32){\line(1,2){3}}
 \end{picture}
\]
\label{figure:for:lemma:forbasesectionUT:dc}
\end{figure}

Suppose that $D_{\mathscr{M} ( - \mathscr{V}_{e})}(v^\prime) < 0$. 
Since $D_{\Mscr} - [v] \geq 0$, this means that 
$D_{\mathscr{M} ( - \mathscr{V}_{e})} (v') =-1$. 
We use Lemma~\ref{lemma:freeatnode:graphversion:nonemptyset}. 
Set $V := \{ x \in \Gamma_{\min} \mid D_{\mathscr{M} ( - \mathscr{V}_{e})} (x) < 0 \}$.
In this case, we have $V= \{v^\prime\}$.
Set $\Gamma_1:= \Gamma_{\min}\setminus \mathrm{relin}(\Delta_q)$. 
Since $( \mathscr{X} , \mathscr{L})$
satisfies condition (ii)
in
Proposition~\ref{prop:pregoodmodel}, we have 
$\deg\left(\rest{(D_{\mathscr{M}}-[v])}{\mathrm{relin}(e)}\right) \leq 1$. 
It follows that
\[
\deg\left(\rest{(D_{\mathscr{M}}-[v])}{\Gamma_{\min}\setminus\mathrm{relin}(e)}\right) 
= \deg\left(D_{\mathscr{M}}-[v]\right)  - 
\deg\left(\rest{(D_{\mathscr{M}}-[v])}{\mathrm{relin}(e)}\right) 
\geq g - 1 \geq 1. 
\]
Since $D_{\mathscr{M} ( - \mathscr{V}_{e})}(v^\prime) = -1$,
we have 
$(D_{\mathscr{M}} - [v])(v') = -1$.
Further, we have 
$v' \in 
\Gamma_{\min}\setminus\mathrm{relin}(e)$.
It follows form the above inequality that 
there exists $u \in \Gamma_{\min}\setminus\mathrm{relin}(e)$ such that $u \neq v^\prime$ and 
$(D_{\mathscr{M}}-[v])(u) \geq 1$. 

Let $\Gamma_{11}^\circ$ be the connected component of $\Gamma_1 \setminus \{v^\prime\}$ that contains $w$, and 
let $\Gamma_{12}^\circ$ be the connected component of $\Gamma_1 \setminus \{v^\prime\}$ that contains $u$. 
Note that $\Gamma_{11}^\circ \neq \Gamma_{12}^\circ$. 
Indeed,
since  $\Delta_q \subset e^\prime$, we have
$\Gamma_{11}^\circ \subset \mathrm{relin}(e)$;
on the other hand, since $u \notin  \mathrm{relin}(e)$,
$\Gamma_{12}^\circ \not\subset  \mathrm{relin}(e)$.
Since the coefficients of $D_{\mathscr{M}(- \mathscr{V}_{e})}$ are nonnegative except at $v'$,
we have
\[
\deg\left(\rest{D_{\mathscr{M}(- \mathscr{V}_{e})}}{\Gamma_{11}^\circ}\right) 
\geq \deg(2 [w]) = 2 \geq 1,   \quad
\deg\left(\rest{D_{\mathscr{M}(- \mathscr{V}_{e})}}{\Gamma_{12}^\circ}\right) 
\geq \deg([u]) = 1. 
\]
Further,
if we set $\Gamma_{1j} := \Gamma_{1j}^\circ \cup \{v^\prime\}$ for $j = 1, 2$,
then 
the valence of $\Gamma_{11} \cup \Gamma_{12}$ at $v^\prime$ is at least $2 = - D_{\mathscr{M} ( - \mathscr{V}_{e})}(v^\prime) + 1$. 
Thus by Lemma~\ref{lemma:freeatnode:graphversion:nonemptyset},
$\Lscr( - \mathscr{V}_{e})$ is free at $q$. 
\QED

{\sl Proof  of Proposition~\ref{prop:existence:sections:stronger}.}\quad 
We construct a model $( \mathscr{X} , \mathscr{L} )$. 
We take an end-vertex $v$ of $e$ as follows: 
if $e$ is of disconnected type, we take $v$ to be the head of $\vec{e}$
(for (1));
if $e$ is of connected type,
we take $v$ such that $v \in e^\prime$
(for (2)).
We apply 
Proposition~\ref{prop:pregoodmodel}, where we take $v$ as $x$. 
Then we get a model $(\Xscr, \Lscr)$ such that: 
$\Xscr$ is a good model; 
$D_{\Mscr} - [v] \geq 0$ on $\Gamma_{\min}$; 
$(\Xscr, \Lscr)$ satisfies 
conditions (i) and (ii) in Proposition~\ref{prop:pregoodmodel}. 
By Lemma~\ref{lemma:forbasesectionUT:dc}(1), 
$\Lscr$ is free at any $p \in \Sing ( \mathscr{X}_s )$,
which implies that
there exists a global section of $\widetilde{s}_{0} \in H^0(\Lscr)$ 
that does not vanish at any $p \in \Sing ( \mathscr{X}_s )$. 
Thus
in assertion (1),
setting $\widetilde{s}_{0}^{(e)}$ to  be $\widetilde{s}_{0}$
proves
the existence of $\widetilde{s}_{0}^{(e)}$;
in assertion (2),
setting $\widetilde{s}_{0}^{(e')}$ to  be $\widetilde{s}_{0}$
proves
that of $\widetilde{s}_{0}^{(e')}$.

Suppose that $e$ is of disconnected type. We show (1). 
By Lemma~\ref{lemma:forbasesectionUT:dc}(2),  
$\Lscr ( - \mathscr{V}_{\vec{e}})$ is free at any $p \in 
\Sing ( \mathscr{X}_s )$,
which implies that
there exists a global section of 
$\widetilde{s}^{(e)}_{1,-} \in H^0( \Lscr( - \mathscr{V}_{\vec{e}}))$ 
that does not vanish at any $p \in \Sing ( \mathscr{X}_s )$.
Let $\widetilde{s}^{(e)}_{1}$ be the image of $\widetilde{s}^{(e)}_{1,-}$ 
by the natural inclusion $\Lscr( - \mathscr{V}_{\vec{e}})
\hookrightarrow \Lscr$.
Then $\widetilde{s}^{(e)}_{1}$ satisfies the required condition.
This completes the proof of (1).

To show (2),
suppose that $e$ is of connected type. 
Then
the same argument as above which uses Lemma~\ref{lemma:forbasesectionUT:dc}(3)
instead of Lemma~\ref{lemma:forbasesectionUT:dc}(2)
shows the existence of an $e^\prime$-unimodularity section.
\QED

We are ready to construct 
a unimodular tropicalization of the minimal skeleton $\Gamma_{\min}$. 
In fact, we construct a tropicalization map
that is not only unimodular but also separate
two points which are in some special positions.
To do this we prove one more lemma.

\begin{Lemma}
\label{lemma:ES:d-c}
Let $e$ be an edge of disconnected type.
Put an orientation on $e$, and let $\vec{e}$ denote
the oriented edge.
Let $\Gamma_1$ be the connected component of $\Gamma \setminus
\mathrm{relin}(e)$ containing the head of $\vec{e}$, 
and let $\Gamma_2$ be the other connected component.
Let $(\mathscr{X} , \mathscr{L})$
be a  
model of $(X,L)$ such that $\mathscr{X}$ is a good model. 
Let
$\widetilde{s}_0$ and $\widetilde{s}_1$ be global sections of $\mathscr{L}$
such that\textup{:}
\begin{enumerate}
\item[(i)]
$\widetilde{s}_0(p) \neq 0$ for any $p \in \Sing(\Xscr_s)$\textup{;} 
\item[(ii)]
$\zero (\widetilde{s}_1) - \mathscr{V}_{\vec{e}}$ is effective on $\Xscr$ 
and is trivial on some open neighborhood of $\Sing ( \mathscr{X}_s )$.
\end{enumerate}
Let $h$ denote the rational function
that is the restriction to $X$
of the rational function
$\widetilde{s}_1 / \widetilde{s}_0$ on $\mathscr{X}$.
Set $\varphi := - \log |h|$.
Let $\lambda_e$ denote the length of $e$.
Then $\varphi ( \Gamma_1 ) = \{  0 \}$
and $\varphi ( \Gamma_2 ) = \{  \lambda_e \}$.
\end{Lemma}

\Proof
We define the rational function $g$ on $\Xscr$ by $g := \widetilde{s}_1 / \widetilde{s}_0$, 
so that $h = \rest{g}{X}$. 
By conditions (i) and (ii), 
there exists an open neighborhood $\mathscr{U}$
of $\Sing ( \mathscr{X}_s )$
such that
$\zero (g) - \mathscr{V}_{\vec{e}} = \left(\zero (\widetilde{s}_1) - \mathscr{V}_{\vec{e}}\right)- \zero(\widetilde{s}_0)$
is trivial on $\mathscr{U}$.

First, we
take any $p \in \Sing ( \mathscr{X}_s)_{\subset \Gamma_1}$.
Let $C_1$ and $C_2$ be the irreducible components
of $\mathscr{X}_s$ such that $p \in C_1 \cap C_2$.
By the definition of $\mathscr{V}_{\vec{e}}$,
we have $\ord_{C_i} ( \mathscr{V}_{\vec{e}}) = 0$ for $i=1,2$.
By Lemma~\ref{lem:forinjectivity}, it follows
that $\varphi ( \Delta_p ) = 0$.
Since $\Gamma_1 = \bigcup_{p \in \Sing ( \mathscr{X}_s)_{\subset \Gamma_1}} \Delta_p$,
this proves $\varphi ( \Gamma_1 ) = \{  0 \}$.

Next, we take any $p' \in \Sing ( \mathscr{X}_s)_{\subset \Gamma_2}$.
Let $C_1'$ and $C_2'$ be the irreducible components
of $\mathscr{X}_s$ such that $p' \in C_1' \cap C_2'$.
By the definition of $\mathscr{V}_{\vec{e}}$,
we have $\ord_{C_i'} ( \mathscr{V}_{\vec{e}}) = \lambda_e$ for $i=1,2$.
By Lemma~\ref{lem:forinjectivity}, it follows
that $\varphi ( \Delta_p ) = \lambda_e$.
Since $\Gamma_2 = \bigcup_{p' \in \Sing ( \mathscr{X}_s)_{\subset \Gamma_2}} \Delta_{p'}$,
this proves $\varphi ( \Gamma_2 ) = \{  \lambda_e \}$.
\QED

\begin{Theorem} \label{theorem:UT:canonical}
Let $X$ be a connected smooth projective curve over $K$
of genus $g \geq 2$, 
and let $L$ be a line bundle over $X$.
Suppose that $\deg (L) \geq 3g-1$.
Then there exist nonzero global sections 
$s_0 , \ldots , s_N \in H^{0} (L)$ such that
the map
$\varphi : X^\an \to \TT\PP^N$
defined by
\[
\varphi := 
( - \log |s_0| : \cdots :
- \log |s_N| )
\]
gives a unimodular tropicalization of the minimal skeleton $\Gamma_{\min}$.
Furthermore, we can take $s_0 , \ldots , s_N$ 
in such a way that $\varphi$ has
the following properties. 
\begin{enumerate}
\item[(i)]
Let $e \in E(\Xscr^{st})$ be of disconnected type. Then $\rest{\varphi}{e}$ is an isometry. 
\item[(ii)]
Let $e \in E(\Xscr^{st})$ be of connected type. 
Let $e^{\prime}$ and $e^{\prime\prime}$ be the $1$-simplices arising 
by dividing $e$ at the middle point.
Then $\rest{\varphi}{e^{\prime}}$ and 
$\rest{\varphi}{e^{\prime\prime}}$ are isometries.  
\item[(iii)]
Let $e \in E(\Xscr^{st})$ be of disconnected type. 
Let $x , y \in \Gamma_{\min}$ be points 
that do not belong to the same connected component
of $\Gamma_{\min} \setminus \mathrm{relin}(e)$.
Then $\varphi (x) \neq \varphi (y)$
and $\varphi (x) , \varphi (y) \notin \varphi ( \mathrm{relin}(e) )$.
\end{enumerate}
\end{Theorem}

\Proof
We write $E(\Xscr^{\st}) = \{e_i\}_{i=1}^{m_1} \amalg  \{e_i\}_{i=m_1+1}^{m_1 + m_2}$, where 
$e_i$ ($1 \leq i \leq m_1$) are the  edges of disconnected type and 
$e_i$ ($m_1 + 1 \leq i \leq m_1 + m_2$) are the edges of connected type. 
For $1 \leq i \leq m_1$, we
give an orientation to $e_i$ and let $\vec{e_i}$ denote the oriented edge. 
Then 
Proposition~\ref{prop:existence:sections:stronger}(1) gives a model 
$(\Xscr_i, \Lscr_i)$ of $(X, L)$, 
a section 
$\widetilde{s}^{(e_i)}_{0}  \in H^0(\Lscr_{i})$
such that $\widetilde{s}^{(e_i)}_{0}(p) \neq 0$ for any $p \in \Sing(\Xscr_s)$, 
and a section
$\widetilde{s}^{(e_i)}_{1}  \in H^0(\Lscr_{i})$
such that
$\zero \left(\widetilde{s}^{(e_i)}_{1} \right) - \mathscr{V}_{\vec{e_i}}$
is effective on $\Xscr$ and is trivial on some open neighborhood of $\Sing ( \mathscr{X}_s )$. 
Note in particular
that 
$\widetilde{s}^{(e_i)}_{0}$ is an $e_i$-base section, and 
$\widetilde{s}^{(e_i)}_{1}$
is an $e_i$-unimodularity section.

For $m_1 + 1 \leq i \leq m_1+ m_2$, 
we divide $e_i$ into two $1$-simplices at the middle point
and let $e_{i}^\prime$ and $e_{i}^{\prime\prime}$ be the $1$-simplices arising from this
subdivision.
Proposition~\ref{prop:existence:sections:stronger}(2)
gives a model 
$(\Xscr_{i}^\prime, \Lscr_{i}^\prime)$ 
(resp. $(\Xscr_{i}^{\prime\prime}, \Lscr_{i}^{\prime\prime})$) of $(X, L)$, 
an $e_{i}$-base section 
$\widetilde{s}^{(e_{i}^\prime)}_{0}  \in H^0(\Lscr_{i}^\prime)$ 
(resp. $\widetilde{s}^{(e_{i}^{\prime\prime})}_{0}  \in H^0(\Lscr_{i}^{\prime\prime})$), 
and an $e_{i}^\prime$-unimodularity section 
$\widetilde{s}^{(e_{i}^\prime)}_{1}  \in H^0(\Lscr_{i}^\prime)$ 
(resp. an $e_{i}^{\prime\prime}$-unimodularity section 
$\widetilde{s}^{(e_{i}^{\prime\prime})}_{1}  \in H^0(\Lscr_{i}^{\prime\prime})$). 

Set $N := 2m_1 + 4 m_2 - 1$.
We consider $\rest{\widetilde{s}^{(e_i)}_0}{X}, \rest{\widetilde{s}^{(e_i)}_{1}}{X} \in H^{0} (L)$ for 
$1 \leq i \leq m_1$  and  
$\rest{\widetilde{s}^{(e_{i}^\prime)}_0}{X}, \rest{\widetilde{s}^{(e_{i}^\prime)}_{1}}{X}, 
\rest{\widetilde{s}^{(e_{i}^{\prime\prime})}_0}{X},
\rest{\widetilde{s}^{(e_{i}^{\prime\prime})}_{1}}{X} \in H^{0} (L) 
$ for $m_1 + 1 \leq i \leq m_1+ m_2$,
and we
denote those $(N+1)$ global sections of $L$ by $s_0, \ldots, s_{N}$. 
Then it is straightforward from Lemmas~\ref{lemma:actuallyunimodular:discon}
and \ref{lemma:actuallyunimodular:conn} that the map $\varphi$ determined by 
$s_0 , \ldots , s_N$ gives a unimodular tropicalization of $\Gamma_{\min}$ 
having 
properties (i) and (ii).

Furthermore,
this tropicalization also has property (iii).
Indeed, let $e$, $x$,  and $y$ be as in (iii),
and
take $i = 1 , \ldots , m_1$ with $e = e_i$.
By the definition of $s_0 , \ldots , s_N$,
there exist $a , b = 1 ,\ldots , N$ such that 
$s_a = \rest{\widetilde{s}^{(e_i)}_0}{X}$ and $s_b = \rest{\widetilde{s}^{(e_i)}_{1}}{X}$.
Set $h := - \log |s_b/ s_a|$.
Then by Lemma~\ref{lemma:ES:d-c},
$h (x) \neq h (y)$. Further, 
by Lemma~\ref{lemma:actuallyunimodular:discon} and Lemma~\ref{lemma:ES:d-c}, 
we have $h (x) , h (y) \notin h ( \mathrm{relin} (e_i))$.
This shows that $\varphi$ has property (iii) of the theorem.
\QED

\setcounter{equation}{0}
\section{Faithful tropicalization of minimal skeleta for $g \geq 2$}
\label{section:FTcan}

Let $X$ be a connected smooth projective curve 
over $K$ of genus $g$.
In this section, we assume that $g \geq 2$. 
Let $\Gamma_{\min}$ be the minimal skeleton of $X^{\an}$. 
As explained in \S\ref{subsec:skeleton:weighted:metric:graph}, 
we endow $\Gamma_{\min}$ with the canonical weight function $\omega$, 
so that $\bar{\Gamma}_{\min} = (\Gamma_{\min}, \omega)$ is the minimal weighted skeleton. 
Let $\mathscr{X}^{\st}$ denote the stable model of $X$ over $R$. 
Also, as explained in \S\ref{subsec:skeleton:weighted:metric:graph}, we endow $\Gamma_{\min}$ with the canonical finite metric structure with the set of vertices $V(\Xscr^{\st}) = V(\bar{\Gamma}_{\min})$ 
and the set of edges $E(\Xscr^{\st}) = E(\bar{\Gamma}_{\min})$. 
We note that $g = g(\bar{\Gamma}_{\min}) \geq g(\Gamma_{\min})$. 

Let $L$ be a line bundle over $X$. 
If $\deg (L) \geq 3g-1$, 
then we have already constructed enough global sections
of $L$ which give a unimodular tropicalization
of $\Gamma_{\min}$. 
However, we have not yet obtained a faithful tropicalization of $\Gamma_{\min}$. 
In view of Theorem~\ref{theorem:UT:canonical}(i)(ii)(iii),
we still need to construct global sections of $L$ 
that separate distinct points $x ,y$ 
of $\Gamma_{\min}$ as below: 
\begin{enumerate}
\item[(i)]
$e \in E ( \mathscr{X}^{\st})$ is an edge of connected type with middle point $w$, 
and $x , y \in \mathrm{relin}(e)$ belong to the different connected components of 
$\mathrm{relin}(e) \setminus \{w\} $; 
\item[(ii)]
$e, f\in E ( \mathscr{X}^{\st})$ are distinct edges, and 
$x \in \mathrm{relin}(e)$ and $y \in f$;
\item[(iii)]
$x , y \in V ( \mathscr{X}^{\st})$ are distinct vertices.
\end{enumerate}

We will construct global sections
which take care of (i) in 
\S\ref{subsection:AS}, 
(ii) in \S\ref{subsection:ES},
and (iii) in \S\ref{subsection:VS}, 
and then we will show 
in \S\ref{subsec:ft:min:skeleton}
that those sections
give a faithful tropicalization of $\Gamma_{\min}$. 
This section will be technically the most difficult in this paper. 

\subsection*{Notation and terminology of \S\ref{section:FTcan}} 
We put together the notation and terminology that are used 
throughout this section.
Let $X$ be a connected smooth projective curve 
over $K$ of genus $g \geq 2$, and let $ L$ be a line bundle over $X$. 
The minimal skeleton $\Gamma_{\min}$ of $X^{\an}$ is endowed with the canonical weight function 
and the canonical finite graph structure as above. 
A good model of $X$ always means a model defined in Definition~\ref{def:goodmodel:new}.  
Given a model $(\Xscr, \Lscr)$ of $(X, L)$ such that $\Xscr$ is a good model, 
we put $\Mscr := \Lscr \otimes \omega_{\Xscr/R}^{\otimes -1}$ as before, and 
for a line bundle $\Nscr$ over $\Xscr$, 
let $D_{\Nscr}$ be the divisor on $\Gamma_{\min}$ defined in \eqref{eqn:def:DM}.
For a node $p \in \Sing(\Xscr_s)$, let $\Delta_p$ denote the canonical $1$-simplex corresponding to $p$ 
in $\Gamma_{\min}$ (cf. Definition~\ref{def:canonical:1:simplex:node}). 

\subsection{Separating points on an edge of connected type} 
\label{subsection:AS}
In this subsection, we construct global sections
which will be used to separate points in 
an edge of connected type.

Let $e \in E ( \mathscr{X}^{\st})$ be en edge of connected type with middle point $w$. 
Let $\mathscr{X}$ be a good model of $X$. By Lemma~\ref{lemma:make:clear}, 
we take the maximal $(-2)$-chains $E$ of $\mathscr{X}_s$ with $\Delta_E = e$. 
Let $E_1 , \ldots , E_{2 \ell - 1}$ be the irreducible components of $E$. 
Since $\mathscr{X}$ is a good model, we have $2 \ell - 1 \geq 3$. 
We give the numbering for 
$E_1 , \ldots , E_{2 \ell - 1}$ as is
illustrated in Figure~\ref{figure:conn:type} in \S\ref{subsection:stepwise}. 
Note that
 $w = [E_\ell]$.
As in Figure~\ref{figure:conn:type} in \S\ref{subsection:stepwise}, 
let $E_0$ and $E_{2 \ell}$ denote the irreducible components of 
$\mathscr{X}_s - E$ with $E _0 \cap E_1 \neq \emptyset$ and 
$E _{2 \ell -1} \cap E_{2 \ell} \neq \emptyset$, and we set 
$\{ p_i \} := E_{i} \cap E_{i+1}$ for $i= 0 , \ldots , 2 \ell - 1$. 

Let $\lambda_i$ denote the multiplicity of $\mathscr{X}$ at $p_i$.
Since $\mathscr{X}$ is a good model, we have $\lambda_{0} = \lambda_{\ell - 1}
= \lambda_{\ell} = \lambda_{2 \ell - 1}$.
It follows that there exist unique fundamental vertical divisors
$\mathscr{A}_{1}$ and  $\mathscr{A}_{2}$ with 
\[
 \Supp(\mathscr{A}_{1}) = E_1 + \cdots + E_{\ell - 1} 
 \quad\text{and}\quad  
 \Supp(\mathscr{A}_{2}) = E_{\ell + 1} + \cdots + 
E_{2 \ell - 1}.  
\]
Then $\Delta_{\Supp ( \mathscr{A}_1 )}$ and $\Delta_{\Supp ( \mathscr{A}_2 )}$
are the two $1$-simplices in $e$ that arise by dividing $e$ at $w$. 

\begin{Definition}
[separating divisor for an edge of connected type, separating divisor for $e$]
\label{def:assymetric:divisor}
Let $e \in E ( \mathscr{X}^{\st})$ be an edge of connected type, and 
let $\mathscr{X}$ be a good model of $X$. 
We call a divisor $\mathscr{A}$ on $\Xscr$ a \emph{separating divisor for $e$} 
if $\mathscr{A} = \mathscr{A}_1$ or $\mathscr{A} = \mathscr{A}_2$. 
\end{Definition}

We remark that the notion of a separating divisor 
is well-defined for $e$, and it
does not depend on the two choices of the numbering $E_1 , \ldots , E_{2 \ell - 1}$
for the irreducible components of~$E$.

\begin{Definition}
[separating section for an edge of connected type, separating section for $e$]
\label{def:asymmetric:section:along:e}
Let $e \in E ( \mathscr{X}^{\st})$ be an edge of connected type. 
Let $( \mathscr{X} , \mathscr{L})$ be a model of $(X,L)$
such that $\mathscr{X}$ is a good model.
Let $\widetilde{s}$ be a nonzero global section of $\mathscr{L}$.
We call $\widetilde{s}$
a \emph{separating section for $e$}
if there exists
a separating divisor $\mathscr{A}$
for $e$
such that
$\zero \left( \widetilde{s} \right) -
\mathscr{A}$ is effective on $\Xscr$ 
and
is trivial on some open neighborhood of
$\Sing ( \mathscr{X}_s )_{\subset e}
\setminus \Supp ( \mathscr{A})$.
\end{Definition}

Let $e^\prime$ and $e^{\prime\prime}$ be the two $1$-simplices that arise
by dividing $e$ at $w$. We take the 
separating divisor $\mathscr{A}$ with $\Delta_{\Supp ( \mathscr{A})} = e^{\prime\prime}$. 
Then $\Sing ( \mathscr{X}_s )_{\subset e}
\setminus \Supp ( \mathscr{A}) = \Sing (\Supp ( \mathscr{A}))_{\subset e^\prime}$.
(For example, if $e^{\prime\prime} =
\Delta_{E_{\ell+1} + \cdots + E_{2\ell-1}}$, then 
$\mathscr{A} = \mathscr{A}_2$ and 
$\Sing ( \mathscr{X}_s )_{\subset e}
\setminus \Supp ( \mathscr{A}) =\{p_{0}, \ldots, p_{\ell-1}\}$.) 

We use the following lemma to separate points 
in $\mathrm{relin}(e)$ in \S\ref{subsec:ft:min:skeleton}. 
Recall that an edge-base section is defined in Definition~\ref{def:edge-base:section}.

\begin{Lemma}
\label{lemma:asymmetric:section:separate:new}
Let $e \in E(\Xscr^{\st})$ be an edge of connected type with middle point $w$. 
Suppose that there exists a model $(\Xscr, \Lscr)$ of $(X, L)$ such that 
$\Xscr$ is a good model and that there exist 
an $e$-base section $\widetilde{s}^{(e)}_0$ 
and a separating section $\widetilde{s}^{(e)}_{1}$ for $e$. 
We set $g := \widetilde{s}^{(e)}_1/\widetilde{s}^{(e)}_{0}$ 
and define a nonzero rational function $h$ on $X$ 
by $h:= \rest{g}{X}$.    
Let $e^\prime$ and $e^{\prime\prime}$ be the two $1$-simplices that
arise by dividing $e$ at $w$. 
Then for any $x \in \mathrm{relin}(e^\prime)$ 
and $y \in \mathrm{relin}(e^{\prime\prime})$,
we have
$- \log |h (x)| \neq - \log |h (y)|$.
\end{Lemma}

\Proof
Since 
$\widetilde{s}^{(e)}_0(p) \neq 0$ for any $p \in \Sing(\Xscr_s)_{\subset e}$, 
we take an open neighborhood $\Uscr$ of $\Sing(\Xscr_s)_{\subset e}$ such that 
$\rest{\widetilde{s}^{(e)}_0}{\Uscr}$ is 
nowhere vanishing
over $\Uscr$.
Then $\zero ( \rest{\widetilde{s}^{(e)}_0}{\Uscr} ) = 0$.
By the definition of a separating section $\widetilde{s}^{(e)}_{1}$ for $e$,
there exists a separating divisor $\mathscr{A}$ for $e$
as in Definition~\ref{def:asymmetric:section:along:e}.
Set $g := \widetilde{s}^{(e)}_{1} / \widetilde{s}^{(e)}_0$.
Then 
$\zero(g) - \mathscr{A}$ is effective over $\Uscr$. 

Note that $\Delta_{\Supp (\mathscr{A} )}$ equals $e^{\prime}$ or $e^{\prime\prime}$.
Without loss of generality, we assume that $\Delta_{\Supp (\mathscr{A} )}
= e^{\prime\prime}$. 
We set 
$\varphi := - \log |h| = \log \left| \rest{g}{X}\right|$. 
Let
$C \in \Irr(\Xscr_s)$ with $[C] \in \relin(e^{\prime\prime})$. Then 
we have $\ord_{C}(\mathscr{A}) > 0$. Further, since $e^{\prime\prime} \subset e$, 
we have $C \cap \Uscr \neq \emptyset$, so that the generic point of $C$ belongs to $\Uscr$. 
Then it follows from Lemma~\ref{lem:forinjectivity}(1) that 
$\varphi ([C]) \geq \ord_{C}(\mathscr{A}) > 0$. 
By Lemma~\ref{lem:forinjectivity}(2),
we have 
$\varphi (y) > 0$ for any $y  \in \mathrm{relin}(e^{\prime\prime})$.
On the other hand,
since $\widetilde{s}^{(e)}_{1}$ does not vanish
at any $p \in \Sing ( \mathscr{X}_s )_{\subset e^\prime} = 
\Sing ( \mathscr{X}_s )_{\subset e} \setminus \Supp(\mathscr{A})$,
$g^{-1}$ is regular on some neighborhood of $\Sing ( \mathscr{X}_s )_{\subset e^\prime}$.
By
Lemma~\ref{lem:forinjectivity},
we see that
$\varphi (x) \leq 0$ for any $x \in e^\prime$.
This proves the lemma.
\QED

In view of Lemma~\ref{lemma:asymmetric:section:separate:new}, 
our next task is to show the existence of a model 
of $(X,L)$ that has an $e$-base section and a separating section for $e$.

\begin{Proposition}
\label{prop:model:AS}
Suppose that $\deg (L) \geq 3 g -1$.
Then for any edge $e \in E ( \mathscr{X}^{\st})$ of connected type,
there exists a model $(\Xscr , \Lscr)$ of $(X,L)$ such that $\Xscr$ is a good model 
and 
such that there exist an $e$-base section and a separating section for $e$ of $\mathscr{L}$.
\end{Proposition}

\Proof
We take an edge $e \in E ( \mathscr{X}^{\st})$ of connected type.
Let $w$ denote the middle point of $e$. 
Applying Proposition~\ref{prop:pregoodmodel} for $x = w$,
we obtain a model $( \mathscr{X} , \mathscr{L} )$ of $(X,L)$
with the following properties: 
\begin{enumerate}
\item[(i)]
$\mathscr{X}$ is a good model of $X$; 
\item[(ii)]
$\mathscr{M}$ is vertically nef (i.e., $D_{\mathscr{M}} \geq 0$),
$D_{\mathscr{M}}(w) \geq 1$,
and $\deg (\rest{D_{\mathscr{M}}}{\mathrm{relin}(e)}) \leq 2$. 
\end{enumerate}
We are going to show that there exist 
an $e$-base section and a separating section for $e$ of $\Lscr$. 

\medskip
{\bf Step 1.}\quad 
We consider the existence of an $e$-base section of $\mathscr{L}$. 
Take any $p \in \Sing ( \mathscr{X}_s)_{\subset e}$
and put $\widetilde{\Gamma} := \Gamma_{\min} \setminus \mathrm{relin} (\Delta_p)$.
Since $p$ is of connected type,
$\widetilde{\Gamma}$ is connected.
By condition (ii) above,
$\rest{D_{\mathscr{M}}}{\widetilde{\Gamma}}$ 
is effective and $\deg\left(\rest{D_{\mathscr{M}}}{\widetilde{\Gamma}}\right) 
\geq \deg \left(D_{\mathscr{M}}\right) - \deg(\rest{D_{\mathscr{M}}}{\mathrm{relin} (e)}) \geq (g+1) -2 \geq 1$.
By Lemma~\ref{lemma:freeatnode:graphversion:emptyset},
$\mathscr{L}$ is free at $p$,
which shows the existence of an $e$-base section.

\medskip
{\bf Step 2.}\quad 
We show that $\mathscr{L}$ has a
separating section for $e$. 
The proof will be somewhat similar to that of Lemma~\ref{lemma:forbasesectionUT:dc}(3). 
By Lemma~\ref{lemma:make:clear}, we take 
the maximal $(-2)$-chains of $\Xscr_s$ with $\Delta_E = e$. 
We use the notation as in the beginning of this subsection. 
In particular,  $E_1, \ldots, E_{2\ell-1}$ are the irreducible components 
of $E$, where numbering is given as in Figure~\ref{figure:conn:type}. 
We denote by $e^\prime$ and $e^{\prime\prime}$ the $1$-simplices arising by dividing $e$ at $w$
such that $e^\prime
= \Delta_{E_1 + \cdots + E_{\ell}}$ and 
$e^{\prime\prime} = \Delta_{E_{\ell+1} + \cdots + E_{2\ell-1}}$. 
We consider the separating divisor $\mathscr{A}_2$ for $e$. 
Recall that $\Sing ( \mathscr{X}_s )_{\subset e}
\setminus \Supp ( \mathscr{A}_2 ) =
\Sing ( \mathscr{X}_s )_{\subset e^\prime} = \{p_0, \ldots, p_\ell\}$. 

We set $\mathscr{L}^\prime := \mathscr{L} ( - \mathscr{A}_2)$.
We show that $\mathscr{L}^\prime$ is free at any
$q \in \Sing ( \mathscr{X}_s )_{\subset e^\prime}$
by using Lemma~\ref{lemma:freeatnode:graphversion:emptyset}
and Lemma~\ref{lemma:freeatnode:graphversion:nonemptyset}.
Take any $q \in \Sing ( \mathscr{X}_s )_{\subset e^\prime}$
and set $\Gamma_1 := \Gamma_{\min} \setminus \mathrm{relin} ( \Delta_q )$.
Since $q$ is of connected type, $\Gamma_1$
is connected.
We set
$\mathscr{M}^\prime := \mathscr{L}^\prime \otimes \omega_{\mathscr{X}/R}^{\otimes -1}
= \mathscr{M} ( - \mathscr{A}_2 )$. 
We set $w_i = [E_i] \in \Gamma_{\min}$ for $i = 0, \ldots, 2 \ell$. Note that 
$w = w_\ell$. 
Then we have 
\begin{equation}
\label{eqn:D:M:prime}
D_{\mathscr{M}^\prime} 
= D_{\mathscr{M}} -D_{\OO_{\mathscr{X}} (\mathscr{A}_2)} 
= (D_{\mathscr{M}} - [w]) + [w_{\ell + 1}] +  
[w_{2 \ell - 1}] - [w_{2 \ell}]. 
\end{equation}
By condition (ii), we have $D_{\mathscr{M}} - [w] \geq 0$, and thus 
$V := \{v \in V(\Xscr^{\st}) \mid D_{\mathscr{M}^\prime} (v) <0\}$
is either the empty set or equal to the one point set $\{w_{2 \ell}\}$. 
We argue it case by case.

\smallskip
{\bf Case 1.}\quad 
Suppose that $V = \emptyset$, i.e., 
$D_{\mathscr{M}^\prime} \geq 0$. 
Since 
\[
\deg\left(\rest{D_{\mathscr{M}^\prime}}{\Gamma_1}\right) 
= \deg\left(\rest{D_{\mathscr{M}}}{\Gamma_1}\right) 
\geq 
\deg \left(D_{\mathscr{M}}\right) -\deg(\rest{D_{\mathscr{M}}}{\mathrm{relin} (e)}) \geq (g+1) -2 \geq 1, 
\]
Lemma~\ref{lemma:freeatnode:graphversion:emptyset} tells us that 
$\mathscr{L}^\prime$ is free at $q$. 

\smallskip
{\bf Case 2.}\quad 
Suppose that $V = \{w_{2 \ell}\}$.  In this case, 
$\rest{D_{\mathscr{M}^\prime}} {\Gamma_{\min} \setminus {\{w_{2 \ell}}\}} \geq 0$ and 
$D_{\mathscr{M}^\prime} (w_{2 \ell}) = -1$ by \eqref{eqn:D:M:prime}.  
Let $\Gamma_{11}^{\circ}$ be the connected component of 
$\Gamma_1 \setminus \{w_{2 \ell}\}$
such that $w_{2 \ell - 1} \in \Gamma_{11}^\circ$. 
Then
$\deg
\left( \rest{D_{\mathscr{M}^\prime}}{\Gamma_{11}^{\circ}} \right) \geq 
D_{\mathscr{M}^\prime} (w_{2 \ell - 1}) \geq 1$. 
Since 
$\deg \left(
\rest{D_{\mathscr{M}}}{\mathrm{relin} (e)}
\right) \leq 2$, we get from \eqref{eqn:D:M:prime} that 
\[
\deg \left(
\rest{D_{\mathscr{M}^\prime}}{\mathrm{relin} (e)}
\right) 
= 
\deg \left(
\rest{D_{\mathscr{M}}}{\mathrm{relin} (e)}\right)  
+ \deg(- [w] + [w_{\ell + 1}] +  [w_{2 \ell - 1}] )
\leq 2 + 1 = 3.
\]
Since $\deg (D_{\mathscr{M}^\prime}) = 
\deg (D_{\mathscr{M}}) \geq g+ 1 \geq 3$ and $D_{\mathscr{M}^\prime} (w_{2\ell}) = -1$,
we have  
\[
\deg \left( \rest{D_{\mathscr{M}^\prime}}{
\Gamma_{\min}\setminus ( \mathrm{relin} (e) \cup \{w_{2\ell}\})} 
 \right) 
 = 
 \deg (D_{\mathscr{M}^\prime}) - 
 \deg \left(
\rest{D_{\mathscr{M}^\prime}}{\mathrm{relin} (e)}
\right) 
 - D_{\mathscr{M}^\prime} (w_{2\ell})
 \geq 1.
\]
Since $\Gamma_{\min}\setminus ( \mathrm{relin} (e) \cup \{w_{2\ell}\}) \subset \Gamma_1 \setminus \{{w_{2\ell}\}}$ and 
$\Gamma_{11}^\circ \subset \mathrm{relin} (e)$, 
we see that there exists a
connected component $\Gamma_{12}^{\circ}$
of $\Gamma_1 \setminus \{w_{2 \ell}\}$ such that
$\Gamma_{12}^{\circ} \neq \Gamma_{11}^{\circ}$
and $\deg
\left( \rest{D_{\mathscr{M}^\prime}}{\Gamma_{12}^{\circ}} \right) \geq 1$.

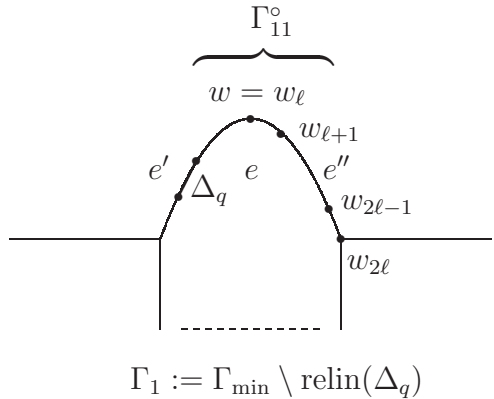
\begin{figure}[!h]
\[
\setlength\unitlength{0.08truecm}
\begin{picture}(100, 65)(0,0)
 \put(5, 25){\line(1,0){25}}
 \qbezier(30, 25)(45, 65)(60, 25)
 \put(60, 25){\line(1,0){25}}
 \put(60, 25){\line(0,-1){15}}
 \put(30, 25){\line(0,-1){15}}
  \put(60, 25){\circle*{1.5}}
  \put(45, 45){\circle*{1.5}}
  \put(33, 32){\circle*{1.5}}
  \put(36, 38){\circle*{1.5}}
  \multiput(33.5, 10)(2,0){12}{\line(1,0){1}}
  \put(25, 0){$\Gamma_1 := \Gamma_{\min}\setminus\relin(\Delta_q)$}
  \put(58, 30){\circle*{1.5}}
  \put(50, 42.5){\circle*{1.5}}
  \put(53, 42){$w_{\ell+1}$}
  \put(61, 20){$w_{2\ell}$}
  \put(35, 32){$\Delta_q$}
  \put(28, 35){$e^\prime$}
  \put(57, 35){$e^{\prime\prime}$}
  \put(44, 35){$e$}
  \put(60, 30){$w_{2\ell-1}$}
  \put(38, 48){$w = w_{\ell}$}
  \put(45, 60){$\Gamma_{11}^\circ$}
  \put(35.5, 50){$\overbrace{\phantom{AAAAAA}}$} 
  \thicklines
  \put(33, 32){\line(1,2){3}}
 \end{picture}
\]
\label{figure:for:prop:model:AS}
\end{figure}

Let $\Gamma_{11} := \Gamma_{11}^{\circ} \cup \{w_{2\ell}\}$ be the closure of 
$\Gamma_{11}^{\circ}$ in $\Gamma_1$,
and
let $\Gamma_{12} := \Gamma_{12}^{\circ} \cup \{w_{2\ell}\}$ be the closure of 
$\Gamma_{12}^{\circ}$ in $\Gamma_1$. 
Since the valence of $\Gamma_{11} \cup
\Gamma_{12}$ at $w_{2\ell}$ is at least $2
= - D_{\mathscr{M}^\prime} (w_{2\ell}) + 1$, one can apply 
Lemma~\ref{lemma:freeatnode:graphversion:nonemptyset}.  
(With the notation of Lemma~\ref{lemma:freeatnode:graphversion:nonemptyset}, 
condition (i) of Lemma~\ref{lemma:freeatnode:graphversion:nonemptyset} is vacuous, 
and for condition (ii), we take $s=2$ and consider $\Gamma_{11} $ and $\Gamma_{12}$.) 
This proves that $\mathscr{L}^\prime$ is free at $q$ also in this case.

Thus $\mathscr{L}^\prime$ is free at $q$ in any case. It follows that 
there exists a global section $s^\prime $ of $\mathscr{L}^\prime  
:= \mathscr{L} ( - \mathscr{A}_2 ) $
such that $s^\prime (q) \neq 0$ for any $q \in 
\Sing ( \mathscr{X}_s )_{\subset e^\prime} = 
\Sing ( \mathscr{X}_s )_{\subset e}
\setminus \Supp ( \mathscr{A}_2 )$.
Let $s$ be the image of $s^\prime$
by the natural homomorphism $\mathscr{L} ( - \mathscr{A}_2 ) 
\hookrightarrow \mathscr{L}$.
Then $s$ is a separating section for $e$.
\QED

\subsection{Separating points in different edges}
\label{subsection:ES}

Let $e, f \in E (\mathscr{X}^{\st})$ be distinct edges of $\Gamma_{\min}$. 
In this subsection, we construct global sections 
that separate a point in $\mathrm{relin}(e)$ and 
a point in $f$.

\begin{Remark}
\label{remark:disconnect-other}
If $e$ is an edge of disconnected type, then we 
have already constructed such global sections. 
Indeed, let $\varphi$ be the tropicalization map
as in Theorem~\ref{theorem:UT:canonical}.
If $e$ is of disconnected type,
then $f$ is contained in a connected component of $\Gamma
\setminus \mathrm{relin}(e)$. 
By property (iii) of 
Theorem~\ref{theorem:UT:canonical},
we obtain $\varphi (f) \cap \varphi ( \mathrm{relin}(e)) =
\emptyset$, and thus $\varphi$ separates a point in $\mathrm{relin}(e)$ and a point in $f$.
\end{Remark}

With Remark~\ref{remark:disconnect-other}, 
in the rest of this subsection, we assume that $e$ is an edge of connected type.  
Suppose that $\mathscr{X}$ is a good model of $X$. 
By Lemma~\ref{lemma:make:clear}, 
we take the maximal $(-2)$-chain $E$  of $\Xscr_s$ with 
$\Delta_E = e$. By the definition of a good model, $E$ has symmetric multiplicities. 
It follows that there exists a unique fundamental vertical divisor 
with support $E$ (cf. \S\ref{subsec:fvd}),  which we denote by $\mathscr{F}_e$. 

\subsubsection{Edges-separating sections}

\begin{Definition}[edges-separating section, $(e, f)$-separating section]
Let $e, f \in E(\Xscr^{\st})$ be distinct edges such that $e$ is of connected type. 
Let $(\mathscr{X} , \mathscr{L})$ be a 
model of $(X,L)$
such that $\mathscr{X}$ is a good model. 
Let $\widetilde{s}$ be a nonzero global section of $\mathscr{L}$. 
We call $\widetilde{s}$ an \emph{$(e,f)$-separating section}
if 
$
\zero
\left(
\widetilde{s}
\right)
-
\mathscr{F}_e
$
is effective on $\Xscr$ and
is trivial on some open neighborhood of 
$\Sing ( \mathscr{X}_s)_{\subset f}
$.
\end{Definition}

We use 
in \S\ref{subsec:ft:min:skeleton}
the following lemma to separate a point in $\relin(e)$ and a point in $f$. 

\begin{Lemma}
\label{lemma:ES:separate:edges}
Let $e, f \in E(\Xscr^{\st})$ be distinct edges such that $e$ is of connected type. 
Assume that there exists a model 
$(\mathscr{X} , \mathscr{L})$ 
such that $\mathscr{X}$ is a good model and such that 
there exist an $e$-base section 
$\widetilde{s}_0$ \textup{(}cf.~Definition~\textup{\ref{def:edge-base:section}}\textup{)} and 
an $(e,f)$-separating section $\widetilde{s}_1$ of $\mathscr{L}$.
We set $g := \widetilde{s}_1 / \widetilde{s}_0$, and define a nonzero rational function 
$h$ on $X$ by $h:= \rest{g}{X}$. 
Then for any $x \in \mathrm{relin}(e)$
and any $y \in f$,
we have $-\log|h(x)| > 0$ and $-\log|h(y)| \leq 0$,
and in particular, $-\log|h(x)| \neq -\log|h(x)|$. 
\end{Lemma}

\Proof
Since $\widetilde{s}_0$ is an
$e$-base section,
$\zero ( \widetilde{s}_0)$
is trivial over some open neighborhood $\mathscr{U}_0$
of $\Sing ( \mathscr{X}_s )_{\subset e}$.
On the other hand,
since $\widetilde{s}_1$ is an $(e,f)$-separating section,
$\zero ( \widetilde{s}_1) - \mathscr{F}_e$ 
is effective on $\mathscr{X}$
and
is trivial
over some neighborhood $\mathscr{U}_1$ of $\Sing ( \mathscr{X}_s )_{\subset f}$.
Since $e \neq f$,
we have $\Sing ( \mathscr{X}_s )_{\subset f} \cap \Supp ( \mathscr{F}_e ) = 
\emptyset$,
so that shrinking $\mathscr{U}_1$ if necessary,
we may furthermore assume that 
$\mathscr{U}_1 \cap \Supp ( \mathscr{F}_e ) = \emptyset$.
Then $\zero ( \widetilde{s}_1)$ is trivial on $\mathscr{U}_1$.
By the definition of $g$,
$\zero (g) - \mathscr{F}_s = (\zero( \widetilde{s}_1) - \mathscr{F}_s) - \zero( \widetilde{s}_0)$ is effective
on $\mathscr{U}_0$,
and $- ( \zero (g) ) = \zero (g^{-1})$ is effective
over $\mathscr{U}_1$.

Take any $C \in \Irr ( \mathscr{X}_s)$ with $[C] \in \mathrm{relin} (e)$.
Since $\zero (g) - \mathscr{F}_e$ is effective
on $\mathscr{U}_0$,
it follows from
Lemma~\ref{lem:forinjectivity}(1) 
that $ - \log |h ([C]) | \geq \ord_C(\mathscr{F}_e) > 0$. 
By Lemma~\ref{lem:forinjectivity}(2),
we see that $- \log |h (x)| > 0$ 
for any $x \in \mathrm{relin} (e)$.
Similarly,
take any $C' \in \Irr ( \mathscr{X}_s)$
with $[C'] \in f$.
Since $- ( \zero (g) ) = \zero (g^{-1})$ is effective
on $\mathscr{U}_1$,
it follows from Lemma~\ref{lem:forinjectivity}(1) 
that $ - \log |h ([C]) | \leq 0$.
By Lemma~\ref{lem:forinjectivity}(2),
we see that $- \log |h (y)| \leq 0$ 
for any $y \in f$.
This concludes the lemma.
\QED

Our task is then to show the existence of a model of $(X, L)$ that 
has an $e$-base section and an $(e,f)$-separating section. 
The goal of this subsection is to prove the following proposition.

\begin{Proposition}
\label{prop:model:ES}
Suppose that $\deg (L) \geq 3g -1$. 
Let $e, f \in E(\Xscr^{\st})$ be distinct edges such that $e$ is of connected type. 
Then there exists a  model $( \mathscr{X} , \mathscr{L} )$ of $(X,L)$
such that $\mathscr{X}$ is a good model  
and such that $\mathscr{L}$ has an $e$-base section and
an $(e , f)$-separating section.
\end{Proposition}

The above proposition follows the proposition below.

\begin{Proposition} \label{prop:sep:edges}
Suppose that $\deg (L) \geq 3g -1$, and let 
$e$ and $f$ be as in Proposition~\ref{prop:model:ES}. 
Then there exists a model $( \mathscr{X} , \mathscr{L})$
of $(X,L)$ with the following properties.
\begin{enumerate}
\item[(i)]
$\mathscr{X}$ is a good model.
\item[(ii)]
For any $p \in \Sing ( \mathscr{X}_s )_{\subset e}$,
$\mathscr{L}$ is free at $p$.
\item[(iii)]
Let $\mathscr{V}_e$ denote the stepwise vertical divisor for $e$.
Set $\mathscr{L}_1 := \mathscr{L} ( - \mathscr{V}_e)$.
Then for any $q \in \Sing ( \mathscr{X}_s )_{\subset f}$,
$\mathscr{L}_1$ is free at $q$.
\end{enumerate}
\end{Proposition}

We prove that Proposition~\ref{prop:sep:edges}  implies 
Proposition~\ref{prop:model:ES}. 
Let $( \mathscr{X} , \mathscr{L} )$ be as in Proposition~\ref{prop:sep:edges}. 
Then since $\mathscr{L}$ is free at any $p \in \Sing (\mathscr{X}_s)_{\subset e}$,
there exists a global section $\widetilde{s}_0$ of $\mathscr{L}$
such that $\widetilde{s}_0 (p) \neq 0$ for any such $p$.
Thus we have an $e$-base section.
Since $\mathscr{L}_1 = \mathscr{L} ( - \mathscr{V}_e)
$ is free at any $q \in \Sing (\mathscr{X}_s)_{\subset f}$,
there exists a global section $\widetilde{s}_1^-$ of $\mathscr{L}_1$
such that $\widetilde{s}_1^- (q) \neq 0$ for any such $q$.
Let $\widetilde{s}_1$ be the image of $\widetilde{s}_1^-$ by the canonical injection 
$\mathscr{L}_1 \hookrightarrow \mathscr{L}$.
Then
$\zero (\widetilde{s}_1) - \mathscr{V}_e$ is effective and is trivial on some open neighborhood of 
$\Sing (\mathscr{X}_s)_{\subset f}$.
Since $\mathscr{V}_e - \mathscr{F}_e$ is effective
and 
is trivial on some open neighborhood of $\Sing (\mathscr{X}_s)_{\subset f}$,
it follows that
$\zero (\widetilde{s}_1) - \mathscr{F}_e  = (\zero (\widetilde{s}_1) - \mathscr{V}_e) 
+ (\mathscr{V}_e - \mathscr{F}_e)$ is effective 
and is trivial on some open neighborhood of 
$\Sing (\mathscr{X}_s)_{\subset f}$.
Thus $\widetilde{s}_1$ is an $(e , f)$-separating section,
and we obtain Proposition~\ref{prop:model:ES}.

\medskip
Thus our goal is to show Proposition~\ref{prop:sep:edges}.
We will prove it according to the position of $e$ and $f$ in $\Gamma_{\min}$.
The construction of a model will be done separately according to the following two cases:
\begin{enumerate}
\item[(A)]
$f$ is of disconnected type, or $f$ is
of connected type and 
$\Gamma_{\min} \setminus ( \mathrm{relin} (e) \cup \mathrm{relin} (f) )$
is connected;
\item[(B)] 
$f$ is of connected type and
$\Gamma_{\min} \setminus ( \mathrm{relin} (e) \cup \mathrm{relin} (f) )$
is not connected.
\end{enumerate}
For case (A), we may take the same kind of models 
obtained by
Proposition~\ref{prop:pregoodmodel}
as before,
but for case (B), we need to construct another kind of model.
The construction will be done below.

\subsubsection{Model for case \textup{(B)}}
We construct a suitable model $(\Xscr, \Lscr)$ for case (B). 
Under the condition of (B), we remark that  
$e$ nor $f$ is a loop and that 
$\Gamma_{\min} \setminus ( \mathrm{relin} (e) \cup \mathrm{relin} (f) )$
has exactly two connected components.
Indeed, since $e$ is of connected type,
$\Gamma_{\min} \setminus \mathrm{relin} (e)$
is connected;
since $
\Gamma_{\min} \setminus ( \mathrm{relin} (e) \cup \mathrm{relin} (f) )
=
( \Gamma_{\min} \setminus \mathrm{relin} (e) ) \setminus \mathrm{relin} (f)$
is not connected, it follows that
this has exactly two connected components.

\begin{Proposition}
\label{prop:separation:special2}
Suppose that $\deg (L) \geq 3g -1$, and let 
$e$ and $f$ be as in Proposition~\textup{\ref{prop:model:ES}}. 
Assume that Case~(B) applies for $e, f$. 
Let $\Gamma_1$ and $\Gamma_2$ be the connected components
of $\Gamma_{\min}\setminus (\mathrm{relin} (e) \cup \mathrm{relin}(f))$.
Then there exists a model $(\mathscr{X} , \mathscr{L})$
of $(X,L)$ 
such that
$\mathscr{X}$ is a good model and 
$\deg \left(
\rest{D_{\mathscr{M}}}{\Gamma_i}
\right) \geq 1$ for $i= 1, 2$. 
\end{Proposition}

To construct a model as in Proposition~\ref{prop:separation:special2},
we will apply Proposition~\ref{prop:pregoodmodel:pre} 
together with the following lemma on the $\Lambda$-metric graph $\Gamma_{\min}$.

\begin{Lemma}
\label{lemma:separation:special1}
Let $e$ and $f$ be edges of $\Gamma_{\min}$ such that 
case~\textup{(B)} applies for $e$ and $f$. Let $\Gamma_1$ and $\Gamma_2$ be 
as in Proposition~\ref{prop:separation:special2}. 
Let $D \in \Div_{\Lambda} ( \Gamma_{\min} )$ be a divisor with $\deg (D) \geq 
g + 1$. Then there exists an effective divisor $E \in \Div_{\Lambda} ( \Gamma_{\min} )$
such that $E \sim D$
and $\deg (\rest{E}{\Gamma_i}) \geq 1$ for both $i =1,2$.
\end{Lemma}

\Proof
Let $w_1$ and $w_2$ be the end vertices of $f$
such that $w_i \in \Gamma_i$. 
\begin{figure}[!h]
\[
\setlength\unitlength{0.08truecm}
\begin{picture}(60, 50)(0,0)
 \put(5, 15){\line(0,1){20}}
 \put(50, 15){\line(0,1){20}}
 \put(5, 35){\line(2,1){10}}
 \put(5, 15){\line(2,-1){10}}
 \put(50, 35){\line(-2,1){10}}
 \put(50, 15){\line(-2,-1){10}}
 \multiput(18, 40)(2,0){10}{\line(1,0){1}}
 \multiput(18, 9)(2,0){10}{\line(1,0){1}}
 \put(26, 45){$\Gamma_1$}
 \put(26, 2){$\Gamma_2$}
  \put(5, 15){\circle*{1.5}}
  \put(5, 35){\circle*{1.5}}
  \put(50, 15){\circle*{1.5}}
  \put(50, 35){\circle*{1.5}}
  \put(52, 13){$w_2$} 
  \put(52, 33){$w_1$} 
  \put(46, 23){$f$} 
  \put(7, 23){$e$} 
 \end{picture}
\]
\label{figure:for:lemma:separation:special1:1}
\end{figure}
First we suppose that $\deg (D) \geq g (\Gamma_{\min}) + 2$.
Then using Riemann's inequality (Proposition~\ref{prop:RIforMG}),
we have an effective divisor $E^\prime$
with $E^\prime  \sim D -  [w_1] - [w_2]$, and 
$E := E^\prime +  [w_1] + [w_2] \in \Div_{\Lambda} ( \Gamma_{\min} )$ is a 
desired effective divisor. Thus we are done.

In the remaining of the proof, 
we assume that $\deg (D) \leq g (\Gamma_{\min}) + 1$. 
Since $g = g ( \bar{\Gamma}_{\min})\geq g ( \Gamma_{\min})$ and we assume that $\deg (D) \geq g + 1$, 
we have $\deg (D) = g (\Gamma_{\min}) + 1$ and $g ( \Gamma_{\min}) = g$. 
Then the following 
claim follows from Remark~\ref{remark:valence:stable}.

\begin{Claim}
\label{claim:lemma:separation:special1:valence:3}
Every vertex of $V(\Xscr^{\st})$ has valency at least $3$. 
\end{Claim}

We note the following claim.

\begin{Claim}
We have $g (\Gamma_{i} ) \geq 1$ for any $i = 1,2$. 
\end{Claim}

To prove the above claim by contradiction,
we assume that $g (\Gamma_{i} ) = 0$ for some $i$. Without loss of generality, 
we assume that $g (\Gamma_{1} ) = 0$. 
Let $v_1$ be the end point of $e$ with $v_1 \in \Gamma_1$.

Suppose that $\Gamma_{1} $ is a singleton.
Then $\Gamma_{1} = \{ w_1 \}$ and $w_1 = v_1$. 
Since $\Gamma_{1}$ is a connected component of 
$\Gamma_{\min} \setminus ( \mathrm{relin} (e) \cup \mathrm{relin} (f) )$,
$e$ and $f$ are the edges of $\Gamma_{\min}$ emanating from $w_1$.
Since $e$ nor $f$ is a loop, this means that 
the valence of $\Gamma_{\min}$ at $w_1$ equals $2$.
However,
that contradicts to Claim~\ref{claim:lemma:separation:special1:valence:3}.

Suppose that 
$\Gamma_{1} $ is a chain.
Then $\Gamma_1$ is a segment connecting $w_1$ and $v_1$.
Let $f'$ be an edge of $\Gamma_{\min}$ emanating from $w_1$
with $f' \neq f$.
Since $\Gamma_{1}$ is a connected component of 
$\Gamma_{\min} \setminus ( \mathrm{relin} (e) \cup \mathrm{relin} (f) )$,
we have $f' \subset \Gamma_1$.
It follows that
the valence of $\Gamma_{\min}$
at $w_1$ equals $2$.
However,
that contradicts to Claim~\ref{claim:lemma:separation:special1:valence:3}.

Thus $\Gamma_1$ is a tree with at least three 
vertices at which $\Gamma_1$ has valence $1$.
Take a valence $1$ point
$u$ of $\Gamma_1$ with $u \neq v_1$ and $u \neq w_1$.
Then
$u$ is a point at which $\Gamma$ also has valence $1$.
However, that is impossible by 
Claim~\ref{claim:lemma:separation:special1:valence:3}.
This concludes that $g (\Gamma_{1} ) \geq 1$,
and thus
we have shown the claim.

\medskip
Let $w$ be the middle point of $f$.
By Proposition~\ref{prop:reduced:Lambda-div},
we take the $w$-reduced divisor $D_{w} \in \Div_{\Lambda} ( \Gamma )$ 
with $D_{w} \sim D$. Since $\deg (D) = g (\Gamma_{\min}) + 1$, 
$D_{w}$ is effective and $D_{w} (w) \geq 1$
(cf. Remark~\ref{rmk:on:reduced:divisors}). 

\begin{Claim}
\label{claim.reduced:divisor:new3}
If $\Supp (D_{w}) \cap \Gamma_i = \emptyset$
for some $i = 1,2$,
then
$D_{w} (w) \geq 2$.
\end{Claim}

We prove the above claim.
Without loss of generality, we assume that 
$\Supp (D_{w}) \cap \Gamma_2 = \emptyset$.
Let $\Gamma'$ be the topological space obtained from $\Gamma$
by contracting $\Gamma_2$ to a singleton,
and let $\alpha : \Gamma \to \Gamma'$ be the natural surjective
map; we
remark that $\{ \alpha (w_2) \} = \alpha ( \Gamma_2 )$.
Let $\beta : \Gamma \setminus \Gamma_2 \to 
\Gamma' \setminus \alpha ( \Gamma_2)$ denote
the restriction 
of $\alpha$.
Since $\beta$ is a homeomorphism
and since $\alpha ( \Gamma_2)$ is a singleton,
we endow $\Gamma'$ with a $\Lambda$-metric 
graph structure such that 
$\beta$
is an isometry.

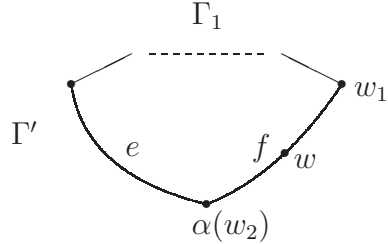
\begin{figure}[!h]
\[
\setlength\unitlength{0.08truecm}
\begin{picture}(60, 40)(0,0)
 \qbezier(10, 25)(12, 10)(32.5, 5)
 \qbezier(55, 25)(45, 10)(32.5, 5)
 \put(10, 25){\line(2,1){10}}
 \put(55, 25){\line(-2,1){10}}
 \put(0, 15){$\Gamma^\prime$}
 \multiput(23,30)(2,0){10}{\line(1,0){1}}
 \put(30, 35){$\Gamma_1$}
  \put(55, 25){\circle*{1.5}}
  \put(32.5, 5){\circle*{1.5}}
  \put(10, 25){\circle*{1.5}}
  \put(45.5, 13.5){\circle*{1.5}}
   \put(57, 23){$w_1$} 
  \put(47, 11){$w$} 
  \put(19, 13){$e$} 
  \put(40, 13){$f$} 
  \put(30, 0){$\alpha(w_2)$}
 \end{picture}
\]
\label{figure:for:lemma:separation:special1:2}
\end{figure}

Let $ D_{w}'$ be the pushout of 
$D_{w}$ by $\alpha$;
in this situation, since $\Supp (D_{w}) \cap \Gamma_2 = \emptyset$,
it is the divisor supported on $\Gamma' \setminus \{ \alpha (w_2)\}$
that coincides with $D_{w} $ via the isometry
$\beta$.
We show that $ D_{w} '$ is an
$\alpha (w)$-reduced divisor.
To argue it by contradiction,
suppose that $ D_{w} '$ is not an
$\alpha (w)$-reduced divisor.
Then there exists a connected compact subset $A$
of $\Gamma' \setminus \alpha (w)$
such that any point in $\partial A$
is saturated for $D_{w} '$.
The set
$\alpha^{-1} ( A )$ is a connected compact
subset of 
$\Gamma \setminus \{ w \}$.
Since $D_{w}' (\alpha (w_2)) = 0$,
$\alpha (w_2) \notin \partial A$.
It follows that $\alpha$ is a homeomorphism
over $\partial A$,
and hence
$\alpha^{-1} ( \partial A ) = \partial 
\alpha^{-1} (  A )$.
Since $\alpha$ is a homeomorphism
over $\partial A$ 
and over $\Supp ( D_{w} ')$,
it follows that
any point in $\partial 
\alpha^{-1} (  A )$
is saturated for
$D_{w}$.
However, this contradicts that $D_w$ is $w$-reduced.
This shows that $D_{w}'$ is $\alpha (w)$-reduced.

Since $g ( \Gamma_2) \geq 1$,
$g ( \Gamma' ) =
 g (\Gamma ) - g ( \Gamma_2)
\leq g (\Gamma ) - 1$,
and hence
$\deg ( D_{w}' ) \geq g ( \Gamma' ) + 2$.
By Remark~\ref{rmk:on:reduced:divisors},
it follows that $D_{w}' (\alpha (w) ) \geq 2$.
Thus $D_{w} (w) \geq 2$. 
We have shown the claim.

\medskip
Let us finish the proof of the lemma. 
If $\Supp (D_{w}) \cap \Gamma_i \neq \emptyset$ for any $i = 1,2$,
then $D_{w}$ is a desired effective divisor,
since $D_{w} \in |D|$. 
Suppose that  $\Supp ( D_{w} ) \cap \Gamma_i = \emptyset$ for some $i = 1, 2$.
Then by Claim~\ref{claim.reduced:divisor:new3}, $D_{w} (w) \geq 2$. 
Then $E^\prime := D_{w} - 2 [w]$ is an effective divisor. 
Since $w$ is the middle point of $f$,
$2 [w] \sim [w_1] + [w_2]$, so that 
$D \sim D_{w} \sim E^\prime + [w_1] + [w_2]$. 
Then $E:=  E^\prime + [w_1] + [w_2]$ is a 
desired effective divisor.
Thus the lemma holds.
\QED

\medskip
\textsl{Proof of Proposition~\ref{prop:separation:special2}.}\quad
We take a divisor $\tilde{D}$ on $X$ such that 
$L \otimes \omega_{X}^{\otimes -1} \cong \OO_{X} ( \tilde{D} )$.
Let $\tau : X^{\an} \to \Gamma_{\min}$ be the retraction map.
Then $\deg
\left(
\tau_{\ast} ( \tilde{D})
\right) \geq (3g-1) - (2g-2) =  g+1$,
and by Lemma~\ref{lemma:separation:special1},
there exists a $D \in 
\left|
\tau_{\ast} ( \tilde{D})
\right|$ such that
$\deg
\left(
\rest{D}{\Gamma_1}
\right)
\geq 
1$
and 
$\deg
\left(
\rest{D}{\Gamma_2}
\right)
\geq 
1$.
By Proposition~\ref{prop:pregoodmodel:pre},
there exist a 
good model $\mathscr{X}$ of $X$
and a line bundle $\mathscr{M}$ over 
$\mathscr{X}$ such that
$D_{\mathscr{M}} = D$.
Set $\mathscr{L} := \mathscr{M} \otimes \omega_{\mathscr{X}/R}$.
Then by construction, $(\mathscr{X} , \mathscr{L})$
is a model that has the required properties
in Proposition~\ref{prop:separation:special2}.
\QED

\subsubsection{Proof of Proposition~\ref{prop:sep:edges}}

Let  $v_1$ and $v_2$ denote the end vertices of $e$ and 
let
$v$ be the middle point of $e$. We remark that if $e$ is a loop, then 
$v_1 = v_2$ and $v$ is the point antipodal to $v_1 (= v_2)$.

Our argument goes according to cases (A) or (B) above.

\medskip
{\bf Case (A): }
Suppose that 
$f$ is of disconnected type, or $f$ is
of connected type and 
$\Gamma_{\min} \setminus ( \mathrm{relin} (e) \cup \mathrm{relin} (f) )$
is connected.

Applying Proposition~\ref{prop:pregoodmodel},
we obtain a model $( \mathscr{X} , \mathscr{L})$
with the following properties: 
\begin{enumerate}
\item[(i)]
$\mathscr{X}$ is a good model;
\item[(ii)]
$D_{\mathscr{M}} - [v_1]$ is effective;
\item[(iii)]
$D_{\mathscr{M}} - [v_1]$ has positive degree over any 
island of $\Gamma_{\min}$. 
\item[(iv)]
$\deg(\rest{D_{\mathscr{M}}}{\relin(e)}) \leq 1$ and $\deg(\rest{D_{\mathscr{M}}}{\relin(f)}) \leq 1$. 
(With comparison with Proposition~\ref{prop:pregoodmodel}(ii), we note that $v_1$ does not lie in the relative interior of $e$ nor $f$.)
\end{enumerate}

We are going to show that $( \mathscr{X} , \mathscr{L})$ enjoys the required properties. 

We take any $p \in \Sing ( \mathscr{X}_s )_{\subset e}$. 
Since $e$ is of connected type and $\Delta_p \subset e$ , we see that $\Gamma_{\min} \setminus 
\relin(\Delta_p)$ is connected. 
Since $\deg(\rest{\mathscr{M}}{\Gamma_{\min} \setminus 
\relin(\Delta_p)}) \geq \deg(\rest{\mathscr{M}}{\Gamma_{\min} \setminus 
\relin(e)}) \geq (g+1) - \deg(\rest{D_{\mathscr{M}}}{\relin(e)}) \geq g > 0$ by condition (iv), 
if follows from Lemma~\ref{lemma:freeatnode:graphversion:emptyset} that 
$\Lscr$ is free at $p$. Thus property (ii) of Proposition~\ref{prop:sep:edges} is verified. 

It remains to verify property (iii) of Proposition~\ref{prop:sep:edges}, 
i.e., $\mathscr{L}_1 = \mathscr{L} (-\mathscr{V}_e)$ is free at any $q  \in \Sing ( \mathscr{X}_s )_{\subset f}$. 
Fixing any $q \in \Sing ( \mathscr{X}_s )_{\subset f}$, we set 
$\widetilde{\Gamma} := \Gamma_{\min} \setminus \mathrm{relin}(\Delta_q)$.
Note that $\widetilde{\Gamma}$ is connected if and only if $f$ is of connected
type. 
If  $\widetilde{\Gamma}$ is not connected,
it has exactly two connected components,
each of which contains some island of $\Gamma_{\min}$. 

We set  
$\mathscr{M}_1 := 
\mathscr{L}_1 \otimes \omega_{\mathscr{X}/R}^{\otimes -1}
=
\mathscr{M}
( - \mathscr{V}_e )$.
Then 
\[
D_{\mathscr{M}_1} = (D_{\mathscr{M}} - [v_1] ) - [v_2] +2 [v]. 
\]
Remark that $D_{\mathscr{M}_1}$ is effective over $\Gamma_{\min} \setminus
\{ v_2 \}$.

Suppose that $D_{\mathscr{M}_1} (v_2) \geq 0$.
Then $D_{\mathscr{M}_1}$ is effective. 
If $\widetilde{\Gamma}$ is connected, then we note from condition (iv) that 
$\deg(\rest{D_{\mathscr{M}_1}}{\widetilde{\Gamma}}) \geq 
(g+1) - \deg(\rest{D_{\mathscr{M}}}{\relin(f)}) \geq g > 0$. 
If $\widetilde{\Gamma}$ is not connected, then  
by condition (iii) above, $D_{\mathscr{M}_1}$ has positive degree
on any connected component of $\widetilde{\Gamma}$.
Thus Lemma~\ref{lemma:freeatnode:graphversion:emptyset}
concludes that $\mathscr{L}_1$ is free at $q$.

Suppose that $D_{\mathscr{M}_1} (v_2) \leq  -1$.
To show that $\mathscr{L}_1$ is free at $q$, we 
use Lemma~\ref{lemma:freeatnode:graphversion:nonemptyset}.
Since $D_{\mathscr{M}} - [v_1]$ is effective, 
we then have $D_{\mathscr{M}_1} (v_2) = -1$
and $D_{\mathscr{M}_1} (x) \geq 0$ for any other $x \in \Gamma$.
If $\widetilde{\Gamma}'$ is a connected component of $\widetilde{\Gamma}$
with $v_2 \notin \widetilde{\Gamma}'$,
then $\rest{D_{\mathscr{M}_1}}{\widetilde{\Gamma}'}$ has positive degree
by condition (iii) above.

Let $\widetilde{\Gamma}_2$ be the connected component of 
$\widetilde{\Gamma}$ with $v_2 \in \widetilde{\Gamma}_2$.
Since $e$ is a connected subspace of $\widetilde{\Gamma}$
and has a common point $v_2$ with $\widetilde{\Gamma}_2$,
we have
$e \subset \widetilde{\Gamma}_2$.
Let $\widetilde{\Gamma}_3^{\circ}$ be the connected component
of $\widetilde{\Gamma}_2 \setminus \{ v_2 \}$ containing $\mathrm{relin} (e)$.
Since $\deg
\left(
\rest{D_{\mathscr{M}_1}}{\mathrm{relin} (e)}
\right) 
\geq
D_{\mathscr{M}_1} (v)
\geq 2$, 
we have $\deg \left(
\rest{D_{\mathscr{M}_1}}{\widetilde{\Gamma}_3^{\circ}}
\right) \geq 2$.

To apply
Lemma~\ref{lemma:freeatnode:graphversion:nonemptyset},
we need to show that
the valence of $\Gamma_3 := \widetilde{\Gamma}_3^{\circ} \cup \{ v_2 \}$
at $v_2$ is at least~$2$.
If $f$ is of disconnected type, then since $e$ is of connected type,
we see that
$\widetilde{\Gamma}_2 \setminus \mathrm{relin}(e)$ is 
still connected.
Suppose that $f$ is of connected type.
Then $\widetilde{\Gamma}$ is connected,
and thus $\widetilde{\Gamma}_2 = \widetilde{\Gamma}$.
By the condition 
of (A), $\widetilde{\Gamma} \setminus \mathrm{relin} (e) 
= \widetilde{\Gamma}_2 \setminus \mathrm{relin} (e) $
is connected. Thus 
for Case (A), $\widetilde{\Gamma}_2 \setminus \mathrm{relin}(e)$ is connected. 
It follows that
the valence of $\widetilde{\Gamma}_3$ at $v_2$ is at least $2$;
otherwise, it equals $1$,
and
since $\widetilde{\Gamma}_3^{\circ}$ is 
a connected component of $\widetilde{\Gamma} \setminus \{ v_2 \}$,
$\widetilde{\Gamma}_2 \setminus \mathrm{relin}(e)$ should not be
connected.
Thus $\widetilde{\Gamma}_3$ at $v_2$ is
at least $2$.

Now, using Lemma~\ref{lemma:freeatnode:graphversion:nonemptyset},
we obtain that
$\mathscr{L}_1$ is free at $q$ for Case (A).

\medskip
{\bf Case (B).}
Suppose that $f$ is of connected type and
$\Gamma_{\min} \setminus ( \mathrm{relin} (e) \cup \mathrm{relin} (f) )$
is not connected.

In this case, since $e$ is not a loop,
$v_1 \neq v_2$. Further, $\Gamma_{\min} \setminus ( \mathrm{relin} (e) \cup \mathrm{relin} (f) )$
has exactly two connected components. We denote these connected components 
by $\Gamma_1$ and $\Gamma_2$, where $v_1 \in \Gamma_1$ and $v_2 \in \Gamma_2$.  
Then by Proposition~\ref{prop:separation:special2},
there exists a model $( \mathscr{X} , \mathscr{L})$
of $(X,L)$ with the following properties: 
\begin{enumerate}
\item[(a)]
$\mathscr{X}$ is a good model;
\item[(b)]
$D_{\mathscr{M}}$ is effective, and
$\rest{D_{\mathscr{M}}}{\Gamma_1}$ 
and $\rest{D_{\mathscr{M}}}{\Gamma_2}$
have positive degree.
\end{enumerate}

We prove that $( \mathscr{X} , \mathscr{L})$ enjoys the required properties. 
Take any $p \in \Sing ( \mathscr{X}_s )_{\subset e}$. 
Since $e$ is of connected type and $\Delta_p \subset e$ , we see that $\Gamma_{\min} \setminus 
\relin(\Delta_p)$ is connected. Since $\Gamma_{\min} \setminus 
\relin(\Delta_p)$ contains $\Gamma_1$ and $\Gamma_2$, 
we have $\deg(\rest{\mathscr{M}}{\Gamma_{\min} \setminus 
\relin(\Delta_p)}) = 
\deg\left(\rest{D_{\mathscr{M}}}{\Gamma_1}\right) 
+ \deg\left(\rest{D_{\mathscr{M}}}{\Gamma_2}\right) 
\geq 2 > 0$ by condition (b) above. 
Then Lemma~\ref{lemma:freeatnode:graphversion:emptyset} shows that 
$\Lscr$ is free at $p$. Thus property (ii) of Proposition~\ref{prop:sep:edges} is verified. 

Let us verify property (iii) of Proposition~\ref{prop:sep:edges}. 
We fix any $q \in \Sing(\Xscr_s)_{\subset f}$. Our goal is to 
show that $\mathscr{L}_1 :=  \mathscr{L}(- \mathscr{V}_e)$ is free at $q$. 
As in Case (A), we set 
$\widetilde{\Gamma} := \Gamma_{\min}\setminus \relin(\Delta_q)$. 
Since $f$ is of connected type, 
$\widetilde{\Gamma}$ is connected. 
By the condition of Case (B), $\widetilde{\Gamma} \setminus 
\mathrm{relin} (e)
$
is not connected.
Let $\widetilde{\Gamma}_1$ 
and $\widetilde{\Gamma}_2$ 
be the connected components 
of $\widetilde{\Gamma} \setminus 
\mathrm{relin} (e)$
with $v_1 \in \widetilde{\Gamma}_1$
and $v_2 \in \widetilde{\Gamma}_2$.
Note that $\Gamma_1 \subset \widetilde{\Gamma}_1$ and 
$\Gamma_2 \subset \widetilde{\Gamma}_2$.  

Set $ \mathscr{M}_1 :=
\mathscr{L}_1 \otimes \omega_{\mathscr{X}/R}^{\otimes -1}
= \mathscr{M}(- \mathscr{V}_e )$ as in Case (A), so that 
$D_{\mathscr{M}_1} = D_{\mathscr{M}} - [v_1] - [v_2] + 2 [v]$.
We note that
\[
V_1 := 
\{ v \in V ( \mathscr{X} ) \mid
D_{\mathscr{M}_1} (v) \leq -1 \} \subseteq \{ v_1 , v_2 \}
.
\]
If $V_1 = \emptyset$,
then $\rest{D_{\mathscr{M}_1}}{\widetilde{\Gamma}}$ is effective
and $\deg\left(\rest{D_{\mathscr{M}_1}}{\widetilde{\Gamma}}\right) \geq 
\deg\left(\rest{D_{\mathscr{M}}}{\Gamma_1}\right) 
+ \deg\left(\rest{D_{\mathscr{M}}}{\Gamma_2}\right) 
\geq 2 > 0$, 
and hence
Lemma~\ref{lemma:freeatnode:graphversion:emptyset}
concludes that
$\mathscr{L}_1$ is free at $q$.

We
suppose that $V_1 \neq \emptyset$.
We want to use
Lemma~\ref{lemma:freeatnode:graphversion:nonemptyset}
to prove that $\mathscr{L}_1$ is free at $q$.
First, 
suppose moreover that
$V_1 = \{ v_1 \}$
or $V_1 = \{ v_2 \}$.
By symmetry, we may and do assume that $V_1 = \{ v_2 \}$.
Then it follows from $D_{\Mscr} \geq 0$ that $D_{\mathscr{M}_1} (v_2) = -1$. 
We set 
$\widetilde{\Gamma}_1^{\prime\circ} := 
\widetilde{\Gamma}_1 \cup \mathrm{relin} (e) $
and $\widetilde{\Gamma}_2^{\prime\circ} :=\widetilde{\Gamma}_2 \setminus \{ v_2 \}$.
We have $\deg 
\left(
\rest{D_{\mathscr{M}_1}}{\widetilde{\Gamma}_1^{\prime\circ}}
\right) \geq 
\deg
\left(
\rest{D_{\mathscr{M}_1}}{\mathrm{relin}(e)}
\right) \geq D_{\mathscr{M}_1} (v) = 2 > 0$.
By property (b) above,
$\deg 
\left(
\rest{D_{\mathscr{M}}}{\widetilde{\Gamma}_2}
\right) \geq 1$.
Since $D_{\mathscr{M}_1} (v_2) = -1$,
we have $D_{\mathscr{M}} (v_2) = 0$.
It follows that $\deg 
\left(
\rest{D_{\mathscr{M}}}{\widetilde{\Gamma}_2^{\prime\circ}}
\right) \geq 1$.
Since $D_{\mathscr{M}_1}$ coincides with $D_{\mathscr{M}}$
on $\widetilde{\Gamma}_2^{\prime\circ}$,
we also have $\deg 
\left(
\rest{D_{\mathscr{M}_1}}{\widetilde{\Gamma}_2^{\prime\circ}}
\right) \geq 1$.
Then there exists a connected component $\widetilde{\Gamma}_2^{\prime\prime\circ}$ of 
$\widetilde{\Gamma}_2^{\prime\circ}$ such that 
 $\deg 
\left(
\rest{D_{\mathscr{M}_1}}{\widetilde{\Gamma}_2^{\prime\prime\circ}}
\right) \geq 1$.
We see that the valence
of
$\widetilde{\Gamma}_1^{\prime\circ} \cup \widetilde{\Gamma}_2^{\prime\prime\circ} \cup \{ v_2 \}$ at $v_2$
is at least $2$.
Then Lemma~\ref{lemma:freeatnode:graphversion:nonemptyset}
concludes that
$\mathscr{L}_1$ is free at $q$.

Finally,
 we consider the case where
$V_1 = \{ v_1 , v_2 \}$.
In this case, 
$\widetilde{\Gamma} \setminus V_1$ is the disjoint union
of 
$\widetilde{\Gamma}_1 \setminus \{ v_1 \}$,
$\widetilde{\Gamma}_2 \setminus \{ v_2 \}$,
and
$\mathrm{relin} (e)$.
Since $D_{\mathscr{M}_1} (v_1) = D_{\mathscr{M}_1} (v_2) = -1$,
we have $D_{\mathscr{M}} (v_1) =D_{\mathscr{M}} (v_2) = 0$.
Since $\deg
\left( \rest{D_{\mathscr{M}}}{\widetilde{\Gamma}_1}
\right) \geq 1$
and $\deg
\left( \rest{D_{\mathscr{M}}}{\widetilde{\Gamma}_2}
\right) \geq 1$,
it follows that
$\deg
\left( \rest{D_{\mathscr{M}}}{\widetilde{\Gamma}_1 \setminus \{ v_1 \}}
\right) \geq 1$
and $\deg
\left( \rest{D_{\mathscr{M}}}{\widetilde{\Gamma}_2 \setminus \{ v_2 \}}
\right) \geq 1$.
Thus 
$\deg
\left( \rest{D_{\mathscr{M}_1}}{\widetilde{\Gamma}_1 \setminus \{ v_1 \}}
\right) \geq 1$
and $\deg
\left( \rest{D_{\mathscr{M}_1}}{\widetilde{\Gamma}_2 \setminus \{ v_2 \}}
\right) \geq 1$.
It follows that there exist connected components 
$\widetilde{\Gamma}_1^{\circ}$ and $\widetilde{\Gamma}_2^{\circ}$
of $\widetilde{\Gamma}_1 \setminus \{ v_1 \}$
and $\widetilde{\Gamma}_2 \setminus \{ v_2 \}$,
respectively,
such that
$\deg
\left( \rest{D_{\mathscr{M}_1}}{\widetilde{\Gamma}_1^{\circ}}
\right) \geq 1$
and $\deg
\left( \rest{D_{\mathscr{M}_1}}{\widetilde{\Gamma}_2^{\circ}}
\right) \geq 1$.
We also have
$\deg
\left( \rest{D_{\mathscr{M}_1}}{\mathrm{relin} (e)}
\right) \geq D_{\mathscr{M}_1} (v) = 2$.
Furthermore,
for each $i=1, 2$, 
the valence of 
$\widetilde{\Gamma}_1^{\circ}
\cup
\widetilde{\Gamma}_2^{\circ} 
\cup
\mathrm{relin}(e)
\cup
V_1$ 
at $v_i$ is at least $2$. 
(Indeed, in this subgraph, $e$ is an edge
emanating from $v_i$,
and 
there exists
in 
$\widetilde{\Gamma}_{i}^{\circ} \cup \{ v_i \}$ 
another edge emanating from $v_i$.)
Thus
Lemma~\ref{lemma:freeatnode:graphversion:nonemptyset} concludes that 
$\mathscr{L}_1$ is free at $q$.
This completes the proof of Proposition~\ref{prop:sep:edges}.
\QED

\subsection{Separating vertices}
\label{subsection:VS}

In this subsection, we consider the separation of vertices. 
Recall that
if $\Xscr$ is a Deligne--Mumford strictly semistable model of $X$ (for example, if $\Xscr$ is a good model), 
then the associated skeleton $S(\Xscr)$ is equal to the minimal skeleton $\Gamma_{\min}$.

\begin{Definition}[vertex-base section, $v_1$-base section; vertices-separating section, 
$(v_1, v_2)$-sep\-a\-rating section]
Let $v_1 , v_2 \in V ( \mathscr{X}^{\st})$ be vertices. 
Let $( \mathscr{X} , \mathscr{L} )$ be a 
model of $(X,L)$ such that $\Xscr$ is a good model.  
Let $C_1$ and $C_2$ denote the irreducible components of $\mathscr{X}_s$
with $[C_1] = v_1$ and $[C_2] = v_2$.
\begin{enumerate}
\item
We call $\widetilde{s} \in H^0 ( \mathscr{L} )$
a \emph{$v_1$-base section} if 
there exists $p \in C_1$ with $\widetilde{s} (p) \neq 0$. 
\item
Suppose that $v_1 \neq v_2$.
We call $\widetilde{s} 
\in H^0 ( \mathscr{L})$ a \emph{$(v_1 , v_2)$-separating section} if 
$\rest{\widetilde{s}}{C_1} = 0$ in $H^0(\rest{\mathscr{L}}{C_1})$ and 
and
$\rest{\widetilde{s}}{C_2} \neq 0$ in $H^0(\rest{\mathscr{L}}{C_2})$.  
\end{enumerate}
\end{Definition}

The following lemma shows that
once one has a model with 
a vertex-base section and a vertices-separating section, 
then one can separate vertices. 

\begin{Lemma}
\label{lemma:vertices-separating}
Let $v_1 , v_2 \in V ( \mathscr{X}^{\st})$ be distinct vertices. 
Suppose that there exists a model 
$( \mathscr{X} , \mathscr{L} )$ of $(X,L)$ such that $\Xscr$ is a good model 
and that there exists a $v_1$-base section $\widetilde{s}_0$ 
and a $(v_1 , v_2)$-separating section $\widetilde{s}_1$. 
We set $g := \widetilde{s}_1 / \widetilde{s}_0$ and define 
a nonzero rational function $h$ on $X$ by $h := \rest{g}{X}$. 
Then $- \log |h (v_1)| > 0$ and $ - \log | h  (v_2) | \leq 0$. 
In particular, we have $- \log |h (v_1)|  \neq 
- \log | h  (v_2) |$. 
\end{Lemma}

\Proof
Let $C_1$ and $C_2$ be the irreducible components of $\mathscr{X}_s$
with $[C_1] = v_1$ and $[C_2] = v_2$.  
Then $\ord_{C_1} (g) > 0$ and $\ord_{C_2}( g )  \leq 0$.
Thus the lemma follows.
\QED

Our task is then to show the existence of a model of $(X,L)$ that 
has an vertex-base section and a vertices-separating section.
The key is the following proposition.

\begin{Proposition} \label{prop:existence:sep:sections}
Let $v_1 , v_2 \in V ( \mathscr{X}^{\st})$ be distinct vertices.
Assume that $v_1$ and $v_2$ do not belong to the same island of $\Gamma_{\min}$
(cf. Definition~\ref{def:island}).
Further, assume that $\deg (L) \geq 3g-1$. 
Then
there exists a model $( \mathscr{X} , \mathscr{L} )$ of $(X,L)$ 
such that $\Xscr$ is a good model and 
such that \textup{(i)} or \textup{(ii)} holds. 
\begin{enumerate}
\item[(i)]
$\mathscr{L}$ has a $v_1$-base section
and a $(v_1 , v_2)$-separating section. 
\item[(ii)]
$\mathscr{L}$ has a $v_2$-base section 
and a $(v_2 , v_1)$-separating section.
\end{enumerate}
\end{Proposition}

The proof of Proposition~\ref{prop:existence:sep:sections}
will require 
a careful analysis 
of $\Gamma_{\min}$. 
Before moving on to the proof, we show the 
following lemma, which gives a sufficient condition in terms of $\Lambda$-metric graphs 
for a model to have a vertices-separated section.

\begin{Lemma} \label{lemma:exist:vert:sep:sections}
Let $v_1 , v_2 \in V ( \mathscr{X}^{\st})$ be distinct vertices.
Let $( \mathscr{X} , \mathscr{L} )$ be model of $(X,L)$ such that $\Xscr$ is a good model. 
Assume that there exists a $v_1' \in \Gamma_{\min} \setminus \{ v_2 \}$
such that $v_1$ and $v_1'$ belong to the same connected component
of $\Gamma_{\min} \setminus \{ v_2 \}$
and such that
$D_{\mathscr{M}} - [v_1 ] - [v_1']$ is effective.
Then there exists a $(v_2 , v_1)$-separating section of $\mathscr{L}$.
\end{Lemma}

\Proof
Since $\mathscr{X}$ dominates $\mathscr{X}^{\st}$,
we note that $v_1 , v_2 \in V( \mathscr{X} )$.
Further,
since $D_{\mathscr{M}} - [v_1 ] - [v_1']$ is effective,
we note that $v_1' \in V( \mathscr{X} )$.
Let $C_1$, $C_2$, and $C_1'$ be the irreducible components 
of $\mathscr{X}_s$ with $[C_1] = v_1$, $[C_2] = v_2$, 
and $[C_1'] = v_1'$. Let 
$\Gamma_1^\circ$ be the connected component of $\Gamma \setminus \{ v_2 \}$
on which $v_1$ and $v_1'$ lie.  
Set $F_1 := \sum_{C \in \Irr ( \mathscr{X}_s ), [C] \in \Gamma_1^\circ} C $.
Fix a point $p_1 \in C_1 (k) \setminus \Sing ( \mathscr{X}_s)$.
Set $\Sigma := F_1 \cap C_2$,
which we regard as a nontrivial effective Cartier divisor on $C_2$.
Since $D_{\mathscr{M}}$ is effective, $\rest{\mathscr{M}}{\mathscr{X}_s}$ is nef.

\begin{figure}[!h]
\[
\setlength\unitlength{0.08truecm}
\begin{picture}(200, 60)(0,0)
 \put(5, 20){\line(1,0){80}}
 \put(70, 20){\line(1,1){15}}
 \put(70, 20){\line(-1,1){15}}
  \put(10, 20){\circle*{1.5}}
  \put(50, 20){\circle*{1.5}}
  \put(70, 20){\circle*{1.5}}
   \put(40, 40){$\Gamma_{\min}$}
  \put(8, 14){$v_1^\prime$}
  \put(48, 14){$v_1$}
  \put(68, 14){$v_2$}
  \put(5, 12){$\underbrace{\phantom{AAAAAAAAAAAAAAAaa}}$} 
  \put(36, 3){$\Gamma_1^\circ$}
 \put(100, 10){\line(1,0){10}}
 \multiput(110,10)(1,0){10}{\line(1,0){0.5}}
 \put(120, 10){\line(1,0){10}}
 \put(135, 25){\line(1,0){10}}
 \multiput(145,25)(1,0){10}{\line(1,0){0.5}}
 \put(155, 25){\line(1,0){10}}
 \put(163, 33){\line(1,0){15}}
 \put(168, 38){\line(1,0){15}}
 \put(173, 43){\line(1,0){15}}
 \put(133, 18){\circle*{1.5}}
 \put(160, 25){\circle*{1.5}}
  \put(115, 23){$C_1^\prime$}
  \put(143, 31){$C_1$}
  \put(183, 50){$C_2$}
  \put(135, 15){$p_1$}
  \put(160, 19){$\Sigma$}
  \put(100, 32){$\overbrace{\phantom{AAAAAAAAAAAAAAA}}$} 
  \put(127, 42){$F_1$}
 \thicklines
  \put(100, 5){\line(1,1){15}}
  \put(123, 8){\line(1,1){20}}
  \put(158, 23){\line(1,1){25}}
  \end{picture}
\]
\label{figure:for:lemma:exist:vert:sep:sections}
\end{figure}

We claim  that $\rest{\mathscr{L}}{F_1} ( - \Sigma )$ is free at $p_1$. Indeed, 
since $\rest{\mathscr{M}}{F_1}$ is nef
and
since
$\deg 
\left(
\rest{\Mscr}{C_1}
\right) \geq 1$
and
$
\deg 
\left(
\rest{\mathscr{M}}{C_1'}
\right) 
\geq 1$,
it follows 
that $\rest{\mathscr{M}}{F_1} ( -p_1 )$ is nef 
and $\deg\left(\rest{\mathscr{M}}{F_1} ( -p_1 )\right) > 0$.
Since $F_1$ is connected, Lemma~\ref{lemma:vanishing0}
gives   
$h^0 \left( 
\left(
\rest{\mathscr{M}}{F_1} ( -p_1 )
\right)^{\otimes -1}
\right) = 0$. 
By the adjunction formula, 
we have
\[
\rest{\mathscr{L} }{F_1} (- \Sigma -p_1) \otimes \omega_{F_1}^{\otimes -1}
=
\rest{\mathscr{L} \otimes \omega_{\mathscr{X}/R}^{\otimes -1}}{F_1} (-p_1)
=
\rest{\mathscr{M}}{F_1} ( -p_1 )
\]
By the Serre duality,
we obtain 
$h^1 \left( \rest{\mathscr{L} }{F_1} (- \Sigma -p_1) \right) =
h^0 \left( 
\left(
\rest{\mathscr{M}}{F_1} ( -p_1 )
\right)^{\otimes -1}
\right) = 0$, and then 
Lemma~\ref{lemma:nonbasepoint} shows that 
$\rest{\mathscr{L}}{F_1} (- \Sigma) $ is free at $p_1$.

It follows that there exists a global section $\eta$
of $\rest{\mathscr{L} }{F_1} (- \Sigma) $ such that $\eta (p_1) \neq 0$.
Note that $\rest{\mathscr{L} }{F_1} (- \Sigma) 
\subset \rest{\mathscr{L} }{F_1} $.
Extending $\eta$ to be zero outside of $F_1$, we obtain a global section
$\eta_1$ of $\rest{\mathscr{L}}{\mathscr{X}_s}$
such that $\rest{\eta_1}{F_1} = \eta$ as a global section of 
$\rest{\mathscr{L}}{F_1}$.
Since $C_2 \notin \Irr (F_1)$, we have
$\rest{\eta_1}{C_2} = 0$.

Since $\rest{\mathscr{M}}{\mathscr{X}_s}$ 
is nef and has positive degree,
$h^0
\left( \rest{\mathscr{M}}{\mathscr{X}_s}^{\otimes -1} \right) = 0$,
and thus
$h^1 \left( \rest{\mathscr{L}}{\mathscr{X}_s} \right) = 0$
by the Serre duality.
By the base-change theorem,
the natural map $H^0 ( \mathscr{L} ) \to 
H^0 \left( \rest{\mathscr{L}}{\mathscr{X}_s} \right)$
is surjective.
It follows that
there exists a global section $\widetilde{s}_1 \in H^0 ( \mathscr{L} )$
such that $\rest{\widetilde{s}_1}{\mathscr{X}_s} = \eta_1$.
Then $\rest{\widetilde{s}_1}{C_1} \neq 0$ and $\rest{\widetilde{s}_1}{C_2} = 0$, so that 
$\widetilde{s}_1$ is a $(v_2 , v_1)$-separating section of $\mathscr{L}$.
\QED

We start the proof of Proposition~\ref{prop:existence:sep:sections}.

\medskip
\textsl{Proof of Proposition~\ref{prop:existence:sep:sections}.}\quad
Let $e \in E(\Xscr^{\st})$ be any edge of disconnected type.  
For a good model $\Xscr$ of $X$, 
we denote by $C_1$ and $C_2$ the irreducible components of $\Xscr_s$ with 
$v_1  = [C_1]$ and $v_2 = [C_2]$. 

\medskip

{\bf Case 1.}
Suppose that $\Gamma_{\min} \setminus \{ v_1 \}$ is connected or 
$\Gamma_{\min} \setminus \{ v_2 \}$ is connected.
Without loss of generality, we may assume that $\Gamma_{\min} \setminus \{ v_2 \}$ is connected. 
By Proposition~\ref{prop:pregoodmodel},
there exists a model $( \mathscr{X} , \mathscr{L})$
of $(X,L)$ with the following properties: 
\begin{enumerate}
\item[(i)]
$\mathscr{X}$ is a good model;
\item[(ii)]
$D_{\mathscr{M}} - [v_1]$ is effective. 
\end{enumerate}

First, suppose that $D_{\mathscr{M}} (v_2) \geq 2$.
Then applying Lemma~\ref{lemma:exist:vert:sep:sections} 
with
$v_2$, $ v_1 $, and $ v_2$ 
in place of $v_1$, $v_2 $, and $ v_1'$,
respectively,
we see that
$\mathscr{L}$ has a $(v_1 , v_2)$-separating section.
Further, since $D_{\mathscr{M}} - [v_1]$ is effective,
Lemma~\ref{lemma:freeatregular:graphversion:effective} shows that
$\mathscr{L}$ is free at any point of $C_1(k) \setminus \Sing ( \mathscr{X}_s )$.
Thus $\mathscr{L}$ has a 
$v_1$-base section,
which proves the proposition.

Next, suppose that 
$D_{\mathscr{M}} (v_2) \leq 1$.
Then,
since $\deg ( D_{\mathscr{M}} - [v_1] ) \geq (g+1) - 1\geq 2$,
we have $\deg \left(
\rest{(D_{\mathscr{M}} - [v_1])}{\Gamma_{\min} \setminus \{ v_2 \}}
\right) \geq 1$.
It follows that
there exists a $v_1' \in \Gamma_{\min}
\setminus \{ v_2 \}$
such that $D_{\mathscr{M}} - [v_1] - [v_1']$ is effective. 
Since $\Gamma_{\min} \setminus \{ v_2 \}$ is assumed to be connected, 
Lemma~\ref{lemma:exist:vert:sep:sections} gives that 
$\mathscr{L}$ has a $(v_2 , v_1)$-separating section.

We need to prove that $\mathscr{L}$ has a $v_2$-base section. 
Since $D_{\mathscr{M}} - [v_1]$ is effective and 
we are assuming that $D_{\mathscr{M}} (v_2) \leq 1$, 
we have $0 \leq D_{\mathscr{M}} (v_2) \leq 1$.
If $D_{\mathscr{M}} (v_2) = 1$, 
then $D_{\mathscr{M}} - [v_2]$ is effective and $\deg ( D_{\mathscr{M}} - [v_2] ) \geq (g+1) - 1\geq 2 > 0$,
so that 
Lemma~\ref{lemma:freeatregular:graphversion:effective} 
gives that $\mathscr{L}$ is free
at any $q \in C_{2}(k) \setminus \Sing ( \mathscr{X}_s)$;
thus $\mathscr{L}$ has a $v_2$-base section.

Suppose that $D_{\mathscr{M}} (v_2) = 0$. 
If the valence of $v_2$ is equal to $1$, then 
the edge $e$ with an end vertex $v_2$ is of disconnected type and 
$v_1$ and $v_2$ belong to different connected components
of $\Gamma_{\min} \setminus \mathrm{relin}(e)$, 
which contradicts the assumption 
that $v_1$ and $v_2$ belong to the same island of $\Gamma_{\min}$. 
Thus the valence of $v_2$ is at least $2$. 
Further, $\Gamma_{\min} \setminus \{ v_2 \}$ is assumed to be connected. 
By Lemma~\ref{lemma:freeatregular:graphversion},
$\mathscr{L}$ is free at any point of $C_2(k) \setminus \Sing (\mathscr{X}_s)$.
Thus $\mathscr{L}$ has a $v_2$-base section.

To sum up, we have shown that there exists 
a model $( \mathscr{X} , \mathscr{L})$ such that 
(i) or (ii) of Proposition~\ref{prop:existence:sep:sections} holds
in Case 1.

\medskip
{\bf Case 2.}
Suppose that $\Gamma_{\min} \setminus \{ v_1 \}$ is not connected and $\Gamma_{\min} \setminus \{ v_2 \}$ 
is not connected.

\smallskip
{\bf Subcase 2.1.}
Suppose that $\Gamma_{\min}$ has an edge $e$ of disconnected type.
Then one of the following holds:
\begin{enumerate}
\item[(I)]
$v_1$ and $\mathrm{relin}(e)$ are contained in
the same connected component of $\Gamma_{\min} \setminus \{ v_2 \}$;
\item[(II)]
$v_2$ and $\mathrm{relin}(e)$ are contained in
the same connected component of $\Gamma_{\min} \setminus \{ v_1 \}$
\end{enumerate} 
Indeed, suppose that (II) does not hold. We take any point $w \in \mathrm{relin}(e)$. Then there exists a continuous map $\gamma: [0, 1] \to \Gamma_{\min}$ such that $\gamma(0) = v_2$, $\gamma(1/2) = v_1$ and 
$\gamma(1) = w$ and such that $\gamma(t) \neq v_2$ for any $t \ne 0$ and  
$\gamma(t) \neq v_1$ for any $t \neq 1/2$. The half path $\gamma: [1/2, 1] \to \Gamma_{\min}$ 
connects $v_1$ and $w$ in $\Gamma_{\min} \setminus \{ v_2 \}$, and thus (I) holds. 

Without loss of generality, we may and do assume that (I) holds. 
Let $\Gamma_1^{\circ}$ be the connected component of $\Gamma_{\min}
\setminus \{ v_2 \}$ with $v_1 \in \Gamma_1^{\circ}$, and 
set $\Gamma_1 := \Gamma_1^{\circ} \cup \{ v_2 \}$, which is
 the closure of 
$\Gamma_1^{\circ}$ in $\Gamma_{\min}$. 
Then $e \subset \Gamma_1$. 

By Proposition~\ref{prop:pregoodmodel},
there exists a model $( \mathscr{X} , \mathscr{L})$ 
of $(X,L)$ such that $\Xscr$ is a good model with the following properties. 
\begin{enumerate}
\item[(i)]
$D_{\mathscr{M}} - [v_1]$ is effective;
\item[(ii)]
$\deg
\left(
\rest{\left(
D_{\mathscr{M}} - [v_1]
\right)}{\Gamma'}
\right) \geq 1$ for any island $\Gamma'$ of $\Gamma_{\min}$.
\end{enumerate}

Since $e \subset \Gamma_1$,
there exists 
an island $\Gamma_1'$ of $\Gamma_{\min}$
satisfying $\Gamma_1' \subset \Gamma_1 \setminus \{ v_2 \}$.
By conditions (i) and (ii) above, 
there exists a $v_1' \in \Gamma_1'$
such that $D_{\mathscr{M}} - [v_1] - [v_1']$ is effective. 
Note that $\Gamma_1 \setminus \{ v_2 \} = \Gamma_1^\circ$ and that 
$v_1, v_1^\prime \in \Gamma_1^\circ$. 
Since $v_1$ and $v_1'$ belong to the 
same 
connected component
$\Gamma_1^\circ$ of $\Gamma_{\min} \setminus \{ v_2 \}$,
Lemma~\ref{lemma:exist:vert:sep:sections}
concludes that $\mathscr{L}$ has a $(v_2 , v_1)$-separating section.

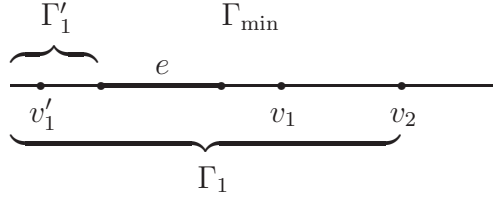
\begin{figure}[!h]
\[
\setlength\unitlength{0.08truecm}
\begin{picture}(90, 40)(0,0)
 \put(5, 20){\line(1,0){80}}
  \put(10, 20){\circle*{1.5}}
  \put(50, 20){\circle*{1.5}}
  \put(70, 20){\circle*{1.5}}
  \put(20, 20){\circle*{1.5}}
  \put(40, 20){\circle*{1.5}}
  \put(40, 30){$\Gamma_{\min}$}
  \put(8, 14){$v_1^\prime$}
  \put(48, 14){$v_1$}
  \put(68, 14){$v_2$}
  \put(29, 22){$e$}
  \put(36, 3){$\Gamma_1$}
  \put(10, 30){$\Gamma_1^\prime$}
  \thicklines
  \put(20, 19.9){\line(1,0){20}}
  \put(20, 20.1){\line(1,0){20}}
  \put(20, 20){\line(1,0){20}}
  \put(5, 12){$\underbrace{\phantom{AAAAAAAAAAAAAAAAa}}$} 
  \put(5, 20){$\overbrace{\phantom{AAAa}}$} 
  \end{picture}
\]
\label{figure:for:prop:existence:sep:sections:case2}
\end{figure}

In the rest of this subcase,
we construct a $v_2$-base section of $\mathscr{L}$. 
Since $D_{\mathscr{M}} - [v_1]$ is effective, we have 
$D_{\mathscr{M}}  (v_2) \geq 0$. 

Suppose that $D_{\mathscr{M}}  (v_2) \geq 1$.
Then $D_{\mathscr{M}}  - [v_2] $ is effective and $\deg ( D_{\mathscr{M}} - [v_2] ) \geq (g+1) - 1\geq 2 > 0$. 
By Lemma~\ref{lemma:freeatregular:graphversion:effective}, 
$\mathscr{L}$ is free
at any $q \in C_{2}(k) \setminus \Sing ( \mathscr{X}_s)$, and 
$\mathscr{L}$ has a $v_2$-base section.

Suppose that $D_{\mathscr{M}}  (v_2) = 0$.
Let $\Gamma_2'$ be the island of $\Gamma_{\min}$ with $v_2 \in \Gamma_2'$.
Since $D_{\mathscr{M}}  (v_2) = 0$, 
it follows from condition (ii) above that
there exists a connected component $\Gamma_2^{\prime\prime\circ}$
of $\Gamma_2' \setminus \{ v_2 \}$
such that
$\deg
\left(
\rest{D_{\mathscr{M}}}{\Gamma_2^{\prime\prime\circ}}
\right) \geq 1$. 
Let $\Gamma_2^{\prime\prime} := \Gamma_2^{\prime\prime\circ} \cup \{v_2\}$ 
be the closure of $\Gamma_2^{\prime\prime}$ in $\Gamma_2'$ (and thus 
the closure of $\Gamma_2^{\prime\prime}$ in $\Gamma_{\min}$). 

We claim that the valence of $\Gamma_2^{\prime\prime}$ at $v_2$
is at least $2$. 
To derive a contradiction, assume that the valence is $1$.
Then there exists a unique edge $e_2$ with end vertex $v_2$ such that $e_2 \subset \Gamma_2''$.
Since the valence is $1$, $e_2$ is not a loop.
Further, $v_2$ and $\Gamma_2^{\prime\prime\circ} \setminus \relin(e_2)$ 
do not belong to the same connected component of $\Gamma_2' \setminus \mathrm{relin}(e_2)$. 
It follows that 
$\Gamma_2' \setminus \mathrm{relin}(e_2)$ is not connected. 
On the other hand, since $e_2 \subset \Gamma_2'$, this contradicts 
$\Gamma_2'$ being an island of $\Gamma_{\min}$.

Now,
wince the valence of $\Gamma_2^{\prime\prime}$ at $v_2$
is at least $2$,
Lemma~\ref{lemma:freeatregular:graphversion}
shows that 
$\mathscr{L}$ is free at any $q \in C_2(k) \setminus
\Sing ( \mathscr{X}_s)$.
Thus $\mathscr{L}$ has a $v_2$-base section.

\medskip
{\bf Subcase~2.2.}
We consider the case
where $\Gamma_{\min}$
does not have an edge of disconnected type. 

In this case, by Lemma~\ref{lemma:goodmodel:VS}, which we prove below, 
we take a model $(\mathscr{X} , \mathscr{L})$ of $(X, L)$ satisfying conditions (i) or (ii) 
in Lemma~\ref{lemma:goodmodel:VS}.

Suppose that this model satisfies condition (i). 
Then by Lemma~\ref{lemma:exist:vert:sep:sections},
$\mathscr{L}$ has a $(v_2 , v_1)$-separating section.
Further, by 
the same argument using Lemma~\ref{lemma:freeatregular:graphversion}
as in Case~1, $\mathscr{L}$ has a 
$v_2$-base section.
Thus Proposition~\ref{prop:existence:sep:sections} holds in this case.

Suppose that $(\mathscr{X} , \mathscr{L})$ satisfies (ii) in Lemma~\ref{lemma:goodmodel:VS}.
Then by Lemma~\ref{lemma:exist:vert:sep:sections},
$\mathscr{L}$ has a $(v_1 , v_2)$-separating section.
Further, 
since $D_{\mathscr{M}} (v_1) \geq 1$, 
the same argument using Lemma~\ref{lemma:freeatregular:graphversion:effective}
as in Case~1 
shows that
$\mathscr{L}$ has a 
$v_1$-base section.
This completes the proof of
Proposition~\ref{prop:existence:sep:sections}.
\QED

\begin{Lemma}
\label{lemma:goodmodel:VS}
Let $v_1, v_2 \in V(\Xscr^{\st})$ be distinct vertices of $
\Gamma_{\min}$. 
Assume that $\deg(L) \geq 3g-1$. 
Suppose that 
$\Gamma_{\min}\setminus\{v_1\}$ and 
$\Gamma_{\min}\setminus\{v_2\}$ are not connected 
and that $\Gamma_{\min}$ does not have an edge of disconnected type. 
Let $\Gamma_1^\circ$ be the connected component of 
$\Gamma_{\min}\setminus\{v_2\}$ with $v_1 \in \Gamma_1^\circ$, 
and set $\Gamma_1 := \Gamma_1^\circ \cup\{v_2\}$,
which is the closure 
of $\Gamma_1^\circ$ in $\Gamma_{\min}$. 
Then there exists a model $(\mathscr{X} , \mathscr{L})$
of $(X,L)$
such that
$D_{\mathscr{M}} - [v_1]$ is effective
and such that
one of the following holds\textup{:} 
\begin{enumerate}
\item[(i)]
there exists a
 $v_1' \in \Gamma_1^{\circ}$ such that
$
D_{\mathscr{M}} - [v_1] - [v_1']
$ is effective;
\item[(ii)]
$D_{\mathscr{M}} (v_2) \geq 2$.
\end{enumerate}
\end{Lemma}

\Proof
Since $\Gamma_1$ contains $v_1$ and $v_2$, 
$\Gamma_1$ is not a singleton. 
We set $\Gamma_2 := \Gamma_{\min} \setminus \Gamma_1^\circ$. 
Since $\Gamma_{\min} \setminus \{ v_2 \}$ is not connected,
$\Gamma_2$ is not a singleton.
We have $\Gamma_1 \cap \Gamma_2 = \{v_2\}$ and 
$\Gamma_1 \cup \Gamma_2 = \Gamma_{\min}$, and
$\Gamma_1$ and $\Gamma_2$ are connected.
We endow $\Gamma_1$ and $\Gamma_2$ with 
a finite graph structure by restricting that of $\Gamma_{\min}$. 

We use the following claim.

\begin{Claim}
\label{claim:goodmodel:VS:divisorongraph}
Let $E \in \Div_{\Lambda} ( \Gamma_{\min} )$ with $\deg (E) \geq g$.
Suppose that for any $E' \in |E|$,
$\rest{E'}{\Gamma_1^{\circ}}$ is trivial.
Then there exists an $E'' \in |E|$ such that
$E'' (v_2) \geq 2$.
\end{Claim}

Let us show the claim in three steps. 

\medskip
{\bf Step 1.}
We prove that 
$\Gamma_1$ does not have an edge of disconnected type. 
To argue by contradiction,
assume that $\Gamma_1$ has  an edge $e$ of disconnected type.
Then 
$\Gamma_1 
\setminus \mathrm{relin}(e) $ has two connected component.
Let $\Gamma_1'$ be the connected component of $\Gamma_1 
\setminus \mathrm{relin}(e) $ with $v_2 \in \Gamma_1'$,
and let $\Gamma_1''$ be the other.
Then 
$\Gamma_1'' \cap \Gamma_2 = \emptyset$ and $\Gamma_1' \cap \Gamma_2
= \{ v_2 \}$,
and we see that
$\Gamma_1''$ and
$\Gamma_1' \cup \Gamma_2$ are connected components of
$\Gamma_{\min} \setminus \mathrm{relin}(e)$.
This means that $e$,  as an edge of $\Gamma_{\min}$,
is an edge of disconnected type.
However, this contradicts the assumption 
of $\Gamma_{\min}$ in the lemma.
Thus $\Gamma_1$ does not have an edge of disconnected type.

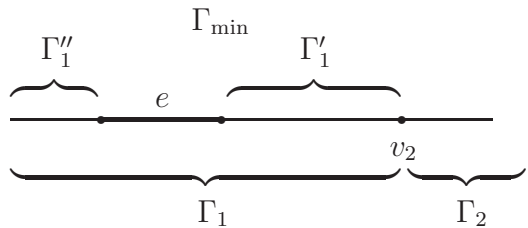
\begin{figure}[!h]
\[
\setlength\unitlength{0.08truecm}
\begin{picture}(90, 45)(0,0)
 \put(5, 20){\line(1,0){80}}
  \put(70, 20){\circle*{1.5}}
  \put(20, 20){\circle*{1.5}}
  \put(40, 20){\circle*{1.5}}
  \put(35, 35){$\Gamma_{\min}$}
  \put(68, 14){$v_2$}
  \put(29, 22){$e$}
  \put(36, 3){$\Gamma_1$}
  \put(79, 3){$\Gamma_2$}
  \put(10, 30){$\Gamma_1^{\prime\prime}$}
  \put(53, 30){$\Gamma_1^{\prime}$}
  \thicklines
  \put(20, 19.9){\line(1,0){20}}
  \put(20, 20.1){\line(1,0){20}}
  \put(20, 20){\line(1,0){20}}
  \put(5, 12){$\underbrace{\phantom{AAAAAAAAAAAAAAAAa}}$} 
  \put(5, 20){$\overbrace{\phantom{AAAa}}$} 
  \put(41, 20){$\overbrace{\phantom{AAAAAAaa}}$} 
  \put(71, 12){$\underbrace{\phantom{AAAAA}}$} 
  \end{picture}
\]
\label{figure:for:lemma:goodmodel:VS}
\end{figure}

{\bf Step 2.}
We prove $g ( \Gamma_1 ) \geq 2$.
To argue by contradiction,
suppose that $g ( \Gamma_1 ) \leq 1$.
Since $\Gamma_1$ is not a singleton and since 
$\Gamma_1$ does not have an edge of disconnected type by Step~1,
it follows that $g(\Gamma_1) \neq 0$. 
Then $g ( \Gamma_1 ) = 1$.
Since $\Gamma_1$ does not have an edge of disconnected type,
we see
that $\Gamma_1$ is a circle.
Thus $\Gamma_1 \setminus \{ v_1 \}$
is connected.
Since 
$\Gamma_2$ is connected and
$( \Gamma_1 \setminus \{ v_1 \}) \cap \Gamma_2 = \{ v_2 \}$,
it follows that
$\Gamma_{\min} \setminus \{ v_1 \}
= (\Gamma_1 \setminus \{ v_1 \}) \cup \Gamma_2$ is connected.
However, this contradicts to the assumption that 
$\Gamma \setminus \{ v_1 \}$ is not connected. 
Thus $g ( \Gamma_1 ) \geq 2$.

\medskip
{\bf Step 3.}
Since $\deg (E) \geq g$,
there exists an $E' \in |E|$ by Riemann's inequality (cf. Proposition~\ref{prop:riemann:inequality}).
By the assumption of the claim, 
$\Supp (E') \subset \Gamma_2$,
and hence we regard $E'$ as a divisor on $\Gamma_2$.
Since 
$g ( \Gamma_1) + g (\Gamma_2 ) = g(\Gamma_{\min}) \leq g(\bar{\Gamma}_{\min}) = g$,
Step~2
gives 
$g ( \Gamma_2 ) \leq g-2$.
Thus $\deg (E' - 2[v_2]) \geq g (\Gamma_2)$.
It follows from
Proposition~\ref{prop:RIforMG}
that there exists an effective
divisor $E'_2 \in \Div_{\Lambda} ( \Gamma_2)$ such that $E'_2 + 2[v_2] \sim E'$
as divisors on $\Gamma_2$.
Regard $E'_2 + 2[v_2]$ as a $\Lambda$-divisor on $\Gamma$,
which we denote by $E''$.
Then,
since
$\Gamma_1 \cup \Gamma_2 = \Gamma_{\min}$ and $\Gamma_1 \cap \Gamma_2
= \{ v_2 \}$,
we have
 $E'' \sim E$ on $\Gamma$ as well, and $E'' (v_2) \geq 2$.
This completes the proof of the claim.

\medskip
We resume the proof of Lemma~\ref{lemma:goodmodel:VS}.
Take a divisor $\tilde{D}_1 \in \Div ( X )$ such that
$L \otimes \omega_{X}^{\otimes -1} \cong \OO_X ( \tilde{D}_1 )$. 
Recall from \eqref{eqn:specialization:map} that $\tau_*: \Div(X) \to \Div_{\Lambda} ( \Gamma_{\min} )$ denotes the specialization map.
Set $E := \tau_{\ast} ( \tilde{D}_1) - [v_1]$.
Then $\deg (E) = \deg( \tilde{D}_1) - 1 \geq (g+1) - 1 = g$.
By Proposition~\ref{prop:riemann:inequality},
$|E| \neq \emptyset$.

We divide our argument into two cases.
First,
suppose that there exists an $E' \in |E|$ 
such that
$\deg \left( \rest{E'}{\Gamma_1^{\circ}} \right) \geq 1$.
Then there exists a $v_1' \in \Gamma_{1}^{\circ}$ 
such that
$E' (v_1') \geq 1$.
We set $D := E' + [v_1]$.
Remark that $D \sim \tau_{\ast} ( \tilde{D}_1)$.
By Proposition~\ref{prop:pregoodmodel:pre},
there exists a 
good model $\mathscr{X}$ and a line bundle $\mathscr{M}$
over $\mathscr{X}$ such that
$D_{\mathscr{M}} = D$.
Set $\mathscr{L} := \mathscr{M} \otimes \omega_{\mathscr{X}/R}$.
Then $D_{\mathscr{M}} - [v_1] = D - [v_1] = E'$, which is effective.
Further,
the above construction immediately shows that
$( \mathscr{X} , \mathscr{L} )$ is a model of $(X , L)$
which satisfies condition (i) of Lemma~\ref{lemma:goodmodel:VS}.
This completes the proof of the lemma
in this case.

Next, suppose that there does not exist an $E' \in |E|$ 
such that
$\deg \left( \rest{E'}{\Gamma_1^{\circ}} \right) \geq 1$.
Then by Claim~\ref{claim:goodmodel:VS:divisorongraph},
there exists an $E'' \in |E|$ such that
$E'' (v_2) \geq 2$.
Set $D := E'' + [v_1]$.
By the same way as above,
we construct a model $(\mathscr{X} , \mathscr{L})$
such that $D_{\mathscr{M}} = D$,
where $\mathscr{M} = \mathscr{L} \otimes \omega_{\mathscr{X}/R}^{\otimes -1}$.
Note that $D_{\mathscr{M}} - [v_1] = E''$ is effective.
Further, it is straightforward to see that
this model satisfies (ii) of Lemma~\ref{lemma:goodmodel:VS}.
Thus we obtain the lemma.
\QED

\subsection{Faithful tropicalization of the minimal skeleton}
\label{subsec:ft:min:skeleton}
We construct a faithful tropicalization of 
the minimal skeleton $\Gamma_{\min}$. 

\begin{Theorem} \label{theorem:FT:canonical}
Let $X$ be a connected smooth projective curve over $K$ of genus $g \geq 2$,
and let $L$ be a line bundle over $X$.
Suppose that $\deg (L) \geq 3g-1$.
Then there exist $s_0 , \ldots , s_N \in H^0 ( X,L)$
such that the map
$\varphi : X^\an \to \TT\PP^N$
defined by $\varphi := ( - \log |s_0| : \cdots : - \log |s_N|)$
gives a faithful tropicalization of $\Gamma_{\min}$.
\end{Theorem}

\Proof
{\bf Step 1.}
By Theorem~\ref{theorem:UT:canonical},
we take
$s_0 , \ldots , s_{N_1} \in H^0 ( X,L)$
such that
the map
$\varphi_1 : X^\an \to \TT\PP^{N_1}$
defined by $\varphi_1 := ( - \log |s_0| : \cdots : - \log |s_{N_1}|)$
satisfies conditions (i)--(iii) in Theorem~\ref{theorem:UT:canonical} 
for $\varphi = \varphi_1$.
In particular, 
$\varphi_1$ is unimodular. 

\medskip
{\bf Step 2.}
We take care of the separation of points in the interior
of an edge.
We set  $\{e_{1} , \ldots , e_{\alpha}\} := \{ e \in E ( \mathscr{X}^{\st}_s ) \mid 
\text{$e$ is of connected type}\}$. 
For each $i = 1 ,\ldots , \alpha$,
Proposition~\ref{prop:model:AS} gives 
a model $( \mathscr{X}_i , \mathscr{L}_i )$
of $(X,L)$ such that $\mathscr{L}_i$ has an $e_i$-base section
$\widetilde{s}^{(e_i)}_{0}$
and a separating  section
$\widetilde{s}^{(e_i)}_{1}$ for $e_i$.
Set $s_{N_1+2i-1} := \rest{\widetilde{s}^{(e_i)}_{0}}{X}$
and 
 $s_{N_1+2i} := \rest{\widetilde{s}^{(e_i)}_{1}}{X}$
for $i = 1 , \ldots , \alpha$,
and set $N_2 := N_1 + 2 \alpha$.
Define $\varphi_2 : X^{\an} \to \TT\RR^{N_2}$
by
$\varphi_2 :=
( - \log |s_0| : \cdots : - \log |s_{N_1}| : 
- \log |s_{N_1+1}| : \cdots : - \log |s_{N_2}|)$.
Since $\varphi_1$ is unimodular, so is $\varphi_2$
(cf. Remark~\ref{remark:unimodular}).

Let $x, y \in \Gamma$ be distinct points
that belong to the interior of some common edge $e$.
If $e$ is of disconnected type, then $\varphi_1 (x) \neq \varphi_1 (y)$
by condition (i) in Theorem~\ref{theorem:UT:canonical},
and hence $\varphi_2 (x) \neq \varphi_2 (y)$.
Suppose that $e$ is of connected type.
We divide $e$ at the middle point
and let $e^\prime$ and $e^{\prime\prime}$ be the $1$-simplices arising from this division.
If $x , y \in e^\prime$ or $x , y \in e^{\prime\prime}$,
then by condition (ii) in Theorem~\ref{theorem:UT:canonical},
we have $\varphi_1  (x) \neq \varphi_1 (y)$ and hence $\varphi_2 (x) \neq \varphi_2 (y)$.
If $x \in \mathrm{relin}(e^{\prime})$ and $y \in \mathrm{relin}(e^{\prime\prime})$,
then by
and Lemma~\ref{lemma:asymmetric:section:separate:new},
we obtain $\varphi_2 (x) \neq \varphi_2 (y)$.
Thus $\varphi_2$ separate two distinct points that belong
to the interior of the same edge.

\medskip
{\bf Step 3.}
We take care of separation of two distinct points 
one of which is not a vertex.
We set 
\[
\{(e_1, f_1), \ldots, (e_\beta, f_\beta)\} 
:= 
\{
(e, f) \in E ( \mathscr{X}^{\st}) \times E ( \mathscr{X}^{\st})
\mid 
\text{$e$ is of connected type, $e \neq f$}. 
\}
\]
For each $i$,
Proposition~\ref{prop:model:ES}
gives a model
a model $( \mathscr{X}_i , \mathscr{L}_i )$
of $(X,L)$ such that $\mathscr{L}_i$ has an $e_i$-base section
$\widetilde{s}^{(e_i)}_{0}$ 
and an $(e_i, f_i)$-separating section
$\widetilde{s}^{(e_i , f_i)}_{1}$.
Set $s_{N_2+2i-1} := \rest{\widetilde{s}^{(e_i)}_{0}}{X}$
and 
 $s_{N_2+2i} := \rest{\widetilde{s}^{(e_i , f_i)}_{1}}{X}$
for $i = 1 , \ldots , \beta$,
and set $N_3 := N_2 + 2 \beta$.
Define $\varphi_3 : X^{\an} \to \TT\RR^{N_3}$
by
$\varphi_3 :=
( - \log |s_0| : \cdots : - \log |s_{N_2}| : 
- \log |s_{N_2+1}| : \cdots : - \log |s_{N_3}|)$.
Since $\varphi_2$ is unimodular, so is $\varphi_3$
(cf. Remark~\ref{remark:unimodular}).

Let $x , y \in \Gamma$ be distinct points.
Assuming that not both of $x ,y$ are vertices, let us 
show that $\varphi_3 (x) \neq \varphi_3 (y)$.
Without loss of generality, we may and do assume that there exists
an edge $e$ with $x \in \mathrm{relin}(e)$.
If $y \in \mathrm{relin}(e)$,
then
 $\varphi_2 (x) \neq \varphi_2 (y)$
by Step~2,
and hence
$\varphi_3 (x) \neq \varphi_3 (y)$.
Suppose that $y \notin \mathrm{relin}(e)$.
Then there exists an edge $f$ with $f \neq e$
such that $y \in f$.
If $e$ is of disconnected type,
then $\varphi_1 (x) \neq \varphi_1 (y)$
as noted in
Remark~\ref{remark:disconnect-other},
and hence $\varphi_3 (x) \neq \varphi_3 (y)$.
Suppose that $e$ is of connected type.
Then
we take $1 \leq i \leq \beta$ with $(e,f) = (e_i , f_i)$. 
By Lemma~\ref{lemma:ES:separate:edges},
$- \log |s_{N_2+2i} / s_{N_2+2i-1} (x)|
\neq - \log |s_{N_2+2i} / s_{N_2+2i-1} (y)|$,
and hence
$\varphi_3 (x) \neq \varphi_3 (y)$.
Thus $\varphi_3 (x) \neq \varphi_3 (y)$
if $x \neq y$ and not both of $x$ and $y$ are vertices.

\medskip
{\bf Step 4.}
We take care of separation of two distinct vertices.
Set 
\[
W := \{ (v,w) \in V ( \mathscr{X}^{\st})
\times V ( \mathscr{X}^{\st})
\mid
\text{$v \neq w$, and $v$ and $w$ belong to the same island of $\Gamma_{\min}$}
\}.
\]
For any $(v,w) \in W$,
Proposition~\ref{prop:existence:sep:sections}
gives a model $( \mathscr{X}^{(v,w)} , \mathscr{L}^{(v,w)})$
such that $\mathscr{L}^{(v,w)}$ has 
a $v$-separating section and a $(v,w)$-separating section
or has a 
$w$-separating section and a $(w,v)$-separating section.
We take an (indexed) subset $\{(v_i, w_i)\}_{i=1}^{\gamma}$
such that 
\begin{itemize}
\item
the model $\mathscr{L}^{(v_i,w_i)}$
has a $v_i$-base section 
$\widetilde{s}^{(v_i,w_i)}_0$
and a $(v_i, w_i)$-separating section
$\widetilde{s}^{(v_i , w_i)}_{1}$;
\item
for any $(v , w) \in W$, there exists a unique $i = 1 ,\ldots , \gamma$ 
such that
$\{ v ,w \} = \{ v_i, w_i \}$.
\end{itemize}
Set $s_{N_3+2i-1} := \rest{\widetilde{s}^{(v_i,w_i)}_{0}}{X}$
and 
 $s_{N_3+2i} := \rest{\widetilde{s}^{(v_i , w_i)}_{1}}{X}$
for $i = 1 , \ldots , \gamma$,
and set $N := N_3 + 2 \gamma$.
We define
$\varphi : X^{\an} \to \TT\RR^{N}$
by $\varphi
:= ( - \log |s_0| : \cdots : - \log |s_{N_3}| : 
- \log |s_{N_3+1}| : \cdots : - \log |s_{N}|)$.
Since $\varphi_3$ is unimodular, so is $\varphi$
(cf. Remark~\ref{remark:unimodular}).

Take any distinct $v, w \in V (\mathscr{X}^{\st})$.
If $v$ and $w$ do not belong to 
some common island,
then by condition (iii) in Theorem~\ref{theorem:UT:canonical} 
we have $\varphi_1 (v) \neq \varphi_1 (w)$,
and hence $\varphi (v) \neq \varphi (w)$.
Suppose that $v$ and $w$ belong to the same island.
Then there exists $1 \leq i \leq \gamma$ with $(v,w) = (v_i , w_i)$ 
or $(w,v) = (v_i , w_i)$.
By Lemma~\ref{lemma:vertices-separating}, 
it follows that
$- \log |s_{N_3+2i} / s_{N_3 +2i -1} (v)|
\neq - \log |s_{N_3+2i} / s_{N_3 +2i -1} (w)|$,
and hence
$\varphi (v) \neq \varphi (w)$.
This completes the proof of the theorem.
\QED

\setcounter{equation}{0}
\section{Faithful tropicalization of minimal skeleta in low genera}
\label{section:FT:lowgenus}
Let $X$ be a connected smooth projective curve of genus $g \geq 0$.
In the previous section, we have shown the existence of a faithful tropicalization
of the minimal skeleton of $X^{\an}$ when $g \geq 2$.
In this section, we prove the existence when $g = 0,1$. 

\begin{Theorem}
\label{thm:FT:low:genus}
Let $X$ be a connected smooth projective curve of genus $g = 0,1$,
and let $\Gamma_{\min}$ be a minimal skeleton of $X^{\an}$.
Let $L$ be a line bundle over $X$.
Suppose that
\[
\deg (L)
\geq
\begin{cases}
1
&
\text{if $g=0$,} 
\\
3
&
\text{if $g=1$.}
\end{cases}
\]
Then there exist $s_0 , \ldots , s_N \in H^{0} (X,L)$
such that the tropicalization map
$\varphi : X^{\an} \to \TT\PP^N$
defined by 
$\varphi =
(- \log |s_0| : \cdots : - \log |s_N|)$
gives a faithful tropicalization of $\Gamma$.
\end{Theorem}

\subsection{Genus $0$ case}
\label{subsection:case:singleton:genus:0}
Assume that $g = 0$. 
Then the existence of a faithful tropicalization 
of a minimal skeleton is trivial.
Indeed, a minimal skeleton is a point. 
Thus if $\deg (L) \geq 1$,
then there exists two non-zero sections $s_0 , s_1 \in H^0 (L)$,
and the tropicalization with respect to those sections gives
a faithful tropicalization of a minimal skeleton.

\subsection{Genus $1$ case}
\label{subsection:case:genus:1}
Assume that $g = 1$. 
Then there exists a unique minimal skeleton of $X^{\an}$, which we denote by 
$\Gamma_{\min}$ as before (cf. \S\ref{subsec:properties:skeleta}).  
We have two possibilities of $\Gamma_{\min}$: a singleton or a circle.
If $\Gamma_{\min}$ is a singleton, then 
the existence of faithful tropicalization is again trivial.
Indeed, if $\deg (L) \geq 3$, then there exist two
nontrivial
sections $s_0 , s_1 \in H^0 (L)$,
and the tropicalization with respect to those sections gives
a faithful tropicalization of $\Gamma_{\min}$.

In the remaining of this subsection, we assume that 
$\Gamma_{\min}$ is a circle. 
We fix a point $v_0 \in \Gamma_{\min, \Lambda}$, and we let 
$e$ denote the loop edge of $\Gamma_{\min}$ such that $\mathrm{relin}
(e) = \Gamma_{\min} \setminus \{ v_0 \}$.
We are going to show the existence of a faithful tropicalization 
of $\Gamma_{\min}$ with 
a similar strategy to that we used for curves of genus at least $2$. 

\begin{Definition}[good model for $g = 1$]
\label{def:goodmodel:g=1}
Let $\mathscr{X}$ be a model of $X$.
We say that $\mathscr{X}$ is a \emph{good model}  
if it satisfies the following conditions. 
\begin{enumerate}
\item[(i)]
The model $\Xscr$ is Deligne--Mumford strictly semistable. 
\item[(ii)]
There exists a $C_0 \in \Irr (\mathscr{X}_s )$
such that $[C_0] = v_0$
and such that $E := \mathscr{X}_s - C_0$
is a $(-2)$-chain with symmetric multiplicities.
\end{enumerate}
\end{Definition}

\begin{Lemma}
\label{lemma:pregoodmodel:g=1}
Let $X$ be a connected smooth projective curve of genus $1$, 
and let $L$ be a line bundle over $X$ with $\deg (L) \geq 3$. 
Then for any $x,y \in \Gamma_{\min, \Lambda}$,
there exists a model $( \mathscr{X} , \mathscr{L})$ of $(X,L)$
such that $\mathscr{X}$ is a good model in Definition~\textup{\ref{def:goodmodel:g=1}} 
and such that
$D_{\mathscr{L}} -  [x] - [y]$ is effective.
\end{Lemma}

\Proof
Recall from \eqref{eqn:specialization:map} that $\tau_*: \Div(X) \to \Div_{\Lambda} ( \Gamma_{\min} )$ denotes 
the specialization map. Take a divisor $\tilde{D}$ on $X$ with $L \cong \OO_X ( \tilde{D} )$. 
Then $\tau_{\ast} ( \tilde{D} ) \in \Div_{\Lambda} ( \Gamma_{\min} )$
and the linear system
$\left| \tau_{\ast} ( \tilde{D} ) \right|$ does not depend
on the choice of $\tilde{D}$ (cf. the proof of Proposition~\ref{prop:pregoodmodel:pre}). 
Since $\deg ( \tau_{\ast} ( \tilde{D} ) -  [x] - [y] ) \geq 1$,
$\left| \tau_{\ast} ( \tilde{D} ) \right| \neq \emptyset$
by Proposition~\ref{prop:riemann:inequality}.
Thus there exists an effective
$E \in \Div_{\Lambda} (\Gamma_{\min})$
such that $E + [x] + [ y] \sim \tau_{\ast} ( \tilde{D} )$.
By the same argument as in the proof of Proposition~\ref{prop:pregoodmodel:pre},
we construct a model
$( \mathscr{X} , \mathscr{L} )$
of $(X,L)$
such that
$\mathscr{X}$ is good in the above sense
and such that 
$D_{\mathscr{L}} = E +  [x] + [y]$.
\QED

Let $(\mathscr{X} , \mathscr{L} )$
be a model of $(X, L)$ such that $\Xscr$ is a good model in Definition~\ref{def:goodmodel:g=1}.
Let $\widetilde{s}$ be a nonzero global section of $\mathscr{L}$.
We call $\widetilde{s}$ a \emph{base section} if $\widetilde{s} (p) \neq 0$ for 
any $p \in \Sing ( \mathscr{X}_s )$.
Further,
let $C_0$ and $E$ be as in Definition~\ref{def:goodmodel:g=1}(ii).
Since $E$ is a maximal $(-2)$-chain with
symmetric multiplicities and $\Delta_E = e$,
we define the stepwise vertical divisor $\mathscr{V}_e$ by the same way as in
\S\ref{subsection:stepwise}.
We call $\widetilde{s}$ a \emph{unimodularity section} if 
$\zero (\widetilde{s}) - \mathscr{V}_e$
is effective on $\Xscr$ and is trivial on some open neighborhood of $\Sing ( \mathscr{X}_s )$.

\begin{Proposition}
\label{prop:model:US:g=1}
Let $X$ and $ L$ be as in Lemma~\textup{\ref{lemma:pregoodmodel:g=1}}. 
Then
there exists a model $(\mathscr{X} , \mathscr{L} )$
such that
$\mathscr{X}$ is a good model in Definition~\ref{def:goodmodel:g=1} and 
such that
$\mathscr{L}$ has a base section and a unimodularity section.
\end{Proposition}

\Proof
Applying Lemma~\ref{lemma:pregoodmodel:g=1} to $x = y = v_0$,
we obtain a model $(\mathscr{X} , \mathscr{L})$ of $(X,L)$
such that $\mathscr{X}$ is a good model and 
$D_{\mathscr{L}} - 2 [v_0]$ is effective.
By the same argument using Lemma~\ref{lemma:freeatnode:graphversion:emptyset}
as in the proof of
Lemma~\ref{lemma:forbasesectionUT:dc}(1),
we see that $\mathscr{L}$ is free at any $p \in \Sing ( \mathscr{X}_s )$,
and then
by the same argument as the proof  of Proposition~\ref{prop:existence:sections:stronger}, we obtain a base section.

We set $\mathscr{L}_1 := \mathscr{L} ( - \mathscr{V}_e )$. 
By the same argument as in the proof of
Lemma~\ref{lemma:forbasesectionUT:dc}(3),
we see that $\mathscr{L}_1$ is free at any $p \in \Sing ( \mathscr{X}_s )$.
By the same argument as the proof  of Proposition~\ref{prop:existence:sections:stronger}, we obtain a unimodularity section of $\mathscr{L}$.
\QED

Let $\mathscr{X}$ be a good model in Definition~\ref{def:goodmodel:g=1}.
Since the maximal $(-2)$-chain $E$ 
has symmetric multiplicity,
we can consider a separating divisor $\mathscr{A}$ for $e$
as in Definition~\ref{def:assymetric:divisor} also in this setting.
Let $\mathscr{L}$ be a model of $L$.
We call a nonzero global section $s$ of $\mathscr{L}$
an \emph{separating section}
if there exists a separating divisor $\mathscr{A}$
(for $e$)
such that
$\zero (s ) - \mathscr{A}$
is trivial on some open neighborhood of 
$\Sing ( \mathscr{X}_s ) \setminus \Supp ( \mathscr{A})$ 
(Here we remark that $\Sing ( \mathscr{X}_s )_{\subset e} = \Sing ( \mathscr{X}_s )$.)  

\begin{Lemma}
Let $X$ and $ L$ be as in Lemma~\textup{\ref{lemma:pregoodmodel:g=1}}. 
Then there exists a model
$(\mathscr{X} , \mathscr{L} )$ of $(X,L)$
such that $\mathscr{X}$ is a good model in Definition~\textup{\ref{def:goodmodel:g=1}} and
such that
$\mathscr{L}$ has a base section and a
separating section.
\end{Lemma}

\Proof
Let $v_1 \in \Gamma_{\min}$ be the point antipodal to $v_0$. In other words, 
let $v_1$ be the middle point of $e$. 
Applying Lemma~\ref{lemma:pregoodmodel:g=1} to $x = v_0$
and $y = v_1$,
we obtain a model $(\mathscr{X} , \mathscr{L})$ of $(X,L)$
such that $\mathscr{X}$ is a good model and 
$D_{\mathscr{L}} -  [v_0] - [v_1]$ is effective.
By the same argument as in Step~1 of the proof
of Proposition~\ref{prop:model:AS},
one sees that $\mathscr{L}$ has a base section.

We set $\mathscr{L}_1 :=
\mathscr{L} ( - \mathscr{A})$.
By the same argument as in Step~2 of 
the proof of Proposition~\ref{prop:model:AS},
one sees that
$\mathscr{L}_1$ is free at any $q \in 
\Sing ( \mathscr{X}_s ) \setminus \Supp ( \mathscr{A})$,
and hence there exists a section
$\widetilde{s}_1 \in H^0 ( \mathscr{L}_1 )$ such that
$\widetilde{s}_1 (q) \neq 0$ for any $q \in 
\Sing ( \mathscr{X}_s ) \setminus \Supp ( \mathscr{A})$.
The same argument as in the last paragraph of the proof of 
Proposition~\ref{prop:model:AS} gives 
the existence of a separating section of $\mathscr{L}$.
\QED

Since we have a model that has a base section and a unimodularity section,
the same argument as the proof of Theorem~\ref{theorem:UT:canonical}
gives us a unimodular tropicalization of $\Gamma_{\min}$.
Furthermore, since
we have a model that has a base section and a separating section,
the same argument as the proof of Theorem~\ref{theorem:FT:canonical}
gives us a faithful tropicalization of $\Gamma_{\min}$.

Together with \S\ref{subsection:case:singleton:genus:0}, 
we have shown Theorem~\ref{thm:FT:low:genus}.

\setcounter{equation}{0}
\section{Faithful tropicalization of an arbitrary skeleton}
\label{sec:ft:general:model}
In this section, we prove the main Theorem~\ref{thm:main:unimodular:faithful}.

\subsection*{Notation and terminology of \S\ref{sec:ft:general:model}} 
Throughout this section, we use the following notation and terminology. 
Let $X$ be a connected smooth projective over $K$ of genus $g \geq 0$,
and let $L$ be a line bundle over $X$. For any compact skeleton $\Gamma$ of $X^{\an}$, 
let $\tau_{\Gamma}: X^{\an} \to \Gamma$ denote 
the retraction map with respect to $\Gamma$.

\subsection{Geodesic paths}
\label{subsection:geodesic:line}
Recall that $X^{\an}\setminus X(K)$ has a canonical metric structure (cf. \cite[Corollary~5.7]{BPR2}). 
In this subsection, we briefly review basic properties of geodesic paths in $X^{\an} \setminus X(K)$.
We refer to \cite[\S5]{BPR2} for details.

Let $x , y \in X^{\an} \setminus X (K)$,
and let $\gamma$ be a path 
in $X^{\an} \setminus X (K)$
from $x$ to $y$,
which we mean a continuous map from a
closed interval $[a , b] \subset \RR$ to $X^{\an} \setminus X (K)$
such that $\gamma (a) = x$ and $\gamma (b) = y$.
We call $\gamma$
a \emph{geodesic path connecting $x$ and $y$}
if $\gamma$ is an isometry from $[a,b]$ to its image 
(cf. \cite[\S5]{BPR2}).
In what follows,
we identify a geodesic path with its image in $X^{\an} \setminus X (K)$.
Remark that two geodesic paths with the same image differ only by
parameterizations:
if $\gamma_1 : [a_1 , b_1] \to X^{\an} \setminus X (K)$ 
and $\gamma_2 : [a_2 , b_2] \to X^{\an} \setminus X (K)$ are 
geodesic  paths connecting $x$ and $y$
such that $\gamma_1 ([a_1 , b_1]) = \gamma_2 ([a_2 , b_2])$,
then there exists a
unique isometry $\phi : [a_1 , b_1] \to 
[a_2 , b_2]$ such that $\gamma_1 = \gamma_2 \circ \phi$.

Before discussing geodesic paths,
we show the following
lemma on the retraction map, which will be frequently used.

\begin{Lemma}
\label{lemma:connectedcomponents}
Let $\Gamma$ be a compact skeleton of $X^{\an}$, and let $x$ be a point in $X^{\an} \setminus \Gamma$. 
Let $A$ be the connected component of $X^{\an} \setminus \Gamma$
with $x \in A$. Then $A$ is the connected component of $X^{\an}
\setminus \{ \tau_{\Gamma} (x) \}$ with $x \in A$.
\end{Lemma}

\Proof
By the definition of $\tau_{\Gamma}$,
$\overline{A} = A \cup \{  \tau_{\Gamma} (x) \}$
and $\tau_{\Gamma} (A) = \{  \tau_{\Gamma} (x) \}$. 
Let $B$ be the connected component of $X^{\an} \setminus  \{  \tau_{\Gamma} (x) \}$
with $B \cap A \neq \emptyset$.
Since $A \subset X^{\an} \setminus  \{  \tau_{\Gamma} (x) \}$
and $A$ is connected, we have $A \subset B$.
Set $C := B \setminus A$.
Then we have $B = A \cup C$ and $A \cap C = \emptyset$.
Since $ \tau_{\Gamma} (x) \notin B$, 
we have $C = B \setminus A = B \setminus (A \cup \{  \tau_{\Gamma} (x) \})
= B \setminus \overline{A}$, which is an open subset of $B$.
On the other hand, since $X^{\an}$ is locally connected, 
$A$ is an open subset of $X^{\an}$ and hence of $B$. 
Since $B$ is connected and $A \neq \emptyset$,
we have $C = \emptyset$, and thus $A = B$. 
This shows the lemma. 
\QED

The following lemma shows basic properties of
geodesic paths,
which are essentially treated 
in \cite{BPR2}.

\begin{Lemma}
\label{lemma:geodesic:path}
Let $x , y \in X^{\an}$.
Suppose that there exists a compact skeleton $\Gamma$ 
of $X^{\an}$ such that $x,y \in \Gamma$. 
We fix a minimal skeleton $\Gamma_{\min}$. 
\begin{enumerate}
\item
Then there exists a geodesic path connecting $x$ and $y$.
\item
Furthermore, suppose that 
$\tau_{\Gamma_{\min}} (x) = \tau_{\Gamma_{\min}} (y)$. 
Then a geodesic path connecting $x$ and $y$
\textup{(}as the image in $X^{\an}$\textup{)}
is unique. 
In particular, if there exists a compact skeleton $\Gamma_1$ of $X^{\an}$ such that
$\tau_{\Gamma_1} (x) = \tau_{\Gamma_1} (y)$, then we have the uniqueness 
of the geodesic path connecting $x$ and $y$. 
\end{enumerate}
\end{Lemma}

\Proof
Assertion (1) is a direct consequence of 
the first paragraph of \cite[\S5.3]{BPR2}.

For assertion (2), put $z := \tau_{\Gamma_{\min}} (x) = \tau_{\Gamma_{\min}} (y)$. 
First, we prove that the uniqueness of a geodesic path
connecting $x$ and $z$.
If $x = z$, then it is trivial, so that we assume that
$x \neq z$.
Then
$x \notin \Gamma_{\min}$.
Let $\Gamma_x^\circ$ 
be the connected component of 
$\Gamma \setminus \Gamma_{\min}$ 
containing $x$,
and let $\Gamma_x$ be the closure of $\Gamma_x^{\circ}$ in $X^{\an}$.
Then we claim that $\Gamma_x$ is tree.
Indeed,
let $A_x$ be the connected component of $X^{\an} \setminus \Gamma_{\min}$ that 
contains $x$. Since $\Gamma_x^\circ \subset X^{\an} \setminus \Gamma_{\min}$ and 
$\Gamma_x^\circ$ is connected, we then have $\Gamma_x^\circ \subset A_x$. 
By the definition of $\tau_{\Gamma_{\min}}$,
we have $\overline{A_x} = A_x \cup \{  z \}$.
Note that,
since $\Gamma$ is connected and $\Gamma_{\min} \subset \Gamma$
is closed, $\Gamma_x^{\circ}$ is not closed.
This means that
 the closure $\Gamma_x$ of $\Gamma_x^\circ$
in $\Gamma$ is equal to $\Gamma_x = \Gamma_x^\circ\cup \{z\}$. 
On the other hand, by \cite[Lemma~3.4(3)]{BPR2}, $A_x$ 
is an open ball.
Then the argument
in \cite[\S5.8]{BPR2}
shows that
$\Gamma_x$ is a tree. 

Let $\gamma$ be a geodesic path connecting $x$ and $z$.
By Lemma~\ref{lemma:connectedcomponents}, $A_x$ is the 
connected component of $X^{\an} \setminus \{ z \}$ with $x \in A_x$.
Since $\gamma \setminus \{ z \}$ is a connected subspace
in $X^{\an} \setminus \{ z \}$ and
$x \in \gamma \setminus \{ z \}$,
it follows that $\gamma \setminus \{ z \} \subset A_x$.
Further,
by \cite[Corollary~5.10]{BPR2}, $\gamma \subset \Gamma$.
It follows that
$\gamma \setminus \{ z \} \subset \Gamma \cap A_x = \Gamma_x^{\circ}$,
and hence $\gamma \subset \Gamma_x$.
Thus $\gamma$ is a geodesic path in $\Gamma_x$.
Since $\Gamma_x$ is a tree,
this proves that $\gamma$ is unique.
This shows that the geodesic path connecting $x$ and $z$ is unique.

If $y \neq z$, then we define $A_y$, $\Gamma_y^{\circ}$,
and $\Gamma_y$ with $x$ replaced with $y$.
Further,
the same argument as above also shows that a geodesic path connecting
$y$ and $z$ is unique.

Next, let us show that the uniqueness of a geodesic path connecting
$x$ and $y$.
We may and do assume that $x \neq z$ and $y \neq z$,
since otherwise we have already shown the assertion.
We remark that in this situation, $A_x$, $\Gamma_x^{\circ}$, $\Gamma_x$,
$A_y$, $\Gamma_y^{\circ}$, and $\Gamma_y$ are defined above.
Let $\gamma$ be a geodesic path connecting $x$ and $y$.
By \cite[Corollary~5.10]{BPR2}, $\gamma \subset \Gamma$.

Suppose that $A_x = A_y$. Then $\Gamma_x = \Gamma_y$. 
We remark that $\Gamma_x \cap \Gamma_{\min} = \{ z \}$.
Then, since $\gamma$ is geodesic in $\Gamma$, 
one sees that $\gamma \subset \Gamma_x$.
Since $\Gamma_x$ is a tree, $\gamma$ is unique.

Suppose that $A_x \neq A_y$. Then $A_x \cap A_y = \emptyset$.
It follows that  
$\Gamma_x^\circ \cap  \Gamma_y^\circ = \emptyset$, and thus 
$\Gamma_x \cap \Gamma_y = \{ z \}$.
This means that $\Gamma_x \cup \Gamma_y$ is a one point sum at $z$ of 
two trees, so that it is also a tree.
Further, $(\Gamma_x \cap \Gamma_y) \cap \Gamma_{\min} = \{ z \}$.
Since $\gamma$ is geodesic in $\Gamma$,
it follows that
$\gamma \subset \Gamma_x \cup \Gamma_y$.
Since $\Gamma_x \cup \Gamma_y$ is a tree, this shows that $\gamma$
is unique.
This completes the proof of the
first assertion of (2). 

\if
Since $\Gamma_x$ and $\Gamma_y$ are trees, 
we see that $\gamma$ is a unique geodesic path connecting $x$ and $y$ lying in $\tilde{\Gamma}$; in fact, 
$\gamma$ is the union of a unique geodesic path connecting $x$ and $z$ in $\Gamma_x$ and 
a unique geodesic path connecting $z$ and $y$ in $\Gamma_y$. 
We have shown that a geodesic path connecting $x$ and $y$ that lie in $\tilde{\Gamma}$ is unique. 

Suppose that $\gamma_1$ and $\gamma_2$ 
are geodesic paths connecting $x$ and $y$. Then we take 
a compact skeleton $\gamma_1 , \gamma_2\Gamma$ 
by \cite[Corollary~5.10]{BPR2}.
Then the above argument gives that $\gamma_1 = \gamma_2$. 
This completes the proof of the former assertion of (2).  
\fi

For the second assertion of (2), if $g \geq 1$, then $\Gamma_{\min}$ is unique and 
is contained in $\Gamma_1$. It follows from 
$\tau_{\Gamma_1} (x) = \tau_{\Gamma_1} (y)$ that 
$\tau_{\Gamma_{\min}} (x) = \tau_{\Gamma_{\min}} (y)$. 
If $g = 0$, then $\Gamma_{\min}$  is a singleton, and we always 
have $\tau_{\Gamma_{\min}} (x) = \tau_{\Gamma_{\min}} (y)$. 
Thus the second assertion follows from the first assertion. 
\QED

Let $x $ and $ y$ be points in $X^{\an}$ as in Lemma~\ref{lemma:geodesic:path}. 
We assume that $\tau_{\Gamma_{\min}}(x) = \tau_{\Gamma_{\min}}(y)$, 
so that there exists a unique geodesic path $\gamma$ connecting $x$ and $y$. 
We call 
the image of $\gamma$ in $X^{\an}$ 
the \emph{geodesic segment connecting $x$ and $y$}. We use the notation 
$[x,y]$ for the geodesic segment connecting $x$ and $y$. 
Further,
if $x \neq y$, then we set
\[
[x , y) := [x,y] \setminus \{ y \}, \quad 
(x , y] := [x,y] \setminus \{ x \}, \quad\text{and}\quad
(x , y) := [x,y] \setminus \{ x, y \}.
\]

We give some properties of geodesic segments. 
The following two lemmas are also essentially treated in \cite{BPR2}. 

\begin{Lemma}
\label{lemma:to:be:used}
Let $\Gamma$ and $\Gamma^\prime$ be compact skeleta with 
$\Gamma \subset \Gamma^\prime$ 
and let $x \in \Gamma^\prime_\Lambda \setminus \Gamma_\Lambda$. 
Then we have the following. 
\begin{enumerate}
\item
\label{lemma:to:be:used:addingsegments}
The union $\Gamma \cup [ \tau_{\Gamma} (x) , x ]$ is a compact skeleton
of $X^{\an}$.
\item
\label{lemma:to:be:used:3}
We have $[\tau_{\Gamma}(x), x] \subset \Gamma^\prime$. Further, we have 
$\Gamma \cap [\tau_{\Gamma}(x), x] = \{\tau_{\Gamma}(x)\}$. 
\item
\label{lemma:to:be:used:4}
For any compact skeleton $\Gamma_1$ with $\Gamma \subset \Gamma_1 \subset \Gamma^\prime$, 
we have $\tau_{\Gamma_1}(x) \in [\tau_{\Gamma}(x), x]$. 
\end{enumerate}
\end{Lemma}

\Proof
We first note that, 
since $x$ and $\tau_{\Gamma} (x)$ is contained in 
a compact skeleton and
$\tau_{\Gamma}(x) = \tau_\Gamma(\tau_{\Gamma}(x))$, 
$[x, \tau_\Gamma(x)]$ makes sense by Lemma~\ref{lemma:geodesic:path}. 

\smallskip
(\ref{lemma:to:be:used:addingsegments})
We write $\Gamma^\prime = S( \mathscr{X}^\prime)$ for some 
strictly semistable model $\Xscr^\prime$ of $X$. 
By \cite[Corollary~5.10]{BPR2},
we have 
$[\tau_{\Gamma} (x) , x ] \subset S( \mathscr{X}^\prime)$.
By successively contracting the irreducible components
$E$ such that $[E] \notin \Gamma \cup [\tau_{\Gamma} (x) , x ]$,
we obtain a strictly semistable model $\mathscr{X}^{\prime\prime}$
such that $\Gamma \cup [\tau_{\Gamma} (x) , x ] = S( \mathscr{X}^{\prime\prime})$.
This shows that $\Gamma \cup [\tau_{\Gamma} (x) , x ]$
is a compact skeleton.

\smallskip
(\ref{lemma:to:be:used:3})
Since $x , \tau_{\Gamma} (x) \in \Gamma'$,
$[\tau_{\Gamma} (x) , x] \subset \Gamma'$ by
\cite[Corollary~5.10]{BPR2}. 
For the second assertion, 
let $A$ be the connected component of $X^{\an} \setminus \Gamma$
with $x \in A$.
By Lemma~\ref{lemma:connectedcomponents}, 
$A$ is the connected component of $X^{\an} \setminus \{ 
\tau_{\Gamma} (x) \}$ with $x \in A$.
Since $(\tau_{\Gamma} (x) , x]$ is a connected subspace
of $X^{\an} \setminus \{ \tau_{\Gamma} (x) \}$ containing $x$, 
we get $(\tau_{\Gamma} (x) , x]
\subset A$.
This proves that
$\Gamma \cap [\tau_{\Gamma}(x), x] = \{\tau_{\Gamma}(x)\}$. 

\smallskip
(\ref{lemma:to:be:used:4})
By the first assertion in (\ref{lemma:to:be:used:3}) applied to $\Gamma \subset \Gamma_1$ and $\tau_{\Gamma_1}(x)$, we have $[\tau_{\Gamma} (x) , \tau_{\Gamma_1} (x)] \subset \Gamma_1$. 
By the second assertion in (\ref{lemma:to:be:used:3}), we then have 
$[\tau_{\Gamma} (x) , \tau_{\Gamma_1} (x)] 
\cap [\tau_{\Gamma_1} (x) , x] = 
\{ \tau_{\Gamma_1} (x) \}
$, which means that
$[\tau_{\Gamma} (x) , \tau_{\Gamma_1} (x)] 
\cup [\tau_{\Gamma_1} (x) , x]$ is a geodesic segment
connecting $\tau_{\Gamma} (x)$ and $x$.
By the uniqueness of a geodesic segment of Lemma~\ref{lemma:geodesic:path}(2),
$[\tau_{\Gamma} (x) , \tau_{\Gamma_1} (x)] 
\cup [\tau_{\Gamma_1} (x) , x]
= [\tau_{\Gamma} (x) , x]$.
Thus $\tau_{\Gamma_1} (x) \in [\tau_{\Gamma} (x) , x]$.
\QED

Lemma~\ref{lemma:to:be:used}(\ref{lemma:to:be:used:addingsegments})
can be also stated as follows. 

\begin{Lemma}
\label{lemma:addingsegments}
Let $\Gamma$ be a compact skeleton of $X^{\an}$.
Let $x$ be a point of $X^{\an} \setminus \Gamma$.
Assume that there exists a compact skeleton $\Gamma'$ 
such that $x \in \Gamma'_{\Lambda}$.
Then $\Gamma \cup [ \tau_{\Gamma} (x) , x ]$ is a compact skeleton
of $X^{\an}$.
\end{Lemma}

\Proof
By \cite[Proposition~3.13(3)]{BPR2},
there exists a compact skeleton $\Gamma^{\prime\prime}$
such that $\Gamma \cup \{ x \} \subset \Gamma^{\prime\prime}$. 
Then the assertion is reduced to Lemma~\ref{lemma:to:be:used}(\ref{lemma:to:be:used:addingsegments}) 
applied to $\Gamma^{\prime\prime}$. 
\QED

Next, we would like to discuss a path one of whose endpoint is a point in $X(K)$. 
Here, we write $[0 , + \infty] = \TT$ by notation.

\begin{Lemma}
\label{lemma:end=canonialend}
Let $\Gamma$ be a compact skeleton,
and
let $P$ be any point in $X(K)$. 
Then the following hold.
\begin{enumerate}
\item
There exists a unique continuous map $\gamma: [0, + \infty] \to X^{\an}$ with 
$\gamma(0) = \tau_{\Gamma} (P)$ and $\gamma(1) = P$ 
such that $\gamma$ 
restricts to an isometry from
$[0 , + \infty )$ to its image.
\item 
There exists a strictly semistable pair
$( \mathscr{X} ; \sigma )$ such that $S( \mathscr{X}) = \Gamma$ and $\sigma (K) = P$. 
\item
We take $\gamma$ as in \textup{(1)} and $( \mathscr{X} ; \sigma )$ as in \textup{(3)}. 
We set 
\[
[\tau_{\Gamma} (P) , P] := \gamma\left([0, + \infty ]\right)
\quad\text{and}\quad  
[\tau_{\Gamma} (P) , P) := [\tau_{\Gamma} (P) , P] \setminus \{ P \} = \gamma\left([0, + \infty)\right).
\] 
Then $\Delta (\sigma) = [\tau_{\Gamma} (P) , P)$ and 
$\overline{\Delta (\sigma)} = [\tau_{\Gamma} (P) , P]$. 
\end{enumerate}
\end{Lemma}

\Proof
First, we prove the uniqueness property in (1). 
Let $\gamma^\prime:  [0, + \infty] \to X^{\an}$ be a continuous map that 
satisfies the same property as $\gamma$. We are going to show that 
$\gamma = \gamma^\prime$.  

We set $I := \gamma([0, + \infty])$ and $J:= \gamma^\prime([0, + \infty])$. 
It suffices to show that $I = J$. To argue by contradiction, assume that 
$I \neq J$.
Then
one easily sees that $I \nsubseteq J$ and $J \nsubseteq I$,
and
there exist $x \in I \setminus J$ and $y \in J \setminus I$
such that the distance $d(\tau_{\Gamma} (P) , x)$
from $\tau_{\Gamma} (P)$ to $x$ and
the distance
$d(\tau_{\Gamma} (P) , y)$ from $\tau_{\Gamma} (P)$ to $y$ are in $\Lambda$.
By the uniqueness of geodesics (cf. Lemma~\ref{lemma:geodesic:path}(2)),
we have
$[\tau_{\Gamma} (P) , x] \subset I$,
$[\tau_{\Gamma} (P) , y] \subset J$.
Further, 
we see that
there exists $z \in [\tau_{\Gamma} (P) , x] \cap 
[\tau_{\Gamma} (P) , y]$ such that 
$[\tau_{\Gamma} (P) , x] \cap 
[\tau_{\Gamma} (P) , y] = [ \tau_{\Gamma} (P) , z]$.
By Lemma~\ref{lemma:addingsegments},
$\Gamma_0 := \Gamma \cup [\tau_{\Gamma} (P) , x]$ is a compact skeleton.
By the choice of $z$, one sees that $z = \tau_{\Gamma_0} (y)$,
which belongs to $\Gamma_{0,\Lambda}$ by Lemma~\ref{lemma:retraction:rational}(1).
Using Lemma~\ref{lemma:addingsegments} again,
we see that $\Gamma_1 := \Gamma_0 \cup [z , y] =
\Gamma \cup 
[\tau_{\Gamma} (P) , x] \cup [\tau_{\Gamma} (P) , y]$
is a compact skeleton.

Let $A$ be the connected component of $X^{\an} \setminus \Gamma_1$ with $P \in A$. 
By Lemma~\ref{lemma:connectedcomponents}, $A$ is the connected component of 
$X^{\an} \setminus \{\tau_{\Gamma_1}(P)\}$, 
and we have $\tau_{\Gamma}(\overline{A}) = \{\tau_{\Gamma_1}(P)\}$ (cf. Lemma~\ref{lemma:retraction:rational}). 
Since $I \setminus [\tau_{\Gamma} (P) , x]$ is a connected subset 
with $I \setminus [\tau_{\Gamma} (P) , x]\subset X^{\an} \setminus \Gamma_1$
and $P \in I \setminus [\tau_{\Gamma} (P) , x] $,
we have $I \setminus [\tau_{\Gamma} (P) , x] \subset A$.
Since $x \in \overline{I \setminus [\tau_{\Gamma} (P) , x]}$,
we have $x \in \overline{A}$. 
It follows that $\tau_{\Gamma_1} (P) = \tau_{\Gamma_1} (x) = x$. 
On the other hand, the same argument shows that $\tau_{\Gamma_1} (P) = y$.
This contradicts with $x \neq y$. 
Thus we obtain the uniqueness property in (1) .

Next we prove (2). 
By Lemma~\ref{lemma:retraction:rational}(2),
we have $\tau_{\Gamma} (P) \in \Gamma_{\Lambda}$.
It follows from
Proposition~\ref{prop:subdivision1:b}
that
there exists a strictly semistable model $\mathscr{X}$ of $X$
such that $S ( \mathscr{X} ) = \Gamma$
and $\tau_{\Gamma} (P) \in V ( \mathscr{X})$.
We take $C_{\tau_{\Gamma} (P)} \in \Irr ( \mathscr{X}_s)$
such that $\tau_{\Gamma} (P) = [ C_{\tau_{\Gamma} (P)} ]$.
By the valuative criterion of properness,
we take a section $\sigma$ of $\mathscr{X} \to \Spec (R)$
such that $\sigma (K) = P$.
Since $\tau_{\Gamma} (P) = [ C_{\tau_{\Gamma} (P)} ]$,
it follows from Lemma~\ref{lemma:retraction:rational}
that
$\red_{\mathscr{X}} (P) = \sigma (k) \in
C_{\tau_{\Gamma} (P)} (k) \setminus \Sing ( \mathscr{X}_s )$.
Thus
$(\mathscr{X} ; \sigma )$
is a strictly semistable pair
as desired.

Finally, we prove the existence of $\gamma$ in (1) and (3). 
Since we have a unique isometry $\Delta(\sigma) \cong \RR_{\geq 0}$ (cf. \S\ref{subsec:skeleta:ass:ssp}), 
we have an isometry $\gamma: [0, + \infty) \to \Delta(\sigma)$. 
To complete the proof of the lemma, we have only 
to show that $\overline{\Delta(\sigma)} = \Delta(\sigma) \cup \{P\}$ and 
that $\gamma$ extends to a homomorphism $\gamma: [0, + \infty] \to \overline{\Delta(\sigma)}$. 
However, this is a fact mentioned in
\cite[Remark~4.12]{GRW}
and in (the proof of)
\cite[Lemma~3.4(2)]{BPR2}. 
Thus there exists $\gamma$ with property described in (1), and 
we have the equalities described in (3). 
\QED

Let $\Gamma$ be a compact skeleton of $X^{\an}$, and let $P$ be a point in $X(K)$. 
By slight abuse of terminology, we call $[\tau_{\Gamma} (P) , P]$ 
a {\em path connecting $\tau_{\Gamma} (P)$ and $P$}. 
By Lemma~\ref{lemma:end=canonialend}, 
$[\tau_{\Gamma} (P) , P]$ exists and is determined only by $\Gamma$ and $P$. 
We also write
\[
(\tau_{\Gamma} (P) , P) := [\tau_{\Gamma} (P) , P]
\setminus \{ \tau_{\Gamma} (P) , P \}.
\]

\subsection{Stepwise vertical divisor associated to a point in $X(K)$}
\label{subsection:setting1}

Let $P$ be a point in $X(K)$, 
and let $\Gamma$ be a compact skeleton of $X^{\an}$.
By Lemma~\ref{lemma:end=canonialend},
there exists a strictly semistable pair 
$(\mathscr{X} ; \sigma)$
such that 
$S(\mathscr{X}) = \Gamma$ and $\sigma (K) = P$,
and we have
$[ \tau_{\Gamma} (P) , P ) = \Delta(\sigma)$
and
$S( \mathscr{X} ; \sigma ) = \Gamma
\cup [ \tau_{\Gamma} (P) , P )$.
Fix a minimal skeleton $\Gamma_{\min} \subset \Gamma$.
Then 
\begin{align*}
& [ \tau_{\Gamma_{\min}} (P) , P ) =  [\tau_{\Gamma_{\min}} (P) ,
\tau_{\Gamma} (P)] \cup
[ \tau_{\Gamma} (P) , P ),
 \\
& [\tau_{\Gamma_{\min}} (P) ,
\tau_{\Gamma} (P)] \cap
[ \tau_{\Gamma} (P) , P ) = \{ \tau_{\Gamma} (P) \}.
\end{align*}

For $P$, $\mathscr{X}$, and $\Gamma_{\min}$ above,
we define a Cartier divisor $\mathscr{W}_P$ on $\mathscr{X}$,
which is a sum of fundamental divisors, 
as follows.
Let $E$ be the chain of curves in $\mathscr{X}_s$
characterized by 
$ \Delta_E
= 
[\tau_{\Gamma_{\min}} (P) ,
\tau_{\Gamma} (P)]
$.
We write the irreducible components of $E$ for
$E_0 , \ldots , E_{r}$
in such a way that
 $[E_0] = \tau_{\min} (P)$ 
and $E_{\alpha} \cap E_{\alpha+1} \neq \emptyset$ for 
$\alpha = 0 , \ldots , r - 1$;
then $[E_r] = \tau_{\Gamma} (P)$.
Note that $\sigma (k) \in E_{r} \setminus \Sing ( \mathscr{X}_s )$.
Let $p_{\alpha}$ denote the node in $\mathscr{X}_s$
with $\{ p_{\alpha} \} = E_{\alpha} \cap E_{\alpha+1}$
for $\alpha = 0 , \ldots , r - 1$.

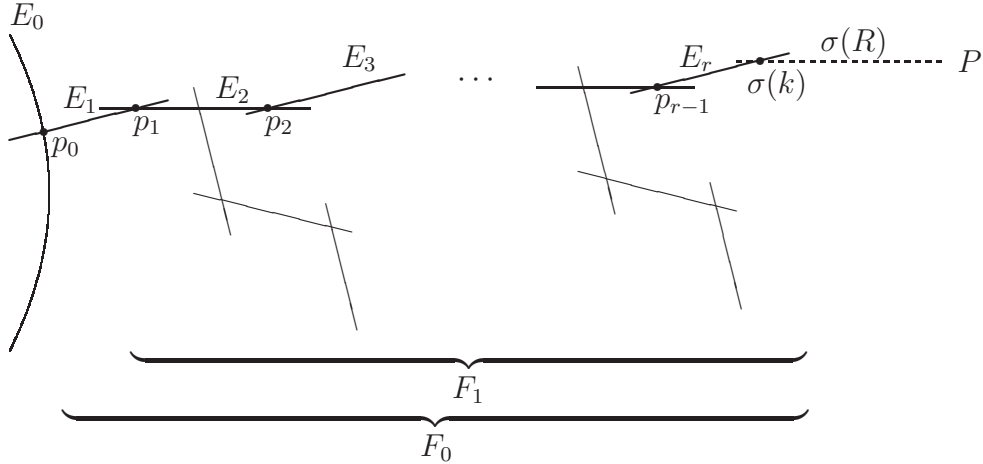
\begin{figure}[!h]
\[
\setlength\unitlength{0.07truecm}
\begin{picture}(180, 85)(0,0)
  \qbezier(0, 20)(15, 50)(0, 80)
   \thicklines
  \put(0,60){\line(4,1){30}}
  \put(17,66){\line(1,0){40}}
  \put(45,65){\line(4,1){30}}
  \put(100,70){\line(1,0){30}}
  \put(118,69){\line(4,1){30}}
  \multiput(138,75)(2, 0){20}{\line(1,0){1}}
  \thinlines
  \put(35,70){\line(1,-4){7}}
  \put(35,50){\line(4,-1){30}}
  \put(60,48){\line(1,-4){6}}
  \put(108,74){\line(1,-4){7}}
  \put(108,54){\line(4,-1){30}}
  \put(133,52){\line(1,-4){6}}
  \put(85, 70){$\cdots$}
  \put(0, 82){$E_0$}
  \put(39, 68){$E_2$}
  \put(10, 66){$E_1$}
  \put(63, 74){$E_{3}$}
  \put(127, 74){$E_r$}
  \put(154, 77){$\sigma(R)$}
  \put(6.5, 61.5){\circle*{1.5}}
  \put(24, 66){\circle*{1.5}}
  \put(49, 66){\circle*{1.5}}
  \put(123, 70){\circle*{1.5}}
  \put(142.5, 75){\circle*{1.5}}
  \put(8, 58){$p_0$}
  \put(23.5, 62){$p_1$}
  \put(48.5, 62){$p_{2}$}
  \put(122.5, 66){$p_{r-1}$}
  \put(140, 69){$\sigma(k)$}
  \put(180,73){$P$}
   \put(23, 20){$\underbrace{\phantom{AAAAAAAAAAAAAAAAAAAAAAAAAAAAA}}$} 
   \put(84, 11){$F_1$}
   \put(10, 9){$\underbrace{\phantom{AAAAAAAAAAAAAAAAAAAAAAAAAAAAAAAA}}$}
 
   \put(78, 0){$F_0$}
\end{picture}
\]
\caption{A part of the configuration of $(\Xscr, \sigma)$  with $\sigma(
K) = P$}
\label{figure:configuration:ends}
\end{figure}

If $r = 0$, then we set $\mathscr{W}_{P} = 0$. 
Assume that $r \geq 1$.
For each $0 \leq \alpha \leq r -1$,
the partial normalization at $p_{\alpha}$ has exactly two connected components.
Let $F_{\alpha}$ be the connected component that contains $E_{r}$.
We regard $F_{\alpha}$ as a curve in $\mathscr{X}_s$.
Then $F_{\alpha}$ satisfies condition (F) in \S\ref{subsec:fvd}.
Let $\mathscr{F}_{\alpha}$ be the fundamental vertical divisor with
support $F_{\alpha}$.
We define a vertical effective
divisor on $\mathscr{X}$ for 
$P$
to be
\begin{align}
\label{align:stepwise:for:ends}
\mathscr{W}_{P}
:=
\sum_{\alpha = 0}^{r - 1} \mathscr{F}_{\alpha}
.
\end{align}
We call this divisor the \emph{stepwise vertical divisor associated to $P$.}

The divisor $\mathscr{W}_P$ depends on the choice of $\Gamma_{\min}$.
However, if $g \geq  1$,
then $\Gamma_{\min}$ is unique 
and hence $\mathscr{W}_P$ is determined by $P$ and
$\mathscr{X}$.

\medskip
We compute the order
of $\mathscr{W}_P$ at each irreducible component of $\mathscr{X}_s$
when $r \geq 1$.
By the definition of $\mathscr{W}_P$,
we have
$\Supp ( \mathscr{W}_P) = F_0$,
which is a tree of 
irreducible components.
For any $C \in \Irr (\Xscr_s - F_0)$, we have $\ord_{C} ( \mathscr{W}_{P} ) = 0$.
For each $j = 0 , \ldots , r-1$, 
let $\lambda_{j}$ denote the multiplicity 
at $p_{j}$, which also equals the length of the canonical $1$-simplex 
$\Delta_{p_j}$ corresponding to $p_{j}$ (cf. \S\ref{subsection:skeleta}); 
then for any $\alpha = 1, \ldots, r$, 
we have $\ord_{E_{\alpha}} ( \mathscr{W}_{P} ) = \sum_{j=0}^{\alpha-1} \lambda_{j}$, 
which also equals to the distance
between $[E_0]$ and $[E_{\alpha}]$. 
For any 
$C \in \Irr ( F_0) \setminus  \Irr ( E)$,
since $F_0$ is a tree of irreducible components,
there exists a unique $\alpha_C = 1 , \ldots , r$
such that $[E_{\alpha_C}]$ is the point in
$\{ [E_\alpha ] \mid 1 \leq \alpha \leq r \}$
that is the nearest to $[C]$;
then we have
$\ord_{C} ( \mathscr{W}_P ) = 
\ord_{E_{\alpha_C}} ( \mathscr{W}_P )
= \sum_{j=0}^{\alpha_C-1} \lambda_{j}$.
To put together, we have the following. 
\begin{equation}
\label{eqn:remark:order:stepwise:leaf}
\ord_{C} \left( \mathscr{W}_P \right)
=
\begin{cases}
0
&
\text{if $C \in \Irr (\Xscr_s - F_0)$,}
\\
\sum_{j=0}^{\alpha-1} \lambda_{j} & 
\text{if $C = E_{\alpha}$ ($\alpha = 1 , \ldots , r$),}
\\
\sum_{j=0}^{\alpha_C-1} \lambda_{j}
&
\text{if $C \in \Irr ( F_0) \setminus  \Irr ( E)$}. 
\end{cases}
\end{equation}

We also compute the degree of the stepwise vertical divisor $\mathscr{W}_P$ over each $C \in \Irr(\Xscr_s)$. 
\begin{equation}
\label{eqn:remark:degree:stepwise:leaf}
\deg \left(\rest{\OO(\mathscr{W}_P}{C}) \right)
=
\begin{cases}
1
&
\text{if $C = E_0$,}
\\
-1 & 
\text{if $C = E_r$,}
\\
0
&
\text{otherwise.} 
\end{cases}
\end{equation}

\subsection{Base sections and $P$-unimodularity sections}
\label{subsec:base:sections:and:P:unimodulatity:sections}
We define key notions of global sections of a model,
which will give global sections of $L$ 
used for a faithful tropicalization of an arbitrary skeleton.

Let $\Gamma$ be a compact skeleton of $X^{\an}$.
We fix a minimal skeleton $\Gamma_{\min} \subset \Gamma$. 
Let $P$ be a point in $X(K)$. 
Let $(\mathscr{X} ; \sigma)$ be a strictly semistable pair
such that $S( \mathscr{X} ) = \Gamma$ and $\sigma (K) = P$.
Let $\mathscr{L}$ be a line bundle over $\mathscr{X}$
such that $( \mathscr{X} , \mathscr{L})$ is a model of $(X,L)$.
Let $B$ be a finite subset of $\mathscr{X}_s (k) \setminus 
( \Sing ( \mathscr{X}_s) \cup \{ \sigma (k) \})$ (We allow $B = \emptyset$).

\begin{Definition}[base section and $P$-unimodularity section]
\label{def:basesection:WBsection:leaf}
Let $\Gamma$, $\Gamma_{\min}$, $P$, $(\mathscr{X} ; \sigma)$, $\Lscr$, and $B$  be 
as above. 
\begin{enumerate}
\item
A nonzero global section $\widetilde{s}_0$ of $\mathscr{L}$ 
is called a \emph{base section with respect to $B$} if
$\widetilde{s} (p) \neq 0$ for any $p \in 
B \cup \Sing ( \mathscr{X}_s ) \cup \{ \sigma (k) \}$.
\item
A nonzero global section $\widetilde{s}_1$ of $\mathscr{L}$ 
is called a \emph{$P$-unimodularity section with
respect to $B$ } 
if $\zero (\widetilde{s}) - \mathscr{W}_{P} - \sigma (R)$ 
is trivial on some open neighborhood of $B \cup \Sing ( \mathscr{X}_s ) \cup \{ \sigma (k) \}$. 
\end{enumerate}
\end{Definition}

The following lemma shows us importance
of a base section and a $P$-unimodularity section.

\begin{Lemma} \label{lemma:faithful:ends2}
Let $\Gamma$, $\Gamma_{\min}$, $P$, $(\mathscr{X} ; \sigma)$, $\Lscr$, and $B$  be 
as above. 
Suppose that there exist 
a base section $\widetilde{s}_{0}$ of $\mathscr{L}$ with respect to $B$
and a $P$-unimodularity section $\widetilde{s}_{1}$ of $\mathscr{L}$ with respect to $B$. 
Set $s_0 := \rest{\widetilde{s}_{0}}{X}$
and  $s_{1} := \rest{\widetilde{s}_{1}}{X}$.
Further, set $h := s_{1} / s_0$,
which is a non-zero rational function on $X$.
Define $\varphi :
X^\an \setminus X (K) \to \RR$
by $\varphi := - \log |h|$.
Then the following hold.
\begin{enumerate}
\item
The function
$\varphi$ gives an isometry from $[\tau_{\Gamma_{\min}} (P) , P) 
$
to $[ 0 , + \infty )$.
\item
The restriction of $\varphi$ to $\Gamma \setminus 
(\tau_{\Gamma_{\min}} (P) , P)$ is locally constant.
\item
Let $P_1$ be a point in $X(K)$. 
Suppose that $\red_{\mathscr{X}} (P_1)  \in B$. 
Then 
the restriction of $\varphi$ to $[\tau_{\Gamma} (P_1) , P_1 )$
is constant.
\end{enumerate}
\end{Lemma}

\Proof
Set $g := \widetilde{s}_{1}
/\widetilde{s}_{0}$,
which is a non-zero rational function on $\mathscr{X}$.
Note that $\rest{g}{X} = h$. 
By the definitions of a base section  with respect to $B$ 
and a $P$-unimodularity section with respect to $B$, 
there exists an open neighborhood $\mathscr{U}$ of 
$B \cup \Sing ( \mathscr{X}_s ) \cup \{ \sigma (k) \}
$
such that $\zero (g) - \mathscr{W}_{P} - \sigma (R)$
is trivial on $\mathscr{U}$. 

\medskip
Let us prove (1).  
Let $E$ be the chain in $\mathscr{X}_s$ such that 
$\Delta_E = [\tau_{\Gamma_{\min}} (P) , \tau_{\Gamma} (P)] $.
Let $\lambda_{E}$ denote the length of
$\Delta_{E} \subset \Gamma$.
(With the notation
in \S\ref{subsection:setting1}, 
$\lambda_{E}$ equals the distance between $[E_0]$
and $[E_{r}]$ in $\Gamma$, 
which is equal to $\sum_{j=0}^{r-1} \lambda_{j}$.) 
Since $\zero (g) - \mathscr{W}_{P} - \sigma (R)$
is trivial on $\mathscr{U}$,
the equality \eqref{eqn:remark:order:stepwise:leaf}
shows that
$\varphi$ gives an isometry from $\Delta_{E}$
to $[0 , \lambda_{E}]$.

Let $\ell$ be a local equation of the Cartier divisor $\sigma (R)$
and we take $\varpi \in K^{\times}$ such that $v_K ( \varpi ) = 
\lambda_{E}$.
Then there exists a rational function 
$u$
on $\mathscr{X}$
such that
$u$
is a unit regular function on some open neighborhood of $\sigma (R)$
and satisfies $g = \varpi u \ell$.
Since $- \log 
\left|
\rest{\ell}{X}
\right|$ gives an isometry from  $\Delta (\sigma)= [ \tau_{\Gamma} (P) , P)$
to $[0 , + \infty )$,
$\varphi = - \log |h|$
gives an isometry from 
$\Delta (\sigma) $ to $[ \lambda_{E} , + \infty)$.
Since 
we have the simplicial decomposition
$[ 
\tau_{\Gamma_{\min}} (P) , P ) = \Delta_{E} \cup \Delta ( \sigma)$,
this concludes that
$\varphi$ gives an isometry from $[ 
\tau_{\Gamma_{\min}} (P) , P )$
to $[ 0 , + \infty)$.
Thus we have (1). 

\medskip
Let us prove (2).
Let $\Gamma^{\circ}$ be a connected component
of $\Gamma \setminus (\tau_{\Gamma_{\min}} (P) , P)$.
Let $D$ be the connected curve in $\mathscr{X}_s$ such that
$\Delta_D$ equals the closure of $\Gamma^{\circ}$ (cf. Definition~\ref{def:canonical:1:simplex:curve}). 
Let $F_0 , \ldots , F_{r-1}$ be the connected curves in $\mathscr{X}_s$
as in the paragraph where we define $\mathscr{W}_P$.
Then one sees that for each $\alpha = 0 , \ldots r-1$,
we have $\Irr (D) \subset  \Irr (F_{\alpha})$
or $\Irr (D) \cap  \Irr (F_{\alpha}) = \emptyset$.
It follows from the definition of $\mathscr{W}_P$ that
there exists $\varpi' \in K^{\times}$
such that for any $C \in \Irr (D)$,
$\ord_{C} (\mathscr{W}_P) = - \log | \varpi' | = : \lambda'$.
Since $\zero (g) - \mathscr{W}_{P} - \sigma (R)$
is trivial on $\mathscr{U}$,
this shows that on some neighborhood of $D \cap \mathscr{U}$,
$\zero (g) = \zero ( \varpi' )$.
Since $D \cap \Sing ( \mathscr{X}_s ) \subset \mathscr{U}$,
it follows that
the restriction of $\varphi$ to $\Delta_D$
equals the constant $\lambda'$.
Thus $\varphi$ is constant over $\Gamma^{\circ}$,
and we have~(2). 

\medskip
Let us prove (3). 
Let $\sigma_1$ be the section of $\mathscr{X} \to \Spec (R)$
such that $\sigma_1 (K) = P_1$.
Since 
$\sigma_1 (k) = \red_{\mathscr{X}} (P_1) \in 
B \subset 
\mathscr{X}_s (k) \setminus 
 \Sing ( \mathscr{X}_s ) $,
$(\mathscr{X} , \sigma_1 )$
is a strictly semistable pair,
and $\Delta ( \sigma_1 ) = [ \tau_{\Gamma} (P_1) , P_1)$ by 
Lemma~\ref{lemma:end=canonialend}.
Recall that $\zero (g) - \mathscr{W}_{P} - \sigma (R)$
is trivial on the neighborhood $\mathscr{U} \supset B$.
Then, since $\sigma_1 (k) \in B$
and $B \cap  ( \Sing (\mathscr{X}_s) \cup \{ \sigma (k) \}) = \emptyset$,
there exist $\varpi_1 \in R \setminus \{ 0 \}$
and a rational function $u$ that is a unit regular function at $\sigma_1 (k)$
such that $g = \varpi_1 u$
on some neighborhood of $\sigma_1 (k)$.
Since $\varphi = - \log |\rest{g}{X}|$,
that shows that $\varphi$ equals the constant $- \log |\varpi_1|_K$
on $\Delta ( \sigma_1)$.
This proves (3).
\QED

We want to show the following proposition, which assures the existence 
of a model $(\Xscr, \Lscr)$ 
that one can apply Lemma~\ref{lemma:faithful:ends2}.  
Here $t(g)$ is the quantity in Theorem~\ref{thm:main:unimodular:faithful}. 

\begin{Proposition}
\label{prop:main:FT:ends}
Let $\Gamma$ be a compact skeleton of $X^{\an}$.
We fix a minimal skeleton
$\Gamma_{\min} \subset \Gamma$.
Let $P$ be a point in $X(K)$.
Assume that $\deg (L) \geq t(g)$.
Then there exists a
model $(\mathscr{X}, \mathscr{L})$ of 
$(X,L)$ with the following properties. 
\begin{enumerate}
\item[(i)]
There 
exists a section $\sigma$
of $\mathscr{X} \to \Spec (R)$ such that 
$( \mathscr{X} ; \sigma )$ is
a strictly semistable pair with $\Gamma
= S ( \mathscr{X} )$ 
and $\sigma (K) = P$.
\item[(ii)]
For any finite subset $B \subset \mathscr{X}_s (k)
\setminus ( \Sing (\mathscr{X}_s) \cup \{ \sigma(k) \})$,
$\mathscr{L}$
has a base section with respect to $B$
and
a $P$-unimodularity section with respect to $B$ 
\textup{(}cf. Definition~\textup{\ref{def:basesection:WBsection:leaf}}\textup{)}.
\end{enumerate}
Further, if $\Xscr^0$ is a strictly semistable model of $X$ with $S(\Xscr^0) = \Gamma$, then 
we may take $\Xscr$ so that $\Xscr$ dominates $\Xscr^0$. 
\end{Proposition}

The proof of Proposition~\ref{prop:main:FT:ends}
will be given in \S\ref{subsection:proof:prop:main:FT:ends}.

\subsection{Good model}
\label{subsec:good:model:classical:point}
In this subsection, we fix a minimal skeleton $\Gamma_{\min}$
of $X^{\an}$. When $g \geq 2$, we endow $\Gamma_{\min}$ with 
a canonical weight function and a canonical finite graph structure 
in~\S\ref{subsec:skeleton:weighted:metric:graph}. 

\begin{Definition}[good model for a point $P \in X(K)$]
\label{def:goodDMmodel:end}
Let $P$ be a point in $X(K)$. 
Let $( \mathscr{X} , \mathscr{L})$ be a 
model of $(X,L)$. 
We call $( \mathscr{X} , \mathscr{L})$ a 
\emph{good model for $P$} if it satisfies the following conditions, 
where $\mathscr{M}:= \mathscr{L} \otimes \omega_{\mathscr{X}/R}^{\otimes -1}$.
\begin{enumerate}
\item[(i)]
The model $\mathscr{X}$ is a 
strictly semistable model of $X$ such that $S ( \mathscr{X})
= \Gamma_{\min}$.
\item[(ii)]
Let $C_P$ be the irreducible component of $\mathscr{X}_s$
with $\red_{\mathscr{X}} (P) \in C_P$
and set $v_P := [ C_P ]$.
Then
$
D_{\mathscr{M}} - [v_P]$ is effective on $\Gamma_{\min}$.  
\item[(iii)]
We have $\deg (D_{\mathscr{M}}) \geq 3$.
Further, if $g \geq 2$,
then
for any island $\Gamma_i$ of $\Gamma_{\min}$,
we have
$\deg
\left(
\rest{(D_{\mathscr{M}}
 - [v_P])}{\Gamma_i}
\right) \geq 1$. 
\end{enumerate}
\end{Definition}

If $g \geq 1$, then condition (i) in Definition~\ref{def:goodDMmodel:end}
is equivalent to $\mathscr{X}$ being a Deligne--Mumford strictly
semistable model of $X$.

We prove a couple of lemmas that will be used
in constructing a base section and a $P$-unimodularity section.

\begin{Lemma}
\label{lemma:vanishing:for:ends:new-1}
Let $( \mathscr{X} , \mathscr{L})$
be a good model of $(X,L)$
for $P$, and set 
$\mathscr{M} := \mathscr{L} \otimes
\omega_{\mathscr{X}/R}^{\otimes -1}$ as before. 
Let $q \in \mathscr{X}_s ( k )$  be a point.
\begin{enumerate}
\item
Suppose that 
$q \notin \Sing ( \mathscr{X}_s)$.
Then we have
$h^{0}
\left(
\left(
\rest{\mathscr{M}}{\mathscr{X}_s} 
\left( - q -  \red_{\mathscr{X}} (P) 
\right)
\right)^{\otimes -1}
\right) = 0$
and
$
h^{0}
\left(
\rest{\mathscr{M}}{\mathscr{X}_s} ( -  \red_{\mathscr{X}} (P) 
)^{\otimes -1}
\right) = 0$.
\item
Suppose that $q \in \Sing ( \mathscr{X}_s)$.
Let $\nu : \widetilde{\mathscr{X}_s} \to \mathscr{X}_s$
be the partial normalization at $q$.
Then 
$
h^{0}
\left(
\nu^{\ast} 
\left(
\rest{\mathscr{M}}{\mathscr{X}_s} ( -  \red_{\mathscr{X}} (P) )
\right)^{\otimes -1}
\right) = 0$.
\end{enumerate}
\end{Lemma}

\Proof
We set
$M := \rest{\mathscr{M}}{\mathscr{X}_s}
(- \red_{\mathscr{X}} (P))$
and $D_{M} := 
\sum_{C \in \Irr ( \mathscr{X}_s)} \deg
\left(
\rest{M}{C}
\right) [C]$.
Then $D_{M} = D_{\mathscr{M}} - [v_P]$. 
Thus 
$D_{M}$ is effective by Definition~\ref{def:goodDMmodel:end}(ii)
and $\deg(D_M) \geq 2$  by Definition~\ref{def:goodDMmodel:end}(iii).

\smallskip
(1)
Suppose that $q \notin \Sing ( \mathscr{X}_s )$. 
We take a unique $C_1 \in \Irr ( \mathscr{X}_s )$
such that $q \in C_1$,
and we set $v_1 := [C_1] \in V ( \mathscr{X} )$.
By Definition~\ref{def:goodDMmodel:end}(iii),
we have $\deg
\left(
D_{M} - [v_1]
\right) \geq 1$.

Let us show the first equality 
$h^0 \left( 
(M(-q))^{\otimes -1}
\right) = 0$. 
If $D_{M} - [v_1]$ is effective, then 
$M (-q)$ is nef and $\deg(M (-q)) > 0$. 
By Lemma~\ref{lemma:vanishing0},
we get 
$h^0 \left( 
(M(-q))^{\otimes -1}
\right) = 0$.
Suppose that $D_{M} - [v_1]$ is not effective.
Then $(D_{M} - [v_1]) (v_1) = -1$,
and $v_1$ is the only point at which $D_{M} - [v_1]$ is negative.
Further,
we note that the valence of $\Gamma_{\min}$ at $v_1$
is at least $2$.
Indeed, suppose that this is not the case.
Then $\{ v_1 \}$ is an island,
and hence by condition (iii) of Definition~\ref{def:goodDMmodel:end},
$D_{M} (v_1) \geq 1$.
This means that 
$D_{M} - [v_1]$ is effective,
which is a contradiction.
This proves that the valence of $\Gamma_{\min}$ at $v_1$
is greater than or equal to $2$.
Thus we see that the assumptions described in 
Remark~\ref{remark:proofs:lemma:freeatnode:graphversion:emptyset}(2) 
are fulfilled,
and hence by this remark,
we obtain $h^{0}
\left(
\left(
M (-q)
\right)^{\otimes -1}
\right) = 0$,
which is
the first equality. 

Since 
$M^{\otimes -1} \subset M (-q)^{\otimes -1}$
the second equality follows from the first one.

\smallskip
(2)
Suppose that $q \in \Sing ( \mathscr{X}_s )$.
Set
$\widetilde{\Gamma} := \Gamma_{\min} \setminus \mathrm{relin} (\Delta_q)$.
Then $D_M$ is effective by Definition~\ref{def:goodDMmodel:end}(ii)
and has positive degree on any connected component of $\widetilde{\Gamma}$
by Definition~\ref{def:goodDMmodel:end}(iii).
Thus, as noted in Remark~\ref{remark:proofs:lemma:freeatnode:graphversion:emptyset}(1),
we obtain
$
h^{0}
\left(
\nu^{\ast}(M)^{\otimes -1}
\right) = 0$.
\QED

\begin{Lemma}
\label{lemma:vanishing:for:ends:new}
With notation and assumptions being in Lemma~\textup{\ref{lemma:vanishing:for:ends:new-1}}, 
we have the following:
\begin{enumerate}
\item
$\rest{\mathscr{L}}{\mathscr{X}_s} \left( 
- \red_{\mathscr{X}} (P) 
\right)$ is free at any $q \in \mathscr{X}_s(k)$.
\item
$\mathscr{L}$ is free at any $q \in \mathscr{X}_s(k)$.
\end{enumerate}
\end{Lemma}

\Proof
Assertion (1) follows from
Lemma~\ref{lemma:nonbasepoint} and Lemma~\ref{lemma:vanishing:for:ends:new-1}.
We prove (2). 
For any $q \in \mathscr{X}_s (k) \setminus \Sing ( \mathscr{X}_s )$,
we have
$h^{0}
\left(
\left(
\rest{\mathscr{M}}{\mathscr{X}_s} 
\left( - q -  \red_{\mathscr{X}} (P) 
\right)
\right)^{\otimes -1}
\right) = 0$
by Lemma~\ref{lemma:vanishing:for:ends:new-1}(1),
and
hence
$h^{0}
\left(
\rest{\mathscr{M}}{\mathscr{X}_s} ( -  q
)^{\otimes -1}
\right) = 0$.
For any
$q \in \Sing ( \mathscr{X}_s )$,
if $\nu : \widetilde{\mathscr{X}_s} \to 
\mathscr{X}_s$ denotes the partial normalization at $q$,
then
$
h^{0}
\left(
\nu^{\ast} 
\left(
\rest{\mathscr{M}}{\mathscr{X}_s} ( -  \red_{\mathscr{X}} (P) )
\right)^{\otimes -1}
\right) = 0$
by Lemma~\ref{lemma:vanishing:for:ends:new-1}(2),
and hence
$h^{0}
\left(
\nu^{\ast} 
\left(
\rest{\mathscr{M}}{\mathscr{X}_s}
\right)^{\otimes -1}
\right) = 0$.
By Lemma~\ref{lemma:nonbasepoint},
it follows that
$\rest{\mathscr{L}}{\mathscr{X}_s}$ is base-point free.
The rest of the argument is the same as 
the proof of Lemma~\ref{lemma:freeatregular:graphversion:effective} when $q \in \mathscr{X}_s (k) \setminus \Sing ( \mathscr{X}_s )$ and 
that
of Lemma~\ref{lemma:freeatnode:graphversion:emptyset} when $q \in \Sing ( \mathscr{X}_s )$. 
\QED

\subsection{Proof of Proposition~\ref{prop:main:FT:ends}}
\label{subsection:proof:prop:main:FT:ends}

Let $\Gamma$ be a compact skeleton of $X^{\an}$.
We fix a minimal skeleton
$\Gamma_{\min} \subset \Gamma$.
For any model
$\mathscr{X}$ of $X$ with $S (\mathscr{X}) = \Gamma$,
by \cite[Theorem~4.11]{BPR2},
there exist a unique model $\mathscr{X}^{\min}$ of $X$ and 
a homomorphism $\mu : \mathscr{X} \to \mathscr{X}^{\min}$ extending the
identity on $X$ such that
$S ( \mathscr{X}^{\min} ) = \Gamma_{\min}$ (namely, $\mathscr{X}^{\min}$ 
is minimal)
and $V( \mathscr{X}^{\min}) = 
V ( \mathscr{X} ) \cap \Gamma_{\min}$.
We use this notation in the sequel.

In this subsection, we prove Proposition~\ref{prop:main:FT:ends}.
First, we construct a model that dominates a good model for a point in $X(K)$.

\begin{Proposition}
\label{prop:good:Gamma-model:ends}
Let $P$ be a point in $X(K)$. 
Assume that
$\deg (L) \geq t(g)$.
Then there exists a
model $(\mathscr{X}, \mathscr{L})$ of 
$(X,L)$ with the following properties. 
\begin{enumerate}
\item[(i)]
There 
exists a section $\sigma$
of $\mathscr{X} \to \Spec (R)$ such that 
$( \mathscr{X} ; \sigma )$ is
a strictly semistable pair with $\Gamma
= S ( \mathscr{X} )$ and 
and $\sigma (K) = P$.
\item[(ii)]
Let
$\mathscr{X}^{\min}$ 
and 
$\mu: \mathscr{X}\to \mathscr{X}^{\min}$ be 
as above.
Then
there exists a line bundle $\mathscr{L}^{\min}$ over $\mathscr{X}^{\min}$
such that $( \mathscr{X}^{\min} , \mathscr{L}^{\min})$ is a good
model for $P$ 
\textup{(}cf. Definition~\textup{\ref{def:goodDMmodel:end}}\textup{)}
and $\mathscr{L} = \mu^{\ast} (\mathscr{L}^{\min})$.
\end{enumerate}
Further, if $\Xscr^0$ is a strictly semistable model of $X$ with $S(\Xscr^0) = \Gamma$, then 
we may take $\Xscr$ so that $\Xscr$ dominates $\Xscr^0$.
\end{Proposition}

\Proof
By Lemma~\ref{lemma:retraction:rational}(1),
we have $\tau_{\Gamma} (P) \in \Gamma_{\Lambda}$. 
Let $\Xscr^0$ be a strictly semistable model of $X$ with $S(\Xscr^0) = \Gamma$. 
By Lemma~\ref{prop:subdivision1:b},
there exists a strictly semistable model
$\mathscr{X}^1$ of $X$ such that $\Gamma = S( \mathscr{X}^1)$
and $V(\mathscr{X}^1) = V( \mathscr{X}^0) \cup \{ \tau_{\Gamma} (P) \}$.
We take a minimal model $\mathscr{X}^{1,\min}$ of $X$ such that
$S(\mathscr{X}^{1 , \min}) = \Gamma_{\min}$ and $V(\mathscr{X}^{1,\min})
= V(\mathscr{X}^{1}) \cap \Gamma_{\min}$.
We note that $\mathscr{X}^{1,\min}$ is dominated by $\mathscr{X}^1$. 
Let $C_P$ be the irreducible component of $\mathscr{X}^{1,\min}_s$ 
with $\red_{\mathscr{X}^{1,\min}}(P) \in C_P$, and we set $v_P := [C_P]$. 

\medskip
{\bf Case 1.}\quad 
Suppose that $g \geq 2$. Then $\mathscr{X}^{1,\min}$
is a Deligne--Mumford strictly semistable model of $X$.
Since $\deg (L) \geq t(g) = 3g-1$,
we use Proposition~\ref{prop:pregoodmodel} with $v_P$
in place of $x$ to obtain 
a model $( \mathscr{X}^2 , \mathscr{L}^2)$ with properties 
in Proposition~\ref{prop:pregoodmodel} such that the identity morphism on $X$ extends to a morphism $\mathscr{X}^2
\to \mathscr{X}^{1,\min}$.
We remark that $S(\mathscr{X}^2) = \Gamma_{\min}$. 
We set $\mathscr{X} := 
\mathscr{X}^1 \times_{\mathscr{X}^{1 , \min}} \mathscr{X}^2$
and let $\mathscr{L}$ be the pullback of $\mathscr{L}^2$ to $\mathscr{X}$.
Then by the construction of $\mathscr{X}^2$ and $\mathscr{X}$,
we have $\mathscr{X}^{\min} = \mathscr{X}^2$. 
Set $\mathscr{L}^{\min} := \mathscr{L}^2$.
Since $\left(\mathscr{X}^2, \Lscr^2\right)$ has the properties in 
Proposition~\ref{prop:pregoodmodel}, we see that 
$\left(\mathscr{X}^{\min}, \mathscr{L}^{\min} \right)$ is a good model 
for $P$ and thus has
property (ii) of Proposition~\ref{prop:good:Gamma-model:ends}. 
Since
$\deg\left( \mathscr{L}^{\min}\otimes \omega_{\mathscr{X}^{\min}/R}^{\otimes -1} \right) 
\geq (3g-1) - (2g-2) = g + 1 \geq 3$, the first condition 
of Definition~\ref{def:goodDMmodel:end}(iii) is satisfied. 

By the valuative criterion of properness, we take a section $\sigma$ of $\mathscr{X} \to \Spec (R)$
such that $\sigma (K) = P$.
By Lemma~\ref{lemma:retraction:rational}(2),
$\sigma (k) \in \mathscr{X}_s (k) \setminus \Sing (\mathscr{X}_s)$,
and thus $( \mathscr{X} ; \sigma)$ is a strictly semistable pair.
By the construction of $\Xscr$, we have $S(\Xscr) = \Gamma$, and thus
$\mathscr{X}$ has 
property (i) in the proposition.

\medskip
{\bf Case 2.}\quad 
We consider the case where
$g = 0$ or $g = 1$. 
First, suppose that $g = 1$ and that the minimal skeleton is not a singleton.
Then a similar argument using
Lemma~\ref{lemma:pregoodmodel:g=1} 
instead of  Proposition~\ref{prop:pregoodmodel}
gives a model with the required properties. 
Here, we remark that
$ \deg
\left(\rest{\Lscr \otimes\omega_{\mathscr{X}^{\min}/R}^{\otimes -1}}{\mathscr{X}^{\min}_s}\right) = \deg(L) \geq t(1) = 3$,
since $\omega_{X}$ is trivial. 

Suppose that 
$g = 1$
and that
the minimal skeleton is a singleton,
or suppose that $g = 0$.
In this case,  $\mathscr{X}^{1 , \min}$ is
a smooth proper model,
and there exists a line bundle
 $\mathscr{L}^{1 , \min}$ over $\mathscr{X}^{1,\min}$ such that 
$\rest{\mathscr{L}^{1,\min}}{X} = L$.
Set $\mathscr{X} := \mathscr{X}^1$ and 
$\mathscr{L} := \mu^{\ast} ( \mathscr{L}^{1,\min})$, where 
$\mu : \mathscr{X} 
\to \mathscr{X}^{1,\min}$ is the morphism extending the identity morphism on $X$.
Then it is straightforward to check that $( \mathscr{X}, \mathscr{L})$
has the required properties.
\QED

\begin{Lemma}
\label{lemma:proof:key:ends}
Let $P$ be a point in $X(K)$.
Let $( \mathscr{X} ; \sigma )$ be a strictly semistable pair 
such that $S ( \mathscr{X} ) = \Gamma$ and $\sigma (K) = P$.
Let $\mathscr{L}$ be a line bundle over $\mathscr{X}$.
Assume that 
there exists a line bundle $\mathscr{L}^{\min}$
over $\mathscr{X}^{\min}$ such that 
$( \mathscr{X}^{\min} , \mathscr{L}^{\min})$ is a good model for $P$
and $\mathscr{L} = \mu^{\ast} ( \mathscr{L}^{\min})$.
Then 
for any finite subset $B \subset \mathscr{X}_s  (k) \setminus 
( \Sing ( \mathscr{X}_s) \cup \{ \sigma(k) \})$,
$\mathscr{L}$ has a base section with respect to $B$
and a $P$-unimodularity section with respect to $B$.
\end{Lemma}

\Proof
Let $B$ be any finite subset of $\mathscr{X}_s (k) \setminus
( \Sing ( \mathscr{X}_s) \cup \{ \sigma(k) \})
$.
First, we construct a base section with respect to $B$.
By Lemma~\ref{lemma:vanishing:for:ends:new}(2),
there exists an $\widetilde{s}_0^{\,\min}
\in H^0 
( \mathscr{L}^{\min} )$
such that 
$\widetilde{s}_0^{\,\min} ( p ) \neq 0$ for any $p \in \mu \left(B \cup \Sing ( \mathscr{X}_s \right)
\cup \{ \sigma (k) \})$.
Set $\widetilde{s}_0 := \mu^{\ast} ( \widetilde{s}_0^{\,\min} ) \in H^0 \left( \mathscr{L} \right)$.
Then
we have $\widetilde{s}_0 (p) \neq 0$ for any 
$p \in B \cup \Sing ( \mathscr{X}_s) \cup \{ \sigma (k) \}$.
This shows that $\widetilde{s}_0$ is a base section
with respect to $B$.

In the following, we construct a $P$-unimodularity section with respect to $B$. 
We use the notation in \S\ref{subsection:setting1}. 
In particular, $\mathscr{W}_P$ is the stepwise vertical divisor 
associated to $P$ defined in \eqref{align:stepwise:for:ends}, 
and $F_0$ is the support of $\mathscr{W}_P$. 
We take a unique
connected curve $D$ in $\mathscr{X}_s$
such that $\mu$ 
restricts to an isomorphism $\phi : D \to \mathscr{X}^{\min}_s$. 
We remark that $D \cap \sigma (R) \neq \emptyset$ if and
only if $\mathscr{W}_P = 0$,
and if this is the case, then $D \cap \sigma (R)$ is a singleton.
Let $p_0$ be the point in $\Xscr_s (k)$ 
such that
\[
\{ p_0 \}
=
\begin{cases}
D \cap F_0
& \text{if $\mathscr{W}_P \neq 0$,}
\\
D \cap \sigma (R)
& \text{otherwise.}
\end{cases}
\] 
We set
$
\mathscr{L}_-
:=
\mathscr{L}
\left(
-
\left(
\mathscr{W}_{P}
+
\sigma (R)
\right)
\right)
$. 
Then $\rest{\OO_{\mathscr{X}} 
\left( \mathscr{W}_{P} + \sigma (R)
\right)}{D} = \OO_{D} (p_0)$.
Since $\phi ( p_0) = \red_{\mathscr{X}^{\min}} ( P )$,
we have $\OO_{D} ( p_0) = \phi^\ast \left( \OO_{\mathscr{X}^{\min}_s} ( 
\red_{\mathscr{X}^{\min}} ( P ) ) \right)$.
It follows that
{\allowdisplaybreaks
\begin{align*}
\rest{\mathscr{L}_-}{D}
&=
\rest{\mathscr{L}
\left(
-
\left(
\mathscr{W}_{P}
+
\sigma (R)
\right)
\right)}{D}
\\
&
=
\rest{
\mu^{\ast}
\left(
\mathscr{L}^{\min}
\right)}{D}
\otimes
\phi^\ast \left( \OO_{\mathscr{X}^{\min}_s} (- 
\red_{\mathscr{X}^{\min}} ( P )
 ) \right)
\\
&=
\phi^\ast
\left(
\rest{\mathscr{L}^{\min}}{\mathscr{X}^{\min}_s} 
(- 
\red_{\mathscr{X}^{\min}} ( P )
 )
\right)
.
\end{align*}}
Since $\phi$ is an isomorphism,
it follows from
Lemma~\ref{lemma:vanishing:for:ends:new}
that $
\rest{\mathscr{L}_-}{D}
$ is basepoint free.
Thus there exists
$\xi_- \in H^0
\left(
\rest{\mathscr{L}_-}{D}
\right)$
such that $\xi_- (q) \neq 0$
for any 
$
q \in ( B  \cup \Sing ( \mathscr{X}_s ) \cup \{ \sigma (k) \})\cap D
$.

Let us define 
$\eta_- \in H^0\left(\rest{\mathscr{L}_-}{\mathscr{X}_s}\right)$ such that
$\eta_- (q) \neq 0$
for any $q \in 
B \cup \Sing ( \mathscr{X}_s ) \cup \{ \sigma (k) \}
$.
If $\mathscr{X}_s = D$,
then we set $\eta_- := \xi_-$,
which satisfies the required condition. 
Suppose that
$\mathscr{X}_s \neq D$.
Since $\mu ( \mathscr{X}_s - D )$ is a 
finite set of points, $\rest{\mathscr{L} 
}{\mathscr{X}_s - D} \cong \OO_{\mathscr{X}_s - D}$.
Further, it follows from \eqref{eqn:remark:degree:stepwise:leaf} 
that
$\rest{\OO_{\mathscr{X}} \left( \mathscr{W}_{P}
+
\sigma(R) \right)}{\mathscr{X}_s - D} \cong \OO_{\mathscr{X}_s - D}$.
Thus 
$\rest{\mathscr{L}_-
}{\mathscr{X}_s - D} \cong \OO_{\mathscr{X}_s - D}$.
It follows that there exists a
(unique) section $\eta_-
\in H^0 \left(
\rest{\mathscr{L}_-
}{\mathscr{X}_s}
\right)$
such that $\rest{\eta_-}{D} = \xi_-$. 
We remark that $\phi^{-1}\left(\mu ( \mathscr{X}_s - D )\right) 
\subset \Sing ( \mathscr{X}_s ) \cap D$. 
Since
$\xi_- (q) \neq 0$
for any 
$
q \in ( B  \cup \Sing ( \mathscr{X}_s ) \cup \{ \sigma (k) \}) \cap D
$,
it follows that
$\eta_-$ is nowhere vanishing over $\mathscr{X}_s - D$.
This shows that
$\eta_- (q) \neq 0$
for any $q \in B  \cup \Sing ( \mathscr{X}_s ) \cup \{ \sigma (k) \}
$.

Put $M := \rest{\mathscr{L}_- \otimes \omega_{\mathscr{X}/R}^{\otimes -1}}{\mathscr{X}_s}$.
Put $\Sigma := D \cap  \left( \mathscr{X}_s -D \right)$.
By the adjunction formula,
we see that
{\allowdisplaybreaks
\begin{align*}
\rest{M}{D} &= 
\rest{\mathscr{L}_-}{D} \otimes 
\omega_{D} ^{\otimes -1}
( -\Sigma )
\\
&=
\phi^{\ast}
\left(
\rest{\mathscr{L}^{\min}  \otimes 
\omega_{\mathscr{X}^{\min}/R} ^{\otimes -1}}{\mathscr{X}^{\min}_s}
(- \red_{\mathscr{X}^{\min}} (P))
\right)
( - \Sigma )
\\
&=
\phi^\ast 
\left(
\rest{\mathscr{M}^{\min}}{\mathscr{X}^{\min}_s}
(- \red_{\mathscr{X}^{\min}} (P))
\right)
( - \Sigma ),
\end{align*}
}%
where
$\mathscr{M}^{\min} :=
\mathscr{L}^{\min}  \otimes 
\omega_{\mathscr{X}^{\min}/R} ^{\otimes -1}$.
Since $(\Xscr^{\min}, \Lscr^{\min})$ is a good model for $P$, 
Lemma~\ref{lemma:vanishing:for:ends:new-1}(1) gives 
that
$h^0
\left(\left(
\rest{\mathscr{M}^{\min}}{\mathscr{X}^{\min}_s}
(- \red_{\mathscr{X}^{\min}} (P))
\right)^{\otimes -1}
\right) = 0$.
Since $\phi$ is an isomorphism, we get 
\[
h^0\left((\rest{M}{D} (\Sigma))^{\otimes -1} \right) 
= h^{0}
\left(
\phi^{\ast}
\left(
\rest{\mathscr{M}^{\min}}{\mathscr{X}^{\min}_s}
(- \red_{\mathscr{X}^{\min}} (P))
\right)^{\otimes -1}
\right) = 0. 
\]

Let us prove $h^0 \left( M^{\otimes -1} \right) = 0$.
If $\mathscr{X}_s = D$, then 
$h^0 \left( M^{\otimes -1} \right) = h^0\left((\rest{M}{D} (\Sigma))^{\otimes -1} \right)  = 0$. 
Suppose that $\mathscr{X}_s \neq D$.
We take any connected component 
$F'$
of $\mathscr{X}_s - D$.
We compute $\rest{M}{F'}$ by using the adjunction formula.
Since
$\rest{\mathscr{L}_-}{F'} \cong \OO_{F'}$, as noted above,
we have $\rest{M}{F'} \cong
\rest{\omega_{\mathscr{X}/R}}{F'}^{\otimes -1}$.
Let $p'$ be the point with $\{ p' \} = D \cap F'$.
Then we have $\rest{\omega_{\mathscr{X}/R}}{F'} = \omega_{F'} (p')$,
and thus $\rest{M}{F'}^{\otimes -1} \cong 
\omega_{F'} (-p')$.
Applying the Riemann--Roch formula on $F'$ for $\OO_{F'} (p')$,
we see that
$h^0 \left( \OO_{F'} (p') \right) - 
h^{0} \left(
\omega_{F'} (-p')
\right) = 2$.
Since $h^0 \left( \OO_{F'} (p') \right) =2$,
it follows that $
h^0
\left(
\rest{M}{F'}^{\otimes -1}
\right)
=
h^{0} \left(
\omega_{F'} (-p')
\right) = 0$.
By Lemma~\ref{lemma:vanishing7},
we obtain $h^0
\left(
M^{\otimes -1}
\right) = 0$.

By the Serre duality,
$
h^1
\left(
\rest{\mathscr{L}_-}{\mathscr{X}_s}
\right) =
h^0
\left(
M^{\otimes -1}
\right) = 0
$.
By the base-change theorem,
the restriction map $H^0 \left( 
\mathscr{L}_-\right)
\to 
H^0 \left( 
\rest{\mathscr{L}_-}{\mathscr{X}_s} \right)
$
is surjective.
Thus
there exists a global section $\widetilde{s}_-
\in H^0
\left(
\mathscr{L}_-
\right)$ such that $\rest{\widetilde{s}_-}{\mathscr{X}_s} = \eta_-$.
For 
any 
$
q \in B  \cup \Sing ( \mathscr{X}_s ) \cup \{ \sigma (k) \}
$,
since
 $\eta_- (q) \neq 0$,
we have $\widetilde{s}_- (q) \neq 0$.
Let $\widetilde{s}$ be the image of $\widetilde{s}_-$
by the natural inclusion $\mathscr{L}_- \hookrightarrow \mathscr{L}$.
Then we see that
$\widetilde{s}$
is
a $P$-unimodularity section
with respect to $B$.
Thus we complete the proof of the lemma.
\QED

\medskip
\textsl{Proof of Proposition~\ref{prop:main:FT:ends}.}\quad
Let us complete the proof of Proposition~\ref{prop:main:FT:ends}.
We take a model $( \mathscr{X}, \mathscr{L})$ as in
Proposition~\ref{prop:good:Gamma-model:ends}.
By condition (i) in Proposition~\ref{prop:good:Gamma-model:ends},
there exists a section $\sigma$
of $\mathscr{X} \to \Spec (R)$ such that
$( \mathscr{X} ; \sigma)$
is a strictly semistable pair such that $\Gamma =
S ( \mathscr{X})$ 
and $\sigma(K) = P$.
Thus $( \mathscr{X}, \mathscr{L})$ has property (i) in 
Proposition~\ref{prop:main:FT:ends}.
By condition (ii) in Proposition~\ref{prop:good:Gamma-model:ends},
$( \mathscr{X} , \mathscr{L})$
satisfies the condition of Lemma~\ref{lemma:proof:key:ends}. 
Then by this lemma,
for any finite subset $B \subset \mathscr{X}_s (k)
\setminus ( \Sing (\mathscr{X}_s) \cup \{ \sigma (k) \})$,
$\mathscr{L}$ has a base section with respect to $B$ and
a $P$-unimodularity section with respect to $B$.
Thus $(\mathscr{X} , \mathscr{L})$ has also property (ii) in 
Proposition~\ref{prop:main:FT:ends}.
\QED

\subsection{Proof of Theorem~\ref{thm:main:unimodular:faithful}}
In this section, we complete the proof of 
the main Theorem~\ref{thm:main:unimodular:faithful},
which is restated below.

\begin{TheoremNoNum}[$=$ Theorem~\ref{thm:main:unimodular:faithful}]
Let $X$ be a connected smooth projective curve over $K$ of genus $g \geq 0$.
Let $\Gamma$ be a skeleton of $X^{\an}$.
Let $L$ be a line bundle over $X$.
Suppose that 
\[
\deg (L) 
\geq
t(g) := 
\begin{cases}
1 & \text{if $g=0$,} 
\\
3 & \text{if $g=1$,} 
\\
3g-1 & \text{if $g\geq 2$.} 
\end{cases}
\]
Then there exist $s_0 , \ldots , s_N \in H^0 ( X,L)$
such that the map
$\varphi : X^\an \to \TT\RR^N$
defined by $\varphi := ( - \log |s_0| : \cdots : - \log |s_N|)$
gives a faithful tropicalization of $\Gamma$.
\end{TheoremNoNum}

The following two lemmas will be used to separate points of a skeleton 
with ends. 

\begin{Lemma}
\label{lemma:faithful:ends}
Let $\Gamma$ be a skeleton.
We write
$\Gamma = S ( \mathscr{X} ;
\sigma_1 , \ldots , \sigma_m)$
for some strictly semistable pair $( \mathscr{X} ;
\sigma_1 , \ldots , \sigma_m)$.
Let $\Gamma_{\min} \subset S( \mathscr{X})$ be a minimal skeleton.
Set $P_i := \sigma_i (K)$ for $i = 1, \ldots , m$. 
Assume that $\deg (L) \geq t (g)$.
Then for any $i = 1 , \ldots , m$,
there exist nonzero global sections $s_0^{(i)}$ and
$s_1^{(i)}$ of $L$ such that 
the function $\varphi^{(i)} : \Gamma \to \RR$
defined by
$\varphi^{(i)} := - \log |s_1^{(i)} / s_0^{(i)}|$
has the following properties\textup{:}
\begin{enumerate}
\item[(i)]
$\varphi^{(i)}$ gives an isometry
$[\tau_{\Gamma_{\min}} (P_i) , P_i ) \to [0 , + \infty)$\textup{;}
\item[(ii)]
for $j \neq i$,
the restriction of $\varphi^{(i)}$ to $[\tau_{\Gamma_{\min}} (P_j) , P_j )
\setminus [\tau_{\Gamma_{\min}} (P_i) , P_i )$ is a constant function,
and
if $(\tau_{\Gamma_{\min}} (P_j) , P_j )
\cap (\tau_{\Gamma_{\min}} (P_i) , P_i ) = \emptyset$,
then the constant is $0$\textup{;}
\item[(iii)]
$\varphi^{(i)} ( \Gamma_{\min} ) = \{ 0 \}$.
\end{enumerate}
\end{Lemma}

\Proof
We fix any $i = 1 , \ldots , m$.
Since $\deg ( L ) \geq t (g)$,
we apply Proposition~\ref{prop:main:FT:ends}
in place of $P$, $\Gamma$, and $\mathscr{X}^0$
with $P_i$, $S(\mathscr{X})$, and $\mathscr{X}$ respectively. 
Then we obtain
a model $(\mathscr{X}^i , \mathscr{L}^i)$ 
of $(X,L)$
that satisfies
the following conditions. 
\begin{enumerate}
\item[(a)]
there exists a section $\sigma^{P_i}$ of $\mathscr{X}^i
\to \Spec (R)$ such that $(\mathscr{X}^i ; \sigma^{P_i})$
is a strictly semistable pair and such that 
$S (\mathscr{X}^i ) = S( \mathscr{X})$ 
and $\sigma^{P_i} (K) = P_i$;
\item[(b)]
the identity on $X$ extends to a
morphism $\mu^i : \mathscr{X}^i \to \mathscr{X}$;
\item[(c)]
for any finite subset $B \subset
\mathscr{X}^i_s (k) \setminus 
(\Sing ( \mathscr{X}^i_s) \cup \{ \sigma^{P_i} (k) \} )$,
$\mathscr{L}^i$ has a base section 
with respect to $B$
(for $(\mathscr{X}^i ; \sigma^{P_i})$ in (a) above)
and a $P_i$-unimodularity section
with respect to $B$.
\end{enumerate}

We set $\mathcal{L} ( \Gamma )
= \{ P_1, \ldots , P_m \}$.
Note that $\sigma^{P_i}(k) = \red_{\mathscr{X}^i} (P_i) 
\in \red_{\mathscr{X}^i} ( \mathcal{L} (\Gamma))$,
and
we set $B^i := \red_{\mathscr{X}^i} ( \mathcal{L} (\Gamma))
\setminus \{ \sigma^{P_i}(k)\}$.
We claim the following.

\begin{Claim}
\label{claim:Bi}
\begin{enumerate}
\item[(1)]
Any point in $B^{i}$ is a non-singular point of $\Xscr^i_s$. 
\item[(2)]
For any $j = 1 , \ldots , m$ with $j \neq i$,
we have $\red_{\mathscr{X}^i} (P_j) \in B^i$.
\end{enumerate}
\end{Claim}

We prove the claim. 
Note that $\red_{\mathscr{X}} ( \mathcal{L} (\Gamma) ) =
\{ \sigma_1(k) , \ldots , \sigma_m(k) \}$. Since 
$(\Xscr; \sigma_1, \ldots, \sigma_m)$ is a strictly semistable pair, 
$\red_{\mathscr{X}} ( \mathcal{L} (\Gamma) )$ is contained in the 
set of non-singular points of $\mathscr{X}_s$. 
Since $\rest{\mu^i}{\mathscr{X}^i_s} : \mathscr{X}^i_s \to \mathscr{X}_s$
is a contraction of $(-2)$-curves 
and $\red_{\mathscr{X}} = \rest{\mu^i}{\mathscr{X}^i_s}
\circ \red_{\mathscr{X}^i}$,
we see that any point in $\red_{\mathscr{X}^i} ( \mathcal{L} (\Gamma) )$
and 
hence any point in $B^i$ is a non-singular point of $\Xscr^i_s$. 
Thus we have (1).

Since $( \mathscr{X} ; \sigma_1 , \ldots , \sigma_m)$
is a strictly semistable pair,
$\rest{\red_{\mathscr{X}}}{\mathcal{L}(\Gamma)}$ is injective.
Since $\red_{\mathscr{X}} = \rest{\mu^i}{\mathscr{X}^i_s}
\circ \red_{\mathscr{X}^i}$, it follows that 
$\rest{\red_{\mathscr{X}^i}}{\mathcal{L}(\Gamma)}$ is injective.
Thus
$\# \rest{\red_{\mathscr{X}^i}}{\mathcal{L}(\Gamma)} = m$,
which shows that $\red_{\mathscr{X}^i} (P_1) ,
\ldots , \red_{\mathscr{X}^i} (P_m)$ are all distinct.
Since $B^i = \red_{\mathscr{X}^i} ( \mathcal{L}(\Gamma) ) \setminus 
\{ \red_{\mathscr{X}^i} (P_i) \}$,
(2) in the claim holds.

\medskip
Since we have Claim~\ref{claim:Bi}(1),
property (c) allows us to take a base section 
$\widetilde{s}_0^{(i)}$ of $\mathscr{L}^{i}$
with respect to $B^i$
and a $P_i$-unimodularity section
$\widetilde{s}_1^{(i)}$ of $\mathscr{L}^{i}$
with respect to $B^i$.
We set $s_0^{(i)} := \rest{\widetilde{s}_0^{(i)}}{X}$
and $s_1^{(i)} := \rest{\widetilde{s}_1^{(i)}}{X}$.

Let us show that $\varphi^{(i)} = - \log |s_1^{(i)} / s_0^{(i)}|$
has the required properties.
It follows from Lemma~\ref{lemma:faithful:ends2}(1)
that it has property (i).
By Lemma~\ref{lemma:faithful:ends2}(2),
$\varphi^{(i)}$ is locally constant on $\Gamma \setminus
(\tau_{\Gamma_{\min}} (P_i) , P_i )$.
Take any $j \neq i$.
Since we have  Claim~\ref{claim:Bi}(2),
Lemma~\ref{lemma:faithful:ends2}(3)
implies that
$\varphi^{(i)}$ is constant on $[ \tau_{\Gamma} (P_j), P_j)$. 
Since $[\tau_{\Gamma_{\min}} (P_j) , P_j) 
=  [\tau_{\Gamma_{\min}} (P_j) , \tau_{\Gamma} (P_j)] \cup 
[\tau_{\Gamma} (P_j) , P_j) \subseteq \Gamma \cup [\tau_{\Gamma} (P_j) , P_j)$, 
it follows that $\varphi^{(i)}$ is locally constant
on $[\tau_{\Gamma_{\min}} (P_j) , P_j)  \setminus ( \tau_{\Gamma_{\min}} (P_i) , P_i )$.
Since $[\tau_{\Gamma_{\min}} (P_j) , P_j)  \setminus ( \tau_{\Gamma_{\min}} (P_i) , P_i )$ is connected, 
this proves that $\varphi^{(i)}$ has the former property in (ii).
Since $\tau_{\Gamma_{\min}} (P_i) \in \Gamma_{\min}$
and $\varphi^{(i)} (\tau_{\Gamma_{\min}} (P_i)) = 0$ by property (i), 
Lemma~\ref{lemma:faithful:ends2}(2) implies that $\varphi^{(i)} ( \Gamma_{\min}) = \{ 0 \}$. 
Thus it has property (iii).
Finally, 
suppose that $(\tau_{\Gamma_{\min}} (P_j) , P_j )
\cap
(\tau_{\Gamma_{\min}} (P_i) , P_i ) = \emptyset$.
Then $[\tau_{\Gamma_{\min}} (P_j) , P_j )$
is contained in the connected component of 
$\Gamma \setminus (\tau_{\Gamma_{\min}} (P_i) , P_i )$
that contains $\Gamma_{\min}$.
Since $\varphi^{(i)}$ is locally constant on 
$\Gamma \setminus (\tau_{\Gamma_{\min}} (P_i) , P_i )$ and 
$\varphi^{(i)} ( \Gamma_{\min}) = \{ 0 \}$,
this shows that 
$\varphi^{(i)}$ equals $0$ on $[\tau_{\Gamma_{\min}} (P_j) , P_j )$.
Thus it also has the latter property in (ii).
\QED

\begin{Lemma}
\label{lemma:remark:toseparate:ends}
We keep the notation in Lemma~\textup{\ref{lemma:faithful:ends}}.
Let $x$ and $y$ be distinct points in $\Gamma$.
\begin{enumerate}
\item
If $x , y \in [\tau_{\Gamma_{\min}} (P_i) , P_i)$ for some $i$,
then $\varphi^{(i)} (x) \neq \varphi^{(i)} (y)$. 
\item
If $x \in (\tau_{\Gamma_{\min}} (P_i) , P_i) \setminus 
 (\tau_{\Gamma_{\min}} (P_j) , P_j)$ 
and $y \in (\tau_{\Gamma_{\min}} (P_j) , P_j)
\setminus  (\tau_{\Gamma_{\min}} (P_i) , P_i)$ for some $i \neq j$,
then $\varphi^{(i)} (x) \neq  \varphi^{(i)} (y)$.
\item
If $x \in (\tau_{\Gamma_{\min}} (P_i) , P_i)$ and
$y \in \Gamma_{\min}$,
then $\varphi^{(i)} (x) \neq \varphi^{(i)} (y)$.
\end{enumerate}
\end{Lemma}

\Proof
Assertion (1) follows from Lemma~\ref{lemma:faithful:ends}(i). 

Assertion (3) follows from Lemma~\ref{lemma:faithful:ends}(i)(iii). 
Indeed, we have $\varphi^{(i)} (x) > 0$ by Lemma~\ref{lemma:faithful:ends}(i), and 
$\varphi^{(i)} (y) = 0$ by Lemma~\ref{lemma:faithful:ends}(iii).  

We show (2).
Suppose that $x \in (\tau_{\Gamma_{\min}} (P_i) , P_i) \setminus 
 (\tau_{\Gamma_{\min}} (P_j) , P_j)$ 
and $y \in (\tau_{\Gamma_{\min}} (P_j) , P_j)
\setminus  (\tau_{\Gamma_{\min}} (P_i) , P_i)$ for some $i \neq j$.
If $\tau_{\Gamma_{\min}} (P_j) \neq \tau_{\Gamma_{\min}} (P_i)$,
then Lemma~\ref{lemma:faithful:ends}(i)(ii)
gives $\varphi^{(i)} (x) > 0 = \varphi (y)$.
Suppose that $\tau_{\Gamma_{\min}} (P_j) = \tau_{\Gamma_{\min}} (P_i)$.
Then there exists $z \in 
[\tau_{\Gamma_{\min}} (P_i) , P_i) \cap 
 [\tau_{\Gamma_{\min}} (P_j) , P_j)$
such that 
$[\tau_{\Gamma_{\min}} (P_i) , z ] = 
[\tau_{\Gamma_{\min}} (P_i) , P_i) \cap 
 [\tau_{\Gamma_{\min}} (P_j) , P_j)$.
Note that $x \in 
(\tau_{\Gamma_{\min}} (P_i) , P_i) \setminus 
 (\tau_{\Gamma_{\min}} (P_j) , P_j)
=
[\tau_{\Gamma_{\min}} (P_i) , P_i) \setminus 
[\tau_{\Gamma_{\min}} (P_i) , z]
$.
Then
by Lemma~\ref{lemma:faithful:ends}(i)(ii), 
we see that
$\varphi^{(i)} (y) = \varphi^{(i)} (z) < \varphi^{(i)} (x)$.
\QED

\medskip
\textsl{Proof of the main Theorem~\ref{thm:main:unimodular:faithful}.}\quad
If $\Gamma$ is minimal, then 
we have already shown a faithful tropicalization of $\Gamma$ 
in Theorem~\ref{theorem:FT:canonical}
and Theorem~\ref{thm:FT:low:genus},
so that we may assume that it is not minimal. 
We write $\Gamma = S( \mathscr{X} ; \sigma_1 , \ldots , \sigma_r )$
for some strictly semistable pair, where we allow 
$r = 0$.
We add sections $\sigma_{r+1} , \ldots , \sigma_{m}$ with 
$m > r$ so that
any $(-1)$-curve in $\mathscr{X}_s$ meets 
at least one of $\sigma_1 (k) , \ldots , \sigma_m (k)$. 
Since
$\Gamma = S(\mathscr{X} ; \sigma_1 , \ldots , \sigma_{r})
\subset S(\mathscr{X} ; \sigma_1 , \ldots , \sigma_{m})$,
it suffices to prove the theorem for 
$S(\mathscr{X} ; \sigma_1 , \ldots , \sigma_{m})$.
Thus
replacing $\Gamma$
if necessary, we assume that $\Gamma = S(\mathscr{X} ; \sigma_1 , \ldots , \sigma_{m})$.

Take any minimal skeleton $\Gamma_{\min} \subset \Gamma$.
We set $P_i := \sigma_i (K)$ for $i=1 , \ldots , m$.
We claim that
\begin{equation}
\label{eqn:completeskeleton}
S(\mathscr{X} ; \sigma_1 , \ldots , \sigma_{m})
=
\Gamma_{\min} \cup \bigcup_{i = 1}^{m} [ \tau_{\Gamma_{\min}} (P_i),
 P_i). 
\end{equation}
Indeed, if $S(\mathscr{X} ; \sigma_1 , \ldots , \sigma_{m}) \supsetneq \Gamma_{\min} \cup \bigcup_{i = 1}^{m} [ \tau_{\Gamma_{\min}} (P_i),
 P_i)$, then there exists  
$E \in \Irr (\mathscr{X}_{s})$
such that $[E]$ has valence $1$ in $S( \mathscr{X} )$
and $[E] \in S (\mathscr{X} ; \sigma_1 , \ldots , \sigma_{m})
\setminus \Gamma_{\min} \cup \bigcup_{i = 1}^{m} [ \tau_{\Gamma_{\min}} (P_i),
 P_i)$. 
Then $E$ is a $(-1)$-curve in $\mathscr{X}_s$
with $\sigma_j (k) \notin E$ for any $j = 1 , \ldots , m$,
which contradicts our assumption of (a newly replaced) $\Gamma$. 
Thus \eqref{eqn:completeskeleton} holds. 

By Theorem~\ref{theorem:FT:canonical}
(when $g \geq 2$)
and by Theorem~\ref{thm:FT:low:genus}
(when $g =0,1$),
there exist $s_0 , \ldots , s_{N'} \in H^0 ( X,L)$
such that
the map
$\varphi' : X^\an \to \TT\RR^{N'}$
defined by $\varphi' := ( - \log |s_0| : \cdots : - \log |s_{N'}|)$
gives
a faithful tropicalization of $\Gamma_{\min}$.

Note that
$\mathcal{L} (\Gamma) = \{ P_1 , \ldots , P_m\}$.
For each $i = 1 , \ldots , m$, we use 
Lemma~\ref{lemma:faithful:ends} to obtain 
global sections $s_0^{(i)}$ and $s_1^{(i)}$ such that
$\varphi^{(i)} := - \log |s_1^{(i)} / s_0^{(i)}|$ enjoys properties
(i)--(iii) of Lemma~\ref{lemma:faithful:ends}. 
Set $s_{N' + 2i-1} := s^{(i)}_0$
and $s_{N' + 2i} := s^{(i)}_{1}$
for $i = 1 , \ldots , m$,
and set $N := N' + 2 m$.
Further, we
set $\varphi :=
( - \log |s_0| : \cdots : - \log |s_{N'}| :
- \log |s_{N'+1}| : \cdots : - \log |s_{N}|)$.

We show that $\varphi$ gives a faithful tropicalization 
of $\Gamma$. 
Since $\varphi'$ is gives a faithful tropicalization
of $\Gamma_{\min}$,
$\varphi$ also gives a faithful tropicalization
of $\Gamma_{\min}$.
By Lemma~\ref{lemma:faithful:ends}(i),
we see that $\varphi$ gives a unimodular tropicalization
over $[\tau_{\Gamma_{\min}} (P_i) , P_i)$
for each $i=1, \ldots , m$. 
Thus by \eqref{eqn:completeskeleton},
$\varphi$ is a unimodular tropicalization of $\Gamma$.

To conclude that $\varphi$ is a faithful tropicalization of 
$\Gamma$, we need to show that $\varphi$ is injective on $\Gamma$.  
Let $x , y \in \Gamma$ be distinct points. If $x, y \in \Gamma_{\min}$, 
then we have shown that $\varphi (x) \neq \varphi (y)$ by 
Theorem~\ref{theorem:FT:canonical} and Theorem~\ref{thm:FT:low:genus}. 
Without loss of generality, we may assume that $x \notin 
\Gamma_{\min}$. 
By (\ref{eqn:completeskeleton}),
we take an $i=1, \ldots , m$ with $x \in (\tau_{\Gamma_{\min}} (P_i) , P_i)$.
\begin{itemize}
\item
If $y \in \Gamma_{\min}$, then
Lemma~\ref{lemma:remark:toseparate:ends}(3) shows that
$\varphi (x) \neq \varphi (y)$.
\item
If $y \in (\tau_{\Gamma_{\min}} (P_i) , P_i)$,
then Lemma~\ref{lemma:remark:toseparate:ends}(1) 
shows that $\varphi (x) \neq \varphi (y)$.
\item
Suppose that $y \notin (\tau_{\Gamma_{\min}} (P_i) , P_i)$
and $y \notin \Gamma_{\min}$.
Then 
by \eqref{eqn:completeskeleton},
there exists $j \neq i$
with $y \in (\tau_{\Gamma_{\min}} (P_j) , P_j ) \setminus (\tau_{\Gamma_{\min}} (P_i) , P_i)$. 
If $x \not\in (\tau_{\Gamma_{\min}} (P_j) , P_j )$, then 
Lemma~\ref{lemma:remark:toseparate:ends}(2)
shows that $\varphi (x) \neq \varphi (y)$. 
If $x \in (\tau_{\Gamma_{\min}} (P_j) , P_j )$, then 
Lemma~\ref{lemma:remark:toseparate:ends}(1) applied to $\varphi^{(j)}$ 
shows that $\varphi (x) \neq \varphi (y)$. 
\end{itemize}
Thus $\varphi$ is injective on $\Gamma$, so that
$\varphi$ is a faithful tropicalization of $\Gamma$. 
This completes the proof of
Theorem~\ref{thm:main:unimodular:faithful}. 
\QED

\setcounter{equation}{0}
\section{Complementary results}
\label{sec:complementary}

\subsection{Theorem~\ref{thm:main:unimodular:faithful} is optimal for curves in low genera}
In this subsection, we prove Theorem~\ref{thm:best:possible}, which shows that the bounds
in Theorem~\ref{thm:main:unimodular:faithful}
for $g=0,1,2$ are optimal.

We begin with the case of genus $0$.
In this case, we have the following obvious result.

\begin{Proposition}
\label{prop:complement:0} 
Let $\PP^1$ be the projective line over $K$. 
Let $L$ be a line bundle over $\PP^1$
of degree $0$, i.e., $L \cong \OO_{\PP^1}$.
Then for any non-zero global sections $s_0 , \ldots , s_N$ of $L$,
the image of the tropicalization map
$\varphi : \PP^{1,\an} \to \TT\PP^N$
defined by $\varphi = ( - \log |s_0| :
\cdots : - \log |s_N| )$
is a singleton.
\end{Proposition}

In particular, if $\Gamma$ is a skeleton of $\PP^{1, \an}$
that is not a singleton,
then 
$\OO_{\PP^1}$ does not admit
a faithful tropicalization of $\Gamma$.

Next we consider the case of genus $1$.
Recall that a smooth projective curve $X$ over $K$ of genus $g \geq 1$ is 
called
a \emph{Mumford--Tate curve} if its minimal skeleton (as a $\Lambda$-metric graph) has genus $g$.

\begin{Proposition}
\label{prop:complement:1} 
Let $X$ be a Mumford--Tate smooth projective curve 
over $K$ of genus $1$.
Let $L$ be a line bundle over $X$ of degree $2$.
Then $L$ does not admit a faithful tropicalization of
the minimal skeleton $\Gamma_{\min}$.
\end{Proposition}

\Proof
Let $s_0, \ldots, s_N \in H^0(L)$ be any non-zero global sections.
It is classically known that the associated morphism
$X \to \PP^N$
decomposes into the composite of
a double covering
$X \to \PP^{1}$
and
a morphism $\PP^1 \to \PP^N$.
This means that the tropicalization $\varphi :
X^{\an} \to \TT\PP^N$
factors though $\PP^{1, \an}$. 
Since 
the image $\varphi ( \Gamma_{\min} )$ is a skeleton of 
$\PP^{1,\an}$,
we see that
the first Betti number of
$\varphi ( \Gamma_{\min} )$
is $0$.
On the other hand, we have $g ( \Gamma_{\min}) = 1$.
This
 shows that
$\Gamma_{\min}$ is not homeomorphic to 
$\varphi ( \Gamma_{\min} )$.
Thus $\varphi$ is not a faithful tropicalization of $\Gamma_{\min}$.
\QED

Finally we show that 
a bicanonical system does not 
necessarily give a faithful tropicalization of a skeleton
when $g=2$. 
We remark
that the bicanonical divisor has degree $4 = t (2) - 1$
with the notation in Theorem~\ref{thm:main:unimodular:faithful}.

\begin{Proposition}
\label{prop:complement:2} 
Let $X$ be a Mumford--Tate smooth projective curve over $K$ 
of genus $2$. 
Then 
the bicanonical line bundle $\omega_X^{\otimes 2}$
does not admit a faithful tropicalization of 
the minimal skeleton $\Gamma_{\min}$.
\end{Proposition}

\Proof
For any global sections $s_0, \ldots, s_N \in H^0(\omega_X^{\otimes 2})$,
the associated morphism 
$
X \to \PP^N
$
decomposes into the composite of
the natural double covering   
$X \to \PP^{1}$ 
given by the quotient of the hyperelliptic involution and 
a morphism $\PP^1 \to \PP^N$. 
This means that the tropicalization $\varphi :
X^{\an} \to \TT\PP^N$ factors through $\PP^{1, \an}$. 
Then, by the same argument as in the proof of 
Proposition~\ref{prop:complement:1}, we see that 
$\varphi$ is not a faithful tropicalization of $\Gamma_{\min}$.
\QED

\subsection{A very ample line bundle that does not admits a faithful tropicalization}
\label{subsec:very:ample:no:ft}
In this subsection, we prove Proposition~\ref{prop:countereg:very:ample}, which shows that  
a {\em very ample} line bundle does not necessarily admit 
a faithful tropicalization. Here we recall the assertion. 

\begin{PropositionNoNum}[$=$ Proposition~\ref{prop:countereg:very:ample}]
Let $d \geq 4$ be any integer. 
Then there exists a smooth projective plane curve $X$ of $\PP^2$ of degree $d$ over $K$ 
such that no skeleton $\Gamma$ of $X^{\an}$ has a faithful tropicalization 
associated to $|\OO_X(1)|$, where  $\OO_X(1)$ is the very ample line bundle over $X$ given by 
the restriction of $\OO_{\PP^2}(1)$. 
\end{PropositionNoNum}

\Proof
We take any positive integers $d_1\geq 2$ and 
$d_2 \geq 2$ with $d_1 + d_2 = d$. 
For $i = 1, 2$, we take a homogeneous polynomial $f_{i} \in k[T_0, T_1, T_2]$ of degree 
$d_i$ such that the curve $C_i$ defined by $f_i$ in $\PP^2_k$ is smooth and such that 
$C_1 \cup C_2$ is a strictly semistable curve over $k$. 

First, we construct a plane curve $X$.
We take any $\widetilde{F}_i \in R[T_0, T_1, T_2]$ which maps to $f_i$ under 
the natural surjection $R [T_0, T_1, T_2] \to k [T_0, T_1, T_2]$
arising from $R \to R/\mathfrak{m} = k$. 
We take a homogeneous polynomial $\widetilde{F} \in R[T_0, T_1, T_2]$ of degree $d$ and 
a nonzero element $\varpi \in R$ with $|\varpi|_K < 1$, and we set 
\[
  \widetilde{G} := \widetilde{F}_1 \widetilde{F}_2 + \varpi \widetilde{F} \in R[T_0, T_1, T_2]. 
\] 
When we regard $\widetilde{G}$ as an element of $K[T_0, T_1, T_2]$, we denote 
it by $G$.
We set 
\[
  X :=  \Proj(K[T_0, T_1, T_2]/(G)). 
\]
By taking $\widetilde{F}$ generally, we may and do assume that $X$ is smooth.

Let $\Gamma$ be any skeleton of $X^{\an}$.
To ease notation,
we set $L := \rest{\OO_{\PP^2}(1)}{X}$. 
We take any $N \geq 1$ and any nonzero global sections $s_0, \ldots, s_N \in H^0(L)$, and we consider 
the associated map 
\[
  \varphi: X^{\an} \to \TT\PP^N, \quad 
  p = (p, |\cdot|) \mapsto \left(-\log|s_0(p)|: \cdots: -\log|s_N(p)|\right).  
\]
Then our goal is 
to show that $\varphi$ is not a faithful tropicalization of $\Gamma$.

We set
\[
\Xscr := \Proj(R[T_0, T_1, T_2]/(\widetilde{G})) 
,
\]
which is naturally a model of $X$.
By the definition of $\widetilde{G}$,
the image of $\widetilde{G}$ by the natural homomorphism
$R [ T_0 , T_1 , T_2] \to k [ T_0 , T_1 , T_2]$
equals $f_1 f_2$.
This means that
the special fiber 
$\mathscr{X}_s$ equals the strictly semistable
semistable curve $C_1 \cup C_2$,
and thus $\mathscr{X}$ is a strictly semistable model of $X$.
For each $i = 1, 2$, since $\#(C_1 \cap C_2) = d_1 d_2 > 2$, 
$C_i$ is not a $(-1)$-curve nor a $(-2)$-curve in $\mathscr{X}_s$.
Thus $\mathscr{X}$ is a Deligne--Mumford strictly semistable model
of $X$.
We remark that $X$ has genus $\frac{(d-1)(d-2)}{2} \geq 3$.
Let $\Gamma_{\min}$ be the minimal skeleton of $X^{\an}$.
Then
for each $i=1,2$,
the Shilov point 
$[C_i] \in X^{\an}$ associated to $C_i$ belongs to $\Gamma_{\min}$.
Since $\Gamma_{\min} \subset \Gamma$,
we note that $[C_1] , [C_2] \in \Gamma$. 

Now, it suffices
to show that $\varphi([C_1]) = \varphi([C_2])$.
We set $\Lscr:= \rest{\OO_{\PP_R^2}(1)}{\Xscr}$. Then $(\Xscr, \Lscr)$ is a model of $(X, L)$. 
By adding the global section $\rest{T_0}{X} \in H^0(L)$, we may and do assume that $s_0 =  \rest{T_0}{X}$. 
We set $\widetilde{s}_0 := \rest{T_0}{\Xscr} \in H^0(\Lscr)$. 
It follows from the exact sequence 
\[
 0 \to \OO_{\PP^2}(-d+1) \to \OO_{\PP^2}(1) \to L  \to 0
\]
and $H^1(\OO_{\PP^2}(-d+1)) = 0$ that the restriction map $H^0(\OO_{\PP^2}(1)) \to H^0(L)$ is surjective. 
Thus for any $1 \leq \ell \leq N$, there exists $(a_{\ell}, b_{\ell}, c_{\ell}) \neq (0, 0, 0) \in K^3$ such that $s_{\ell} = \rest{\left(a_{\ell} T_0 + b_{\ell} T_1 + c_{\ell} T_2\right)}{X}$. We take $\alpha_{\ell} \in K$  such that $a_{\ell}^\prime := a_{\ell}/\alpha_{\ell}, b_{\ell}^\prime := b_{\ell}/\alpha_{\ell}, c_{\ell}^\prime := c_{\ell}/\alpha_{\ell}$ belong to $R$ and such that one of $a_{\ell}^\prime, b_{\ell}^\prime, c_{\ell}^\prime$ is a unit of $R$. 
We set $\widetilde{s}_{\ell} = \rest{\left(a_{\ell}^\prime T_0 + b_{\ell}^\prime T_1 + c_{\ell}^\prime T_2\right)}{\Xscr} \in H^0(\Lscr)$ 
and $g_{\ell}:= \widetilde{s}_{\ell}/\widetilde{s}_0$. 
Then $g_{\ell}$ is a nonzero rational function on $\Xscr$. 

Recall that $C_1$ is a curve in $\PP^2_k$ of degree $d_1 \geq 2$,
and in particular it is not linear.
Since
$\rest{(a_{\ell}^\prime T_0 + b_{\ell}^\prime T_1 + c_{\ell}^\prime T_2)}{\PP^2_k} \neq 0$, 
it follows that $\rest{(a_{\ell}^\prime T_0 + b_{\ell}^\prime T_1 + c_{\ell}^\prime T_2)}{C_1} \not\equiv 0$. This means that $\widetilde{s}_{\ell}$ does not vanish at the generic point $\xi_1$ of $C_1$. By the same argument, $\widetilde{s}_{0}$ does not vanish at $\xi_1$. Thus 
$g_{\ell}$ is a unit regular function of $\OO_{\Xscr, \xi_1}$, and 
we have $-\log|g_{\ell}|([C_1]) = 0$. It follows that 
\[
\varphi([C_1]) = (0: -\log |\alpha_{1}|_K: \cdots: -\log |\alpha_{N}|_K). 
\]
By the same argument, we obtain $\varphi([C_2]) = (0: -\log |\alpha_{1}|_K: \cdots: -\log |\alpha_{N}|_K)$.
This proves that
$\varphi([C_1]) = \varphi([C_2])$,
and thus
this completes the proof. 
\QED

\subsection{Comparison with \cite{KY}}
\label{subsection:comparison}
In this subsection,
we clarify the difference between 
the work in \cite{KY}
and the current work
with an emphasis of an advantage of
Theorem~\ref{thm:main:unimodular:faithful}.

The paper \cite{KY} discusses faithful tropicalization of
projective varieties under some conditions of definability over discrete valuation 
rings (of equicharacteristic zero, cf. Remark~\ref{remark:equicharacteristic}). 
It has an obvious advantage of 
treating smooth projective varieties of arbitrary dimension. 
However, we will find that the current paper has an advantage 
when we restrict our attentions to curves, besides the definability
assumption over discrete valuation rings.

First, let us recall what the result
in \cite{KY} says for curves.
The following assertion has essentially been proved in \cite{KY}.

\begin{Theorem}[cf.~{\cite[Theorem~1.1]{KY}}]
\label{thm:KY}
Assume that $R$ has equicharacteristic zero \textup{(}cf. Remark~\textup{\ref{remark:equicharacteristic})}.  
Let $X$ be a smooth projective curve over $K$, 
and let $\Gamma$ be a skeleton of $X^{\an}$.
Let $L$ be a line bundle over $X$.
Assume that there exists a discrete valuation ring $R'$
dominated by $R$
and
a strictly semistable curve $\mathscr{X}' \to \Spec (R')$
such that $\Gamma$ is associated to the model
$\mathscr{X}' \otimes_{R'} R$.
Suppose that there exist
a relatively ample line bundle $\mathscr{N}$ over $\mathscr{X}'$
and an integer $m \geq 2$
such that $\rest{\mathscr{N}^{\otimes m} \otimes \omega_{\mathscr{X}^\prime/R^\prime}}{X}
\cong L$. Then $L$ admits a faithful tropicalization of $\Gamma$.
\end{Theorem}

\Proof
We briefly explain how \cite[Theorem~1.1]{KY} leads to Theorem~\ref{thm:KY}.
Let $K'$ be the fraction field of $R'$. We put an absolute value $|\cdot|_{K'}$ on $K'$ 
by restricting the absolute value $|\cdot|_K$ on $K$.
Set $X' := \mathscr{X}' \otimes_{R'} K'$.
Then the Berkovich analytic space $(X')^{\an}$ over $K'$
and the skeleton $S(\mathscr{X}')$
associated to a strictly semistable model $\mathscr{X}' \to \Spec(R^\prime)$ 
are defined similarly. 
Let $\alpha : X \to X'$ be the natural morphism, and 
the canonical map $\alpha^{\an} : X^{\an} \to (X')^{\an}$ 
induces an isometry 
$\Gamma = S ( \mathscr{X}' \otimes_{R'} R ) \cong S ( \mathscr{X}' )$. 

Set $L' := \rest{\mathscr{N}^{\otimes m} \otimes \omega_{\mathscr{X}'/R'}}{X'}$.
By \cite[Theorem~1.1]{KY},
there exist global sections
$s_0 , \ldots , s_N$ 
of $L'$ such that the associated morphism 
$\psi : (X')^{\an} \to \TT\PP^N$
gives a faithful tropicalization of $S ( \mathscr{X}')$. 
Since $\alpha^{\ast} \left( 
L'
\right) = L$,
we 
regard $s_0 , \ldots , s_N$ as global sections
of $L$,
and we also have the associated tropicalization map
$\varphi : X^{\an} \to \TT\PP^N$.
Further, we have
$\varphi = \psi \circ \alpha^{\an}$.
Since $\alpha^{\an}$ gives the isometry between the skeleta,
this proves that
$\varphi$ also gives a faithful tropicalization of $\Gamma$.
\QED

\begin{Remark}
\label{remark:equicharacteristic}
In fact, 
the assumption in Theorem~\ref{thm:KY} that $R$ has equicharacteristic zero 
is unnecessary, 
as we now explain.
Suppose that $R$ has positive residue characteristic, that is, 
$\mathrm{char} (k) > 0$.
Note that
for any smooth projective curve $Z$ over any field $k$
and for any ample line bundle $N$ over $Z$,
$N^{\otimes m} \otimes \omega_Z$ is basepoint free for any $m \geq 2$.
This means the assumption ``$m \geq \phi (d)$" in 
\cite[Theorem~1.1]{KY} is fulfilled.
Further, for curves, 
the assertions of vanishing of cohomologies
and basepoint-freeness as in \cite[\S3]{KY},
which are the technical keys to prove \cite[Theorem~1.1]{KY},
hold true even in positive characteristic.
Then
one can construct global sections of $\rest{\mathscr{N}^{\otimes m} \otimes \omega_{\mathscr{X}'/R'}}{X}$ that give a faithful tropicalization of $\Gamma$
by essentially the same arguments in \cite{KY}.
Thus the same conclusion of Theorem~\ref{thm:KY} holds true 
without the assumption that $R$ has equicharacteristic zero. 
\end{Remark}

To clarify the crucial difference between the results of the two papers,
we revisit Question~\ref{question:first}.
Recall that Theorem~\ref{thm:main:unimodular:faithful} gives us a
concrete answer to Question~\ref{question:first} that
$d (g) = t (g)$, where $t (g)$ is given in Theorem~\ref{thm:main:unimodular:faithful}. 
Now, we want to ask whether or not
Theorem~\ref{thm:KY} makes some contribution to 
Question~\ref{question:first}.
The following proposition is crucial to this issue.

For a smooth projective curve $X$ over $K$ and a complete discrete valuation ring $R^\prime$ 
dominated by $R$, we say that $X$ is \emph{definable over $R'$}
if there exists a proper curve $\mathscr{X} \to \Spec (R')$
over $R'$ such that $\mathscr{X} \otimes_{R'} R$ is a model of $X$.

\begin{Proposition}
\label{prop:cannotbecoverd}
Let $g$ and $d$  be any positive integers with $g \geq 2$. 
Then there exists a smooth projective curve
$X$ over $K$ of genus $g$ 
that is definable over some complete discrete valuation ring
dominated by $R$
and has the following property\textup{:} 
for any complete discrete valuation ring $R'$
dominated by $R$,
for any Deligne--Mumford strictly semistable model $\mathscr{X}^\prime \to \Spec (R')$
of $X$ such that $\mathscr{X}^\prime$ is regular, 
and for any line bundle $\mathscr{N}$ over $\mathscr{X}^\prime$,
if $\mathscr{N}$ is relatively ample,
then $\deg \left( \rest{\mathscr{N}}{X} \right) \geq d$. 
\end{Proposition}

\Proof
Note that
for each prime number $p$,
\cite[Theorem~29.1]{matsumura} gives a complete discrete valuation ring
$W_p$ with residue field $k$ such that $p$ is a uniformizer in $W_p$.
Set
\begin{align*}
V :=
\begin{cases}
k [[ t ]] &
\text{if $R$ has equicharacteristic,}
\\
W_p
&
\text{if $R$ has mixed characteristic and $\mathrm{char} (k) = p$.}
\end{cases}
\end{align*}
Then $V$ is a discrete valuation ring
with residue field $k$.
Let $\varpi$ be a uniformizer of $V$.
Then there exists 
an injective homomorphism 
$V \to R$ of local rings 
that induces the identity between the residue fields $k$;
see \cite[Theorem~29.2]{matsumura} when $R$ has mixed characteristic.

Since $g \geq 2$,
there exists a stable curve
$D$ over $V / ( \varpi) = k$
of genus $g$ that has
at least two nodes, which we denote by $q_1$ and $q_2$.
Using the deformation theory of stable curves,
we obtain a generically smooth,  strictly semistable and stable curve
$\mathscr{X}^0 \to \Spec (V)$ 
of genus $g$ with special fiber
$D$
such that the completion of the local ring of $\mathscr{X}^0$
at $q_1$ is isomorphic to $V [[x,y]] / (xy - \varpi)$
and that at $q_2$ is isomorphic to
$V [[x,y]]/ (xy - \varpi^{d})$.
We set $X := \mathscr{X}^0 \otimes_{V} K$. 
Then $X$ is a smooth projective curve over $K$ of genus $g$.
We remark that $\mathscr{X}^{\st} := \mathscr{X}^0 \otimes_{V} R$
is the stable model of $X$ over $R$. 
Let $\Gamma_{\min}$ denote the minimal skeleton of $X^{\an}$.
For $i=1,2$, let $e_i \in E(\Xscr^{\st})$ be the edge of $\Gamma_{\min}$ corresponding to $q_i$,
and let $\lambda_i$ denote the length of $e_i$. Then $\lambda_2 = d  \lambda_1$.

Let us prove that $X$ has the required properties. 
We take any complete discrete valuation ring $R'$ 
with a Deligne--Mumford strictly semistable model 
$\mathscr{X}^\prime \to \Spec (R')$
of $X$ such that $\mathscr{X}^\prime$ is regular. 
Let 
$\mathscr{X}^{\prime} \to \mathscr{X}^{\prime \st}$
be the contraction of $(-2)$-curves,
where $\mathscr{X}^{\prime \st}$ is the stable model of $\mathscr{X}^{\prime}$.
Since $\mathscr{X}^{\prime \st} \otimes_{R'} R$ is 
a stable curve with geometric generic fiber
$X$, we have $\mathscr{X}^{\prime \st} \otimes_{R'} R = \mathscr{X}^{\st}$
by the uniqueness of the stable model.
In particular,
$\mathscr{X}^{\prime \st} \otimes_{R'} k = D$,
where $k$ is regarded as an $R'$-algebra by 
the composite $R' \to R \to R/ \mathfrak{m} = k$.
The morphism $\mathscr{X}^{\prime} \to \mathscr{X}^{\prime \st}$
induces a morphism $\alpha :
\mathscr{X}^{\prime} \otimes_{R'} k \to D$,
which is the contraction of $(-2)$-curves in 
the Deligne--Mumford semistable curve
$\mathscr{X}^{\prime} \otimes_{R'} k$ over $k$.
For $i= 1,2$, let $n_i$ denote the number of nodes that
$\alpha$ maps to $q_i$.
Let $\varpi'$ denote the uniformizer of $R'$.
Since $\mathscr{X}^\prime$ is 
strictly semistable and
regular, we note that $n_i ( - \log | \varpi' |_K) = 
\lambda_i$ for $i = 1,2$.
Since $\lambda_2 = d \lambda_1$,
it follows that $n_2 = d  n_1 \geq d$.

We claim
that $\# \Irr ( \mathscr{X}^{\prime} \otimes_{R'} k ) \geq d$.
Since the assertion is trivial if $d = 1$,
we may and do assume that $d \geq 2$.
Then
we have $n_2 \geq 2$,
and hence
there exists a maximal $(-2)$-chain $E$ in $\mathscr{X}^{\prime} \otimes_{R'} k$
such that $\Delta_{E} = e_2$.
Note that $\# \Irr (E) = n_2 -1 \geq d -1$.
Since $\mathscr{X}^{\prime} \otimes_{R'} k$ has an irreducible components that is
not a $(-2)$-curve,
it follows that
$
\# \Irr ( \mathscr{X}^{\prime} \otimes_{R'} k )
\geq 
\# \Irr ( E )  + 1
\geq
d
$, as desired.

Now let $\mathscr{N}$ be any line bundle over $\mathscr{X}^\prime$ that is 
relatively ample. 
Then, since $\mathscr{X}^{\prime} \otimes_{R'} k$ has at least $d$ irreducible components, 
we have 
$\deg
\left(
\rest{\mathscr{N}}{X}
\right)
=
\deg \left( 
\rest{\mathscr{N}}{\mathscr{X}^{\prime} \otimes_{R'} k}
\right) \geq d$. 
Thus the proposition holds.
\QED

Suppose that a smooth projective curve $X$ over $K$ is definable over 
a complete discrete valuation ring $R^\prime$ dominated by $R$, and we take 
a proper curve $\mathscr{X} \to \Spec (R')$ over $R'$ 
such that $\mathscr{X} \otimes_{R'} R$ is a model of $X$. 
Let $L$ be a line bundle over $X$, and let $\Gamma$ be a skeleton of $X^{\an}$. 
We say that $L$ is \emph{definable over $R'$}
if there exists
a line bundle $\mathscr{L}$ over $\mathscr{X}$
such that $\mathscr{L} \otimes_{R'} R$ is
a model of $L$.
Further, we say that $\Gamma$ is \emph{definable over
$R'$}
if there exists a strictly semistable curve
$\mathscr{X} \to \Spec (R')$
over $R'$ 
such that $\mathscr{X} \otimes_{R'} R$ is a model of $X$
and such that $\Gamma$ is the skeleton associated to this model.

We fix any integer $g \geq 2$.
Take any integer $s \geq 2g - 2$.
We are going to show that there exist a smooth projective curve 
$X$
over
$K$ of genus $g$ and a line bundle $L$ over $X$
with $\deg (L) = s$
such that $X$, $L$, and the minimal skeleton 
of $X$ are definable over some complete discrete valuation ring
dominated by $R$,
but such that 
there do {\em not} exist
a complete discrete valuation ring
$R'$ dominated by $R$,
a Deligne--Mumford strictly semistable curve 
$\mathscr{X}' \to \Spec (R')$ with $\mathscr{X}' \otimes_{R'} R = X$,
a relatively ample line bundle over $\mathscr{N}$ over $\mathscr{X}'$,
and an integer $m \geq 2$
such that
$\rest{\mathscr{N}^{\otimes m} \otimes \omega_{\mathscr{X}'/R'}}{X} 
\cong L$.
Indeed, suppose that this is proved.
Then,
the assumptions of Theorem~\ref{thm:KY} are not fulfilled.
Since $s \geq 2g-2$ is taken arbitrarily, this suggests that the existence of
a universal bound
$d (g)$ as in Question~\ref{question:first} cannot be deduced from Theorem~\ref{thm:KY},
so that 
Theorem~\ref{thm:KY} cannot contribute to Question~\ref{question:first}. 

The construction is done by
Proposition~\ref{prop:cannotbecoverd}.
Indeed, 
applying this proposition for
$d := \lceil
(s- 2g + 2)/2
\rceil + 1$,
where
for a real number $x$, 
$\lceil x \rceil$ denotes the 
smallest integer with $\lceil x \rceil \geq x$,
we get
a smooth projective curve $X$ over $K$ of genus $g$
that is definable over some complete discrete
valuation ring dominated by $R$
and enjoys
the properties in the proposition.
We note that
the minimal skeleton of $X$ is also definable over the complete
valuation ring.
Since $X$ is definable over 
some complete discrete
valuation ring dominated by $R$,
there exists
a line bundle $L$ over $X$
with $\deg (L) = s$
that is definable over some complete discrete valuation ring
dominated by $R$.

However, 
for any complete discrete valuation ring $R'$ dominated by $R$,
for
any Deligne--Mumford strictly semistable 
curve $\mathscr{X}' \to \Spec (R')$ with $\mathscr{X}' 
\otimes_{R'} R = X$,
for any
relatively ample line bundle $\mathscr{N}$
over $\mathscr{X}'$,
and for
any integer $m \geq 2$,
we have
$\rest{\mathscr{N}^{\otimes m} \otimes \omega_{\mathscr{X}'/R'}}{X}
\ncong L$.
Indeed,
for relatively ample $\mathscr{N}$, since $X$ enjoys the properties of 
Proposition~\ref{prop:cannotbecoverd},
we have
\[
\deg ( \rest{\mathscr{N}}{X}) 
\geq d > (s- 2g + 2)/2 
,
\]
and hence
\[
\deg \left( 
\rest{\mathscr{N}^{\otimes m} \otimes \omega_{\mathscr{X}'/R'}}{X}
\right)
>
m(s- 2g + 2)/2 + 2g - 2
\geq s
= \deg (L)
.
\]
Thus we get the conclusion.

\setcounter{equation}{0}
\section{Limit of tropicalizations by polynomials of a bounded degree}
\label{sec:limit:trop}

Let $\mathbb{A}^N $ be the $N$-dimensional affine space with 
affine coordinate functions $z_1 , \ldots , z_N$, and let  
$Y \subset \mathbb{A}^N$ be a closed subvariety.  
In \cite{payne}, Payne considers 
the inverse system consisting of all affine embeddings of $Y$
given by polynomials in $z_1, \ldots, z_N$ and shows that the inverse limit 
is homeomorphic to $Y^{\an}$. 
In his construction, there are no restrictions for 
affine embeddings of $Y$. 

In this section, we consider a question whether 
$Y^{\an}$ is homeomorphic to  the inverse limit of 
the affine embeddings 
given by polynomials of degree at most some effective bound. 
When $Y$ is a suitable affine curve, 
we will answer this question in the affirmative, 
as an application of 
Theorem~\ref{thm:main:unimodular:faithful}. 

\subsection{Statement of the result}
\label{subsec:limit:trop:statement}

Let us briefly recall Payne's construction
of limit tropicalizations. 
For any $n \geq 1$, let $\mathbb{A}^n := \Spec 
( K [x_1 , \ldots , x_n] )$ denote the $n$-dimensional affine
space with affine coordinate functions $x_1 , \ldots , x_n$, on which 
the algebraic torus $\GG_m^n :=
\Spec \left(
K [x_1^{\pm 1} , \ldots , x_n^{\pm 1}]
\right)$ acts naturally.
Using the standard coordinates of $\mathbb{A}^n$, 
we have a  map 
\[
\mathbb{A}^{n} (K) \to \TT^n, \quad 
p \mapsto (- \log |x_1 (p)|_K , \ldots ,  - \log |x_n (p)|_K ), 
\]
which extends to a map
\[
\Trop\colon 
\mathbb{A}^{n,\an} \to \TT^n, \quad 
p \mapsto (- \log |x_1 (p)| , \ldots ,  - \log |x_n (p)| )
,
\]
called the standard tropicalization map.
On $\TT^n$ acts $\RR^n = \Trop(\GG_m^{n, \an})$ naturally, which is compatible with 
the map $\Trop$. 

By an {\em equivariant morphism} $\varphi : \mathbb{A}^{m} 
\to \mathbb{A}^{n}$, we mean a morphism that is equivariant by
the torus action with respect to some group homomorphism $\GG_m^m \to \GG_m^n$.
Then $\varphi$ induces a well-defined continuous
map $\Trop(\varphi)\colon \TT^m \to \TT^n$ that is compatible with the tropicalization maps, 
i.e., $\Trop \circ \varphi = \Trop(\varphi) \circ \Trop$. We call 
$\Trop(\varphi)$ the tropicalization of the map $\varphi$. 

A closed embedding $\iota: Y \hookrightarrow \mathbb{A}^{n}$ is called an
\emph{affine embedding}.
Let $I$ denote the set of all affine embeddings of $Y$. 
Then $I$ has a structure of directed set as follows.
Let $\iota_1 : Y \hookrightarrow \mathbb{A}^{n_1}$ and 
$\iota_2 : Y \hookrightarrow \mathbb{A}^{n_2}$ belong to $I$. 
Then we declare that $\iota_1 \leq \iota_2$ if 
there exists an equivariant morphism $\varphi : \mathbb{A}^{n_2} \to \mathbb{A}^{n_1}$
with $\varphi \circ \iota_2 = \iota_1$. 
Furthermore,  
for $\iota_1 : Y \hookrightarrow \mathbb{A}^{n_1}$ and 
$\iota_2 : Y \hookrightarrow \mathbb{A}^{n_2}$, we set $\iota_3 = (\iota_1, \iota_2): Y \hookrightarrow \mathbb{A}^{n_1+n_2}$. Since the natural projections ${\rm pr}_1\colon \mathbb{A}^{n_1+n_2} \to \mathbb{A}^{n_1}$ and 
${\rm pr}_2\colon \mathbb{A}^{n_1+n_2} \to \mathbb{A}^{n_2}$ satisfy 
${\rm pr}_1 \circ \iota_3 = \iota_1$ and ${\rm pr}_2 \circ \iota_3 = \iota_2$, we have 
$\iota_1 \leq \iota_3$ and $\iota_2 \leq \iota_3$.

For an affine embedding $\iota: Y \hookrightarrow \mathbb{A}^{n}$, 
we define the 
\emph{tropicalization $\Trop(Y,\iota)$ of $Y$ with respect to $\iota$}
to be 
the closure $\overline{
\{\Trop(\iota(y)) \mid y \in Y (K) \}}
$ of $\iota (Y (K) )$ in $\TT^n$.
By \cite[Proposition~2.2]{payne},
$
\Trop(Y,\iota) = 
\Trop (\iota^{\an} ( Y^{\an} ))$,
where $\iota^{\an} : Y^{\an} \hookrightarrow \mathbb{A}^{n, \an}$
is the map between the Berkovich spaces associated to $\iota$.
Let 
\[
\pi_\iota : Y^{\an} \to \Trop(Y,\iota)
\] denote the restriction
of $\Trop \circ \iota^{\an}$,
which is surjective and continuous. 
Let $\iota_1 : Y \hookrightarrow \mathbb{A}^{n_1}$ 
and $\iota_2 : Y \hookrightarrow \mathbb{A}^{n_2} $ be elements
of $I$.
Suppose that $\iota_1 \geq \iota_2$.
Then by the definition, there exists an equivariant morphism
$\varphi : \mathbb{A}^{n_2} \to \mathbb{A}^{n_1}$ such that
$\iota_1 = \varphi \circ \iota_2$.
The map $\Trop(\varphi): \TT^{n_2} \to \TT^{n_1}$ induces the map 
$\Trop(Y,\iota_2) \to \Trop(Y,\iota_1)$ by restriction.
Remark that
$\Trop(Y,\iota_2) \to \Trop(Y,\iota_1)$ is 
determined only by $\iota_1$ and $\iota_2$, and is independent of the choice of $\varphi$.

To sum up, $\left( \Trop(Y, \iota ) \right)_{\iota \in I}$ 
constitutes an inverse system 
in the category of topological spaces,
and we have the inverse limit $\varprojlim_{\iota \in I} \Trop(Y, \iota)$.
Let 
\begin{equation}
\label{eqn:limit:tropical:map}
\mathrm{LTrop} : Y^{\an} \to \varprojlim_{\iota \in I} \Trop(Y, \iota)
\end{equation}
to be the continuous map induced from $( \pi_{\iota} )_{\iota \in I}$ by the
universal property of the inverse limit.
Then Payne proves the following theorem.

\begin{Theorem}[{\cite[Theorem~1.1]{payne}}]
\label{thm:limit:trop}
Analytification is the limit of all tropicalizations, 
i.e., 
the map $\mathrm{LTrop}$
in \eqref{eqn:limit:tropical:map} is a homeomorphism. 
\end{Theorem}

Since
$\Trop(Y,\iota)$ equals the closure $\overline{
\{\Trop(\iota(y)) \mid y \in Y (K) \}}
$ of $\iota (Y (K) )$ in $\TT^n$,
Theorem~\ref{thm:limit:trop} suggests that 
one can describe
$Y^{\an}$ in terms of tropical
geometry even without knowing 
the definition of the Berkovich space $Y^{\an}$ associated to $Y$.

Our aim here is to prove that 
$Y^{\an}$ is still homeomorphic to the inverse limit 
of tropicalizations of an {\em effectively bounded degree}.  
Recall that $Y$ is a closed subvariety of $\mathbb{A}^N = \Spec(K[z_1, \ldots, z_N])$. 
We say that an affine embedding $\iota : Y \hookrightarrow \mathbb{A}^n$ has 
\emph{degree at most $D$} if 
there exist $f_1 , \ldots , f_n \in K [z_1 , \ldots , z_N]$ such that $\deg (f_j) \leq D$ for $j = 1, \ldots, n$ 
and $\iota = \left(
\rest{f_1}{Y}, \ldots, \rest{f_n}{Y} \right)$. 
Let $I_{\leq D}$ denote the set of all affine embeddings of $Y$ with degree at most $D$. 
Then since $I_{\leq D} \subset I$, it
is an ordered set, and it is actually a directed set.
Indeed, for $\iota_1 : Y \hookrightarrow \mathbb{A}^{n_1}$ and 
$\iota_2 : Y \hookrightarrow \mathbb{A}^{n_2}$, 
the map
$(\iota_1, \iota_2):  Y \hookrightarrow \mathbb{A}^{n_1 + n_2}$ belongs $I_{\leq D}$ and satisfies 
$\iota_1 \leq (\iota_1, \iota_2)$ and $\iota_2 \leq (\iota_1, \iota_2)$. 
Thus $\left( \Trop(Y, \iota ) \right)_{\iota \in I_{\leq D}}$ 
again constitutes an inverse system in the category of topological spaces,
and we have the inverse limit 
$\varprojlim_{\iota \in I_{\leq D}} \Trop(Y, \iota)$.
We denote by
\begin{equation}
\label{eqn:limit:tropical:map:d}
\mathrm{LTrop}_{\leq D} : Y^{\an} \to \varprojlim_{\iota \in I_{\leq D}} \Trop(Y, \iota)
\end{equation}
the continuous map induced from $( \pi_{\iota} )_{\iota \in I_{\leq D}}$
by the universality of the inverse limit.

We regard $\mathbb{A}^N$ as an open subset of $\PP^N$.
Let $X$ be the closure of $Y$ in $\PP^N$.
Then the degree of $X$ in $\PP^N$ depend only on $Y$
and not on the choice of the open embedding
$\mathbb{A}^N \hookrightarrow \mathbb{P}^N$.
We call $d$ the \emph{degree of $Y$}. 
Further, we say that \emph{$Y$ has smooth compactification in $\PP^N$}
if $X$ is smooth.
This notion is also well defined for $Y$.

The following is 
the main result of this section.

\begin{TheoremNoNum}[$=$ Theorem~\ref{thm:limit:trop:d:intro}]
Let $Y \subset \mathbb{A}^N$ be a closed connected subvariety.  
Assume that $\dim (Y) = 1$ and $Y$ has smooth compactification.
Let $d$ be the degree of $Y$.
Set
\begin{equation}
\label{eqn:thm:limit:trop:d}
D := 
\begin{cases}
\max
\left\{
\left\lceil
\frac{3d^2 - 9 d + 4}{2d}
\right\rceil
,
1
\right\}
& 
\parbox{8cm}{%
\flushleft %
if 
$N \leq 2$, or
if
$N \geq 3$ and
$Y$ is contained in
some affine plane in $\Aff^N$,}
\\
\max \{ d-2 , 1 \}
& \text{otherwise,}
\end{cases}
\end{equation}
where for a real number $x$, $\lceil x \rceil$ denotes the 
smallest integer with $\lceil x \rceil \geq x$. 
Then the map $\mathrm{LTrop}_{\leq D}$ in
\eqref{eqn:limit:tropical:map:d} is a homeomorphism. 
\end{TheoremNoNum}

In Theorem~\ref{thm:limit:trop}, 
we take many polynomials 
$f_1, \ldots, f_n \in K[z_1, \ldots, z_N]$ (including a generating set for 
the coordinate ring of $Y$) and consider the embedding $\iota\colon Y \hookrightarrow \Aff^n$. 
Intuitively and roughly speaking, Theorem~\ref{thm:limit:trop} says that $\Trop(Y, \iota)$ approximates 
$Y^{\an}$ well, and if we take more and more polynomials $f_i$, then $\Trop(Y, \iota)$ approximates 
$Y^{\an}$ more and more. 
Theorem~\ref{thm:limit:trop:d:intro} asserts that
under suitable conditions, we can
do the same thing with only affine embeddings
of degree at most $D$, 
which is  effectively given
in terms of the degree of $X$.

In the proof of the injectivity of $\mathrm{LTrop}_{\leq D}$,
we will use 
Theorem~\ref{thm:main:unimodular:faithful} together with
some results in \S\ref{sec:ft:general:model}.
The proof will be given in the subsequent subsections.

\begin{Example}[Linear embedding of an elliptic curve]
\label{eg:limit:trop:linear}
Suppose that $F \in K[Z_0, Z_1, Z_2]$ is a nonsingular irreducible homogeneous polynomial of degree $d = 3$.  
We put $f(z_1, z_2) := F(1, z_1, z_2) \in K[z_1, z_2]$ and let $Y$ be the 
curve in $\mathbb{A}^2 =
\Spec ( K[z_1, z_2] )$ defined by $f$. Thus  
$Y$ is (an affine part of) an elliptic curve. In this case, $D := 
\max
\left\{
\left\lceil
\frac{3d^2 - 9 d + 4}{2d}
\right\rceil
,
1
\right\} = 1$,
and
Theorem~\ref{thm:limit:trop:d:intro} says that the map 
\begin{equation}
\label{eqn:limit:tropical:map:1}
\mathrm{LTrop}_{\leq 1}
\colon\;
Y^{\an} \to \varprojlim_{\iota \in I_{\leq 1}} \Trop(Y, \iota).
\end{equation}
is a homeomorphism,
namely, the analytification $Y^{\an}$ 
of $Y$ is the inverse limit of tropicalizations $\Trop(Y, \iota)$ 
of linear embeddings. 
\end{Example}

\begin{Example}[Linear embedding is not possible in general]
\label{eg:limit:trop:linear:impossible}
In general, it is false that the analytification of an affine curve is the inverse limit of tropicalizations 
of {\em linear} embeddings,
i.e., 
the map \eqref{eqn:limit:tropical:map:1} 
is not a homeomorphism in general. To see this, 
we take any $d \geq 4$ and take a homogeneous polynomial $F \in K[Z_0, Z_1, Z_2]$ as in the proof of 
Proposition~\ref{prop:countereg:very:ample} in \S\ref{subsec:very:ample:no:ft}. (Here we use $Z_0, Z_1, Z_2$ for the coordinate functions in place of $T_0, T_1, T_2$.) 
Let $X \subset \PP^2$ be the connected smooth projective curve defined by $F$, and 
let $Y \subset \Aff^2$ be the connected smooth affine curve defined by $f(z_1, z_2) := F(1, z_1, z_2)$. 
We take $[C_1], [C_2] \in X^{\an}$ as in the proof 
of Proposition~\ref{prop:countereg:very:ample} in \S\ref{subsec:very:ample:no:ft}. Since $Y^{\an} \setminus Y(K) = X^{\an} \setminus X(K)$, we have $[C_1], [C_2] \in Y^{\an}$.  For any $\iota \in I_{\leq 1}$, the proof of Proposition~\ref{prop:countereg:very:ample} 
shows that $\iota([C_1]) = \iota([C_2])$. Thus the images of $[C_1]$ and $[C_2]$ in $\varprojlim_{\iota \in I_{\leq 1}} \Trop(Y, \iota)$ coincide with each other, 
whence 
$\mathrm{LTrop}_{\leq 1}$ is not a homeomorphism. 
\end{Example}

\begin{Example}[Quadratic embedding of a non-hyperelliptic curve of genus $3$]
\label{eg:limit:trop:quad}
Let $X$ be a connected smooth non-hyperelliptic curve of genus $3$. 
Then the canonical map embeds $X$ as a degree $d = 4$ curve in $\PP^2$. 
Let $Y$ be the restriction of $X$ to $\Aff^2$. Since 
$\lceil \frac{3d^2 - 9 d + 4}{2d} \rceil  = 2$ when $d=4$, 
Theorem~\ref{thm:limit:trop:d:intro} says that 
the analytification $Y^{\an}$ is the limit of tropicalizations 
of quadratic embeddings. 
\end{Example}

It may be interesting to compare Examples~\ref{eg:limit:trop:linear}, \ref{eg:limit:trop:linear:impossible} and \ref{eg:limit:trop:quad} with Cueto--Markwig \cite{CM}, in which they study an algorithmic side of tropicalizations of plane curves, and with Wagner \cite{Wa}. 

\subsection{Polynomial of bounded degree that separates two points}
In order to show that $\mathrm{LTrop}_{\leq D}$ is injective,
we need to find, for any distinct $x , y \in
Y^{\an}$, a polynomial $f (z_1 , \ldots , z_N)$
of degree at most $D$
such that $- \log |f (x)| \neq - \log |f (y)|$.
In this section, we construct such an $f$.
The key to the construction is Proposition~\ref{prop:key:injective}.
In the proof of this proposition,
we use not only
Theorem~\ref{thm:main:unimodular:faithful}
but also
some global sections constructed in \S\ref{sec:ft:general:model}
together with the lemma below.
Here,
we
recall that for any compact skeleton $\Gamma$ of $X^{\an}$, $\tau_{\Gamma}: X^{\an} \to \Gamma$ denotes 
the retraction map with respect to $\Gamma$. 

\begin{Lemma}
\label{lemma:construction:Gamma1}
Let $x$ and $y$ be distinct points in $X^{\an}$ such that 
$x$ does not belong to any compact skeleton.
We fix a minimal skeleton $\Gamma_{\min}$ of $X^{\an}$. 
Then there exists a compact skeleton $\Gamma$
that contains $\Gamma_{\min}$ and satisfies the following conditions\textup{:}
\begin{align*}
& \Gamma = \Gamma_{\min} \cup
[ \tau_{\Gamma_{\min}} (x) , \tau_{\Gamma} (x) ], \quad 
\Gamma_{\min} \cap [\tau_{\Gamma_{\min}} (x) , \tau_{\Gamma} (x)] 
= \{ \tau_{\Gamma_{\min}} (x) \} \\
& 
\tau_{\Gamma} (x) \not\in \Gamma_{\min}, \quad 
\tau_{\Gamma} (x) \neq \tau_{\Gamma} (y).
\end{align*}
\end{Lemma}

\Proof
By \cite[Theorem~5.2]{BPR2}, 
it follows from $x \neq y$ that there exists a compact skeleton 
$\Gamma_0$ containing $\Gamma_{\min}$ such that $\tau_{\Gamma_0}(x) \neq \tau_{\Gamma_0}(y)$. Since 
$x$ does not belong to any compact skeleton, we have $\tau_{\Gamma_{0}} (x) \neq x$, 
and by \cite[Theorem~5.2]{BPR2} there exists a compact skeleton
$\Gamma_1$ with $\Gamma_{0} \subset \Gamma_{1}$
such that $\tau_{\Gamma_0} (x) \neq \tau_{\Gamma_{1}} (x)$. 
Since $\Gamma_{0} \subset \Gamma_{1}$, we have 
$\tau_{\Gamma_{0}} (\tau_{\Gamma_1} (x)) = \tau_{\Gamma_{0}} (x) \neq \tau_{\Gamma_1} (x)$, 
so that $\tau_{\Gamma_1} (x) \not\in \Gamma_{0}$, and in particular, $\tau_{\Gamma_1} (x) \not\in \Gamma_{\min}$. 

By Lemma~\ref{lemma:retraction:rational}, 
we have $\tau_{\Gamma_1} (x) \in \Gamma_{1 , \Lambda}$. 
We set $\Gamma := \Gamma_{\min} \cup [ \tau_{\Gamma_{\min}} (x) , \tau_{\Gamma_1} (x) ]$.
By Lemma~\ref{lemma:addingsegments},
$\Gamma$ is a compact skeleton.
Since $\tau_{\Gamma_{\min}} (x) , \tau_{\Gamma_1} (x)  \in \Gamma_1$, 
Lemma~\ref{lemma:to:be:used}(\ref{lemma:to:be:used:3}) gives us 
$[ \tau_{\Gamma_{\min}} (x) , \tau_{\Gamma_1} (x) ] \subset \Gamma_1$. 
Thus $\Gamma \subset \Gamma_1$. 

We have so far taken a compact skeleton $\Gamma_0$ with $\tau_{\Gamma_0}(x) \neq \tau_{\Gamma_0}(y)$, 
then we have taken a bigger compact skeleton $\Gamma_1$ with $\tau_{\Gamma_0}(x) \neq \tau_{\Gamma_1}(x)$, and then we have defined a compact skeleton $\Gamma$ as a part of this bigger skeleton $\Gamma_1$. 
We are going to show that $\Gamma$ has required properties. 

\smallskip
{\bf Step 1.}\quad
We show $\tau_{\Gamma} (x) \not\in \Gamma_{\min}$
and $\Gamma = \Gamma_{\min} \cup
[\tau_{\Gamma_{\min}} (x) , \tau_{\Gamma} (x) ]$. 
It follows from 
$ \tau_{\Gamma_1} (x) \in \Gamma$ and $\Gamma \subset \Gamma_1$ 
that $\tau_{\Gamma} (x) = \tau_{\Gamma} ( \tau_{\Gamma_1} (x) )
= \tau_{\Gamma_1} (x)$. 
Since $\tau_{\Gamma_1} (x) \not\in \Gamma_{\min}$, 
we have $\tau_{\Gamma} (x) \not\in \Gamma_{\min}$
and
$\Gamma = \Gamma_{\min} \cup
[\tau_{\Gamma_{\min}} (x) , \tau_{\Gamma} (x) ]$.

\medskip
{\bf Step 2.}\quad
We have $\Gamma_{\min} \cap [ \tau_{\Gamma_{\min}} (x) , \tau_{\Gamma} (x)  ]
= \{ \tau_{\Gamma_{\min}} (x) \}$, 
using $\tau_{\Gamma_{\min}} (\tau_{\Gamma} (x)) = \tau_{\Gamma_{\min}} (x)$ 
and the latter assertion of Lemma~\ref{lemma:to:be:used}(\ref{lemma:to:be:used:3}). 

\smallskip
{\bf Step 3.}\quad
We show the last property $\tau_{\Gamma} (x) \neq \tau_{\Gamma} (y)$. 
Noting that $\tau_{\Gamma}(x) \;(= \tau_{\Gamma_1}(x)) \not\in \Gamma_0$, 
we take the connected component $B$ of $X^{\an}\setminus \Gamma_0$ with $\tau_{\Gamma}(x) \in B$. 
By Lemma~\ref{lemma:connectedcomponents}, $B$ is the connected component 
of $X^{\an}\setminus \{\tau_{\Gamma_0}(x)\}$ with $\tau_{\Gamma}(x) \in B$, 
and we have $\tau_{\Gamma_0}(B) = \{\tau_{\Gamma_0}(x)\}$ (cf. Lemma~\ref{lemma:retraction:rational}).  

To prove $\tau_{\Gamma} (x) \neq \tau_{\Gamma} (y)$ 
by contradiction, we assume that $\tau_{\Gamma} (x) = \tau_{\Gamma} (y)$. 
Then we can show that $y \not\in \Gamma$. Indeed, since $\tau_{\Gamma_0} (y) \neq \tau_{\Gamma_0} (x)$ and 
$\Gamma_0 \subset \Gamma_1$, we have $\tau_{\Gamma_1} (y) \neq \tau_{\Gamma_1} (x)$. 
From the obvious equality $\tau_{\Gamma_1} (x)  = \tau_{\Gamma_1} (\tau_{\Gamma_1} (x) )$, 
it follows that
$\tau_{\Gamma_1} (y)
\neq
\tau_{\Gamma_1} (\tau_{\Gamma_1} (x))$,
and hence
$y \neq \tau_{\Gamma_1} (x) = \tau_{\Gamma} (x) = \tau_{\Gamma} (y)$.
This shows $y \notin \Gamma$.
Thus we take the connected component
$A$
of 
$X^{\an} \setminus \Gamma$ with 
$y \in A$. Since we have assumed that
$\tau_{\Gamma} (y) = \tau_{\Gamma} (x)$,
Lemma~\ref{lemma:connectedcomponents} tells us that 
$A$ is the connected component of $X^{\an} \setminus \{ \tau_{\Gamma} (x)\}
=X^{\an} \setminus \{ \tau_{\Gamma} (y)
\}
$ with $y \in A$.

We note that $\tau_{\Gamma_0} (x) \notin \overline{A}$. Indeed, 
since $[\tau_{\Gamma_{\min}} (x) , \tau_{\Gamma} (x) )$ is 
a connected subspace of $X^{\an} \setminus \{ \tau_{\Gamma} (x) \}$, 
$\tau_{\Gamma_{\min}} (x) \in \Gamma_{\min} \subset \Gamma$, and $\Gamma \cap A = \emptyset$, 
we have
$[\tau_{\Gamma_{\min}} (x) , \tau_{\Gamma} (x) ) \cap A = \emptyset$.
Since $\overline{A} = A \cup \{ \tau_{\Gamma} (x) \}$
and $\tau_{\Gamma} (x) \neq \tau_{\Gamma_{\min}} (x) $,
we have $[\tau_{\Gamma_{\min}} (x) , \tau_{\Gamma} (x) ) \cap 
\overline{A} = \emptyset$.
On the other hand, since 
$\tau_{\Gamma_0} (x) \neq \tau_{\Gamma_1}(x)$, 
Lemma~\ref{lemma:to:be:used}(\ref{lemma:to:be:used:4}) tells us that 
$\tau_{\Gamma_0} (x) \in [\tau_{\Gamma_{\min}}(x), \tau_{\Gamma_1}(x)) 
= [\tau_{\Gamma_{\min}}(x), \tau_{\Gamma}(x))$. 
We obtain $\tau_{\Gamma_0} (x) \notin \overline{A}$.

Since $\overline{A}$ is connected and $\tau_{\Gamma_0} (x) \notin \overline{A}$,
it follows that $\overline{A}$ is contained in some connected component
of $X^{\an} \setminus \{ \tau_{\Gamma_0} (x) \}$.
Since $B$ is a connected component of $X^{\an} \setminus \{ \tau_{\Gamma_0} (x) \}$ 
and $ \tau_{\Gamma}(x) \in \overline{A} \cap B$, 
we have $\overline{A} \subset B$.
Then
$\tau_{\Gamma_0} ( \overline{A}) \subset
\tau_{\Gamma_0} (B) = \{ \tau_{\Gamma_0} (x) \}$, 
which implies that
$\tau_{\Gamma_0} (y) = \tau_{\Gamma_0} (x)$.
However, that contradicts the choice of $\Gamma_0$. 
Thus $\tau_{\Gamma} (x) \neq \tau_{\Gamma} (y)$.
\QED

\begin{Remark}[four types of points in $X^{\an}$]
\label{rmk:four:types}
By Berkovich's classification theorem (cf. \cite{Be1}), 
one categorizes the points of $X^{\an}$ into four types.
The \emph{type I} points  are the classical points of $X^{\an}$ (cf. \S\ref{subsec:Berko:skeleta}),
i.e., the 
points in $X(K)$.
The {\em type II} points are
the Shilov points (cf. \S\ref{subsec:Berko:skeleta}), 
i.e.,
the points in
$\mathrm{Sh}(X^{\an}) := \bigcup_{\Gamma} \Gamma_{\Lambda}$, where 
$\Gamma$ runs through all the compact skeleta of $X^{\an}$. 
The {\em type III} points are those
in $X^{\an}$ that belong to some compact skeleta but are not 
Shilov points. 
Thus the set of type III points is equal to 
$\bigcup_{\Gamma} \Gamma \setminus  \bigcup_{\Gamma} \Gamma_{\Lambda}$, 
where $\Gamma$ runs through all the compact skeleta of $X^{\an}$. 
The {\em type IV} points are those 
in $X^{\an}$ that are neither of type I, II, or III. 
\end{Remark}

To prove Theorem~\ref{thm:limit:trop:d:intro}, we need to separate all types of points.
We have separated points in a skeleton to prove Theorem~\ref{thm:main:unimodular:faithful}, namely,
with the terminology 
in Remark~\ref{rmk:four:types},  
we have constructed global sections that separate type II and type III points. 
To separate two points one of which is type I or type IV,
the following is the key proposition. 
In the proof, we use a base section and a $P$-unimodularity section for a suitable $P \in X(K)$ 
and with respect to a suitable $B$,
which we consider in \S\ref{sec:ft:general:model}. 

Let $t (g)$ be as in Theorem~\ref{thm:main:unimodular:faithful}.

\begin{Proposition}
\label{prop:key:injective}
Let $X$ be a connected smooth projective curve over $K$ of genus $g$,
and let $L$ be a line bundle over $X$.
Assume that $\deg (L) \geq t (g)$.
Then for any $x , y \in X^{\an}$,
there exist nonzero
 global sections 
$s_0$ and $ s_1$ of $L$
such that $s_0 (x) \neq 0$ if $x \in X (K)$,
$s_0 (y) \neq 0$
if $y \in X (K)$,
and
$- \log |s_1 / s_0 (x)| \neq - \log |s_1 / s_0 (y)|$.
\end{Proposition}

\begin{Remark}
\label{remark:evaluation:globalsection}
In Proposition~\ref{prop:key:injective},
since $s_0$ and $s_1$ are nonzero sections,
$s_1 / s_0$ is a nonzero rational function on $X$.
It follows that if $x = (x, |\cdot|) \in X^\an \setminus X (K)$,
then 
$- \log |s_1 / s_0 (x)|$ makes sense as a real number.
\end{Remark}

\medskip

\Proof
First, we prove the assertion if $x$ or $y$ is a point in $X(K)$. 
By symmetry, we assume that $x \in X(K)$. 
It follows from $\deg (L) \geq t (g)$ and $t(g) \geq 2 g +1$ that 
$L$ is very ample.

Suppose that
$y \in X(K)$.
Since $L$ is very ample,
there exist global sections $s_0$
 and $s_1$
of $L$ such that $s_0 (x) \neq 0$, $s_0 (y) \neq 0$,
$s_1 (x) = 0$, and $s_1 (y) \neq 0$.
This proves
$- \log |s_1 / s_0 (x)| = + \infty \neq - \log |s_1 / s_0 (y)|$.

Suppose that $y = (y, |\cdot|) \not\in X(K)$. 
Since $L$ is very ample,
there exist nonzero global sections $s_0$
 and $s_1$
of $L$ such that $s_0 (x) \neq 0$ and
$s_1 (x) = 0$.
Then Remark~\ref{remark:evaluation:globalsection}
concludes that
 $- \log |s_1 / s_0 (x)| = + \infty \neq - \log |s_1 / s_0 (y)|$.

Thus in the following we may assume that neither $x$ nor $y$ belongs to $X(K)$. 
If there exists a skeleton $\Gamma$
such that $x,y \in \Gamma$, then the assertion follows immediately
from Theorem~\ref{thm:main:unimodular:faithful}.
Therefore, we may assume that 
at least one of $x,y$ does not belong to any skeleton.
Without loss of generality, we assume that 
$x$ does not belong to any skeleton.

We fix a minimal skeleton $\Gamma_{\min}$ in $X^{\an}$. 
We take a compact skeleton $\Gamma$ that contains $\Gamma_{\min}$ as in
Lemma~\ref{lemma:construction:Gamma1}, namely,
$\Gamma = \Gamma_{\min} \cup [\tau_{\Gamma_{\min}} (x), \tau_{\Gamma}(x)]$.
Since
$\tau_{\Gamma_{\min}} (x), \tau_{\Gamma}(x), \tau_{\Gamma}(y)$ are $\Lambda$-rational points (cf. Lemma~\ref{lemma:retraction:rational}), 
\cite[Theorem~4.11]{BPR2} gives us
a strictly semistable model $\mathscr{X}^0$ of $X$ 
such that
$S ( \mathscr{X}^0) = \Gamma$ and $\tau_{\min} (x)$,  $\tau_{\Gamma}(x)$, $\tau_{\Gamma}(y) 
\in V ( \mathscr{X}^0)$.
To ease notation, we put $w := \tau_{\Gamma} (x)$.
Since $w \in V(\mathscr{X}^0) $,
there exists $C^0_w \in \Irr ( \mathscr{X}^0_s) $ such that $[C^0_w] = w$.
Since $\tau_{\Gamma} (x) =  [C^0_w] \neq x$,
it follows from Lemma~\ref{lemma:retraction:rational}
that $\red_{\mathscr{X}^0} (x) \in C^0_{w} (k)
\setminus \Sing ( \mathscr{X}^0_s)$.
We take a section $\sigma^0$ of $\mathscr{X}^0 \to \Spec (R)$ such that
$\sigma^0 (k) \in C^0_w \setminus \Sing ( \mathscr{X}^0_s )$
and $\sigma^0 (k) \neq \red_{\mathscr{X}^0} (x)$.
Then $( \mathscr{X}^0 ; \sigma^0)$ is a strictly semistable pair.
Set $P := \sigma^0 (K) \in X (K)$.
Since $\sigma^0 (k) \in C^0_w \setminus \Sing ( \mathscr{X}^0_s )$,
we have
$[w , P)=\Delta ( \sigma^0)$ 
and
 $S ( \mathscr{X}^0;   \sigma^0) = 
\Gamma \cup [w , P)$
(cf. Lemma~\ref{lemma:end=canonialend}(3)).

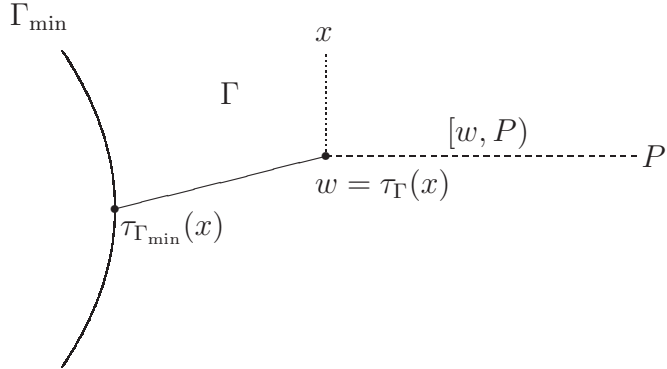
\begin{figure}[!h]
\[
\setlength\unitlength{0.07truecm}
\begin{picture}(150, 90)(0,0)
  \qbezier(10, 5)(30, 35)(10, 65)
  \multiput(60,45)(2, 0){30}{\line(1,0){1}}
  \multiput(60,45)(0, 1){20}{\line(0,1){0.2}}
  \put(20, 35){\circle*{1.5}}
  \put(20,35){\line(4,1){40}}
  \put(60, 45){\circle*{1.5}}
  \put(21, 30){$\tau_{\Gamma_{\min}}(x)$}
  \put(58, 38){$w = \tau_{\Gamma}(x)$}
  \put(40, 55){$\Gamma$}
  \put(83, 48){$[w, P)$}
  \put(0, 70){$\Gamma_{\min}$}
  \put(58, 67){$x$}
  \put(120, 43){$P$}
 \end{picture}
\]
\label{figure:for:prop:key:injective}
\end{figure}

Let $A$ be the connected component of $X^{\an} \setminus \Gamma$ 
with $x \in A$.

\begin{Claim}
\label{claim:disjoint:prop:key:injective}
We have $A \cap [w , P) = \emptyset$.
\end{Claim}

Indeed, to argue by contradiction,
suppose that $A \cap [w , P) \neq \emptyset$. 
Since $w = \tau_{\Gamma} (x) \in \Gamma$ and $A \cap \Gamma = \emptyset$, 
we have  $A \cap (w , P) \neq \emptyset$.
Since $A$ is a connected component of $X^{\an} \setminus \{ w \}$ by Lemma~\ref{lemma:connectedcomponents}, 
it follows that $(w , P) \subset A$, 
and thus $[w , P ] \subset \overline{A} = A \cup \{ w \}$.
In particular, $P \in A$.
This shows that
\[
\{ \sigma^0 (k) \} = \{ \red_{\mathscr{X}^0} (P) \}
\subset \red_{\mathscr{X}^0} (A) \supset \{ \red_{\mathscr{X}^0} (x) \}
.
\]
By Lemma~\ref{lemma:retraction:rational}, $\red_{\mathscr{X}^0} (A)$ is a singleton, so that 
$\sigma^0 (k) = \red_{\mathscr{X}^0} (x)$. 
This contradicts the choice of $\sigma^0$.
Thus 
Claim~\ref{claim:disjoint:prop:key:injective} holds.  

\smallskip
Now, using Proposition~\ref{prop:main:FT:ends},
we obtain 
a strictly semistable model $(\mathscr{X} , \mathscr{L})$ of $(X,L)$
and a section $\sigma$ of $\mathscr{X} \to \Spec (R)$
such that $\Xscr$ dominates $\Xscr^0$ and 
$(\mathscr{X} ; \sigma)$ is a strictly semistable pair 
with the following properties:
\begin{enumerate}
\item[(i)]
$S(\mathscr{X} ) = \Gamma$, $V(\Xscr) \supset V(\Xscr^0)$, 
and $\sigma (K) = P$;
\item[(ii)]
We set 
\[
B :=
\begin{cases}
\{ \red_{\mathscr{X}} (x) , \red_{\mathscr{X}} (y) \}
&
\text{if $\red_{\mathscr{X}} (y) \in \mathscr{X}_s (k) \setminus (\Sing ( \mathscr{X}_s) 
\cup
\{ \sigma (k) \})$,}
\\
\{ \red_{\mathscr{X}} (x) \}
&
\text{otherwise}: 
\end{cases}
\]
Then
$\mathscr{L}$ has a base section $\widetilde{s}_0$
with respect to $B$; 
\item[(iii)]
$\mathscr{L}$ has 
a $P$-unimodularity section $\widetilde{s}_1$
with respect to $B$.
\end{enumerate}
Here we remark that property $V(\Xscr) \supset V(\Xscr^0)$ in (i) holds, since 
$\Xscr$ dominates $\Xscr^0$. 
We also remark that $\red_{\mathscr{X}} (x) \in
 \mathscr{X}_s (k) \setminus (\Sing ( \mathscr{X}_s) 
\cup
\{ \sigma^0 (k) \})$,
since $\red_{\mathscr{X}^0} (x) \in \mathscr{X}^0_s(k)
\setminus \Sing ( \mathscr{X}^0_s) $ 
and $\red_{\mathscr{X}^0} (x)  \neq \sigma^0 (k)$.

\smallskip
For $i=0,1$, set $s_i := \rest{\widetilde{s}_i}{X}$.
We show that these $s_0, s_1$ have required properties. 
Since we assume that neither $x$ nor $y$ belongs to $X(K)$, 
we have $s(x) \neq 0$ and $s(y) \neq 0$ for any nonzero global section $s$ of $L$, 
whence $s_0(x) \neq 0$ and $s_0(y) \neq 0$.  
It remains to show that $- \log | s_1 / s_0 (x)|
\neq - \log | s_1 / s_0 (y)|$. 

\smallskip
{\bf Step 1.} In this step, we consider $- \log | s_1 / s_0 (x)|$:
the  
goal in this step is to
show that 
\begin{equation}
\label{eqn:compari:x:w:1}
- \log | (s_1 / s_0) (x)| = - \log | (s_1 / s_0) (w)|,
\end{equation}
where we recall that $w := \tau_{\Gamma}(x)$. 

We use the notation in Remark~\ref{rmk:four:types}: 
$\mathrm{Sh}(X^\an) $ denotes the set of Shilov points of $X^{\an}$. 
To prove (\ref{eqn:compari:x:w:1}),
we first show that 
for any $u \in A \cap \mathrm{Sh}(X^\an)$,  
\begin{equation}
\label{eqn:compari:x:w:2}
- \log | s_1 / s_0 (u)|
= - \log | s_1 / s_0 (w)|.
\end{equation}
We take any $u \in A \cap \mathrm{Sh}(X^\an)$.
By the definition of $\tau_{\Gamma}$, we have $\tau_{\Gamma}(A) = \{w\}$,
Then since $\tau_{\Gamma}(u) = w = \tau_{\Gamma}(w)$, let 
$[u, w]$ denotes the geodesic segment (see \S\ref{subsection:geodesic:line}). 
Since $u \in \mathrm{Sh}(X^\an)$, 
Lemma~\ref{lemma:to:be:used}
shows that
$\Gamma \cup [w , u]$ is a compact skeleton and
$\Gamma \cap [w , u] = \{ w \}$.
Noting \cite[Theorem~4.11]{BPR2},
we take
a strictly semistable model $\mathscr{X}'$ 
such that 
$S ( \mathscr{X}^\prime ) = \Gamma \cup [w , u]$
and $V( \mathscr{X}' ) = V ( \mathscr{X} ) \cup \{ u \}$;
we remark that
by this equality, the identity on $X$ extends to a morphism
$\mu : \mathscr{X}^\prime \to \mathscr{X}$.
Let $\sigma^\prime$ be the section of $\mathscr{X}^\prime \to \Spec (R)$ with
$\sigma = \mu \circ \sigma^\prime$.
Then 
$( \mathscr{X}^\prime ;  \sigma^\prime 
)
$ is a strictly 
semistable pair, and we have 
$
S ( \mathscr{X}^\prime ; 
\sigma^\prime 
) 
=
\Gamma \cup \Delta(\sigma^\prime)$. 
We are going to show that, for
$( \mathscr{X}^\prime ;  \sigma^\prime 
)$, 
$\mu^{\ast} (\widetilde{s}_0)$ is a base section
of $\mu^{\ast} ( \mathscr{L})$
with respect to the empty set $\emptyset$
and $\mu^{\ast} (\widetilde{s}_1)$ is a $P$-unimodularity
section of $\mu^{\ast} ( \mathscr{L})$
with respect to $\emptyset$.

By Lemma~\ref{lemma:connectedcomponents},
$A$ is a connected component of $X^{\an} \setminus \{ w \}$, and we have 
$u \in A$.
Since
$(w,u]$ is a connected subspace of $X^{\an} \setminus \{ w \}$ 
with $u \in (w,u]$,
it follows that $(w,u] \subset A$.
Thus
$\red_{\mathscr{X}} (( w ,u ]) 
 \subset \red_{\mathscr{X}} (A) = \{ \red_{\mathscr{X}} (x) \}$
(cf. Lemma~\ref{lemma:retraction:rational}),
which is contained in
$B$. 
Here, by the definition of $\mathscr{X}'$,
we note that $\mu : \mathscr{X}' \to \mathscr{X}$ is 
not an isomorphism
over $\red_{\mathscr{X}} (x)$ 
but an isomorphism except over 
$\red_{\mathscr{X}} (x)$.
Since 
$\red_{\mathscr{X}} (x) \in B$
and
$\widetilde{s}_0$
is a base section of $\mathscr{L}$
with respect to $B$,
it follows 
that $\mu^{\ast} ( \widetilde{s}_0 )$ is a base section of $\mu^{\ast} 
( \mathscr{L})$
with respect to $\emptyset$.

To see that
$\mu^{\ast} ( \widetilde{s}_1 )$ is a $P$-unimodularity section
of $\mu^{\ast} 
( \mathscr{L})$ with respect to $\emptyset$,
let $\mathscr{W}_P$ be the 
stepwise vertical divisor associated to $P$
defined in (\ref{align:stepwise:for:ends}).
Since $\widetilde{s}_1$ is a $P$-unimodularity section
of $\mathscr{L}$ with respect to $B$,
there exists an open neighborhood of $\mathscr{U} \subset \mathscr{X}$
of $B \cup \Sing ( \mathscr{X}_s ) \cup \{ \sigma (k) \}$
such that
$\zero (\widetilde{s}_1) - \mathscr{W}_P - \sigma (R)$
is trivial on $\mathscr{U}$.
The pullback of this Cartier divisor by $\mu$ equals
$\zero (\mu^{\ast} (\widetilde{s}_1)) - 
\mu^{\ast} ( \mathscr{W}_P ) - \mu^{\ast} ( \sigma (R) )$,
which is trivial on
$\mu^{-1} (\mathscr{U})$.
We remark that $\mu^{-1} (\mathscr{U})$ is an open neighborhood of 
$\Sing ( \mathscr{X}'_s) \cup \{ \sigma' (k) \}$.
Let $C_w$ be the irreducible component of $\mathscr{X}_s$
with $[C_w] = w$.
Since $\tau_{\Gamma} (x) = w$,
Lemma~\ref{lemma:retraction:rational}
tells us 
that $\red_{\mathscr{X}} (x) \in C_{w} (k) \setminus \Sing (\mathscr{X}_s)$.
We take $\varpi_w \in R$ such that $- \log |\varpi_w|
= \ord_{C_w} (\mathscr{W}_P)$.
Then on some neighborhood of $\red_{\mathscr{X}} (x)$,
the Cartier divisor $\mathscr{W}_P$ is defined by $\varpi_w$.
Thus on some neighborhood of $\mu^{-1} ( \red_{\mathscr{X}} (x) )$,
$\mu^{\ast} ( \mathscr{W}_P ) $ is defined by $\varpi_w$.
Further, $\red_{\mathscr{X}} (x) \neq \sigma (k)$.
By the definition of the stepwise vertical divisor associated to $P$,
it follows that
$\mu^{\ast} ( \mathscr{W}_P ) $ equals the 
the stepwise vertical divisor 
$\mathscr{W}_P' $ on $\mathscr{X}'$ associated to $P$.
Furthermore, since 
$\mu$ is isomorphism over $\sigma (R)$,
we have
$\mu^{\ast} ( \sigma (R) ) = \sigma' (R)$.
Thus 
$\zero (\mu^{\ast} (\widetilde{s}_1)) - 
\mathscr{W}_P' - \sigma' (R)$
is trivial over $\mu^{-1} ( \mathscr{U})$.
This proves that $\mu^{\ast} (\widetilde{s}_1)$ is a
$P$-unimodularity section of $\mu^{\ast} ( \mathscr{L})$
with respect to $\emptyset$.

Note by Claim~\ref{claim:disjoint:prop:key:injective}
that
$[w , u] \cap [w , P) = \{ w \}$.
Then by Lemma~\ref{lemma:faithful:ends2}(2), it follows that
$- \log | s_1 / s_0 (u)| = - \log | s_1 / s_0 (w)|$.
We have shown \eqref{eqn:compari:x:w:2}. 
 
Recall that, since $X^{\an}$ is a locally connected space, 
$A$ is open in $X^{\an}$. Since $\mathrm{Sh} ( X^{\an})$ is dense in $X^{\an}$, 
$A \cap \mathrm{Sh} ( X^{\an}) $ is dense in $A$.
Since $- \log | s_1 / s_0 |$ is a continuous function (with values in $\TT$) on $A$, 
the equality \eqref{eqn:compari:x:w:2} holds for any $u \in A$
as well as for points in $A \cap  \mathrm{Sh} ( X^{\an})$. 
Since $x \in A$, this gives \eqref{eqn:compari:x:w:1}. 

\smallskip
{\bf Step 2.}\quad 
In this step, we consider $- \log | s_1 / s_0 (y)|$. 
Noting that $\tau_{\Gamma_{\min}}(x) = \tau_{\Gamma_{\min}}(\tau_{\Gamma}(x))$ and 
$w:= \tau_{\Gamma}(x)$, we consider the geometric segment 
$[\tau_{\Gamma_{\min}} (x) , w]$ (see \S\ref{subsection:geodesic:line}).  
We first treat the case where $\tau_{\Gamma} (y) \in [\tau_{\Gamma_{\min}} (x) , w]$, 
and then the case where $\tau_{\Gamma} (y) \not\in [\tau_{\Gamma_{\min}} (x) , w]$. 

\smallskip
{\bf Case 1.}
Suppose that $\tau_{\Gamma} (y) \in [\tau_{\Gamma_{\min}} (x) , w]$.
By Lemma~\ref{lemma:faithful:ends}(i),
the function $- \log | s_1 / s_0| : [\tau_{\Gamma_{\min}} (x) , w]
\to \RR$ is an affine function with 
$- \log | s_1 / s_0 ( \tau_{\Gamma_{\min}} (x))| = 0$
and $- \log | s_1 / s_0 ( w)| > 0$. 
By the definition of $\Gamma$, we have $\tau_\Gamma(y) \neq w$, 
so that 
$
- \log | s_1 / s_0 (\tau_{\Gamma_1} (y)  )| < 
- \log | s_1 / s_0 (w)|
$.
Thus
by (\ref{eqn:compari:x:w:1}), 
\begin{align}
\label{eqn:compari:x:y}
- \log | s_1 / s_0 (x)| > 
- \log | s_1 / s_0 (\tau_{\Gamma} (y))|.
\end{align}
If $y = \tau_{\Gamma} (y)$, 
then \eqref{eqn:compari:x:y} concludes that
$- \log | s_1 / s_0 (x)| \neq 
- \log | s_1 / s_0 (y)|$.

Suppose now that $y \neq \tau_{\Gamma} (y)$.
Let $A^\prime$ be the connected component of 
$X^{\an} \setminus \Gamma$
with $y \in A^\prime$.
By Lemma~\ref{lemma:connectedcomponents},
$A^\prime$ is the connected component
$X^{\an} \setminus \{ \tau_{\Gamma} (y) \}$
with $y \in A^\prime$.
Then by the same argument as in Step 1, 
 $ - \log |s_1 / s_0 ( y )|
=
- \log |s_1 / s_0 ( \tau_{\Gamma} 
(y))|$, 
where  we use $\red_{\mathscr{X}} (y) \in B$
if $\red_{\mathscr{X}} (y) \in \mathscr{X}_s (k)
\setminus (
\Sing ( \mathscr{X}_s ) \cup
\{ \sigma^0 (k) \} )$.
Together with (\ref{eqn:compari:x:w:1}),
this concludes that $- \log  |s_1 / s_0 (y)| \neq - \log |s_1 / s_0 ( x )|$.

\smallskip
{\bf Case 2.}
Suppose that $\tau_{\Gamma} (y) \notin [\tau_{\Gamma_{\min}} (x) , w]$.
Since
$\Gamma = \Gamma_{\min} \cup [\tau_{\Gamma_{\min}} (x), w]$
and $\tau_{\Gamma} (y) \in \Gamma$,  
$\tau_{\Gamma} (y) \in \Gamma_{\min}$. 
By Lemma~\ref{lemma:faithful:ends}(iii),
we have $- \log |s_1 / s_0 (\tau_{\Gamma} (y))| = 0$.
If $y = \tau_{\Gamma} (y)$, then 
$- \log |s_1 / s_0 (y)| = 0$;
since 
$- \log |s_1 / s_0 (x)| 
= - \log |s_1 / s_0 (w)| > 0$
by \eqref{eqn:compari:x:w:1}, 
this concludes
$- \log |s_1 / s_0 (x)| \neq - \log |s_1 / s_0 (y)| $.

Suppose that $y \neq \tau_{\Gamma} (y)$.
Let $A^{\prime}$ be the connected component of 
$X^{\an} \setminus 
\Gamma$ with $y \in A^{\prime}$.
Then by the same argument as in Case 1 (or Step 1), 
one shows that the restriction of $- \log |s_1/s_0|$ 
to $A^{\prime}$ equals the constant function $0$. 
Since $y \in A^{\prime}$, 
this concludes $- \log |s_1 / s_0 (x)| \neq - \log |s_1 / s_0 (y)| $.
\QED

\subsection{Proof of Theorem~\ref{thm:limit:trop:d:intro}}
\label{subsec:proof:thm:limit:trop:d:intro}
Let $Y \subset \mathbb{A}^N = \Spec(K[z_1, \ldots, z_N])$ be a (closed) connected smooth affine curve.
From here on,
we regard $\mathbb{A}^N$ as an open subscheme of $\PP^N$ via the embedding
$(z_1 , \ldots , z_N ) \mapsto (1 : z_1 : \cdots : z_N)$.
Let $\OO_{\PP^N} (1)$ denote the tautological line bundle over  $\PP^N$.
Since $Y \subset \mathbb{A}^{N} \subset \PP^N$,
we regard $Y$ as a subscheme of $\PP^N$.
Let $X$ be the closure of $Y$ in $\PP^N$. 
Recall that 
the degree $d$ of $Y$ is defined as the 
the degree of $X$ in $\PP^N$. 
Assume that $Y$ has smooth compactification in $\PP^N$,
i.e.,
$X$ is smooth.
Let $g$ denote the genus of $X$.
We also call $g$ the genus of $Y$.
For any integer $m$, set $\OO_X (m) := \rest{\OO_{\PP^N}(m)}{X}$,
which is a very ample line bundle over $X$.

\medskip
We say that $Y$ is \emph{degenerate} if it is contained in some 
hyperplane of 
$\mathbb{A}^N$
and that 
$X$ is \emph{degenerate} if it is contained in some 
hyperplane of 
$\mathbb{P}^N$.
Then $Y$ is degenerate if and only if $X$ is degenerate. 

If $X$ is non-degenerate, then
$N \leq d$.
Indeed, 
if $X$ is non-degenerate,
then
the restriction $H^0 
( \PP^N , \OO_{\PP^N}(1) ) \to H^0 ( X , \OO_X (1) )$ is injective,
and since $h^{0} (X , \OO_X (1)) \leq \deg ( \OO_X (1) ) + 1 = d+1$
and $h^0 ( \PP^N , \OO_{\PP^N}(1) ) = N+1$,
we have $d \geq N$.

\begin{Proposition}
\label{prop:injective:affine}
Let $Y \subset \mathbb{A}^N = \Spec(K[z_1, \ldots, z_N])$ be a \textup{(}closed\textup{)} connected smooth affine curve
 that has smooth compactification. 
Let $d$ be the degree of $Y$, and let $g$ be the genus of $Y$. 
Assume that $Y$ is non-degenerate.
Let $\ell$ be an integer.
Suppose that
\begin{align}
\label{align:forell}
\ell \geq 
\begin{cases}
\lceil t (g)/d \rceil
&
\text{if $N = 1, 2$,}
\\
\max
\{ \lceil t (g)/d \rceil, d + 1 - N \}
& \text{if $N \geq 3$.}
\end{cases}
\end{align}
Then for any $x ,y \in Y^{\an}$ with $x \neq y$,
there exists a polynomial $f ( z_1 , \ldots , z_N ) \in
K[z_1 , \ldots , z_N]$ of degree at most $\ell$ such that 
$- \log |f (x)| \neq - \log |f (y)|$.
\end{Proposition}

\Proof
Let $X$ be the closure of $Y$ in $\PP^N$. Then 
$X$ is a connected smooth projective curve 
over $K$ of genus $g$ and degree $d$. 
We denote by $\jmath: X \to \PP^N$ the natural embedding. 
Let $Z_0, \ldots, Z_N$ denote the homogeneous coordinate functions of 
$\PP^N$. 

We note that the restriction
$
\jmath^*: 
H^0 ( \PP^N , \OO_{\PP^N}(\ell) )
\to H^0 (X , \OO_{X} (\ell))$ is
surjective. 
Indeed, 
if $N = 1$, then $X = \PP^1$ and the surjectivity is obvious. 
If $N = 2$, then 
the cokernel of $\jmath^*$ is a subspace
of $H^{1} ( \PP^2 , \OO_{\PP^2}(\ell - d) )$. 
Since $H^{1} ( \PP^2 , \OO_{\PP^2}(\ell - d) ) = 0$, 
we see that $\jmath^*$ is surjective. 
If $N \geq 3$,
then since $\ell \geq d + 1 - N$,
\cite[Theorem]{GLP} tells us that
$\jmath^*$ is surjective. (Here we use the assumption that $Y$ is non-degenerate.)

We have $\deg (\OO_X(\ell)) = d \ell$. Since $\ell$ is an integer with $\ell \geq t (g)$, 
Proposition~\ref{prop:key:injective} gives us 
nonzero  $s_0 , s_1 \in H^0 (X, \OO_X(\ell))$ such that
$s_0 (x) \neq 0$, $s_0 (y) \neq 0$,
and
$- \log |s_1 / s_0 (x)| \neq - \log |s_1 / s_0 (y)|$.
Since $\jmath^*: H^0 ( \PP^N , \OO_{\PP^N}(\ell) ) \to H^{0} 
(X , \OO_X (\ell))$ is surjective, 
there exist $F_0 (Z_0, \ldots , Z_N), F_1(Z_0, \ldots , Z_N)
\in H^0 ( \PP^N , \OO_{\PP^N}(\ell))$
such that $\rest{F_0}{X} = s_0$ and $\rest{F_1}{X} = s_1$. 

For $i = 0, 1$, 
we set $f_i (z_1 , \ldots , z_N) := F_i (1 , z_1 , \ldots , z_N) \in K[z_1, \ldots, z_N]$. 
Since $F_i$ is a homogeneous polynomial of degree $\ell$, 
$f_i \in K[z_1, \ldots, z_N]$ is a polynomial of degree at most $\ell$. 
Then we have 
\[
  - \log | (f_1 / f_0) (x)| = - \log |s_1 / s_0 (x)| 
  \neq - \log |s_1 / s_0 (y)| = - \log | (f_1 / f_0) (y)|. 
\]
It follows that
$- \log |f_0 (x)| \neq - \log |f_0 (y)|$ or $- \log |f_1 (x)| \neq - \log |f_1 (y)|$.  
This completes the proof.
\QED

\begin{Remark}
\label{remark:m0e0}
Let $Y$, $d$, and $g$  be as in Proposition~\ref{prop:injective:affine}. 
Let $X$ be the closure of $Y$ in $\PP^N$, so that 
$X$ is a connected smooth projective curve in $\PP^N$ of genus $g$ and degree $d$. 
Then we have the following sufficient condition
for $\ell$ to satisfy (\ref{align:forell}).
\begin{enumerate}
\item
Suppose that $N=1$. Then $Y = \Aff^1$, $d = 1$ and $g = 0$. 
In this case any integer $\ell \geq 1$ satisfies 
\eqref{align:forell}. 
\item
Suppose that $N=2$.
Since $X$ is smooth, we have $g = (d-1)(d-2)/2$. 
If $d = 1, 2$, then $g = 0$ and $t(g)= 1$. 
If $d = 3$, then $g=1$ and $t(g) = 3$. If 
$d \geq 4$, then $t(g) = 3g-1 =  (3 d^2 -9d + 4)/2$.  
Thus if $\ell$ is an integer with 
\[
\ell
\geq 
\begin{cases}
1 & \text{if $d = 1, 2 , 3$}, 
\\ 
\left\lceil
\frac{3d^2 - 9 d + 4}{2d}
\right\rceil
& \text{if $d \geq 4$}, 
\end{cases}
\]
then $\ell$ satisfies 
(\ref{align:forell}).
\item
Suppose that $N \geq 3$. Suppose also that 
$Y$ is non-degenerate. Then $d \geq N$. Dividing 
$d-1$ by $N-1$, we write 
\[
d - 1 = m_0 (N-1) + \epsilon_0,
\]
where $m_0$ and $\epsilon_0$ are integers
with $m_0 \geq 1$ and $0 \leq \epsilon_0 \leq N-2$. 
Let $\pi(d, N)$ be Castelnuovo's number, which is defined by 
\[
 \pi(d, N) := \frac{( m_0 + 1) (\epsilon_0 + d - 1)}{2} - (d-1). 
\]
Then by Castelnuovo's bound theorem (see \cite[Chap.~III~\S2]{ACGH}),
we have
\begin{align}
\label{align:AGM}
g \leq 
\pi(d, N). 
\end{align}
(In \cite[Chap.~III~\S2]{ACGH}, Castelnuovo's number $\pi(d, N)$ is defined to be 
$\frac{m_0(m_0-1)}{2}(N-1) + m_0 \epsilon_0$, and one sees that 
$\frac{( m_0 + 1) (\epsilon_0 + d - 1)}{2} - (d-1) =  \frac{m_0(m_0-1)}{2}(N-1) + m_0 \epsilon_0$.)

Since $d \geq N \geq 3$,
we obtain
\[
t (g)/d
\leq 
\begin{cases}
1 & \text{if $g=0,1$}
\\
\frac{3( m_0 + 1) (\epsilon_0 + d - 1)}{2d} - \frac{3d-2}{d}
& \text{otherwise.}
\end{cases}
\]
Thus if 
$\ell$ is an integer with
\[
\ell
\geq
\max
\left\{
\left\lceil 
\frac{3( m_0 + 1) (\epsilon_0 + d - 1)}{2d} - \frac{3d-2}{d}
\right\rceil
,
d+1-N
\right\}
,
\]
then $\ell$ satisfies 
\eqref{align:forell}.
\end{enumerate}
\end{Remark}

\medskip
We start the proof of Theorem~\ref{thm:limit:trop:d:intro}. 
Recall that $D$ is an integer given by 
\eqref{eqn:thm:limit:trop:d}. 
Our goal is to prove that the map
\[
\mathrm{LTrop}_{\leq D} : Y^{\an} \to \varprojlim_{\iota \in I_{\leq D}} \Trop(Y, \iota)
\]
in
\eqref{eqn:limit:tropical:map:d} is an homeomorphism.

\subsubsection*{Injectivity}
 
The injectivity essentially amounts to the following proposition.

\begin{Proposition}
\label{prop:injective:affine:key}
Let $Y \subset \mathbb{A}^N = \Spec(K[z_1, \ldots, z_N])$ be a \textup{(}closed\textup{)} connected smooth affine curve
that has smooth compactification. 
Let $d$ be the degree of $Y$, and let $g$ be the genus of $Y$. 
Further, let $D$
be as in \eqref{eqn:thm:limit:trop:d}.
Then for any $x ,y \in Y^{\an}$ with $x \neq y$,
there exists a polynomial $f ( z_1 , \ldots , z_N ) \in
K[z_1 , \ldots , z_N]$ of degree at most $D$ such that 
$- \log |f (x)| \neq - \log |f (y)|$.
\end{Proposition}

\Proof
First, we assume that $Y$ is non-degenerate.
If $N = 1$, then the assertion immediately follows from 
Proposition~\ref{prop:injective:affine} and Remark~\ref{remark:m0e0}(1). 
If $N = 2$, then 
Proposition~\ref{prop:injective:affine} and 
Remark~\ref{remark:m0e0}(2) shows the assertion.
Thus we may and do assume that $N \geq 3$.
Recall that $D := \max \{ d-2 , 1\}$. 
By non-degeneracy assumption of $Y$, we have $d \geq N \geq 3$, 
and thus $D = d-2$. 
With the notation in Remark~\ref{remark:m0e0}(3), 
we are going to show that 
\begin{equation}
\label{eqn:prop:injective:affine:key}
 D \geq 
\max
\left\{
\left\lceil 
\frac{3( m_0 + 1) (\epsilon_0 + d - 1)}{2d} - \frac{3d-2}{d}
\right\rceil
,
d + 1 - N
\right\}. 
\end{equation}
Indeed, by Proposition~\ref{prop:injective:affine}
and Remark~\ref{remark:m0e0}(3),
the assertion follows from this inequality.

Let us prove (\ref{eqn:prop:injective:affine:key}).
By the inequality of arithmetic and geometric means,
we have
\[
(N - 1) (m_0 + 1) (\epsilon_0 + d -1)
\leq 
\left( 
\frac{(N-1) (m_0 + 1) + (\epsilon_0 + d - 1)}{2}
\right)^2
.
\]
Since 
$(N-1) (m_0 + 1) + (\epsilon_0 + d - 1)
= 2(d-1) + N-1$ and $3 \leq N \leq d$, 
we have 
\begin{align*}
 (m_0 + 1) (\epsilon_0 + d -1)
& \leq
\frac{(d-1)^2}{N-1} + (d-1) + \frac{N -1}{4} \\
& \leq 
\frac{(d-1)^2}{2} + \frac{5 (d -1)}{4}  
= \frac{2d^2+d - 3}{4}.
\end{align*}
It follows that
\[
\frac{3( m_0 + 1) (\epsilon_0 + d - 1)}{2d} - \frac{3d-2}{d}
\leq
\frac{3(2d^2+d - 3)}{8d}
- \frac{3d-2}{d}
=
\frac{6d^2 - 21 d + 7}{8d}
.
\]
Further,
since
$
(d - 2) -
\frac{6d^2 - 21 d + 7}{8d}
=
\frac{2d^2 + 5 d - 7}{8d} \geq 0
$, 
we have
\[
d - 2 
\geq
\left\lceil
\frac{3( m_0 + 1) (\epsilon_0 + d - 1)}{2d} - \frac{3d-2}{d}
\right\rceil
.
\]
Since
$d - 2 \geq d + 1 - N$,
we obtain \eqref{eqn:prop:injective:affine:key}. 
This completes the proof of the proposition when $Y$ is non-degenerate.

It remains to consider the case where $Y$ is degenerate.
Let $W$ be the smallest linear subspace of $\PP^N$ such that $X \subset W$.
Let $r$ denote the dimension of $W$.
Since $Y \subset W$ and $Y = X \cap \mathbb{A}^N$,
$W \cap \mathbb{A}^{N}$ is an affine subspace of $\mathbb{A}^N$ of dimension $r$
and $Y$ is a closed subvariety of this affine space.
We fix an isomorphism 
$\phi : W \cap \mathbb{A}^{N} \to \mathbb{A}^r$.
Let $w_1 , \ldots , w_r$ denote the affine coordinates of $\mathbb{A}^r$.
We regard $\mathbb{A}^r$ as an open subscheme of $\PP^r$ via
$(w_1 , \ldots , w_r) \mapsto (1 : w_1 : \cdots : w_r)$.
Then $\phi$ extends to a unique isomorphism $W \to \PP^r$.
Note that the closure of $\phi (Y)$ in $\PP^r$ is smooth and non-degenerate
and
has degree $d$.
Further, note that the integer
$D'$ defined by the formula (\ref{eqn:thm:limit:trop:d}) for $\phi (Y) \subset \PP^r$ equals $D$.
It follows from what we have shown above that 
there exists a polynomial $h$ on $w_1 , \ldots , w_r$
of degree at most $D$ such that
$- \log |h (\phi (x))| \neq - \log |h (\phi (y))|$.
Since $W \cap \mathbb{A}^{N}$ is an affine subspace of $\mathbb{A}^N$,
$\phi^{\ast} (w_1) , \ldots , \phi^{\ast} (w_r)$ are 
the restrictions of some polynomials 
of degree $1$ on $z_1 , \ldots , z_N$.
It follows that
there exists a polynomial $f$ on $z_1 , \ldots , z_r$ of degree at most $D$
such that $\rest{f}{Y} = \phi^{\ast} ( h)$.
Then we have $- \log |f (x)| \neq - \log |f (y)|$. 
Thus the proof is complete.
\QED

We prove that  $\mathrm{LTrop}_{\leq D}$ is injective. 
We take any $x, y \in Y^{\an}$ with $x \neq y$. 
It suffices to show that there exists $\iota \in I_{\leq D}$ such that 
$\pi_{\iota} (x) \neq \pi_{\iota}(y)$. 
Take any $\iota_0 \in I_{\leq D}$,
and we write $\iota_0 = ( \rest{f_1}{Y} , \ldots , \rest{f_n}{Y}  )$ for 
some $f_1 , \ldots , f_n \in K [z_1 , \ldots , z_N ]$
of degree at most $D$.
By Proposition~\ref{prop:injective:affine:key},
there exists $f (z_1 , \ldots , z_N )$ such that
$- \log | f ( x) | \neq - \log | f ( y) |$.
Set
$\iota := (\rest{f}{Y} , \rest{f_1}{Y} , \ldots , \rest{f_n}{Y} )$.
Then $\iota$ is an affine embedding of $Y$ of degree at most $D$,
and $\pi_{\iota} (x) \neq \pi_{\iota} (y) $.
This shows that $\mathrm{LTrop}_{\leq D}$ is injective.

\subsubsection*{Surjectivity}
We prove that $\mathrm{LTrop}_{\leq D}$ is surjective.
Take any $(b_\iota)_{\iota \in I_{\leq d}} \in \varprojlim_{\iota \in I_{\leq D}} \Trop(Y, \iota)$. 
For each $\iota \in I_{\leq d}$,
we set $A_\iota := \pi_\iota^{-1}(b_\iota)$. Since $\pi_\iota\colon Y^{\an} \to \trop(Y, \iota)$ is surjective, 
we have $A_\iota \neq \emptyset$
for any $\iota \in I_{\leq d}$. 
We claim that $\{A_{\iota}\}_{\iota \in I_{\leq D}}$ has a finite intersection property, i.e., 
for any finitely many elements $\iota_1, \ldots, \iota_k \in I_{\leq D}$, 
we have 
$
A_{\iota_1}\cap\cdots\cap A_{\iota_k} \neq \emptyset
$.
For $\iota_j: Y^{\an} \hookrightarrow \Aff^{n_j}$, we set 
$\iota = (\iota_1, \ldots, \iota_k): Y^{\an} \hookrightarrow \Aff^{n_1+\cdots+n_k}$, 
and we write $\pr_j: \Aff^{n_1+\cdots+n_k} \to \Aff^{n_j}$ for the natural projection. 
Since $\iota_j = \pr_j \circ \iota$, we have $\iota_j \leq \iota$, and 
the induced map $\Trop(Y, \iota) \to \Trop(Y, \iota_j)$ is the restriction to $\Trop(Y, \iota)$ of 
the the natural projection $\Trop(\pr_j): \TT^{n_1+\cdots+n_k} \to \TT^{n_j}$. 
By the definition of the inverse limit, we have $\Trop(\pr_j)(b_{\iota}) = b_{\iota_j}$. 
It follows that $A_{\iota} \subset A_{\iota_j}$, and we obtain 
$A_{\iota_1}\cap\cdots\cap A_{\iota_k} \supset A_{\iota} \neq \emptyset$. 

Following Bourbaki \cite{bourbakiGT}, we say that 
a continuous map between topological spaces
is \emph{proper} if it is universally closed.
We remark that a fiber of a proper map between
Hausdorff spaces is compact (see, for example, \cite[\S10.2, ~Theorem~1]{bourbakiGT}).

We fix any $\iota_0 \in I_{\leq D}$. 
Since $\pi_{\iota_0}: Y^{\an} \to \Trop(Y, \iota_0)$ 
is proper (see \cite[Lemma~2.1, Proposition~2.2]{payne}) and 
since $Y^{\an}$ and $\Trop(Y, \iota_0)$ are Hausdorff spaces, 
$A_{\iota_0}$ is compact.
Consider the family $\{A_{\iota} \cap A_{\iota_0} \}_{\iota \in I_{\leq D}}$
of subspaces of $A_{\iota_0}$.
Since $\{A_{\iota} \}_{\iota \in I_{\leq D}}$ has a finite intersection property, 
so does $\{A_{\iota} \cap A_{\iota_0} \}_{\iota \in I_{\leq D}}$.
Since $A_{\iota} \cap A_{\iota_0}$ is compact for any $\iota \in I_{\leq D}$,
it follows that $\bigcap_{\iota \in I_{\leq D}} A_{\iota} =
\bigcap_{\iota \in I_{\leq D}} ( A_{\iota} \cap A_{\iota_0} )
\neq
\emptyset$. 
We take a $y \in \bigcap_{\iota \in I_{\leq D}} A_{\iota}$. Then we have 
$\pi_\iota(y) = b_\iota$ for any $\iota \in I_{\leq D}$, and thus 
$\mathrm{LTrop}_{\leq D} (y)
=
\left(\varprojlim_{\iota \in I_{\leq D}}\pi_\iota\right) (y) = (b_\iota)_{\iota \in I_{\leq d}}$. 
Thus $\mathrm{LTrop}_{\leq D}$
is surjective. 

\subsubsection*{Homeomorphism}
We have shown that the map $\mathrm{LTrop}_{\leq D}$ in \eqref{eqn:limit:tropical:map:d}  is continuous and bijective. To conclude that it is a homeomorphism, it suffices to show that it is a closed map. 
Since each tropicalization map
$\pi_\iota: Y^{\an} \to \Trop (Y, \iota )$ is proper,
it follows from \cite[\S10.2, Corollary~4]{bourbakiGT} that
the map 
\[
\prod_{\iota \in I_{\leq D}} \pi_{\iota}: Y^{\an} \to \prod_{\iota \in I_{\leq D}} \Trop (Y, \iota ), 
\quad y \mapsto (\pi_\iota(y))_{\iota \in I_{\leq D}}
\]
is proper. In particular, $\prod_{\iota \in I_{\leq D}} \pi_{\iota}$ is a closed map. 
Since the image of the map $\prod_{\iota \in I_{\leq D}} \pi_{\iota}$ is contained in 
$\varprojlim_{\iota \in I_{\leq D} } \Trop (Y, \iota )$ and 
the topology of
$\varprojlim_{\iota \in I_{\leq D} } \Trop (Y, \iota )$ is 
the subspace topology of $\prod_{\iota \in I_{\leq D}} \Trop (Y, \iota)$, 
we see that $\mathrm{LTrop}_{\leq D}$ is a closed map. 
Thus we complete the proof.

\begin{Remark}
\label{remark:morethanD}
It follows from the above proof that
$\mathrm{LTrop}_{\leq D'}$ is a 
surjective closed continuous map for any positive integer $D'$;
if $D' \geq D$, then $\mathrm{LTrop}_{\leq D'}$ is injective, in addition,
and hence it is
a homeomorphism.
\end{Remark}

We finish this paper by remarking that in the course of the proof of 
Theorem~\ref{thm:limit:trop:d:intro}, we have shown the following as a corollary of Theorem~\ref{thm:main:unimodular:faithful}. 
Note that with the notation below, we have 
$Y^{\an} \setminus Y(K) = X^{\an}\setminus X(K)$ and thus 
any skeleton $\Gamma$ of $X^{\an}$ is contained in $Y^{\an}$. 

\begin{Proposition}
\label{prop:last:remark}
Let $Y \subset \mathbb{A}^N = \Spec(K[z_1, \ldots, z_N])$ be a connected smooth affine curve 
that has smooth compactification $X$. 
Let $d$ be the degree of $Y$. 
Define $D$ by \eqref{eqn:thm:limit:trop:d}. 
Then for any skeleton $\Gamma$ of $X^{\an}$, 
there exist nonzero polynomials $f_1, \ldots f_m \in K[z_1, \ldots, z_N]$ of degree 
at most $D$ such that 
\[
 \psi: Y^{\an} \to \RR^m, \quad
 p = (p, |\cdot|) \mapsto (-\log|f_1(p)|, \ldots, -\log|f_m(p)|)
\]
gives a faithful tropicalization of $\Gamma$. 
\end{Proposition}


\end{document}